\documentclass[11pt, oneside]{article}
\usepackage[width=16cm,height=25cm]{geometry}
\usepackage{graphicx}
\usepackage{color}
\usepackage{amssymb}
\usepackage{amsmath}
\usepackage{amsthm}
\usepackage{bm}
\usepackage{array}
\usepackage[utf8]{inputenc}
\usepackage[english]{babel}
\usepackage{lineno}
\usepackage{listings}
\usepackage{textcomp}
\usepackage[ruled]{algorithm2e}
\usepackage{cancel}
\usepackage{url}
\usepackage{aliases}

\newcommand{\resa}[1]{#1}
\newcommand{\resb}[1]{#1}

\newtheorem{theorem}{Theorem}

\newtheorem{definition}{Definition}

\title{Constraint preserving discontinuous Galerkin method for ideal compressible MHD \resb{on 2-D Cartesian grids}}
\author
{
Praveen Chandrashekar\footnote{Corresponding author} \footnote{TIFR Center for Applicable Mathematics, Bangalore, India. Email: \texttt{praveen@tifrbng.res.in}} \ and \
Rakesh Kumar\footnote{TIFR Center for Applicable Mathematics, Bangalore, India. Email: \texttt{rakesh@tifrbng.res.in}}
}
\date{}							

\begin{document}
\maketitle

\lstset{
   language=bash,
   keywordstyle=\bfseries\ttfamily\color[rgb]{0,0,1},
   identifierstyle=\ttfamily,
   commentstyle=\color[rgb]{0.133,0.545,0.133},
   stringstyle=\ttfamily\color[rgb]{0.627,0.126,0.941},
   showstringspaces=false,
   basicstyle=\small,
   numberstyle=\footnotesize,
   numbers=none,
   stepnumber=1,
   numbersep=10pt,
   tabsize=2,
   breaklines=true,
   prebreak = \raisebox{0ex}[0ex][0ex]{\ensuremath{\hookleftarrow}},
   breakatwhitespace=false,
   aboveskip={0.1\baselineskip},
    columns=fixed,
    upquote=true,
    extendedchars=true,
    backgroundcolor=\color[rgb]{0.9,0.9,0.9}
}

\begin{abstract}
We propose a constraint preserving discontinuous Galerkin method for ideal compressible MHD \resb{in two dimensions and using Cartesian grids,} which automatically maintains the global divergence-free property. The approximation of the magnetic field is achieved using Raviart-Thomas polynomials and the DG scheme is based on evolving certain moments of these polynomials which automatically guarantees divergence-free property. We also develop HLL-type multi-dimensional Riemann solvers to estimate the electric field at vertices which are consistent with the 1-D Riemann solvers. When limiters are used, the divergence-free property may be lost and it is recovered by a divergence-free reconstruction step.   We show the performance of the method on a range of test cases up to fourth order of accuracy.
\end{abstract}
\section{Introduction}
The equations governing ideal, compressible MHD are a mathematical model for plasma and form a system of non-linear hyperbolic conservation laws. While it is natural to try to use Godunov-type numerical methods which have been very successful for other non-linear hyperbolic conservation laws~\cite{Toro2009}, the MHD equations have an additional feature in the form of a constraint on the magnetic field $\B$, i.e., the divergence of  $\B$ must be zero, which may not be satisfied by standard schemes. The non satisfaction of this constraint can yield wrong solutions and the methods can also be unstable~\cite{Toth2000}. Hence various strategies have been developed over the years to deal with this issue. Projection-based methods~\cite{Brackbill1980} use standard schemes to update the solution and \'a posteriori correct the magnetic field to make the divergence to be zero by solving an elliptic equation. Hyperbolic divergence cleaning methods have been developed in~\cite{Dedner2002} by introducing an extra Lagrange multiplier or pressure variable. Constrained transport methods~\cite{Evans1988}, \cite{Gardiner2005} are designed to automatically keep some discrete measure of the divergence to be invariant. A key idea in most of these methods is the staggered storage of variables with the magnetic field components being located on the faces and the remaining hydrodynamic variables being located in cell centers. Divergence-free reconstruction of magnetic field have been developed in conjunction with approximate Riemann solvers~\cite{Balsara1999}, \cite{Balsara2001}, \cite{Balsara2004}, \cite{Balsara2009} which also preserve a discrete divergence constraint. Another class of methods~\cite{Powell1999}, \cite{Janhunen2000}, \cite{Bouchut2007}, \cite{Winters2016}, \cite{Chandrashekar2016} aims to  construct a stable scheme without explicitly making the divergence to be zero and are based on Godunov's symmetrized version of the MHD model~\cite{Godunov1972}; moreover these methods are not conservative since the symmetrized model has source terms.

Many of the ideas first developed in a finite volume setting have been extended to discontinuous Galerkin methods which provide a good framework for constructing high order accurate schemes. A locally divergence-free basis was used in combination with Rusanov fluxes in~\cite{Cockburn2004}. A similar approach based on Godunov's symmetrized MHD model has also been developed~\cite{Guillet2019}, whereas entropy stable DG schemes using SBP type operators have been developed in~\cite{Derigs2018a},~\cite{Bohm2018},~\cite{Liu2018a}. The divergence-free reconstruction idea has been combined in a DG scheme in~\cite{Balsara2017a} for induction equation and for Maxwell's equations in~\cite{Hazra2019}. DG schemes which automatically preserve the divergence condition have been developed in~\cite{Li2011}, \cite{Li2012} using central DG idea and in \cite{Fu2018}, \cite{Chandrashekar2019} using Godunov approach. Maintaining the positivity of solutions is very important and recent work shows a close link between this property and a discrete divergence-free condition~\cite{Wu2018a}, see also~\cite{Janhunen2000}. DG schemes which in practice are positive have been developed in standard formulations under the assumption that the Lax-Friedrich scheme is positive~\cite{Cheng2013}. A provably positive DG scheme has been developed in~\cite{Wu2018} based on Godunov's symmetric MHD form and locally divergence-free basis, but the solutions are not guaranteed to be globally divergence-free and the method is not conservative due to the use of Godunov's symmetrized MHD model.

In the present work, we develop a DG scheme based on tensor product polynomials and in particular using Raviart-Thomas polynomials for the approximation of the magnetic field. This builds on the initial work done for induction equation in~\cite{Chandrashekar2019} and is similar in spirit to~\cite{Fu2018}, which developed up to third order schemes using Brezzi-Douglas Marin (BDM) polynomials (see \cite{Brezzi1991}, Section III.3.2) for the magnetic field.
\begin{enumerate}
\item We develop arbitrarily high order DG schemes for ideal MHD which automatically preserve the divergence constraint and the solutions are globally divergence-free.
\item The magnetic field is approximated using Raviart-Thomas polynomials which have tensor product structure. The degrees of freedom are evolved with a hybrid scheme defined both on faces and cells.
\item We develop multi-dimensional HLL and HLLC Riemann solvers which are consistent with their 1-D counterparts.
\item We couple the DG scheme with divergence-free reconstruction method when a TVD-type limiter is applied since the limiter can destroy the divergence-free property.
\end{enumerate}
Because of the DG foundations, we can achieve arbitrarily high order of accuracy with this approach, at least for smooth solutions. The discretization and evolution of the magnetic field has a hybrid nature in the sense that the scheme is defined both on the faces and inside the cells. The method requires numerical fluxes both on the cell faces and the cell corners. On cell faces, a 1-D Riemann problem is present which can be solved approximately, e.g., using HLL-type of schemes. At the cell corners, multiple states meet defining a multi-dimensional Riemann problem. A HLL-type solver can also be formulated for such problems~\cite{Balsara2010},~\cite{Balsara2014a}, and in particular when using DG methods, it is important that this solver should be consistent with the 1-D Riemann solver. In case of discontinuous solutions, some form of TVD-type limiting strategy is required to control spurious numerical oscillations. But such a limiter applied on the magnetic field can destroy the divergence-free property of the solutions; we then perform a local divergence-free reconstruction of the solution following the ideas in~\cite{Hazra2019} which are extended to the case of Raviart-Thomas polynomials. While we cannot prove positivity of solutions in the framework of divergence-free schemes, we show that a heuristic application of scaling limiters can lead to stable computations, but this topic is still an open problem in the context of constraint preserving schemes which rely on Riemann solvers.

The rest of the paper is organized as follows. In Section~(\ref{sec:model}) we list the MHD equations and introduce suitable notation necessary in the paper. Section(\ref{sec:spaces}) explains the structure of the approximating polynomial spaces and resulting degrees of freedom. Section~(\ref{sec:rtrecon}) shows how to construct the magnetic field inside the cell given the degrees of freedom of the Raviart-Thomas polynomials. The DG scheme is explained in Section~(\ref{sec:dgscheme}) for both the hydrodynamic and magnetic variables, and we also discusses constraint satisfaction on the magnetic field divergence by the numerical scheme. The computation of the numerical fluxes is explained in Section~(\ref{sec:flux}) and the limiting procedure in Section~(\ref{sec:lim}). We then present an extensive set of numerical results in Section~(\ref{sec:res}).
\section{Ideal MHD equations}
\label{sec:model}
In the following, we consider only the two dimensional case and it is then convenient to arrange the variables in the following way to deal with the divergence constraint. Let $\rho,p$ be the density and pressure of the gas, $\tote$ be the total energy per unit volume and $\vel = (v_x,v_y,v_z)$ be the gas velocity. The components of the magnetic field are $\bthree = (B_x,B_y,B_z)$. Define
\[
\conh = [\rho, \ \rho \vel, \ \tote, \ B_z ]^\top, \qquad \B = (\Bx, \By)
\]
then the 2-D ideal MHD equations can be written as a system of conservation laws
\begin{equation}
\df{\conh}{t} + \nabla\cdot\flh(\conh,\B) = 0, \qquad \df{\Bx}{t} + \df{E_z}{y} = 0, \qquad \df{\By}{t} - \df{E_z}{x} = 0
\label{eq:mhdeqn}
\end{equation}
where $E_z$ is the electric field in the $z$ direction given by
\[
E_z = v_y B_x - v_x B_y
\]
and the fluxes $\flh = (\flh_x, \flh_y)$ are of the form
\[
\flh_x = \begin{bmatrix}
\rho v_x \\
P + \rho v_x^2 - B_x^2 \\
\rho v_x v_y - B_x B_y \\
\rho v_x v_z - B_x B_z \\
(\tote + P) v_x - B_x (\vel \cdot \bthree)\\
v_x B_z - v_z B_x
\end{bmatrix}, \qquad
\flh_y = \begin{bmatrix}
\rho v_y \\
\rho v_x v_y - B_x B_y \\
P + \rho v_y^2 - B_y^2 \\
\rho v_y v_z -  B_y B_z \\
(\tote + P) v_y - B_y (\vel \cdot \bthree) \\
v_y B_z - v_z B_y
\end{bmatrix}
\]
where the total pressure $P$ and energy $\tote$ are given by
\[
P = p + \frac{1}{2}|\bthree|^2, \qquad \tote = \frac{p}{\gamma-1} + \half \rho |\vel|^2 + \frac{1}{2} |\bthree|^2
\]
Since magnetic monopoles do not exist, the magnetic field $\bthree$ must have zero divergence. In fact if the divergence is zero at the initial time, then under the action of the induction equation, it remains zero at future times also, and hence is referred to as an involution constraint. Since we consider only 2-D problems in this work, the divergence-free  condition is equivalent to the 2-D divergence of $\B$ being zero, i.e.,
\[
\nabla\cdot\B = \df{B_x}{x} + \df{B_y}{y} = 0
\]
In the above discussion, we have written the equations in the form~(\ref{eq:mhdeqn}) which is suitable for the implementation of the divergence-free scheme in 2-D. For the computation of the numerical fluxes, we have to consider all the equations together in conservation form which can be written as
\begin{equation}
\df{\con}{t} + \df{\flx}{x} + \df{\fly}{y} = 0
\label{eq:mhdfull}
\end{equation}
where
\[
\con = \begin{bmatrix}
\rho \\
\rho v_x \\
\rho v_y \\
\rho v_z \\
\tote \\
B_x \\
B_y \\
B_z \end{bmatrix}, \qquad
\flx = \begin{bmatrix}
\rho v_x \\
P + \rho v_x^2 - B_x^2 \\
\rho v_x v_y - B_x B_y \\
\rho v_x v_z - B_x B_z \\
(\tote + P) v_x - B_x (\vel \cdot \bthree)\\
0 \\
-E_z \\
v_x B_z - v_z B_x
\end{bmatrix}, \qquad
\fly = \begin{bmatrix}
\rho v_y \\
\rho v_x v_y - B_x B_y \\
P + \rho v_y^2 - B_y^2 \\
\rho v_y v_z - B_y B_z \\
(\tote + P) v_y - B_y (\vel \cdot \bthree) \\
E_z \\
0 \\
v_y B_z - v_z B_y
\end{bmatrix}
\]
Let $\mathcal{A}_x = \flx'(\con)$ and $\mathcal{A}_y = \fly'(\con)$ be the flux Jacobians. The Jacobian matrices have real eigenvalues given by
\[
\lambda(\mathcal{A}_d) = \{ v_d - c_{fd}, \ v_d - c_{sd}, \ v_d-c_a, \ v_d, \ 0, \ v_d+c_a, \ v_d + c_{sd}, \ v_d+c_{fd} \}, \qquad d=x,y
\]
where $c_{sd}$, $c_{fd}$ are the slow and fast magnetosonic speeds and $c_a$ is the Alfven wave speed. The Alfven wave speed is given by
\[
c_a = \frac{|B_d|}{\sqrt{\rho}}
\]
and the magnetosonic speeds are given by
\[
c_{sd} = \sqrt{\half \left[ a^2 + |\bm{b}|^2 - \sqrt{(a^2 + |\bm{b}|^2)^2 - 4 a^2 b_d^2} \right]}, \quad c_{fd} = \sqrt{\half \left[ a^2 + |\bm{b}|^2 + \sqrt{(a^2 + |\bm{b}|^2)^2 - 4 a^2 b_d^2} \right]}
\]
where
\[
a = \sqrt{\frac{\gamma p}{\rho}}, \qquad \bm{b} = \frac{|\bthree|}{\sqrt{\rho}}
\]
with $a$ being the sound speed.
\section{Approximation spaces}
\label{sec:spaces}
We map each cell to the reference cell $[-\half,+\half] \times [-\half,+\half]$ with coordinates $(\xi,\eta)$.  Define the tensor product polynomials by
\[
\tpoly_{r,s} = \textrm{span}\{ \xi^i \eta^j : 0 \le i \le r, \ 0 \le j \le s \}
\]
As basis functions for polynomials, we will first construct one dimensional orthogonal polynomials given by
\[
\phi_0(\xi) = 1, \quad
\phi_1(\xi) = \xi, \quad
\phi_2(\xi) = \xi^2 - \frac{1}{12}, \quad
\phi_3(\xi) = \xi^3 - \frac{3}{20} \xi, \quad
\phi_4(\xi) = \xi^4 - \frac{3}{14} \xi^2 + \frac{3}{560}
\]
whose mass matrix  is diagonal with entries given by $m_i = \intod \phi_i^2(\xi)\ud\xi$. Let $k \ge 0$ be the degree of approximation. The hydrodynamic variables $\conh$ are approximated in each cell in the space $\tpoly_{k,k}$, which can be written as
\begin{equation}
\conh(\xi,\eta) = \sum_{i=0}^k \sum_{j=0}^k \conh_{ij} \phi_i(\xi) \phi_j(\eta) \in \tpoly_{k,k}
\label{eq:hydrosol}
\end{equation}
Note that this approximation is in general discontinuous across the cell faces.

Let us  approximate the normal component of $\B$ on each face by one dimensional polynomials of degree $k$. On the vertical faces of cells, we will approximate the $x$-component of $\B$ by
\begin{equation}
b_x(\eta) = \sum_{j=0}^k a_j \phi_j(\eta) \in \poly_k(\eta)
\label{eq:bx}
\end{equation}
while on the horizontal faces, the  $y$-component is approximated by
\begin{equation}
b_y(\xi) = \sum_{j=0}^k b_j \phi_j(\xi) \in \poly_k(\xi)
\label{eq:by}
\end{equation}
For $k \ge 1$, let us also define certain {\em cell moments} which we will use in determining the magnetic field inside the cells. These moments are defined as
\[
\alpha_{ij} = \alpha_{ij}(B_x) := \frac{1}{m_{ij}} \inttd B_x(\xi,\eta) \phi_i(\xi) \phi_j(\eta)\ud\xi\ud\eta, \qquad 0 \le i \le k-1, \quad 0 \le j \le k
\]
\[
\beta_{ij} = \beta_{ij}(B_y) := \frac{1}{m_{ij}} \inttd B_y(\xi,\eta) \phi_i(\xi) \phi_j(\eta)\ud\xi\ud\eta, \qquad 0 \le i \le k, \quad 0 \le j \le k-1
\]
where $m_{ij}$ is given by
\[
m_{ij} = \inttd [\phi_i(\xi) \phi_j(\eta)]^2 \ud\xi\ud\eta = m_i m_j
\]
The test functions used to define the $\alpha$ moments belong to $\tpoly_{k-1,k}$ and those used to define the $\beta$ moments belong to $\tpoly_{k,k-1}$. In fact, these quantities correspond to the degrees of freedom used to define the Raviart-Thomas polynomials~\cite{Raviart1977}, which provide $\hdiv$ conforming approximation of vector fields. In our numerical approach, {\em the quantities $b_x$, $b_y$, $\alpha$, $\beta$ form the discretization of the magnetic field} and these quantities will be evolved forward in time by a DG scheme. Using this information we will reconstruct the magnetic field $\B$ inside the cells. We note that $\alpha_{00}$ and $\beta_{00}$ are the mean values of $B_x$ and $B_y$ in the cells and $a_0$, $b_0$ are the mean values of the normal component of $\B$ on the corresponding faces.

\resb
{
The next section describes how to construct the vector field $\B$ from the information contained in the face and cell moments by solving a local reconstruction problem. We are in particular interested in obtaining approximations which are globally divergence-free vector fields.
\begin{definition}[Globally divergence-free]
We will say that a vector field $\B$ defined on a mesh is globally divergence-free if
\begin{enumerate}
\item $\nabla\cdot\B = 0$ in each cell $K$
\item $\B \cdot \un $ is continuous at each face $F$
\end{enumerate}
\end{definition}
}
\section{RT reconstruction problem}
\label{sec:rtrecon}
\begin{figure}
 \centering
 \includegraphics[width=0.4\textwidth]{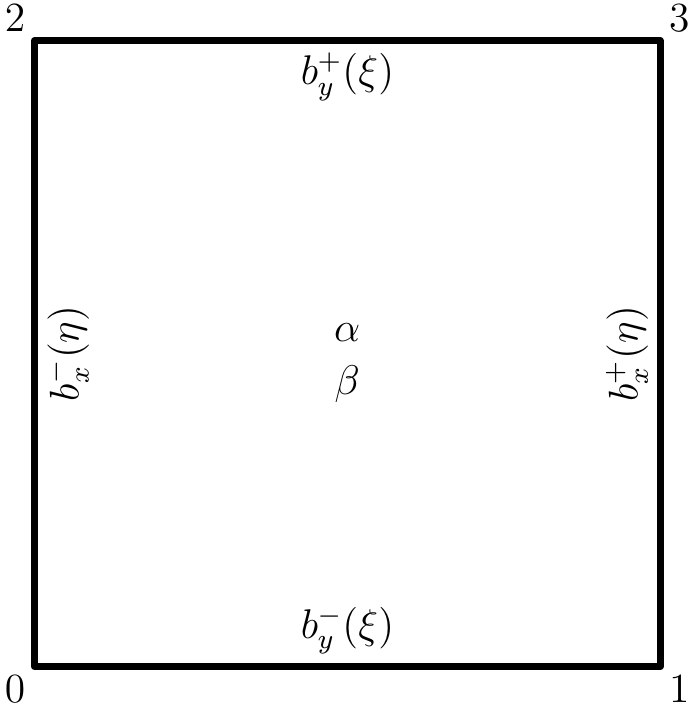}
 \caption{Location of dofs for $\B$}
 \label{fig:bvars}
\end{figure}
Consider a cell as shown in Figure~\ref{fig:bvars}. We are given normal components of $\B$ on the faces in the form of polynomials $b_x^\pm(\eta) \in \poly_k$ and $b_y^\pm(\xi) \in \poly_k$, and also the set of cell moments
\[
\{ \alpha_{ij}, \ 0 \le i \le k-1, \ 0 \le j \le k\}, \qquad \{ \beta_{ij}, \ 0 \le i \le k, \ 0 \le j \le k-1\}
\]
Using this information, we want to construct the magnetic field vector inside the cell. \\

\noindent
{\em RT reconstruction problem}: Find $B_x \in \tpoly_{k+1,k}$ and $B_y \in \tpoly_{k,k+1}$ such that
\[
B_x(\pm \shalf, \eta) = b_x^\pm(\eta), \quad \eta \in [-\shalf,\shalf], \qquad B_y(\xi,\pm\shalf) = b_y^\pm(\xi), \quad \xi \in [-\shalf,\shalf]
\]
\[
\frac{1}{m_{ij}} \inttd B_x(\xi,\eta) \phi_i(\xi) \phi_j(\eta)\ud\xi\ud\eta = \alpha_{ij}, \qquad 0 \le i \le k-1, \quad 0 \le j \le k
\]
\[
\frac{1}{m_{ij}} \inttd B_y(\xi,\eta) \phi_i(\xi) \phi_j(\eta)\ud\xi\ud\eta = \beta_{ij}, \qquad 0 \le i \le k, \quad 0 \le j \le k-1
\]
Since $\textrm{dim }\tpoly_{k+1,k} = \textrm{dim }\tpoly_{k,k+1} = (k+1)(k+2)$, we have $2(k+1)(k+2)$ coefficients to be determined. On the faces we are given $4(k+1)$ pieces of information in terms of the polynomials $b_x^\pm$, $b_y^\pm$, and inside the cell, we have $2k(k+1)$ pieces of information in terms of the cell moments $\alpha,\beta$, and hence we have as many equations as the number of unknowns.
\begin{theorem}
(1) The RT reconstruction problem has a unique solution. (2) If the data $b_x, b_y, \alpha, \beta$ correspond to a divergence-free vector field, then the reconstructed field is also divergence-free.
\end{theorem}
For a proof of the above theorem, we refer the reader to~\cite{Brezzi1991},~\cite{Chandrashekar2019}.  Note that this reconstruction is very local to each cell; it uses data in the cell and on its faces. \resb{We remark that this reconstruction is different from the reconstruction performed in finite volume methods to recover the solution from cell averages. In the present case, we have the full information available in $(b_x,b_y,\alpha,\beta)$ and it is converted into a spatial polynomial $(B_x,B_y)$ by the RT reconstruction step. The two sets od data $(b_x,b_y,\alpha,\beta)$ and $(B_x,B_y)$ contain the same information and are two different ways to represent the magnetic field.}


\paragraph{Properties of $\B$}
The vector field $\B$ obtained from this reconstruction process satisfies certain conditions. Firstly, note that $\B \cdot \un$ has a unique value at a face common to two cells, i.e., the normal component of $\B$ is continuous at all the faces. Secondly, if the data $b_x^\pm, b_y^\pm, \alpha_{ij}, \beta_{ij}$ comes from a divergence-free vector field, then the reconstructed field is also divergence-free~\cite{Chandrashekar2019}. \resb{Hence, the vector field $\B$ will be globally divergence-free also}. The initial data must be generated carefully in order to ensure divergence-free condition. An initial divergence-free vector field has a corresponding stream function which we interpolate to a continuous space of $\tpoly_{k+1,k+1}$ polynomials. The polynomials $b_x, b_y$ on the faces are set equal to the curl of this interpolated stream function and the cell moments $\alpha,\beta$ are obtained by computing the integrals exactly with a Gauss-Legendre rule where we again use the curl of the interpolated stream function, and is explained in Appendix~\ref{sec:ic}. \resb{The DG scheme to be explained in later sections will then ensure that the solutions remain divergence-free at future times also.}

We will now give the solution of the reconstruction problem at different orders. To do this, we first write the reconstructed magnetic field components as a tensor product of orthogonal 1-D polynomials
\begin{equation}
B_x(\xi,\eta) = \sum_{i=0}^{k+1} \sum_{j=0}^k a_{ij} \phi_i(\xi) \phi_j(\eta) \in \tpoly_{k+1,k}, \qquad B_y(\xi,\eta) = \sum_{i=0}^{k} \sum_{j=0}^{k+1} b_{ij} \phi_i(\xi) \phi_j(\eta) \in \tpoly_{k,k+1}
\label{eq:BxBy}
\end{equation}
\paragraph{Reconstruction for $k=0$}
In this case we have constant approximation on the faces for the normal components
\[
b_x^\pm(\eta) = a_0^\pm,  \qquad b_y^\pm(\xi) = b_0^\pm
\]
and the vector field inside the cell is of the form
\[
\Bx(\xi,\eta) = a_{00} + a_{10} \phi_1(\xi), \qquad \By(\xi,\eta) = b_{00} + b_{01} \phi_1(\eta)
\]
which has a dimension of four. The solution of the reconstruction problem is given by
\begin{equation*}
a_{00} = \half(a_0^- + a_0^+), \quad b_{00} = \half (b_0^- + b_0^+), \quad
a_{10} = a_0^+ - a_0^-, \quad b_{01} = b_0^+ - b_0^-
\end{equation*}
Note that no cell moments are present at this order and the reconstruction is determined by the face solution alone.
\paragraph{Reconstruction for $k=1$}
In this case we have linear approximation on the faces and in addition, there are four cell moments. The solution of the reconstruction problem is given in Table~\ref{tab:rec1}.

\begin{table}
\boxed{
\begin{minipage}{0.48\textwidth}
\begin{equation*}
\begin{aligned}
&a_{ij} = \alpha_{ij}, \quad i=0, \ 0\le j \le 1  \\
&a_{10} = a_0^+ - a_0^- \\
&a_{20} = 3(a_0^- + a_0^+ - 2 \alpha_{00}) \\
&a_{11} = a_1^+ - a_1^-  \\
&a_{21} = 3 (a_1^- + a_1^+ - 2 \alpha_{01})
\end{aligned}
\end{equation*}
\end{minipage}
\begin{minipage}{0.48\textwidth}
\begin{equation*}
\begin{aligned}
&b_{ij} = \beta_{ij}, \quad 0 \le i \le 1, \ j=0 \\
&b_{01} = b_0^+ - b_0^-  \\
&b_{02} = 3(b_0^- + b_0^+ - 2 \beta_{00})  \\
&b_{11} = b_1^+ - b_1^-  \\
&b_{12} = 3(b_1^- + b_1^+ - 2 \beta_{10})
\end{aligned}
\end{equation*}
\end{minipage}
}
\caption{Solution of RT reconstruction problem for $k=1$}
\label{tab:rec1}
\end{table}
\paragraph{Reconstruction for $k=2$}
In this case we have quadratic approximation on the faces and in addition, there are 12 cell moments. The solution of the reconstruction problem is given in Table~\ref{tab:rec2}.
\begin{table}
\boxed{
\begin{minipage}{0.5\textwidth}
\begin{equation*}
\begin{aligned}
& a_{ij} = \alpha_{ij}, \quad 0\le i\le 1, \ 0\le j\le 2 \\
& a_{20}=3(a_0^{-}+a_0^{+}-2\alpha_{00})\\
& a_{30}= 10(a_0^{+}-a_0^{-} - \alpha_{10}) \\
& a_{21}=3(a_1^{-}+a_1^{+}-2\alpha_{01})\\
& a_{31}= 10(a_1^{+} - a_1^{-} - \alpha_{11}) \\
& a_{22}=3(a_2^{-}+a_2^{+}-2\alpha_{02}) \\
& a_{32}= 10(a_2^{+} - a_2^{-} - \alpha_{12})
\end{aligned}
\end{equation*}
\end{minipage}
\begin{minipage}{0.5\textwidth}
\begin{equation*}
\begin{aligned}
& b_{ij}  =\beta_{ij}, \quad 0\le i\le 2, \ 0\le j\le 1  \\
& b_{02}=3 (b_0^{-}+b_0^{+}-2\beta_{00}) \\
& b_{03}= 10(b_0^{+} - b_{0}^{-} - \beta_{01})\\
& b_{12} = 3(b_1^{-}+b_1^{+}-2\beta_{10}) \\
& b_{13} = 10(b_1^{+} - b_1^{-} - \beta_{11}) \\
& b_{22}=3(b_2^{-}+b_2^{+}-2\beta_{20}) \\
& b_{23}= 10(b_2^{+} - b_2^{-} - \beta_{21})
\end{aligned}
\end{equation*}
\end{minipage}
}
\caption{Solution of RT reconstruction problem for $k=2$}
\label{tab:rec2}
\end{table}
\paragraph{Reconstruction for $k=3$}
In this case we have cubic approximation on the faces and in addition, there are 24 cell moments. The solution of the reconstruction problem is given in Table~\ref{tab:rec3}.
\begin{table}
\boxed{
\begin{minipage}{0.5\textwidth}
\begin{equation*}
\begin{aligned}
& a_{ij} = \alpha_{ij}, \quad 0\le i\le 2, \ 0\le j\le 3  \\
& a_{30} = 10(a_0^{+}-a_0^{-} - \alpha_{10}) \\
& a_{40} = \frac{35}{3}(3 a_0^{-} + 3 a_0^{+} -6 \alpha_{00} -\alpha_{20})\\
& a_{31} = 10(a_1^{+} - a_1^{-} - \alpha_{11}) \\
& a_{41} = \frac{35}{3}(3 a_1^{-} + 3 a_1^{+} -6 \alpha_{01} -\alpha_{21})\\
& a_{32} = 10(a_2^{+} - a_2^{-} - \alpha_{12})\\
& a_{42} = \frac{35}{3}(3 a_2^{-} + 3 a_2^{+} -6 \alpha_{02} -\alpha_{22})\\
& a_{33} = 10(a_3^{+} - a_3^{-} - \alpha_{13})\\
& a_{43} = \frac{35}{3}(3 a_3^{-} + 3 a_3^{+} -6 \alpha_{03} -\alpha_{23})
\end{aligned}
\end{equation*}
\end{minipage}
\begin{minipage}{0.5\textwidth}
\begin{equation*}
\begin{aligned}
& b_{ij}  =\beta_{ij}, \quad 0\le i\le 3, \ 0\le j\le 2   \\
& b_{03} = 10(b_0^{+} - b_{0}^{-} - \beta_{01})\\
& b_{04} = \frac{35}{3}(3 b_0^{-} + 3 b_0^{+} -6 \beta_{00} -\beta_{02})\\
& b_{13} = 10(b_1^{+} - b_1^{-} - \beta_{11}) \\
& b_{14} = \frac{35}{3}(3 b_1^{-} + 3 b_1^{+} -6 \beta_{10} -\beta_{12})\\
& b_{23} = 10(b_2^{+} - b_2^{-} - \beta_{21})\\
& b_{24} = \frac{35}{3}(3 b_2^{-} + 3 b_2^{+} -6 \beta_{20} -\beta_{22})\\
& b_{33} = 10(b_3^{+} - b_3^{-} - \beta_{31}) \\
& b_{34} = \frac{35}{3}(3 b_3^{-} + 3 b_3^{+} -6 \beta_{30} -\beta_{32})
\end{aligned}
\end{equation*}
\end{minipage}
}
\caption{Solution of RT reconstruction problem for $k=3$}
\label{tab:rec3}
\end{table}

\section{Numerical scheme}
\label{sec:dgscheme}
The basic unknowns in our scheme are the polynomials $b_x,b_y$ approximating the normal component of $\B$ on the cell faces, the cell moments $\alpha,\beta$, and the polynomials approximating the hydrodynamic variables and $B_z$ inside the cells which are grouped inside the set $\conh$. We will devise DG schemes to evolve all these quantities forward in time. We perform spatial discretization using DG scheme and then solve the resulting set of ODE using a Runge-Kutta scheme for time integration.
\begin{figure}
   \centering
   \includegraphics[width=0.6\textwidth]{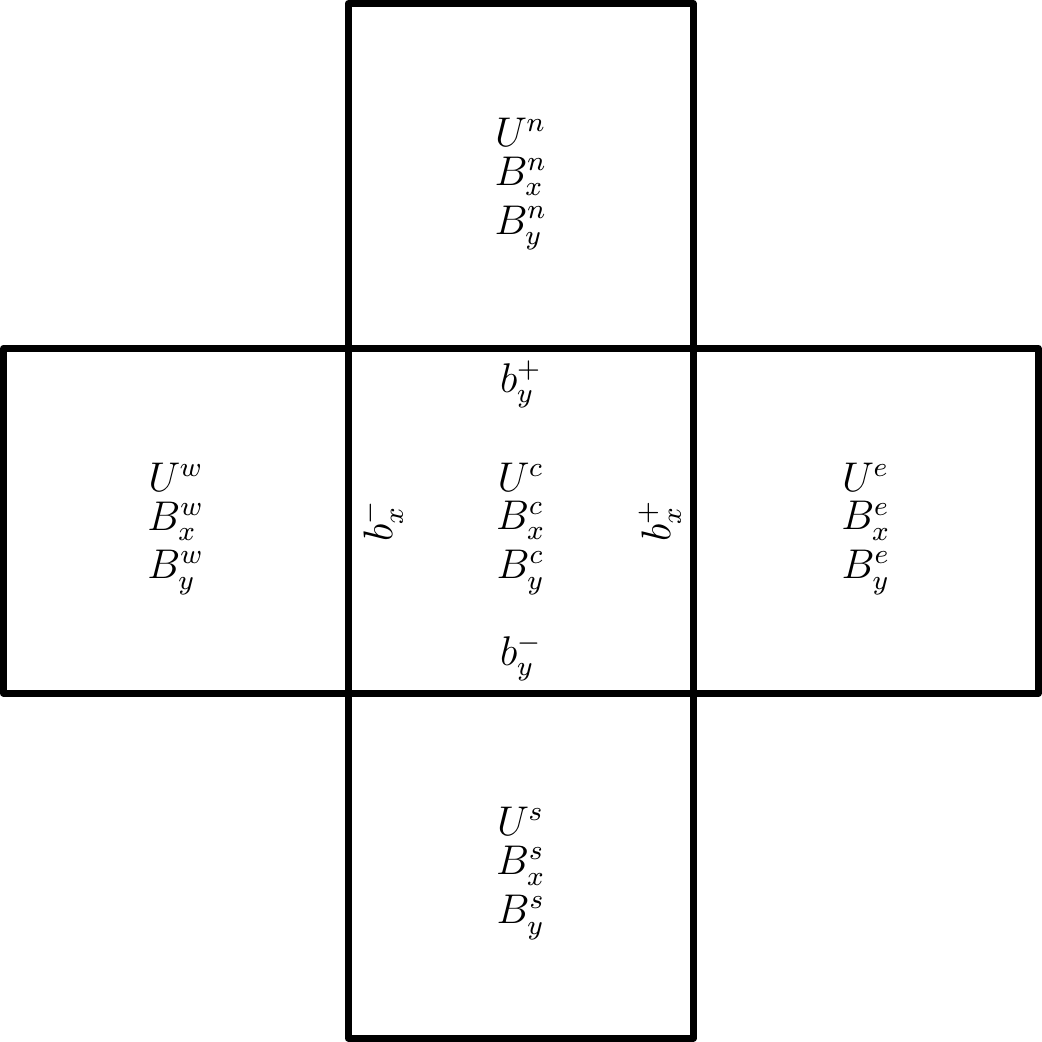}
   \caption{Stencil and variables for DG scheme}
   \label{fig:stencil}
\end{figure}
\subsection{Discontinuous Galerkin method for $\B$ on the faces}
The normal component of $\B$ has been approximated on the faces of our mesh and we want to construct a numerical scheme to evolve these values forward in time. If we observe the equation governing $B_x$, we see that it evolves in time only due to the $y$ derivative of the electric field $E_z$. Restricting ourselves to a vertical face, we see that we have a one dimensional PDE for $B_x$ which can discretized using a 1-D DG scheme applied on the face. Multiplying by a test function $\phi_i(\eta)$ and integrating by parts on a vertical face yields
\[
\intod \df{b_x}{t} \phi_i \ud\eta - \frac{1}{\dy} \intod \hEz \dd{\phi_i}{\eta} \ud\eta + \frac{1}{\dy} [ \tEz \phi_i ] = 0, \qquad 0 \le i \le k
\]
where $\hEz$ is obtained from a 1-D Riemann solver and $\tEz$ is obtained from a multi-D Riemann solver. The face integral is computed using $(k+1)$-point Gauss-Legendre quadrature which results in the semi-discrete scheme
\begin{equation}
m_i \dd{a_i}{t} - \frac{1}{\Delta y} \sum_q \hEz(\eta_q) \dd{\phi_i}{\eta}(\eta_q) \omega_q + \frac{1}{\Delta y}[ \tEz(\shalf) \phi_i(\shalf) - \tEz(-\shalf) \phi_i(-\shalf)] = 0
\label{eq:bxdg}
\end{equation}
Similarly on the horizontal faces, using a test function $\phi_i(\xi)$, the DG scheme for $B_y$ is given by
\[
\intod \df{b_y}{t} \phi_i \ud\xi + \frac{1}{\dx} \intod \hEz \dd{\phi_i}{\xi} \ud\xi - \frac{1}{\dx}[ \tEz \phi_i ] = 0, \qquad 0 \le i \le k
\]
Using $(k+1)$-point Gauss-Legendre quadrature on the face, we obtain the semi-discrete scheme
\begin{equation}
m_i \dd{b_i}{t} + \frac{1}{\Delta x} \sum_q \hEz(\xi_q) \dd{\phi_i}{\xi}(\xi_q) \omega_q - \frac{1}{\Delta x}[ \tEz(\shalf) \phi_i(\shalf) - \tEz(-\shalf) \phi_i(-\shalf)] = 0
\label{eq:bydg}
\end{equation}
\subsection{Discontinuous Galerkin method for $\B$ in the cells}
For $k \ge 1$ we have additional cell moments that are required to reconstruct the magnetic field inside the cells. We can derive evolution equations for these moments using the induction equation and using integration by parts to transfer derivatives onto the test functions. This leads to the following set of semi-discrete equations,
\begin{eqnarray*}
m_{ij} \dd{\alpha_{ij}}{t} &=& \inttd \df{B_x}{t}\phi_i(\xi)\phi_j(\eta)\ud\xi\ud\eta = - \frac{1}{\dy}\inttd \df{E_z}{\eta}\phi_i(\xi)\phi_j(\eta)\ud\xi\ud\eta \\
&=& - \frac{1}{\dy}\intod [\hEz(\xi,\shalf)\phi_i(\xi)\phi_j(\shalf) - \hEz(\xi,-\shalf)\phi_i(\xi)\phi_j(-\shalf)] \ud\xi \\
&& + \frac{1}{\dy} \inttd E_z(\xi,\eta) \phi_i(\xi) \phi_j'(\eta) \ud\xi\ud\eta, \quad 0 \le i \le k-1, \quad 0 \le j \le k
\end{eqnarray*}
and
\begin{eqnarray*}
m_{ij} \dd{\beta_{ij}}{t} &=& \inttd \df{B_y}{t}\phi_i(\xi)\phi_j(\eta)\ud\xi\ud\eta = \frac{1}{\dx}\inttd \df{E_z}{\xi}\phi_i(\xi)\phi_j(\eta)\ud\xi\ud\eta \\
&=& \frac{1}{\dx}\intod [\hEz(\shalf,\eta)\phi_i(\shalf)\phi_j(\eta) - \hEz(-\shalf,\eta)\phi_i(-\shalf)\phi_j(\eta)] \ud\eta \\
&& - \frac{1}{\dx} \inttd E_z(\xi,\eta) \phi_i'(\xi) \phi_j(\eta) \ud\xi\ud\eta, \quad 0 \le i \le k, \quad 0 \le j \le k-1
\end{eqnarray*}
Note that the numerical fluxes $\hEz$ required in the face integrals are obtained from a 1-D Riemann solver. We observe that this is not a Galerkin method because the equation for $\alpha,\beta$ has test functions which are different from the RT polynomials. For example, in the $\alpha$ equation, we use test functions from $\tpoly_{k-1,k}$ whereas $B_x \in \tpoly_{k+1,k}$.
\subsection{Discontinuous Galerkin method for $\conh$ inside cells}
The hydrodynamic variables and $B_z$ which are grouped into the variable $\conh$ are approximated by $\tpoly_{k,k}$ polynomials inside each cell. We will apply a standard DG scheme to the first equation in~(\ref{eq:mhdeqn}); multiplying this equation by a test function $\Phi_i(\xi,\eta)$ and performing an integration by parts over one cell, we get
\[
\begin{aligned}
\inttd \df{\conh^c}{t} \Phi_i(\xi,\eta) \ud\xi \ud\eta - & \inttd\left[ \frac{1}{\dx} \flh_x \df{\Phi_i}{\xi} + \frac{1}{\dy} \flh_y \df{\Phi_i}{\eta} \right] \ud\xi\ud\eta \\
+ & \frac{1}{\dx} \intod \nflh_x^+ \Phi_i(\shalf,\eta)\ud\eta - \frac{1}{\dx} \intod \nflh_x^- \Phi_i(-\shalf,\eta)\ud\eta \\
+ & \frac{1}{\dy} \intod \nflh_y^+ \Phi_i(\xi,\shalf)\ud\xi - \frac{1}{\dy} \intod \nflh_y^- \Phi_i(\xi,-\shalf)\ud\xi = 0
\end{aligned}
\]
where the test functions $\{ \Phi_i, \ i=0,1,\ldots,(k+1)^2-1\}$, are the tensor product basis functions of $\tpoly_{k,k}$ arranged as a one dimensional sequence, $\nflh_x^-$, $\nflh_x^+$ are the numerical fluxes on the left and right faces obtained from a 1-D Riemann solver, and, $\nflh_y^-$, $\nflh_y^+$ are the numerical fluxes on bottom and top faces obtained from a 1-D Riemann solver. The integral inside the cell is evaluated using a tensor product of $(k+1)$-point Gauss-Legendre quadrature while the face integrals are evaluated using $(k+1)$-point Gauss-Legendre quadrature. The fluxes used in the above DG scheme are computed from the solution variables in the following way,
\[
\flh_x = \flh_x(\conh^c, B_{x}^c, B_{y}^c), \qquad \flh_y = \flh_y(\conh^c, B_{x}^c, B_{y}^c)
\]
\[
\nflh_x^+ = \nflh_x((\conh^c, b_x^+, B_{y}^c), (\conh^e, b_x^+, B_{y}^e)), \qquad \nflh_x^- = \nflh_x((\conh^w, b_x^-, B_{y}^w), (\conh^c, b_x^-, B_{y}^c))
\]
\[
\nflh_y^+ = \nflh_y((\conh^c, B_{x}^c,b_y^+), (\conh^n, B_{x}^n,b_y^+)), \qquad \nflh_y^- = \nflh_y((\conh^s, B_{x}^s,b_y^-), (\conh^c, B_{x}^c,b_y^-))
\]
and Figure~\ref{fig:stencil} shows the notation used for the arguments in the flux functions. On the faces, we make use of the normal component $b_x, b_y$ that is already available on the face, and the remaining component is obtained from the RT reconstruction $\B$ inside the cell.
\subsection{Constraints on the magnetic field}
\label{sec:con}
We have completely specified the spatial discretization for all the variables. We are moreover interested in ensuring that the magnetic field remains divergence-free at all times if the initial condition was divergence-free.  The continuity of the normal component of $\B$ is ensured since this is directly approximated in terms of the 1-D polynomials $b_x, b_y$. To be globally divergence-free, the vector field $\B$ must have zero divergence.
\begin{theorem}
The DG scheme satisfies
\[
\dd{}{t}\int_K (\nabla\cdot\B) \phi \ud x \ud y = 0, \qquad \forall \phi \in \tpoly_{k,k}
\]
and since $\nabla\cdot\B \in \tpoly_{k,k}$ this implies that $\nabla\cdot\B$ is constant with respect to time. If $\nabla\cdot\B = 0$ everywhere at the initial time, then this is true at any future time also.
\end{theorem}
The proof can be found in~\cite{Chandrashekar2019} and so we do not repeat it here. In many applications, shocks or other discontinuities may be present or they can develop even from smooth initial data. In these situations,  some form of limiter is absolutely necessary in order to control the numerical oscillations and keep the computations stable. However, if a limiter is used in a post-processing step which is how limiters are applied in DG schemes, then the limited solution may not be divergence-free. We will address this issue subsequently in the paper.
\section{Numerical fluxes}
\label{sec:flux}
\begin{figure}
\begin{center}
\begin{tabular}{cc}
\includegraphics[width=0.45\textwidth]{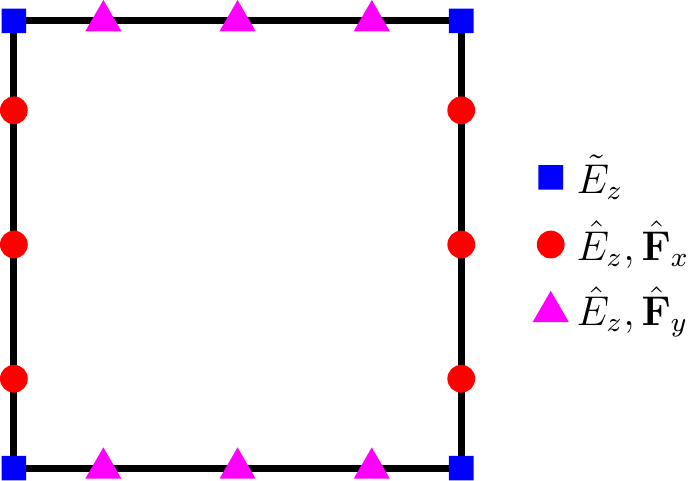} &
\includegraphics[width=0.4\textwidth]{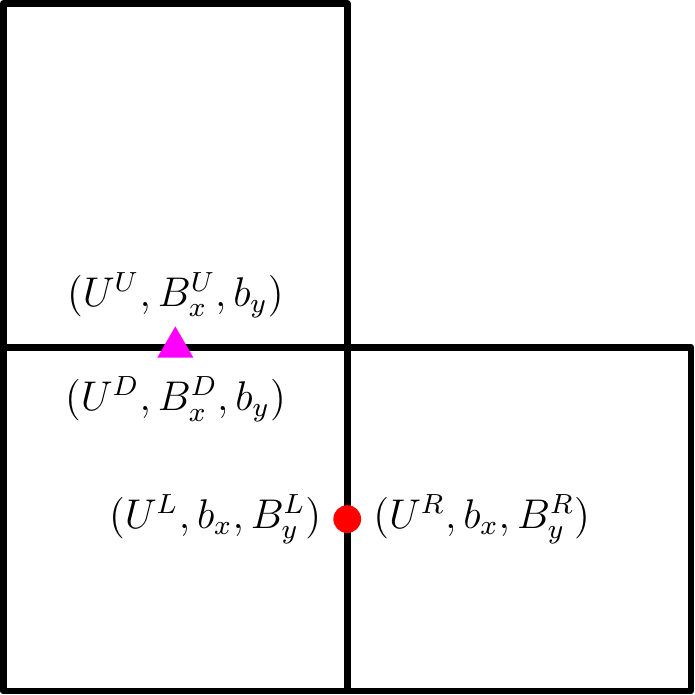} \\
(a) & (b)
\end{tabular}
\caption{(a) Face quadrature points and numerical fluxes. (b) 1-D Riemann problems at a vertical and horizontal face of a cell}
\label{fig:riem1d}
\end{center}
\end{figure}
A major component of the DG scheme is the  specification of numerical fluxes required on the faces and vertices of the cells. We use Gauss-Legendre quadrature on the faces and the required numerical fluxes are shown in figure~(\ref{fig:riem1d}a). The electric field $E_z$ is required at the cell vertices and the face quadrature points. The numerical fluxes $\nflh_x, \nflh_y$ are not required at the cell vertices but only at the face quadrature points which are interior to each face. These fluxes are determined by approximate solution of 1-D and 2-D Riemann problems. On the cell faces, we have a 1-D Riemann problem since the solution is possibly discontinuous. For example at a vertical cell face, see Figure~\ref{fig:riem1d}b, we have the left state $(\conh^L, B_x, B_y^L)$ and a right state $(\conh^R, B_x, B_y^R)$. Note that $B_x$ which is the normal component on a vertical face, has the same value on both sides since this component is directly approximated on the face. The tangential component $B_y$ is obtained by the RT reconstruction in the two cells adjacent to the vertical face and can be discontinuous. We now have the two conserved state variables $\con^L = \con(\conh^L, B_x, B_y^L)$ and $\con^R = \con(\conh^R, B_x, B_y^R)$ and let $\nflx$ denote the numerical flux obtained by solving the 1-D MHD Riemann problem corresponding to these two states. Note that we have to solve the Riemann problem for the full MHD system~(\ref{eq:mhdfull}) to obtain this flux. From this flux, we can obtain the fluxes required for our DG scheme as follows.
\[
\nflh_x = \begin{bmatrix}
(\nflx)_1 \\
(\nflx)_2 \\
(\nflx)_3 \\
(\nflx)_4 \\
(\nflx)_5 \\
(\nflx)_8
\end{bmatrix}, \qquad \hEz = -(\nflx)_7
\]
Similarly, at any horizontal face, see Figure~\ref{fig:riem1d}b, we have the bottom state $\con^D = \con(\conh^D, B_x^D, B_y)$ and top state $\con^U = \con(\conh^U, B_x^U, B_y)$, where we now see that the normal component $B_y$ is continuous. The solution of the 1-D MHD Riemann problem with the two states $\con^D$, $\con^U$ yields the numerical flux $\nfly$ from which we obtain the fluxes required for our DG scheme
\[
\nflh_y = \begin{bmatrix}
(\nfly)_1 \\
(\nfly)_2 \\
(\nfly)_3 \\
(\nfly)_4 \\
(\nfly)_5 \\
(\nfly)_8
\end{bmatrix}, \qquad \hEz = (\nflx)_6
\]
\begin{figure}
\begin{center}
\begin{tabular}{cc}
\includegraphics[width=0.465\textwidth]{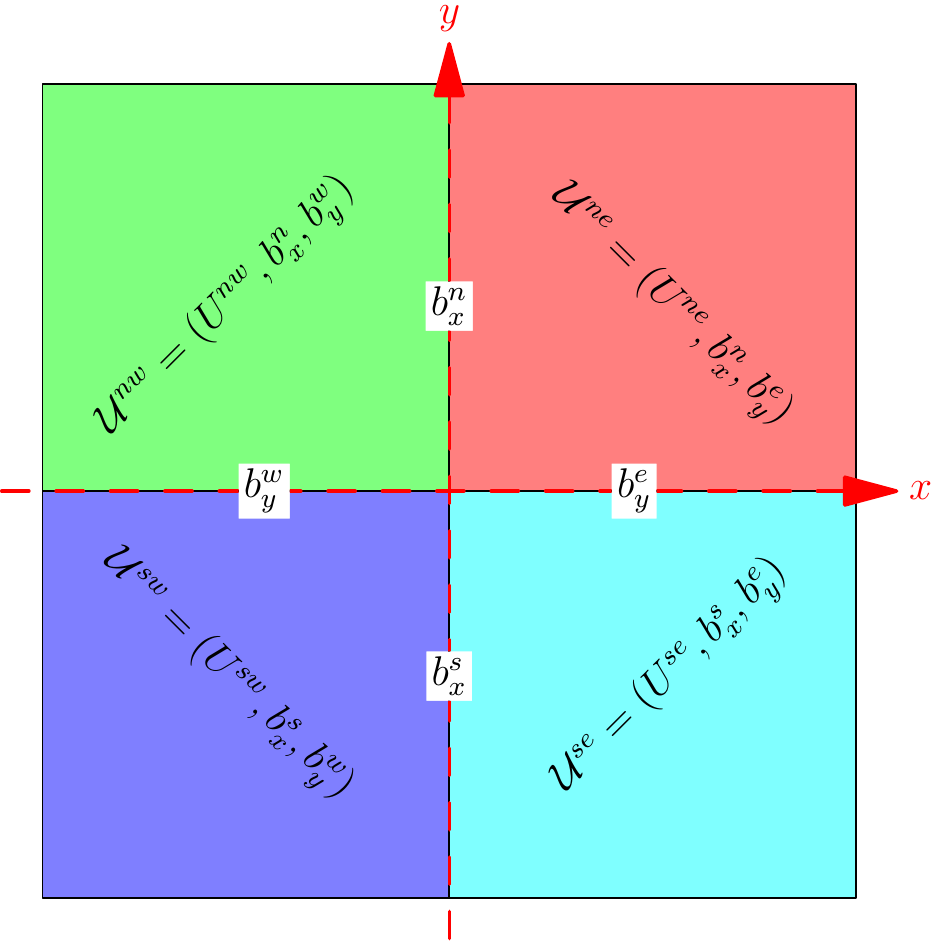} &
\includegraphics[width=0.49\textwidth]{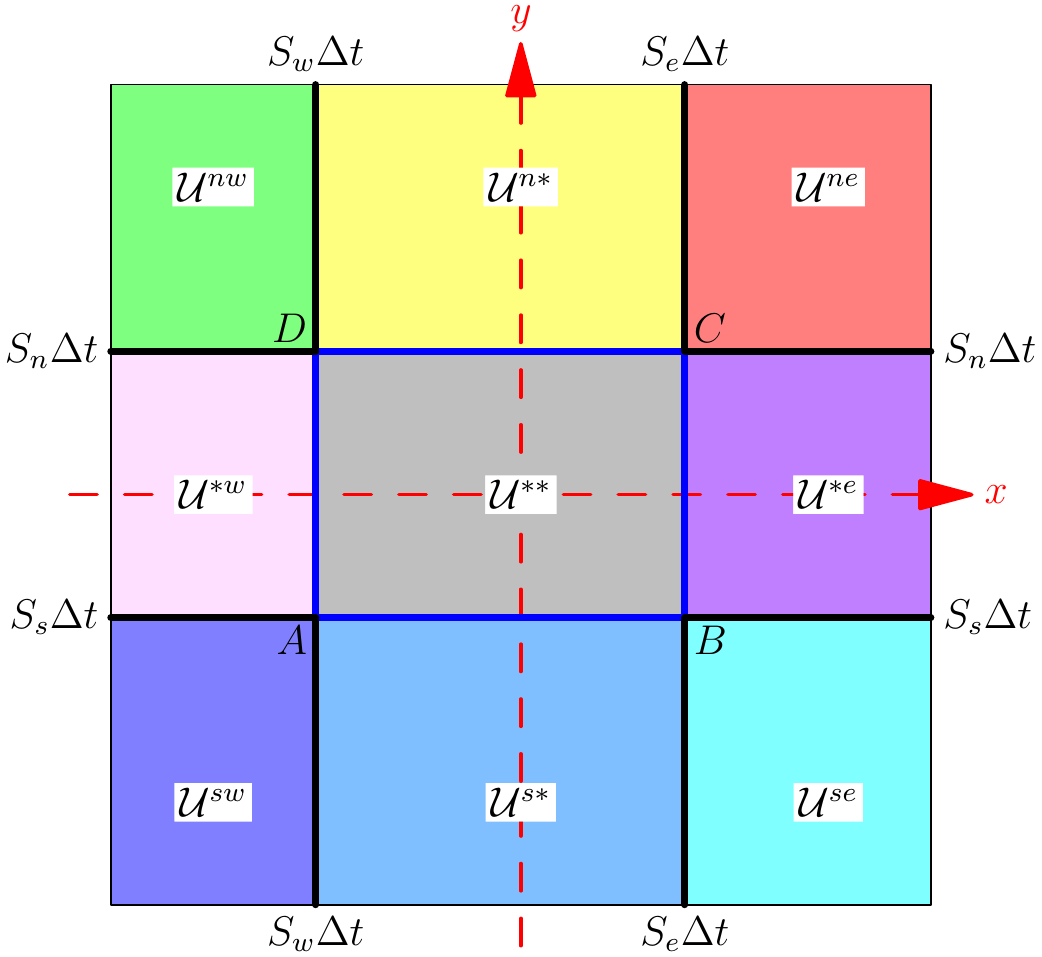} \\
(a) & (b)
\end{tabular}
\caption{Wave model for 2-D Riemann problem: (a) initial condition, (b) waves and solution at time $\Delta t$. The space-time view of the Riemann fan corresponding to Figure~\ref{fig:riem2d}b is shown in Figure~\ref{fig:riem2dc}.}
\label{fig:riem2d}
\end{center}
\end{figure}

\begin{figure}
\begin{center}
\includegraphics[width=0.45\textwidth]{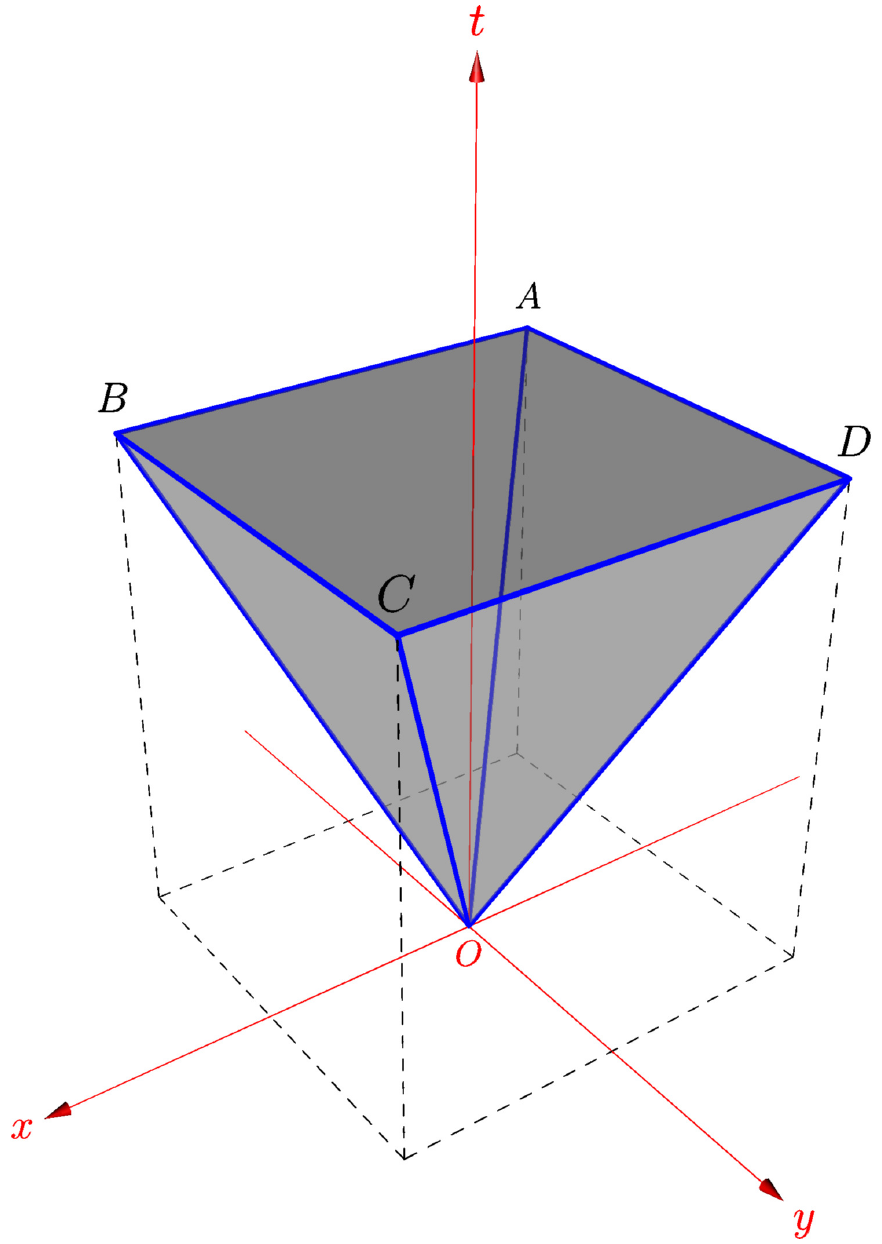}
\caption{Riemann fan in space-time for the 2-D Riemann problem. Figure~\ref{fig:riem2d}b shows the view of the Riemann fan when looking down on the $t$-axis.}
\label{fig:riem2dc}
\end{center}
\end{figure}
Around each vertex, there are four states that define a 2-D Riemann problem as shown in figure~(\ref{fig:riem2d}a). For example at South-West position, the hydrodynamic variables $\conh^{sw}$ are evaluated from the South-West cell solution which is a two dimensional polynomial, the $x$ component of $\B$ is obtained from the 1-D polynomial $b_x^s$ on the South face and the $y$-component of $\B$ is obtained from the 1-D polynomial $b_y^w$ on the West face. The remaining three states are determined in a similar manner. The solution of this 2-D Riemann problem yields the electric field $\tEz$ which is required to update the magnetic field variables stored on the faces in terms of the polynomials $b_x, b_y$. In the following sections, we consider the 1-D Riemann problem in the $x$-direction with initial data
\[
\con(x,0) = \begin{cases}
\con^L, & x < 0 \\
\con^R, & x > 0 \end{cases}
\]
and explain how the $x$-component of the flux is computed from approximate Riemann solvers.
\subsection{Local Lax-Friedrich flux}
The local Lax-Friedrich flux or the Rusanov flux~\cite{Rusanov1962}  is very simple and robust. For the $x$ direction, the numerical flux is given by
\begin{equation}
\nflx = \half[ \flx(\con^L)  + \flx(\con^R)] - \half \alpha_x^{LR} (\con^R - \con^L), \qquad \alpha_x^{LR} = \max\{ \alpha_x(\con^L), \alpha_x(\con^R)\}
\label{eq:lfx}
\end{equation}
where $\alpha_x(\con)$ is the maximum eigenvalue of the $x$ directional flux Jacobian. For the MHD system, the maximum wave speeds along the two directions are given by
\[
\alpha_x = |v_x| + c_{fx}, \qquad \alpha_y = |v_y| + c_{fy}
\]
The electric field that is obtained from this flux is
\begin{equation}
\hEz(\con^L,\con^R) = -(\nflx)_7 = \half( E_z^L + E_z^R) + \half \alpha_x^{LR} (B_y^R - B_y^L)
\label{eq:ezx}
\end{equation}
Finally, we need to specify the electric field $\tEz$ at the vertices of the cells. At any vertex, we have four states that come together giving rise to a 2-D Riemann problem as shown in Figure~\ref{fig:riem2d}a. The electric field at the vertex is estimated as~\cite{Balsara2014a},~\cite{Fu2018}
\begin{equation}
\begin{aligned}
\tEz = \frac{1}{4}(E_z^{sw} + E_z^{se} + E_z^{nw} + E_z^{ne}) - & \half \tilde{\alpha}_y \left( \frac{B_x^{nw} + B_x^{ne}}{2} - \frac{B_x^{sw} + B_x^{se}}{2} \right) \\
+ & \half \tilde{\alpha}_x \left( \frac{B_y^{ne} + B_y^{se}}{2} - \frac{B_y^{nw} + B_y^{sw}}{2} \right)
\end{aligned}
\label{eq:ezv}
\end{equation}
where
\[
\tilde{\alpha}_d = \max\{ \alpha_d(\con^{sw}), \alpha_d(\con^{se}), \alpha_d(\con^{nw}), \alpha_d(\con^{ne}) \}, \qquad d = x, y
\]
Note that since the normal component of $\B$ is continuous across the cell faces, we actually have $B_x^{nw} = B_x^{ne}$, $B_x^{sw} = B_x^{se}$, $B_y^{sw} = B_y^{nw}$ and $B_y^{se} = B_y^{ne}$.
\paragraph{Consistency with 1-D solver}
Now consider a situation where the four states actually form a 1-D Riemann problem, e.g., $\con^{sw} = \con^{nw} = \con^L$ and $\con^{se} = \con^{ne} = \con^R$. Then we have $\tilde{\alpha}_x = \alpha_x^{LR}$, and the electric field at the vertex given by equation~(\ref{eq:ezv}) reduces to
\[
\tEz = \half( E_z^L + E_z^R) + \half \alpha_x^{LR} (B_y^R - B_y^L) = \hEz(\con^L,\con^R)
\]
which coincides with the electric field given by the 1-D Riemann solver in equation~(\ref{eq:ezx}). Hence the estimate~(\ref{eq:ezv}) of the electric field at the vertices has the important continuity property that it reduces to the estimate obtained from the 1-D Riemann solver if the 2-D Riemann data corresponds to a 1-D Riemann data. Moreover, this consistency property is essential to maintain one dimensional solution structures aligned with the grid when solving the problem using the two dimensional scheme.
\subsection{HLL Riemann solver in 1-D}
In the HLL solver~\cite{Harten1983a}, we consider only the slowest and fastest waves in the solution of the Riemann problem. Let us denote these speeds by $S_L = \lminx(\con^L,\con^R)$ and $S_R = \lmaxx(\con^L,\con^R)$ with $S_L < S_R$, and there is an intermediate state $\con^*$ between these two waves. The intermediate state is obtained by satisfying the conservation law over the Riemann fan leading to
\[
\con^* = \frac{S_R \con^R - S_L \con^L - (\fl^R_x - \fl^L_x)}{S_R - S_L}
\]
and the flux is obtained by satisfying the conservation law over one half of the Riemann fan
\[
\flx^* = \flx^L + S_L (\con^* - \con^L) = \flx^R + S_R (\con^* - \con^R) = \frac{S_R \flx^L - S_L \flx^R + S_L S_R (\con^R - \con^L)}{S_R - S_L}
\]
The intermediate state and the above flux are required in the transonic case where $S_L < 0 < S_R$. The numerical flux in the general case is given by
\[
\nflx = \begin{cases}
\flx^L & S_L > 0 \\
\flx^R & S_R < 0 \\
\flx^* & \textrm{otherwise}
\end{cases}
\]
The electric field is obtained from the seventh component of the numerical flux and is given by
\[
\hEz(\con^L,\con^R) = - (\nflx)_7 = \begin{cases}
E_z^L & S_L > 0 \\
E_z^R & S_R < 0 \\
\frac{S_R E_z^L - S_L E_z^R - S_L S_R (B_y^R - B_y^L)}{S_R - S_L} & \textrm{otherwise}
\end{cases}
\]
Moreover, since the Riemann data satisfies $B_x^L = B_x^R$, the HLL solver automatically gives $B_x^* = B_x^L = B_x^R$ and $(\nflx)_6=0$. The wave speed estimates are taken as
\[
\lminx(\con^L,\con^R) = \min\{ v_x^L - c_{fx}^L, \bar{v}_x - \bar{c}_{fx} \}, \quad \lmaxx(\con^L,\con^R) = \max\{ v_x^R + c_{fx}^R, \bar{v}_x + \bar{c}_{fx} \}
\]
where the quantities with an overbar are based on Roe average state if it is physically admissible or based on average of primitive variables, otherwise.
\subsection{HLL solver in 2-D}
Consider a 2-D Riemann problem with data
\[
\con = \begin{cases}
\con^{sw} & x < 0, y < 0 \\
\con^{nw} & x < 0, y > 0 \\
\con^{se} & x > 0, y < 0 \\
\con^{ne} & x > 0, y > 0
\end{cases}
\]
which is illustrated in Figure~\ref{fig:riem2d}a. The four states give rise to four 1-D Riemann solutions and a strongly interacting state in the middle as shown in Figure~\ref{fig:riem2d}b. Similar to Balsara~\cite{Balsara2014a}, we have assumed that the waves are bounded by four wave speeds which are defined as
\[
S_w = \min\{ \lminx(\con^{sw},\con^{se}), \lminx(\con^{nw}, \con^{ne}) \}, \quad S_e = \max\{ \lmaxx(\con^{sw},\con^{se}), \lmaxx(\con^{nw}, \con^{ne}) \}
\]
\[
S_s = \min\{ \lminy(\con^{sw}, \con^{nw}), \lminy(\con^{se}, \con^{ne}) \}, \quad S_n = \max\{ \lmaxy(\con^{sw}, \con^{nw}), \lmaxy(\con^{se}, \con^{ne}) \}
\]
and $\con^{*w}, \con^{*e}, \con^{n*}, \con^{s*}$ are the intermediate states obtained from the 1-D HLL solution. The strongly interacting state in the middle is given by
\[
\begin{aligned}
\con^{**} = \frac{1}{2(S_e - S_w)(S_n - S_s)} \bigg[ & 2 S_e S_n \con^{ne} - 2 S_n S_w \con^{nw} + 2 S_s S_w \con^{sw}  - 2 S_s S_e \con^{se} \\
& - S_n(\flx^{ne} - \flx^{nw}) + S_s (\flx^{se} - \flx^{sw}) - (S_n - S_s)(\flx^{*e} - \flx^{*w}) \\
& - S_e (\fly^{ne} - \fly^{se}) + S_w (\fly^{nw} - \fly^{sw}) - (S_e - S_w)(\fly^{n*} - \fly^{s*}) \bigg]
\end{aligned}
\]
where the transverse fluxes $\flx^{*w}$, $\flx^{*e}$, $\fly^{s*}$, $\fly^{n*}$ are yet to be specified. In particular, the magnetic field components are given by
\[
\begin{aligned}
B_x^{**} =  \frac{1}{2(S_e - S_w)(S_n - S_s)} \bigg[ & 2 S_e S_n B_x^{ne} - 2 S_n S_w B_x^{nw} + 2 S_s S_w B_x^{sw}  - 2 S_s S_e B_x^{se} \\
& - S_e (E_z^{ne} - E_z^{se}) + S_w (E_z^{nw} - E_z^{sw}) - (S_e - S_w)(E_z^{n*} - E_z^{s*}) \bigg]
\end{aligned}
\]
\[
\begin{aligned}
B_y^{**} = \frac{1}{2(S_e - S_w)(S_n - S_s)} \bigg[ & 2 S_e S_n B_y^{ne} - 2 S_n S_w B_y^{nw} + 2 S_s S_w B_y^{sw}  - 2 S_s S_e B_y^{se} \\
& + S_n(E_z^{ne} - E_z^{nw}) - S_s (E_z^{se} - E_z^{sw}) + (S_n - S_s)(E_z^{*e} - E_z^{*w}) \bigg]
\end{aligned}
\]
To determine the fluxes $\flx^{**}, \fly^{**}$ at the vertex, we follow~\cite{Vides2015} and write down the jump conditions across the slanted sides of the space-time Riemann fan, see Figure~\ref{fig:riem2dc}, which are given by
\begin{eqnarray*}
\fly^{**} &=& \fly^{n*} - S_n (\con^{n*} - \con^{**}) \\
\fly^{**} &=& \fly^{s*} - S_s (\con^{s*} - \con^{**}) \\
\flx^{**} &=& \flx^{*e} - S_e (\con^{*e} - \con^{**}) \\
\flx^{**} &=& \flx^{*w} - S_w (\con^{*w} - \con^{**})
\end{eqnarray*}
This is an over-determined set of equations which can be solved by a least-squares method following the ideas in~\cite{Vides2015}. Our interest is only in the estimation of the electric field $E_z^{**}$ and we will not concern ourselves in computing all components of $\flx^{**}, \fly^{**}$. The electric field occurs in both the flux components; the sixth component of the first two equations and the seventh component of the last two equations contain the electric field and these equations are given by
\begin{equation}
\begin{aligned}
E_z^{**} = \ & E_z^{n*} - S_n (B_x^{n*} - B_x^{**}) \\
E_z^{**} = \ & E_z^{s*} - S_s (B_x^{s*} - B_x^{**}) \\
E_z^{**} = \ & E_z^{*e} + S_e (B_y^{*e} - B_y^{**}) \\
E_z^{**} = \ & E_z^{*w} + S_w (B_y^{*w} - B_y^{**})
\end{aligned}
\label{eq:e2djump}
\end{equation}
This is still an over-determined set of equations since there is only one unknown $E_z^{**}$ but four equations. The least-squares solution of this set of equations is just the average of the above four equations
\begin{equation}
\begin{aligned}
E_z^{**} = \frac{1}{4}(E_z^{n*} + E_z^{s*} + E_z^{*e} + E_z^{*w}) - & \frac{1}{4} S_n (B_x^{n*} - B_x^{**}) - \frac{1}{4} S_s (B_x^{s*} - B_x^{**}) \\
+ & \frac{1}{4} S_e (B_y^{*e} - B_y^{**}) + \frac{1}{4} S_w (B_y^{*w} - B_y^{**})
\end{aligned}
\label{eq:Ezss}
\end{equation}
Hence the electric field is given by
\[
\tEz = \begin{cases}
E_z^{*w} & S_w > 0 \\
E_z^{*e} & S_e < 0 \\
E_z^{s*} & S_s > 0 \\
E_z^{n*} & S_n < 0 \\
E_z^{**} & \textrm{otherwise}
\end{cases}
\]
Note that $E_z^{*w}$, etc. are the electric fields obtained from the 1-D HLL solver. The above formula is implemented in computer code using \texttt{if-else-if} statements. Also, since the normal components of magnetic field are continuous we actually have $B_x^{n*} = B_x^{ne} = B_x^{nw}$, $B_x^{s*} = B_x^{se} = B_x^{sw}$, $B_y^{*e} = B_y^{ne} = B_y^{se}$ and $B_y^{*w} = B_y^{nw} = B_y^{sw}$.
\paragraph{Remark}
We do not have to specify all the components of the transverse fluxes $\flx^{*w}$, $\flx^{*e}$, $\fly^{s*}$, $\fly^{n*}$ since our scheme requires only knowledge of the electric field from the 2-D Riemann problem. We only require the sixth component from $\fly^{s*}$, $\fly^{n*}$ and the seventh component from $\flx^{*w}$, $\flx^{*e}$, which corresponds to the electric field. These electric fields are obtained from the 1-D HLL fluxes.
\paragraph{Consistency with 1-D solver}
Suppose that the 2-D Riemann data has jumps only along the $x$ direction, so that $\con^{nw} = \con^{sw} = \con^L$ and $\con^{ne} = \con^{se} = \con^R$. Adopting the notation of the 1-D Riemann solver, we set $S_w = S_L$ and $S_e = S_R$. We will consider the transonic case, since the fully supersonic case is trivial. There is a common value of $B_x$ in all the four states so that $B_x^{n*} = B_x^{s*} = B_x^{**}$ and
\[
B_y^{*w} = B_y^L, \quad B_y^{*e} = B_y^R, \quad B_y^{**} = \frac{S_R B_y^R - S_L B_y^L + E_z^R - E_z^L}{S_R - S_L}
\]
\[
E_z^{n*} = E_z^{s*} = \hEz := \hEz(\con^L,\con^R), \qquad E_z^{*w} = E_z^L, \qquad E_z^{*e} = E_z^R
\]
Now the electric field $\tEz$ from the 2-D Riemann solver is given by
\begin{eqnarray*}
\tEz &=& \frac{1}{4}(E_z^{n*} + E_z^{s*} + E_z^{*e} + E_z^{*w}) - \frac{1}{4} S_n \cancel{(B_x^{n*} - B_x^{**})} - \frac{1}{4} S_s \cancel{(B_x^{s*} - B_x^{**})} \\
&& + \frac{1}{4} S_e (B_y^{*e} - B_y^{**}) + \frac{1}{4} S_w (B_y^{*w} - B_y^{**}) \\
&=& \frac{1}{4}(\hEz + \hEz + E_z^R + E_z^L) +  \frac{1}{4} S_R (B_y^{R} - B_y^{**}) + \frac{1}{4} S_L (B_y^{L} - B_y^{**}) \\
&=& \half \hEz + \frac{1}{4}(E_z^R + E_z^L) +  \frac{1}{4} S_R (B_y^{R} - B_y^{**}) + \frac{1}{4} S_L (B_y^{L} - B_y^{**}) \\
&=& \half \hEz + \half \frac{S_R E_z^L - S_L E_z^R - S_L S_R (B_y^R - B_y^L)}{S_R - S_L} \\
&=& \hEz(\con^L,\con^R)
\end{eqnarray*}
which follows after a little bit of algebraic manipulations. This shows that as far as the electric field is concerned, the 2-D Riemann solver reduces to the 1-D Riemann solver in case the 2-D Riemann data has jumps only along one direction.  The estimate \eqref{eq:Ezss} consists of the average of the four estimates from the four 1-D Riemann problems and some additional jump terms.  We see that the extra jump terms which arise from the strongly interacting state are necessary to achieve consistency with the 1-D solver. These extra terms are not present in other Riemann solvers, e.g., see equation (7)-(9) of~\cite{Balsara1999}. However, those former works are related to finite volume schemes where such consistency property is not strictly necessary. 
\subsection{HLLC Riemann solver in 1-D}
The HLLC Riemann solver~\cite{Toro1994} includes a middle contact wave in addition to the slowest and fastest waves and thus contains two intermediate states $\con^{*L}, \con^{*R}$. In the case of MHD, there are several variants of the solver~\cite{Gurski2004}, \cite{Li2005a}. To simplify the notation, we denote the velocity components as $(u,v,w)$ in this section. Following Batten~\cite{Batten1997}, we take the speed of the middle wave from the HLL intermediate state
\[
S_M = \frac{(\rho u)^*}{\rho^*} = \frac{(S_R - u_R) \rho_R u_R - (S_L - u_L) \rho_L u_L - (P_R - P_L)}{(S_R - u_R)\rho_R - (S_L - u_L) \rho_L}
\]
The intermediate states have the form
\[
\con^{*L} = \begin{bmatrix}
\rho^{*}_L \\
\rho^{*}_L S_M \\
\rho^{*}_L v^{*}_L \\
\rho^{*}_L w^{*}_L \\
\tote^{*}_L \\
B_x^* \\
B_y^* \\
B_z^*
\end{bmatrix}, \qquad
\con^{*R} = \begin{bmatrix}
\rho^{*}_R \\
\rho^{*}_R S_M \\
\rho^{*}_R v^{*}_R \\
\rho^{*}_R w^{*}_R \\
\tote^{*}_R \\
B_x^* \\
B_y^* \\
B_z^*
\end{bmatrix}
\]
with the common intermediate value of the magnetic field being equal to the HLL intermediate state. The intermediate density is obtained from the jump conditions across the left and right waves
\[
\rho^*_\alpha = \rho_\alpha \frac{S_\alpha - u_\alpha}{S_\alpha - S_M}, \qquad \alpha = L, R
\]
Li~\cite{Li2005a} proposes to take the intermediate velocities from the jump conditions across the left and right waves,
\[
v^*_\alpha = v_\alpha + \frac{B_x^\alpha B_y^\alpha - B_x^* B_y^*}{\rho_\alpha (S_\alpha - u_\alpha)}, \qquad w^*_\alpha = w_\alpha + \frac{B_x^\alpha B_z^\alpha - B_x^* B_z^*}{\rho_\alpha (S_\alpha - u_\alpha)}, \qquad \alpha = L, R
\]
and defines the energy as
\[
\tote^{*}_\alpha = \frac{(S_\alpha - u_\alpha) \tote_\alpha - P_\alpha u_\alpha + P^* S_M + B_x^\alpha (\vel_\alpha \cdot \B_\alpha) - B_x^*(\vel^* \cdot \B^*)}{S_\alpha - S_M}, \qquad \alpha = L,R
\]
where $P^*$ is the common intermediate total pressure given by
\[
P^* = P_L + \rho_L(S_L - u_L)(S_M - u_L) = P_R + \rho_R (S_R - u_R) (S_M - u_R)
\]
and $\vel^*, \B^*$ are the values from the HLL intermediate state. Once the two intermediate states are determined, the numerical flux is given by
\[
\nflx = \begin{cases}
\flx^L & S_L > 0 \\
\flx^R & S_R < 0 \\
\flx^L + S_L(\con^{*L} - \con^L) & S_L \le 0 \le S_M \\
\flx^R + S_R(\con^{*R} - \con^R) & S_M \le 0 \le S_R
\end{cases}
\]
The electric field is obtained from the seventh component of the flux
\[
\hEz(\con^L,\con^R) = - (\nflx)_7 = \begin{cases}
E_z^L & S_L > 0 \\
E_z^R & S_R < 0 \\
E_z^L - S_L(B_y^* - B_y^L) & S_L \le 0 \le S_M \\
E_z^R - S_R(B_y^* - B_y^R) & S_M \le 0 \le S_R
\end{cases}
\]
But due to the definition of $B_y^*$ from the HLL intermediate state, the two intermediate values of the electric field are identical and equal to the HLL estimate of the electric field.
\subsection{HLLC Riemann solver in 2-D}
We will assume the same type of wave modeling in 2-D as we adopted in case of HLL solver, except that there is an intermediate wave in all the four 1-D Riemann problems. We do not endow any more sub-structure in the strongly interacting state. The jump conditions across the 2-D Riemann fan given in~(\ref{eq:e2djump}) are still valid since we have a common electric field and magnetic field in the intermediate states of the 1-D Riemann fan. The definition of the electric field at the vertex thus matches with that obtained from the HLL solver. Consequently, the consistency with the 1-D Riemann solvers also follows in same way as it was proved for the HLL solver.
\section{Limiting procedure}
\label{sec:lim}
The basic unknowns in our scheme are the normal component of $b_x$, $b_y$ on the faces, the additional cell moments $\alpha$, $\beta$ stored inside the cells, and the remaining variables $\conh$ inside the cell. As long as we don't apply any limiter on $b_x, b_y, \alpha, \beta$, our algorithm is guaranteed to preserve the initial divergence. However, for computing discontinuous solutions, some form of limiter is absolutely necessary to control spurious numerical oscillations. If any form of limiter is applied which a posteriori modifies these solution variables, then the divergence will not be preserved and some correction has to be applied to recover divergence-free property. We now detail the steps in our limiting strategy including a divergence-free reconstruction step. The choice of which variable set is limited is an important one for systems of conservation laws, and it is found in many studies that applying the limiter to characteristic variables as opposed to conserved variables, gives better control on oscillations and leads to more accurate solutions~\cite{Cockburn1989}. Hence, our limiting strategy will be based on limiting the set of characteristic variables.
\subsection{TVD-type limiter}

\subsection*{Step 1}
Using the basic solution variables for the magnetic field, we perform the RT reconstruction to obtain the polynomial $\B(\xi,\eta)$ representation in each cell. We now have the complete solution in the cell and we proceed to apply a TVD/TVB limiter to this solution polynomial. The basic idea in this process is to check if the linear part of the solution in a cell is smooth relative to the variation of the cell averages around the cell~\cite{Cockburn1990}. Consider the cell indexed by $(i,j)$. Observe that in the expansion~(\ref{eq:hydrosol}),~(\ref{eq:BxBy}), the components $\conh_{10}, a_{10},b_{10}$ give a measure of the $x$ derivative and $\conh_{01},a_{01},b_{01}$ give a measure of the $y$ derivatives. Form the vector of slopes in the two directions
\[
\con_{i,j}^x = \begin{bmatrix}
(\conh_{10})_{i,j} \\
(a_{10})_{i,j} \\
(b_{10})_{i,j} \end{bmatrix}, \qquad \con_{i,j}^y = \begin{bmatrix}
(\conh_{01})_{i,j} \\
(a_{01})_{i,j} \\
(b_{01})_{i,j} \end{bmatrix}
\]
and the differences of the cell averages
\[
\con_{i,j}^{x-} = \begin{bmatrix}
(\conh_{00})_{i,j} - (\conh_{00})_{i-1,j}\\
(a_{00})_{i,j} - (a_{00})_{i-1,j}\\
(b_{00})_{i,j} - (b_{00})_{i-1,j}\end{bmatrix}, \quad \con_{i,j}^{x+} = \begin{bmatrix}
(\conh_{00})_{i+1,j} - (\conh_{00})_{i,j}\\
(a_{00})_{i+1,j} - (a_{00})_{i,j}\\
(b_{00})_{i+1,j} - (b_{00})_{i,j}\end{bmatrix}
\]
\[
\con_{i,j}^{y-} = \begin{bmatrix}
(\conh_{00})_{i,j} - (\conh_{00})_{i,j-1}\\
(a_{00})_{i,j} - (a_{00})_{i,j-1}\\
(b_{00})_{i,j} - (b_{00})_{i,j-1}\end{bmatrix}, \quad \con_{i,j}^{y+} = \begin{bmatrix}
(\conh_{00})_{i,j+1} - (\conh_{00})_{i,j}\\
(a_{00})_{i,j+1} - (a_{00})_{i,j}\\
(b_{00})_{i,j+1} - (b_{00})_{i,j}\end{bmatrix}
\]
The values $a_{00},a_{10},a_{01},b_{00},b_{10},b_{01}$ are known from the polynomial $\B(\xi,\eta)$. Let $\Rx, \Lx, \Ry, \Ly$ be the matrix of right and left eigenvectors based on the cell average state. We now convert the above slopes to characteristic variables
\[
\cvar_{i,j}^x = \Lx \con_{i,j}^x, \quad \cvar_{i,j}^{x\pm} = \Lx \con_{i,j}^{x\pm}, \quad \cvar_{i,j}^y = \Ly \con_{i,j}^y, \quad \cvar_{i,j}^{y\pm} = \Ly \con_{i,j}^{y\pm}
\]
Now apply the minmod limiter on each of the characteristic variables
\[
(\cvar_{i,j}^x)^m = \minmod\left[ \cvar_{i,j}^x, \beta \cvar_{i,j}^{x-}, \beta \cvar_{i,j}^{x+}, M_{i,j}^x\Delta x^2 \right]
\]
\[
(\cvar_{i,j}^y)^m = \minmod\left[ \cvar_{i,j}^y, \beta \cvar_{i,j}^{y-}, \beta \cvar_{i,j}^{y+}, M_{i,j}^y\Delta y^2 \right]
\]
where the minmod function is defined as
\[
\minmod[a,b,c,\delta] = \begin{cases}
a & |a| < \delta \\
s \min(|a|,|b|,|c|) & s = \sign{a} = \sign{b} = \sign{c} \\
0 & \textrm{otherwise}
\end{cases}
\]
If $(\cvar_{i,j}^x)^m = (\cvar_{i,j}^x)$ and $(\cvar_{i,j}^y)^m = (\cvar_{i,j}^y)$ then we do not change the solution in this cell. Otherwise, convert them back to conserved variables by multiplying with the matrices $\Rx$, $\Ry$
\[
(\con_{i,j}^x)^m = \Rx (\cvar_{i,j}^x)^m, \qquad (\con_{i,j}^y)^m = \Ry (\cvar_{i,j}^y)^m
\]
and reset the first order components of $\conh$ to the limited values, i.e., $\conh_{10} = (\con_{i,j}^x)^m$ and $\conh_{01} = (\con_{i,j}^y)^m$, while setting all higher modes to zero. Similarly, we also reset the modes $a_{10}, a_{01}, b_{01}, b_{10}$ of the magnetic field and kill the higher modes of $\B(\xi,\eta)$ if the limiter is active in the current cell. The magnetic field at this stage will not be divergence-free and we will correct this in a later step.
\subsection*{Step 2}
We now loop over the all the faces in the mesh and apply a limiter to the solution polynomials $(b_x, b_y)$ which reside on the faces by making use of the cell solution $\B$ that has already been limited in the previous step. Consider a vertical face on which we have the polynomial $b_x(\eta) = \sum_{j=0}^k a_j \phi_j(\eta)$. We  also have the limited solution polynomials $\B(\xi,\eta)$ in the two adjacent cells of this face and we can evaluate them on the face. From the left cell we obtain\footnote{Here the indices $i,j$ denote the solution modes and not the cell indices.}
\[
B_x^L(\shalf,\eta) = \sum_{i=0}^{k+1}\sum_{j=0}^k a_{ij}^L \phi_i(\shalf) \phi_j(\eta) = \sum_{j=0}^k a_{j}^L \phi_j(\eta), \quad a_j^L = \sum_{i=0}^{k+1} a_{ij}^L \phi_i(\shalf)
\]
and similarly from the right cell we obtain
\[
B_x^R(-\shalf,\eta) = \sum_{i=0}^{k+1}\sum_{j=0}^k a_{ij}^R \phi_i(-\shalf) \phi_j(\eta) = \sum_{j=0}^k a_{j}^R \phi_j(\eta), \quad a_j^R = \sum_{i=0}^{k+1} a_{ij}^R \phi_i(-\shalf)
\]
We now compare the three solutions at the face via a minmod function to decide on the smoothness and modify the solution on the face as follows
\[
a_j \leftarrow \minmod\left( a_j, \beta a_j^L, \beta a_j^R \right), \qquad j=1,\ldots,k
\]
where $\beta \in [1,2]$ can be chosen to be greater than unity in order to allow larger slope similar to the MC limiter. Note that we do not modify the mean value on the face which corresponds to $a_0$ and this is important to perform divergence-free reconstruction in the next step. A similar procedure is applied to limit the face solution polynomials $b_y(\xi)$ located on the horizontal faces in the mesh.
\subsection*{Step 3}
The final step will restore the divergence-free condition on the magnetic field in those cells where the limiter has been active. This involves using the limited facial solution polynomials $b_x^\pm(\eta)$, $b_y^\pm(\xi)$ to reconstruct a divergence-free vector field $\B(\xi,\eta)$ inside the cell. This is achieved by determining the cell moments $\alpha_{ij}, \beta_{ij}$ in a divergence-free manner, after which the RT reconstruction can be performed. The procedure at different orders is explained in Appendix~(\ref{sec:dfr}). At second and third order, the reconstruction can be performed using the limited facial solution alone while at fourth order, we need an additional information from inside the cell, which is taken as the quantity $\omega = b_{01} - a_{10}$ and approximates the curl of $\B$. Note that the quantities $a_{10}$, $b_{01}$ are available to us after the cell solution $\B(\xi,\eta)$ has been limited in Step~1 and we already have a limited estimate of the quantity $\omega$.
\subsection{Positivity limiter}
The density and pressure have to remain positive since otherwise the problem is ill-posed and the computations would break down. High order positivity preserving schemes are usually built on the basis of a first order positive scheme. In a DG scheme, the solution is discontinuous and can be scaled in each cell to make it positive~\cite{Zhang2010b}, \cite{Zhang2010}. However, a constraint preserving DG scheme like the one proposed in this work and also the scheme in~\cite{Fu2018}, do not have a fully discontinuous solution since the normal component of $\B$ has to be continuous. Hence the local scaling limiter idea cannot be applied since scaling the $\B$ field in one cell changes the field in the neighbouring cells. The staggered storage of magnetic field variables also complicates the construction and analysis of positive schemes. Moreover, to the best of our knowledge, there is no first order, provably positive, divergence-free scheme based on Godunov-type approach  available in the literature. However, if we do not demand strictly divergence-free solutions, then provably positive schemes can be developed as in~\cite{Cheng2013}, which uses a standard DG approach, and as in~\cite{Wu2018} where locally divergence-free basis is used for $\B$ which is fully discontinuous across the cell faces and hence can be scaled in a local manner to achieve positivity.

Since there is no rigorous theory of positivity preservation in the framework of constraint preserving, high order DG schemes applied to the MHD system at present, we take a heuristic approach to ensure positivity property whose success can only be judged from numerical experiments. The approach we take here is to ensure that the solution is positive at all the quadrature points where the solution is used to compute the quadratures and fluxes involved in the DG scheme. The set of points $S$ includes the $(k+1)^2$ GL quadrature points inside the cell, the $4(k+1)$ GL quadrature points on the faces and the four corner points. The positivity is achieved by scaling the solution in each cell by following the ideas in~\cite{Zhang2010},~\cite{Cheng2013}. We apply this scaling to the hydrodynamic variables stored in $\conh$ and the RT polynomial $\B$. These polynomials are used to compute all the cell and face integrals in the DG scheme. However, {\em we do not scale the basic degrees of  freedom of the magnetic field which are $b_x, b_y, \alpha, \beta$, which allows us to maintain the divergence-free property of the magnetic field}. The scaling limiter relies on the fact that the average values on the cells and faces are already positive. In some difficult problems like a strong blast wave with low plasma beta (i.e., where thermodynamic pressure $p$ is much smaller than the magnetic pressure $\half |\B|^2$), see Section~(\ref{sec:blast}), the average pressure may also be negative in a few cells, in which case we have to set the pressure to a small positive value. While this is not an ideal solution to this problem, it does allow us to maintain stability of the computations as shown in the results section. The positivity limiter is applied after the TVB limiter and so the $\B$ field has already been limited and is not so badly behaved.
\section{Numerical results}
\label{sec:res}

We have explained the semi-discrete version of the DG scheme in the previous sections which leads to a system of coupled ODE. Starting from the specified initial condition, these ODE are integrated forward in time using Runge-Kutta schemes. For $k=1$ and $k=2$, we use the second and third order strong stability preserving RK schemes~\cite{Shu1988}, respectively, while for $k=3$ we use the 5-stage,  fourth order SSPRK scheme~\cite{Kraaijevanger1991}, \cite{Spiteri2002}. The time step is computed as
\[
\Delta t = \frac{\cfl}{\max\left(\frac{|v_x|+c_{fx}}{\dx} + \frac{|v_y|+c_{fy}}{\dy}\right)}
\]
where the wave speeds are based on cell average values and the maximum is taken over all the cells in the grid. We also use a shock indicator as explained in~\cite{Fu2017} in most of these test cases and the limiter is applied only in those cells which are marked by the indicator. Unless stated otherwise, in all the test cases we use $\cfl = 0.95/(2k+1)$, where $k$ is the degree of the approximating polynomials. The initial condition of the magnetic field variables must be set carefully to ensure that it is divergence-free~\cite{Chandrashekar2019}, and this is explained in Appendix~(\ref{sec:ic}). A high level view of the algorithm in given in Algorithm~(\ref{algo:mhd}).

\begin{algorithm}
\SetAlgoLined
\caption{Constraint preserving scheme for ideal compressible MHD}
\label{algo:mhd}
 Allocate memory for all variables\;
 Set initial condition for $\conh,b_x,b_y,\alpha,\beta$\;
 Loop over cells and reconstruct $B_x, B_y$\;
 Set time counter $t=0$\;
 \While{$t < T$}{
  Copy current solution into old solution\;
  Compute time step $\Delta t$\;
  \For{each RK stage}{
  Loop over vertices and compute vertex flux\;
  Loop over faces and compute all face integrals\;
  Loop over cells and compute all cell integrals\;
  Update solution to next stage\;
  Loop over cells and do RT reconstruction $(b_x,b_y,\alpha,\beta) \to \B$\;
  Loop over cells and apply limiter on $\conh,\B$\;
  Loop over faces and limit solution $b_x$, $b_y$\;
  Loop over cells and perform divergence-free reconstruction\;
  Apply positivity limiter\;
  }
  $t = t + \Delta t$\;
 }
\end{algorithm}
\subsection{Alfven wave}

\begin{figure}
\begin{center}
\begin{tabular}{ccc}
 \includegraphics[width=0.33\textwidth]{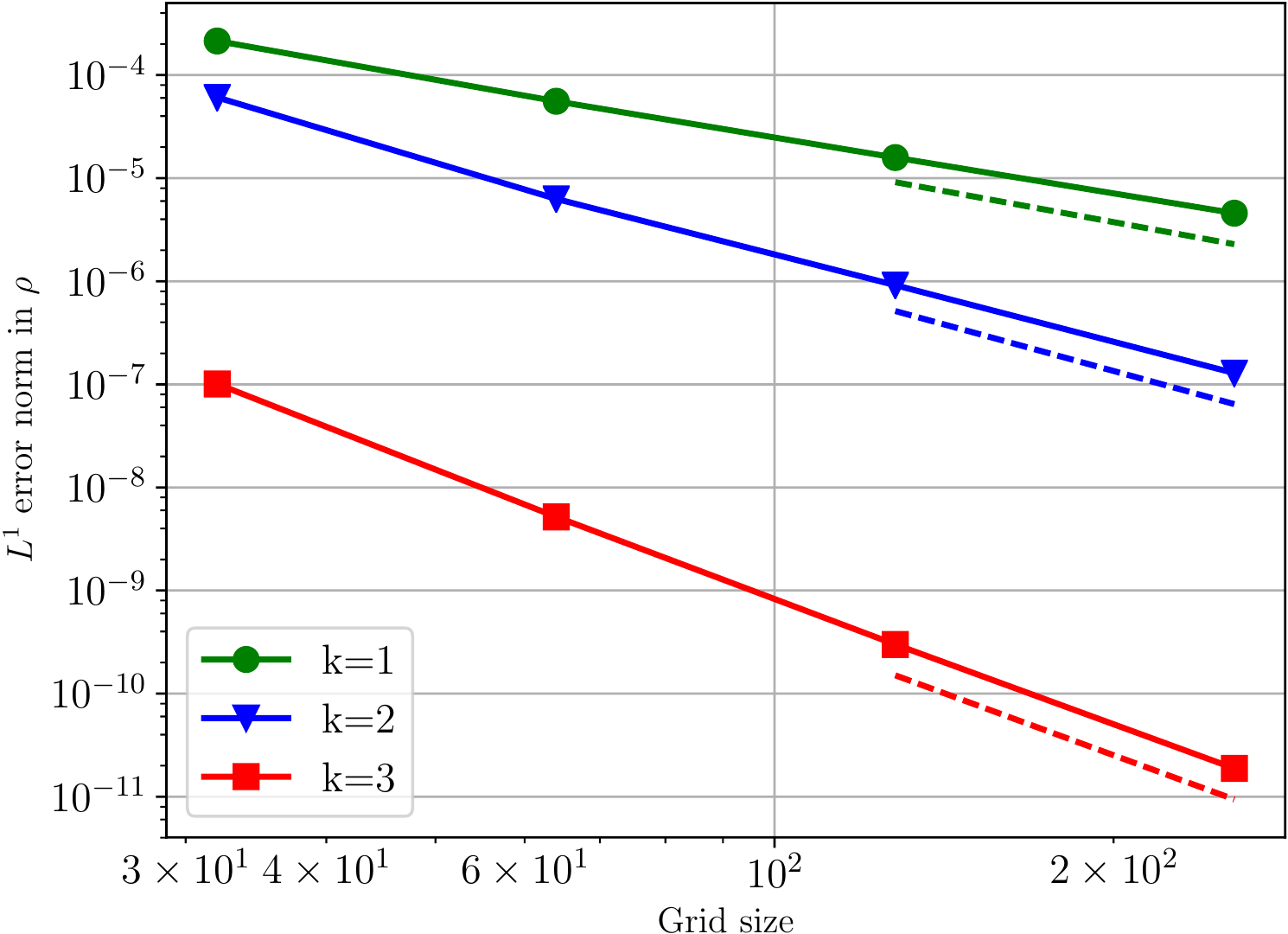} &
\includegraphics[width=0.33\textwidth]{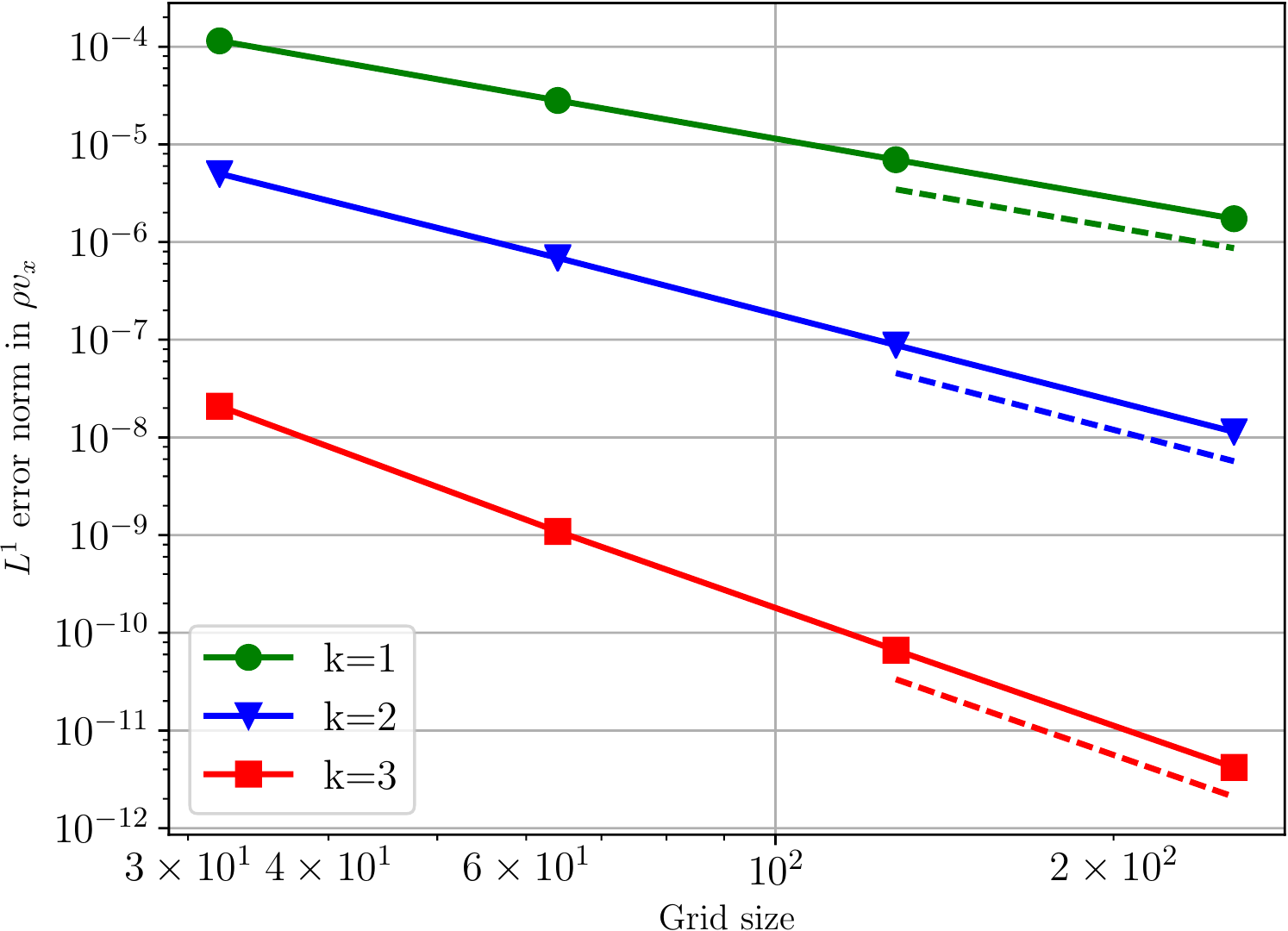} &
\includegraphics[width=0.33\textwidth]{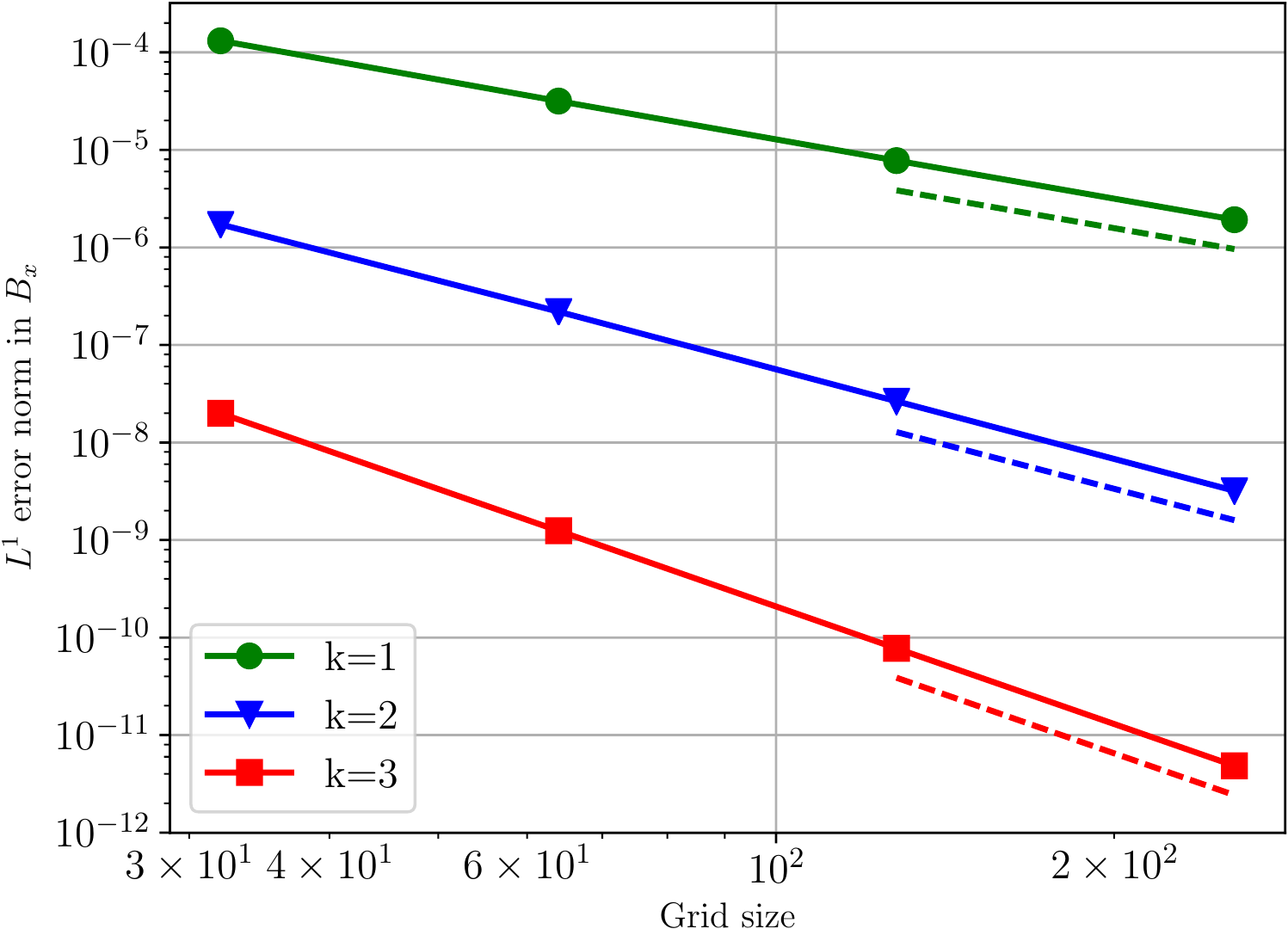}
\end{tabular}
\end{center}
\caption{Comparison of convergence results obtained from  HLLC flux for the smooth Alfven wave problem using degree $k=1,2,3$. The dashed line shows second, third and fourth order rates.}
\label{fig:alfven}
\end{figure}

This test case  involves the propagation of circularly polarized wave over the rectangular domain $[0, 1/\cos \alpha] \times [0,1/\sin \alpha]$ and it is used to test accuracy and convergence of numerical algorithms~\cite{Toth2000}. The parameter $\alpha$ is the angle of wave propagation relative to 
$x-$axis and it is taken to be $\pi/6$. The initial condition is given by
\begin{align*}
 \rho=1, \qquad \vel =v_{\perp}(-\sin \alpha, \cos \alpha, 0), \qquad p=0.1 \\
 B_x = B_{\parallel} \cos \alpha - B_{\perp} \sin \alpha, \qquad B_y = B_{\parallel} \sin \alpha + B_{\perp} \cos \alpha, \qquad B_z= v_z
\end{align*}
where 
\begin{equation*}
 B_{\parallel} =1, \qquad B_{\perp}= v_{\perp}= 0.1 \sin(2\pi(x\cos \alpha+y \sin \alpha))
\end{equation*}
 We have taken periodic boundary conditions in  both directions.
The numerical solution is computed up to time $T=1$ with $\gamma=5/3$.  In Figure~\ref{fig:alfven}, we have compared the error and convergence rates obtained using HLLC flux for $k=1,2,3$. In case of all proposed schemes the optimal rates of 2, 3 and 4 have been achieved for all the variables, i.e., with degree $k$ polynomials, the error is $O(h^{k+1})$. We have observed similar results with HLL and LxF fluxes (not shown here).

\subsection{Smooth magnetic vortex}
\begin{figure}
\begin{center}
\begin{tabular}{ccc}
\includegraphics[width=0.33\textwidth]{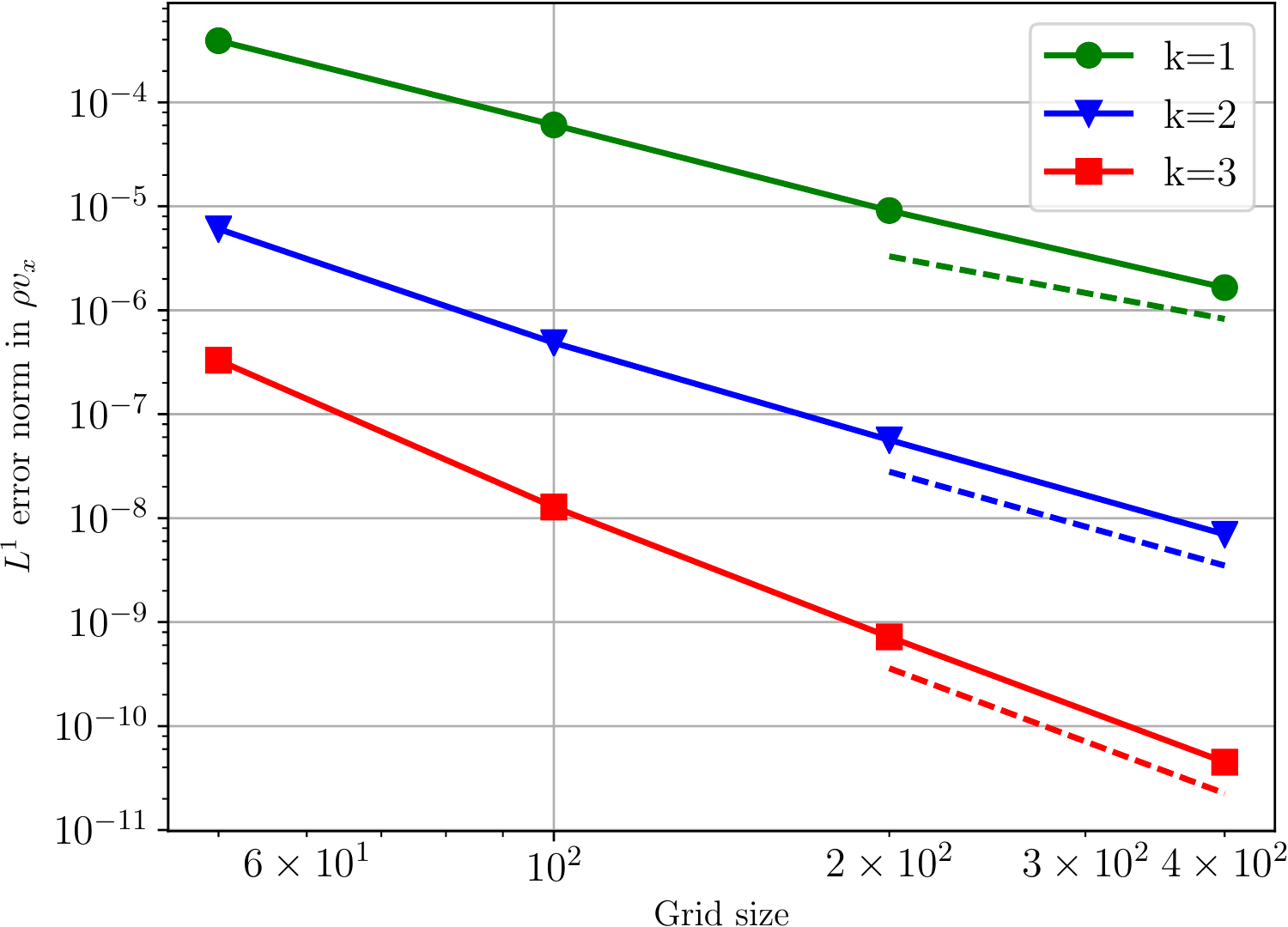} &
\includegraphics[width=0.33\textwidth]{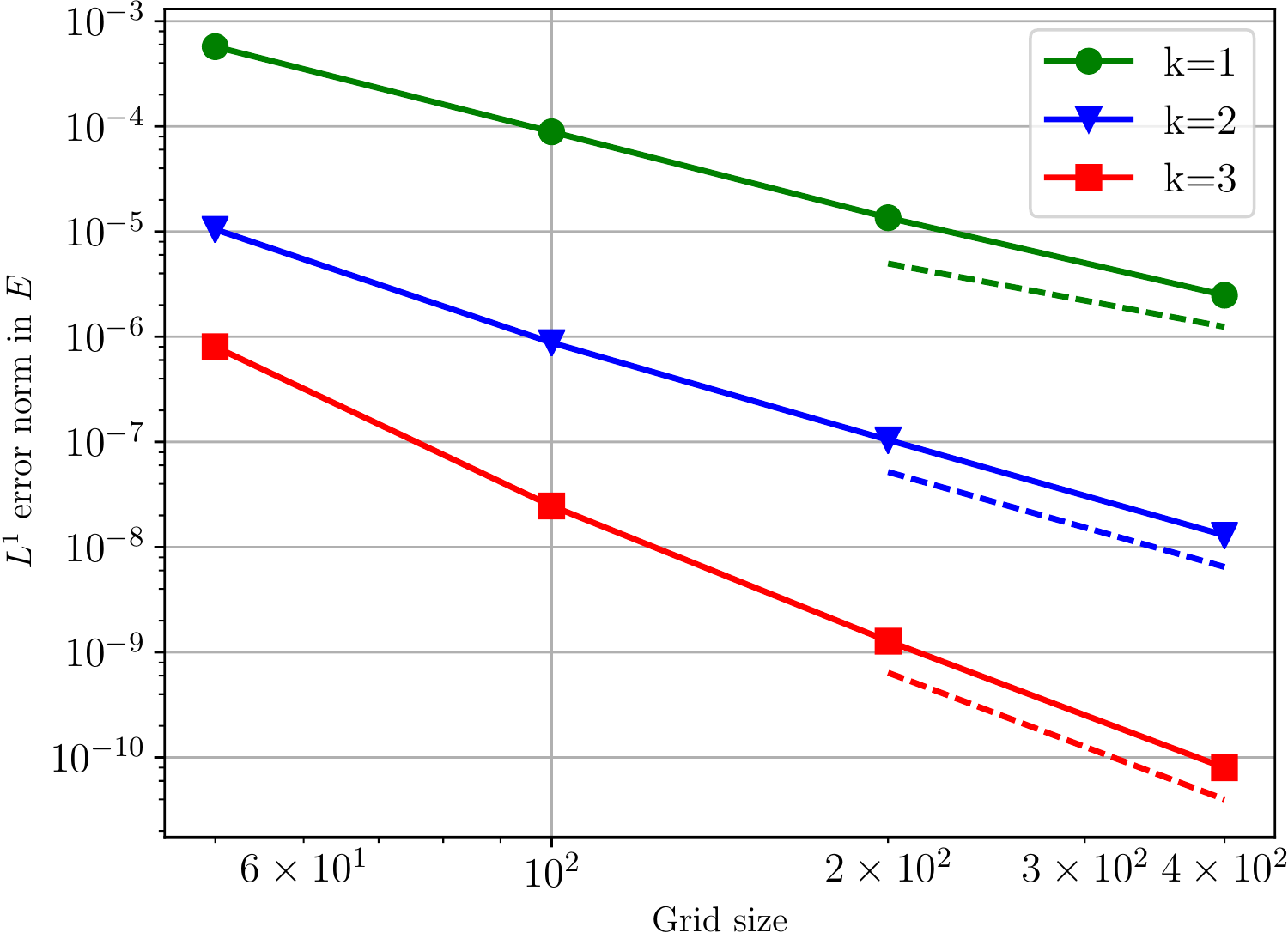} &
\includegraphics[width=0.33\textwidth]{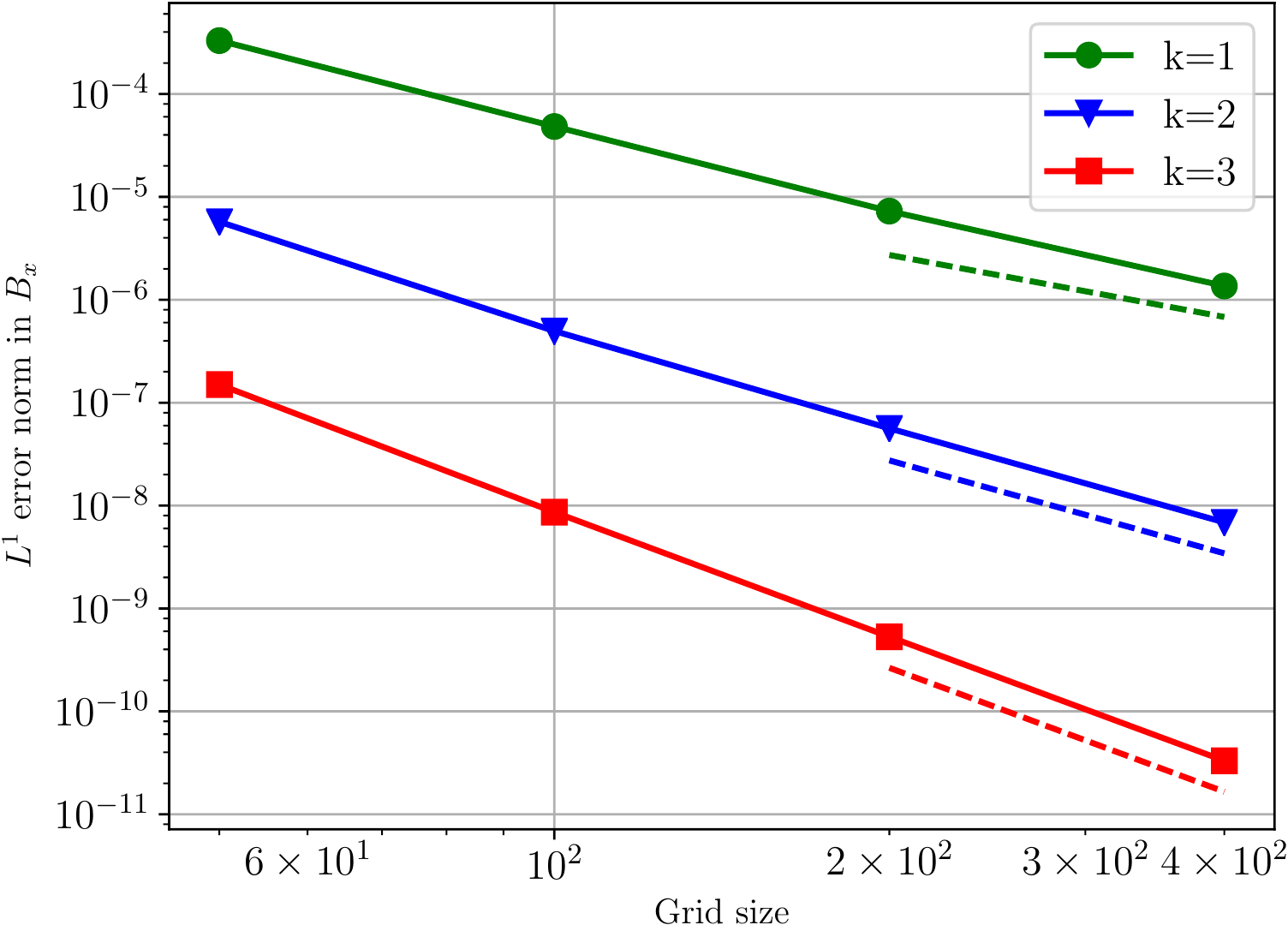}
\end{tabular}
\end{center}
\caption{Comparison of convergence results obtained from  HLLC flux for the smooth vortex  problem using degree $k=1,2,3$. The dashed line shows second, third and fourth order rates.}
\label{fig:vortex}
\end{figure}

This test case involves the propagation of a smooth, constant density  vortex along an oblique direction to the computational mesh and it is truly multidimensional in nature~\cite{Balsara2004}. The problem is initialized over the computational domain $[-10,10]\times [-10,10]$, with periodic boundary conditions in both directions.  The initial unperturbed primitive variables are given by
\begin{equation*}
\rho = 1, \qquad p = 1, \qquad \vel = (1,1,0), \qquad \bthree = (0,0,0)
\end{equation*}
and $\gamma=5/3$.
A vortex is initialized at the origin by adding the fluctuations in velocity and magnetic field which are given by
\begin{align*}
  \delta v_x &= -\frac{\kappa}{2\pi} y\exp(0.5(1-r^2)), \qquad \delta v_y =  \frac{\kappa}{2\pi} x\exp(0.5(1-r^2)), \qquad\delta v_z =0 \nonumber\\
  \delta B_x &=- \frac{\mu}{2\pi}y\exp(0.5(1-r^2)), \qquad \delta B_y = \frac{\mu}{2\pi}x\exp(0.5(1-r^2)), \qquad \delta B_z=0
\end{align*}
and the perturbation in pressure is given by
\begin{equation*}
  \delta p = \left[ \frac{1}{8\pi}\left(\frac{\mu}{2\pi}\right)^2 (1-r^2) - \frac{1}{2} \left(\frac{\kappa}{2\pi}\right)^2 \right] \exp(1-r^2)
\end{equation*}
We have set the parameters  $\kappa = 1$ and $\mu = 1$ in the initial condition. The smooth vortex returns to its initial position after some fixed time-period and facilitates us to measure the accuracy and convergence of the numerical algorithms. The numerical simulations are performed up to time $T=20$ using CFL=0.95. The convergence of the error for some of the variables with respect to grid refinement and for HLLC flux are shown in Figure~\ref{fig:vortex}, which indicates that the optimal rates of 2, 3 and 4 have been achieved for $k=1,2,3$, respectively. 

\begin{figure}
\begin{center}
\begin{tabular}{cc}
 \includegraphics[width=0.48\textwidth]{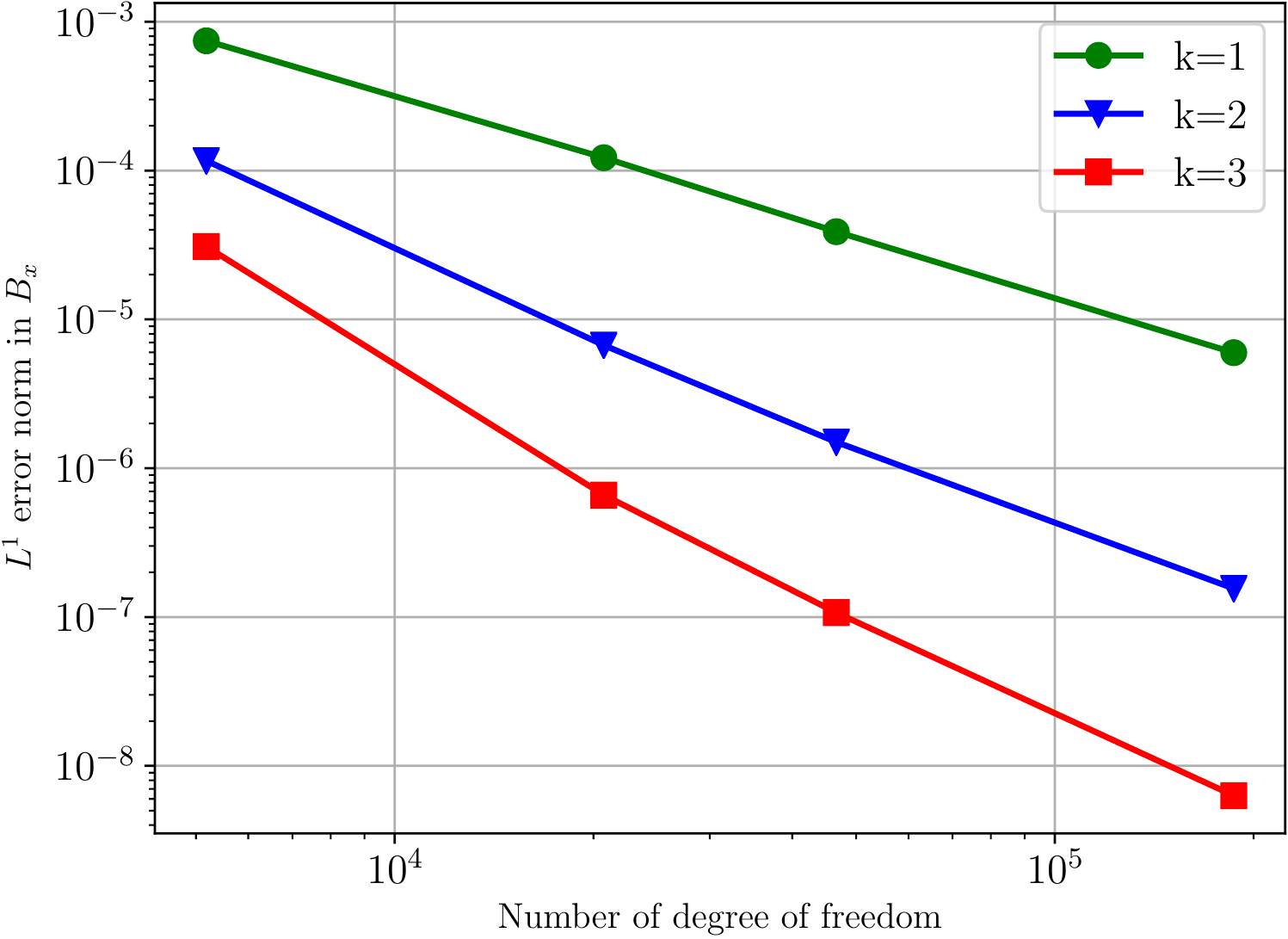} &
\includegraphics[width=0.48\textwidth]{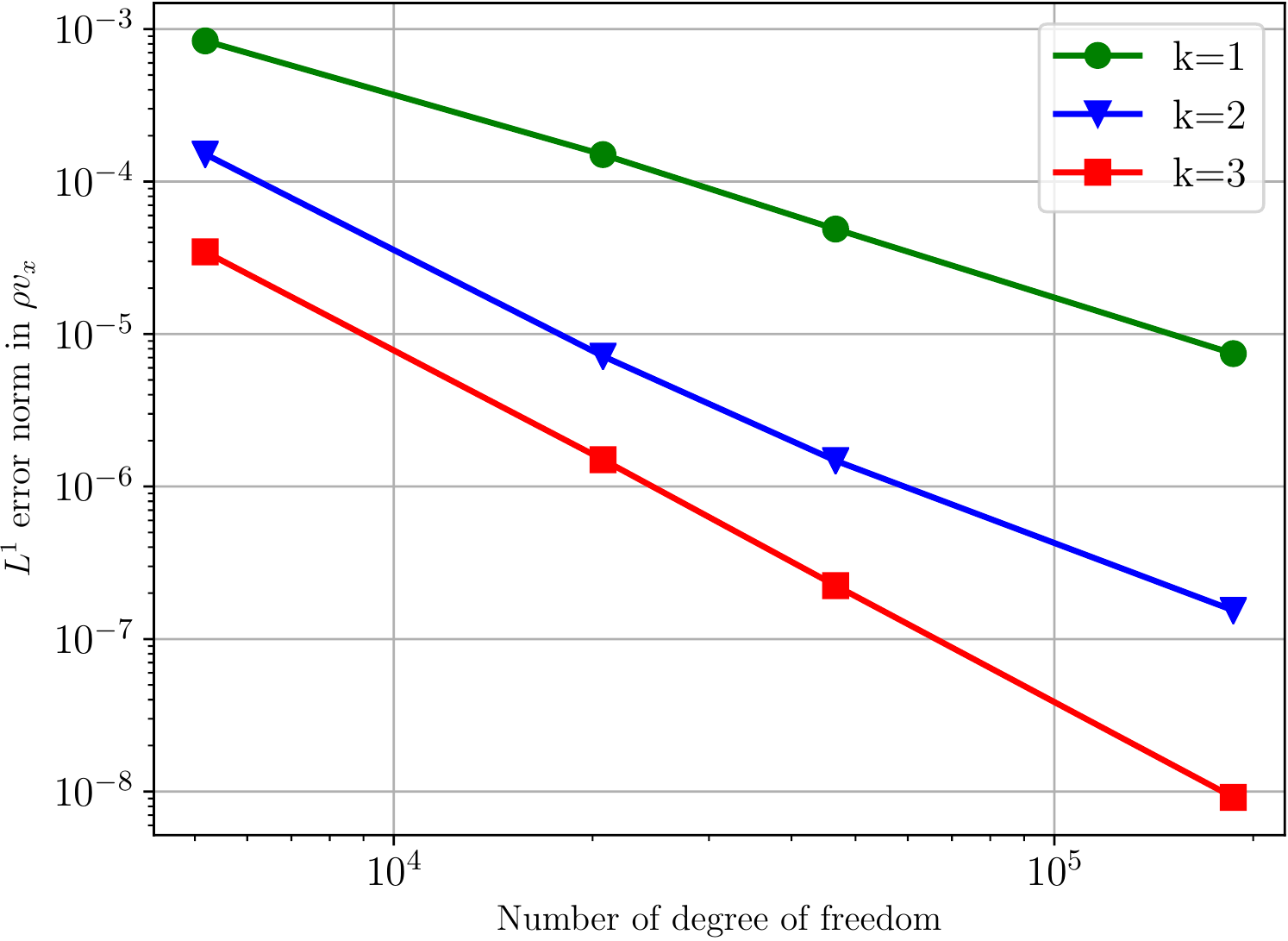} \\
(a)  & (b) 
\end{tabular}
\end{center}
\caption{Comparison of convergence results obtained from  HLLC flux for the smooth vortex  problem using degree for different degree as a function of number of degrees of freedom.}
\label{fig:vortex2}
\end{figure}

\resb
{
To study the usefulness of using high order methods, we compute the vortex problem on different meshes and polynomial degree so that the total number of degrees of freedom are roughly matched. The error norm as a function of the number of degrees of freedom are shown in Figure~\ref{fig:vortex2}.  We observe that to obtain same error level, a low order method (small degree $k$) requires more degrees of freedom than a high order method (large degree $k$).
}


\subsection{Brio-Wu shock tube}
\begin{figure}
\begin{center}
\begin{tabular}{cc}
\includegraphics[width=0.48\textwidth]{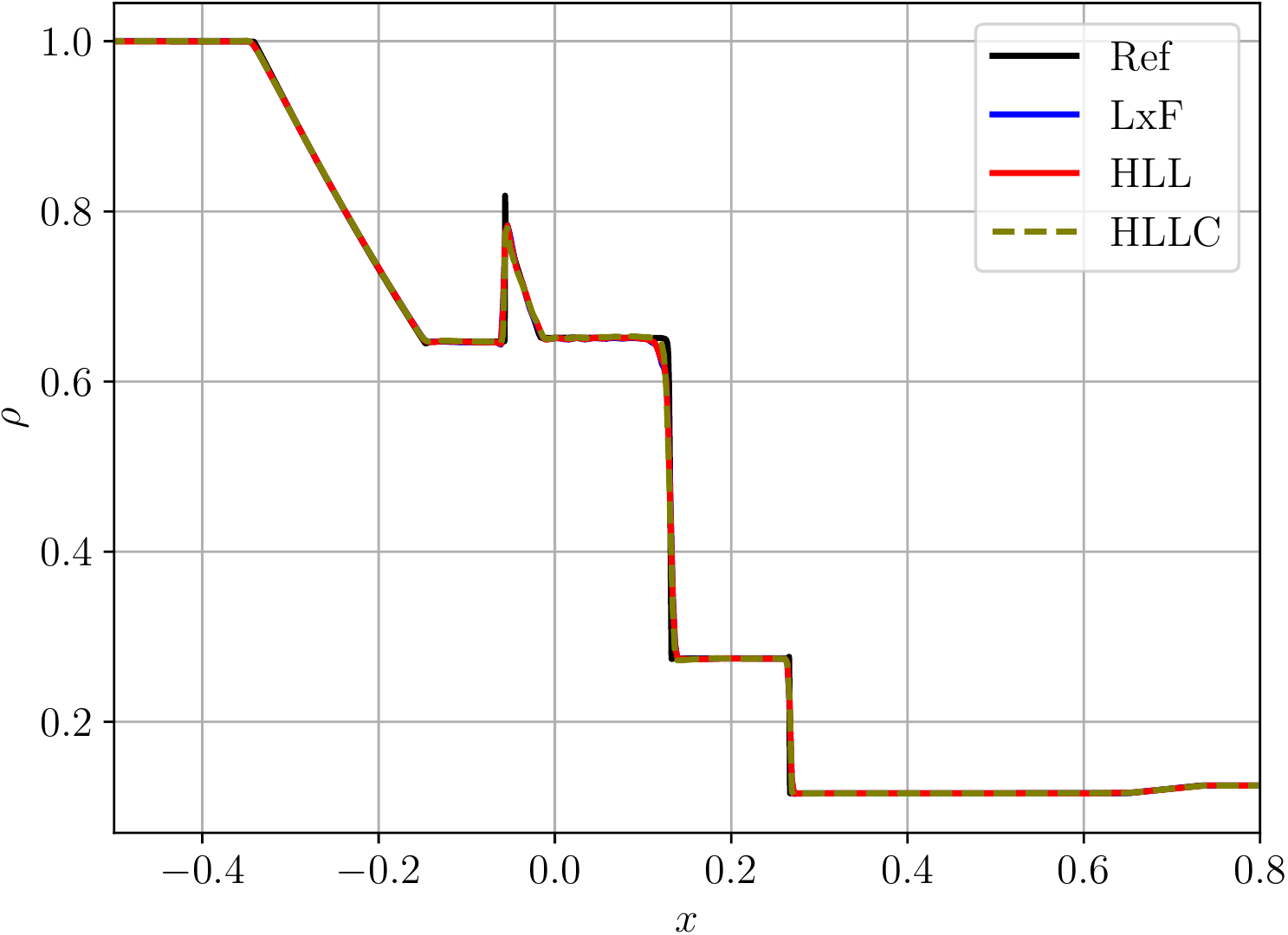} &
\includegraphics[width=0.48\textwidth]{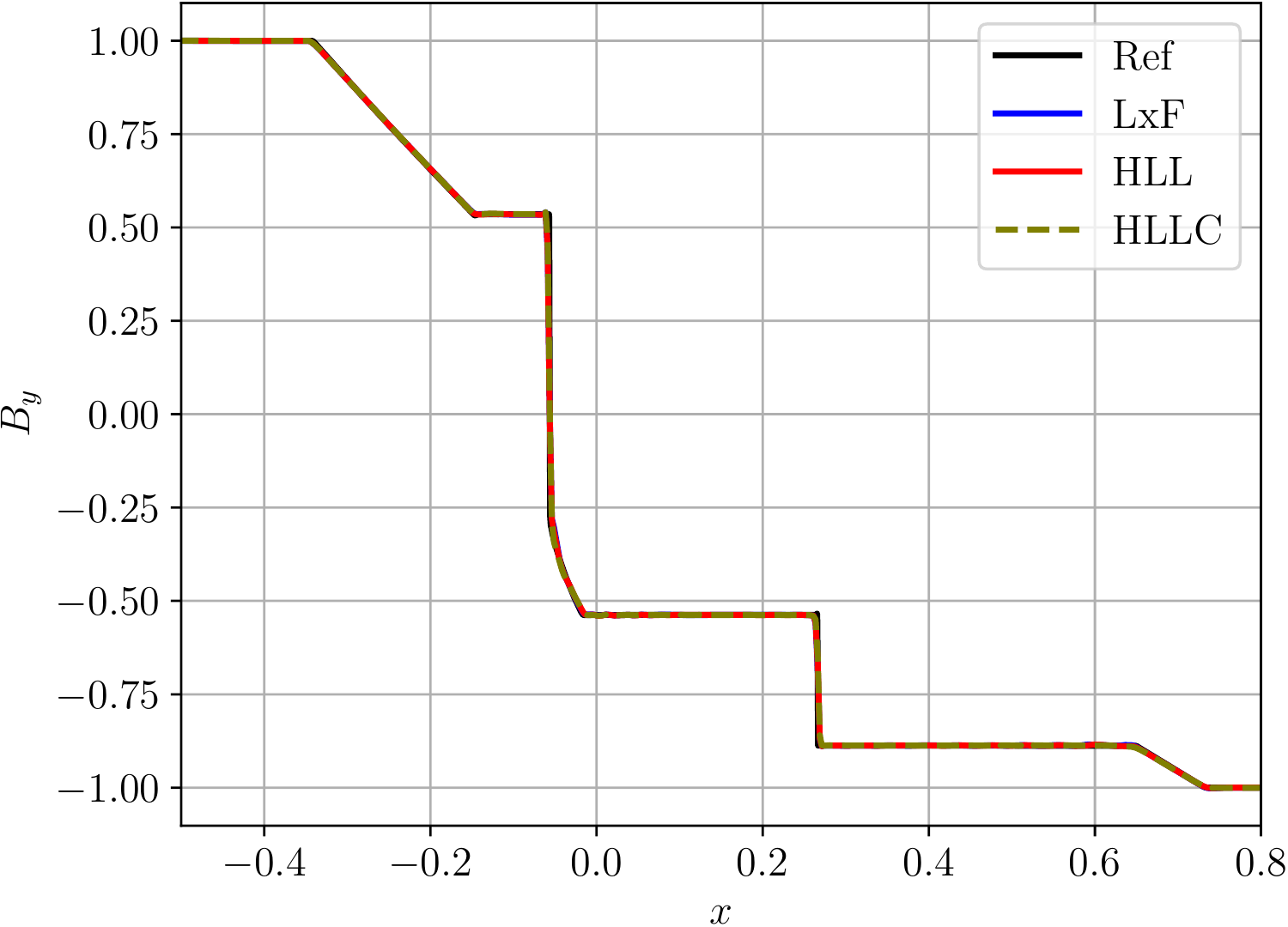} \\
\includegraphics[width=0.48\textwidth]{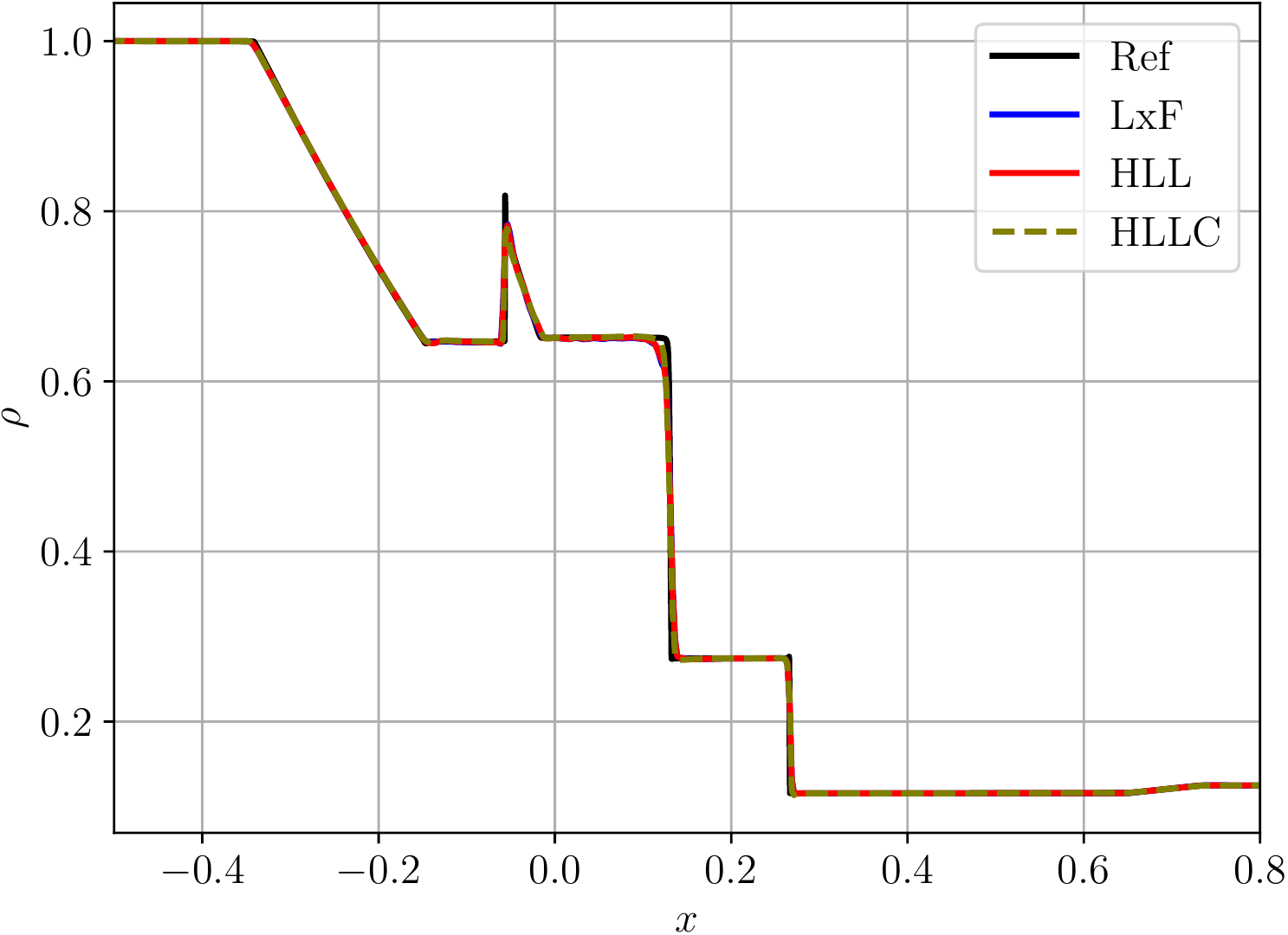} &
\includegraphics[width=0.48\textwidth]{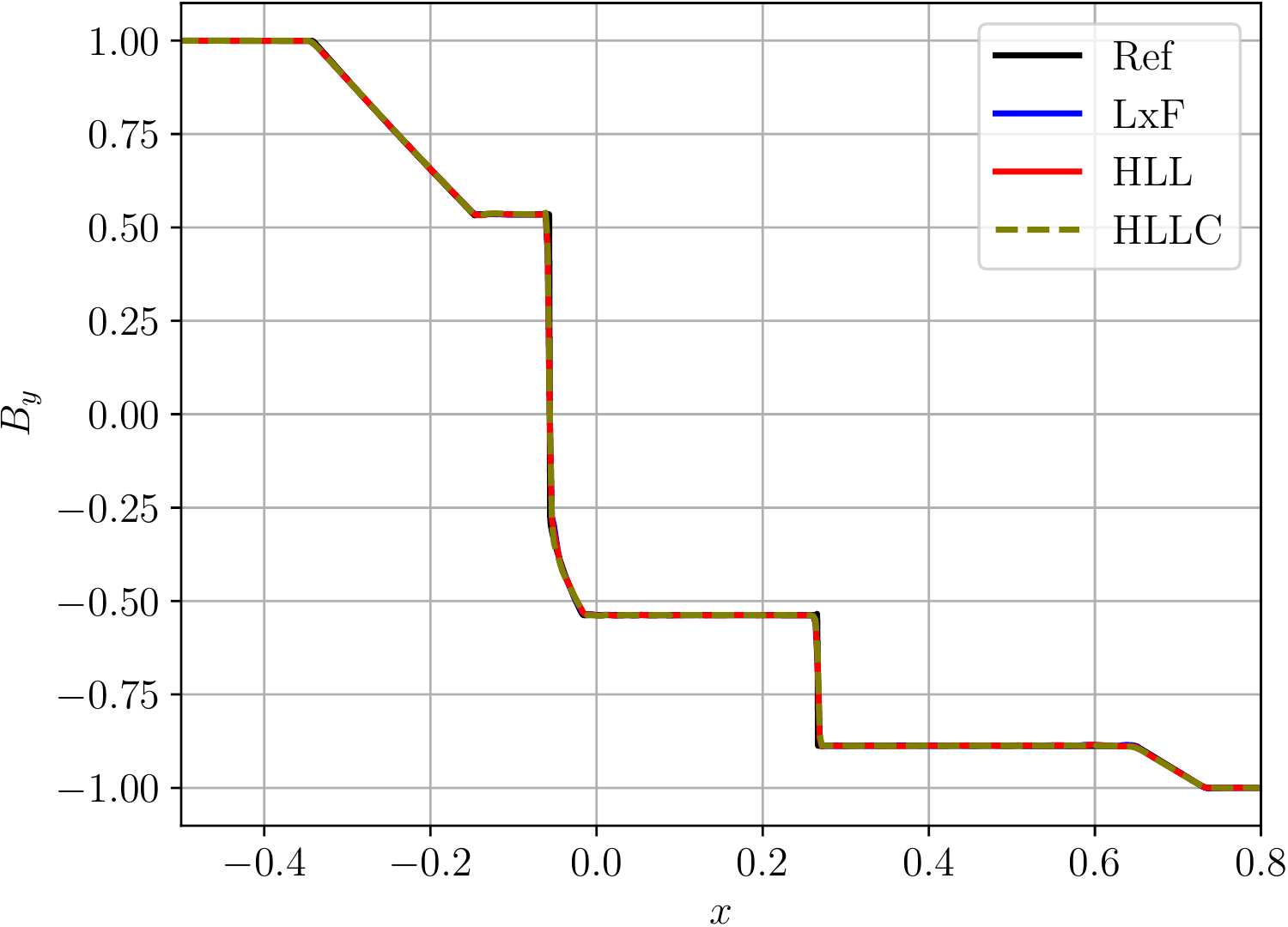} \\
\includegraphics[width=0.48\textwidth]{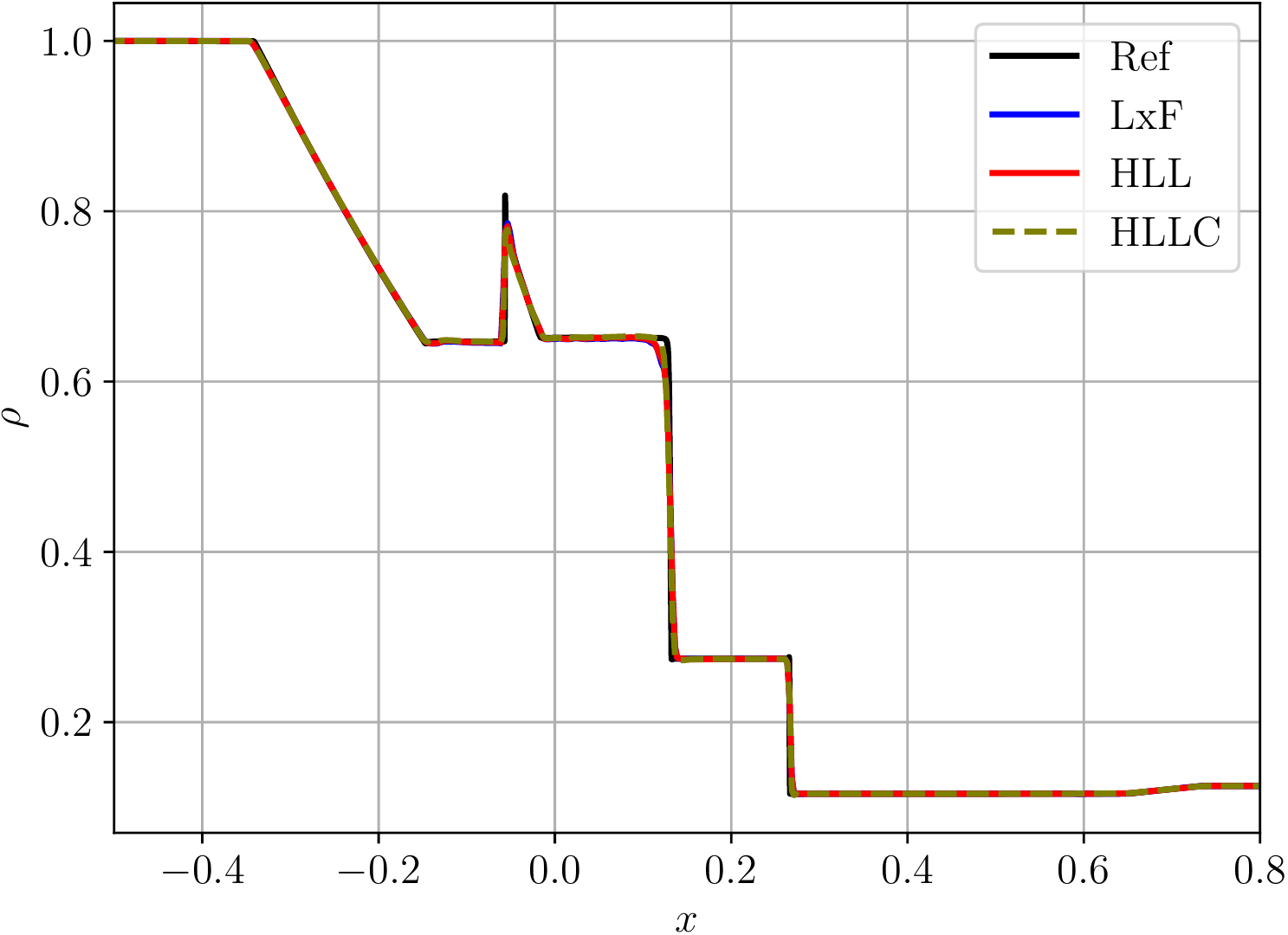} &
\includegraphics[width=0.48\textwidth]{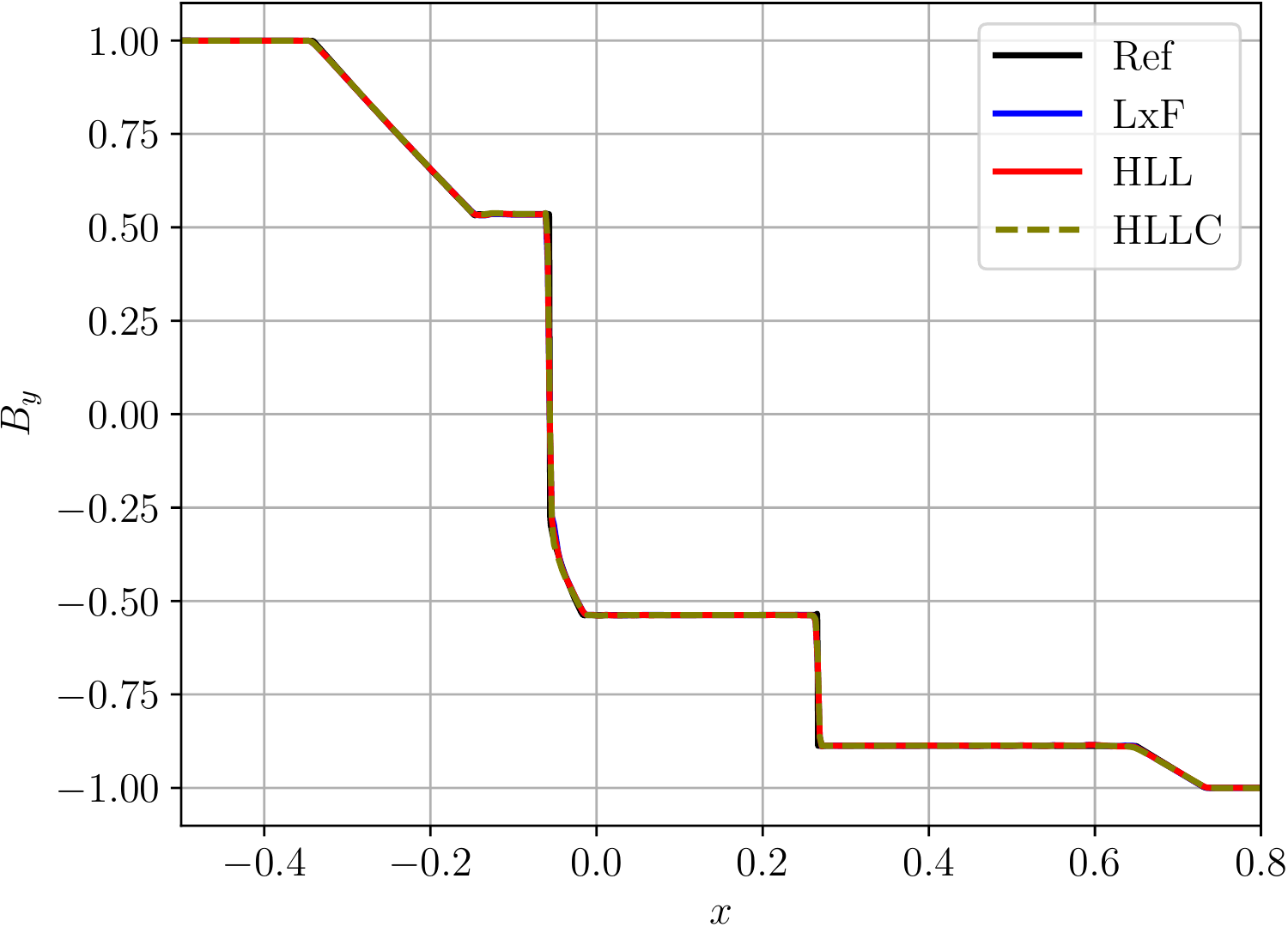}
\end{tabular}
\end{center}
\caption{Brio-Wu test case with 800 cells. Comparison of $\rho$ and $B_y$ obtained with different fluxes. Top row: $k=1$, middle row: $k=2$, bottom row: $k=3$}
\label{fig:bw1}
\end{figure}

\begin{figure}
\begin{center}
\begin{tabular}{ccc}
\includegraphics[width=0.33\textwidth]{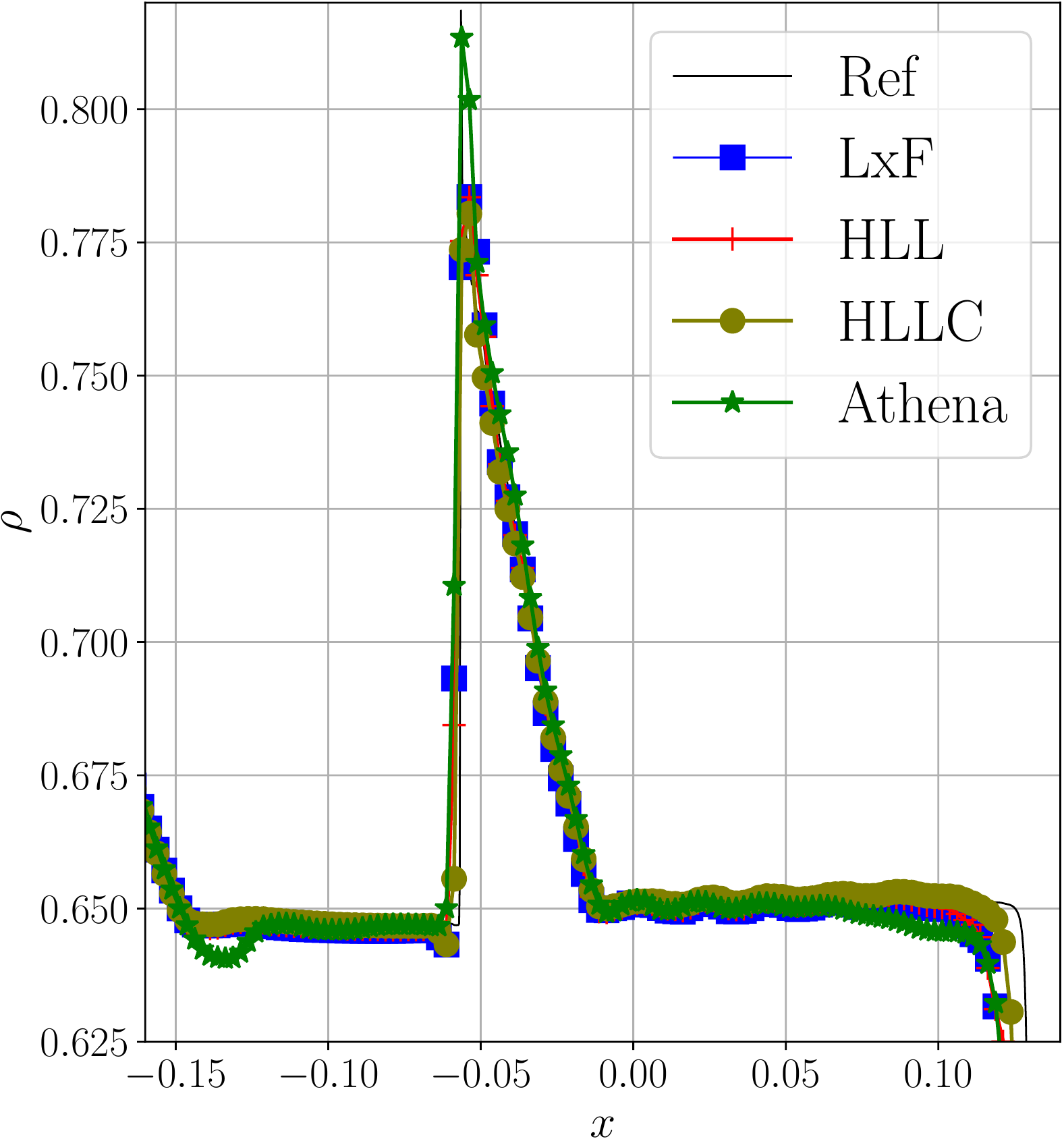} &
\includegraphics[width=0.33\textwidth]{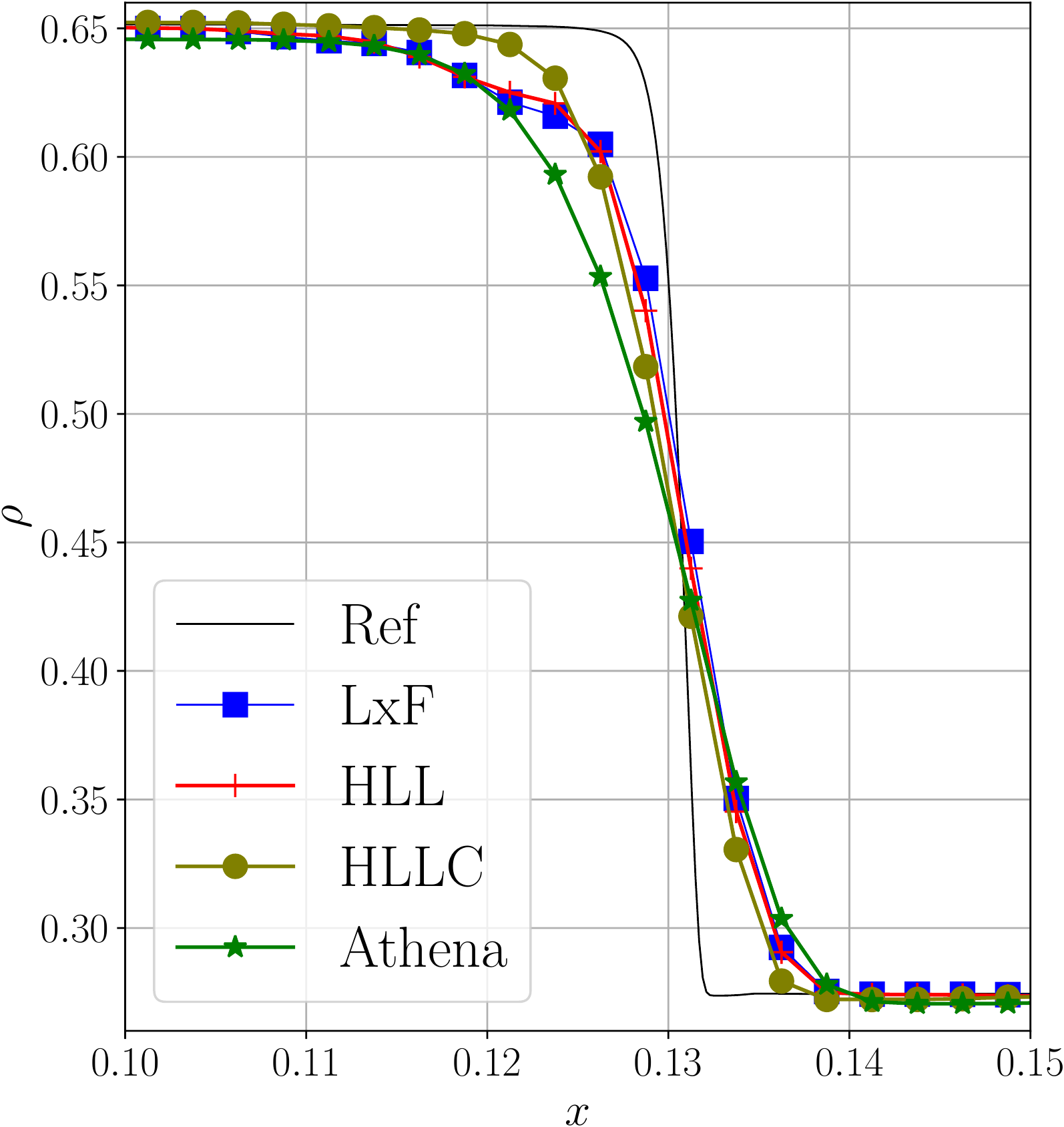} &
\includegraphics[width=0.33\textwidth]{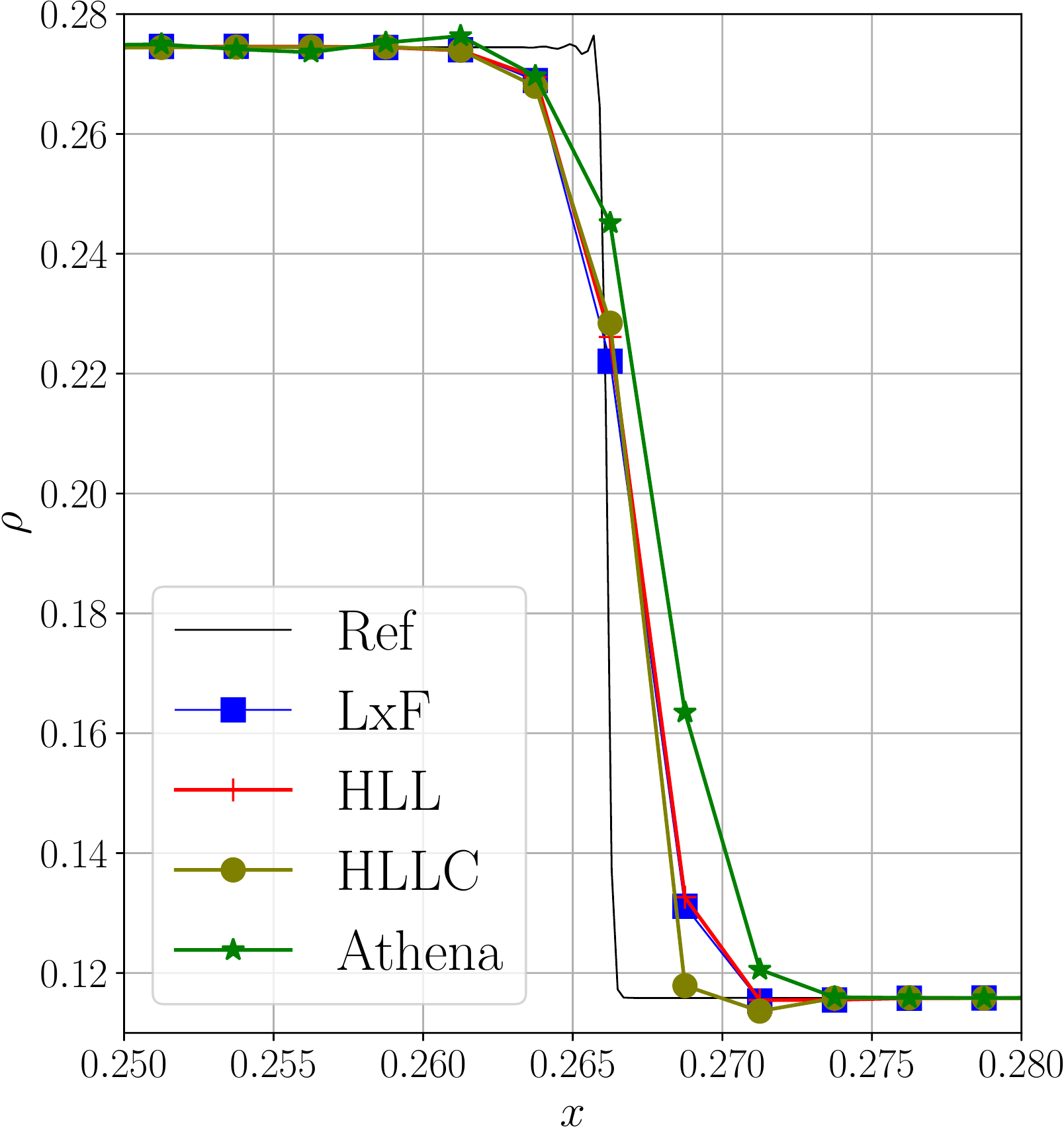} \\
(a) & (b) & (c)
\end{tabular}
\end{center}
\caption{Brio-Wu test case with 800 cells and degree $k=1$. Comparison of $\rho$  obtained with different fluxes and Athena code (a) around compound wave, (b) around contact wave, (c) around shock region.}
\label{fig:bw2}
\end{figure}

\begin{figure}
\begin{center}
\begin{tabular}{ccc}
\includegraphics[width=0.33\textwidth]{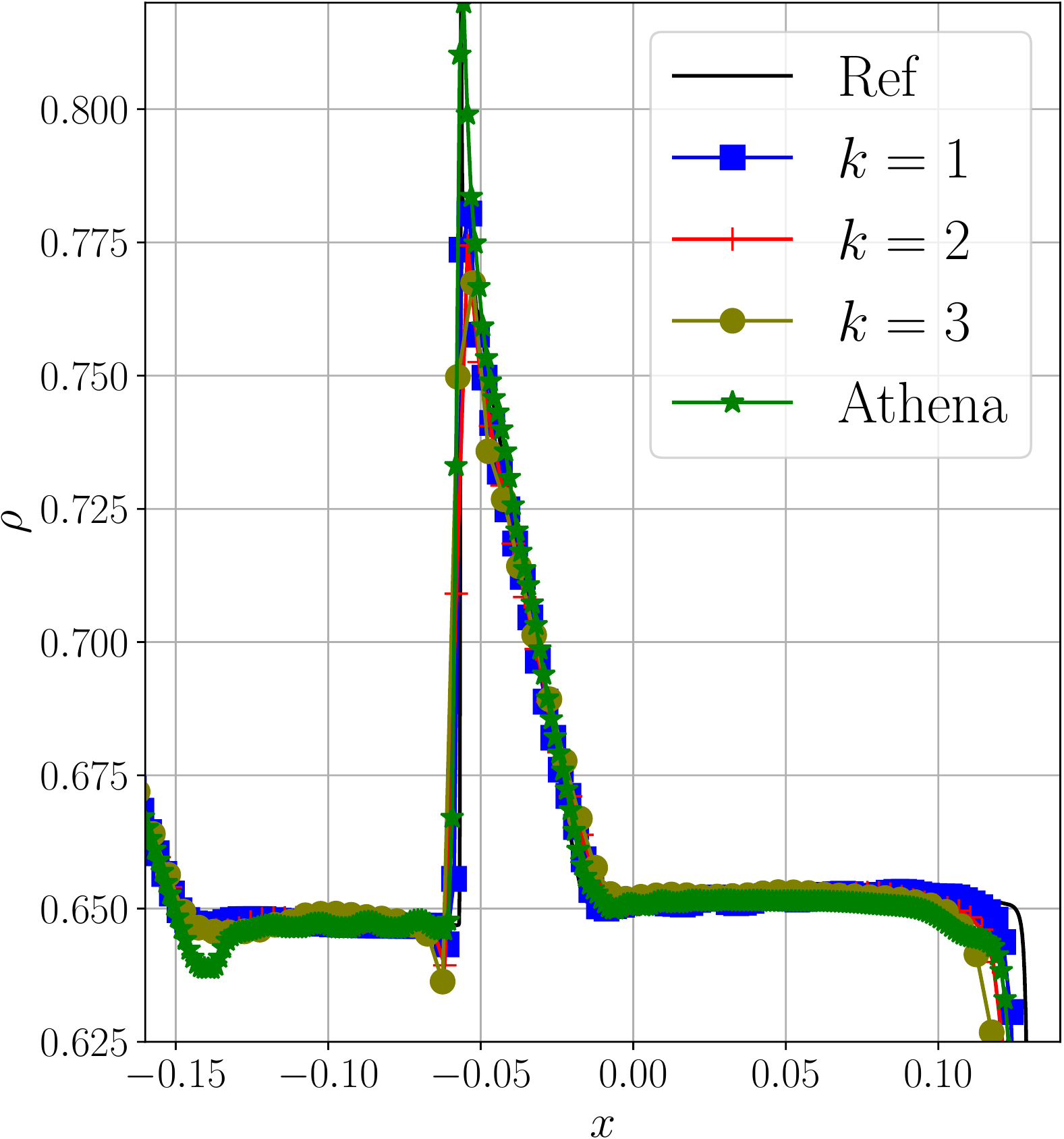} &
\includegraphics[width=0.33\textwidth]{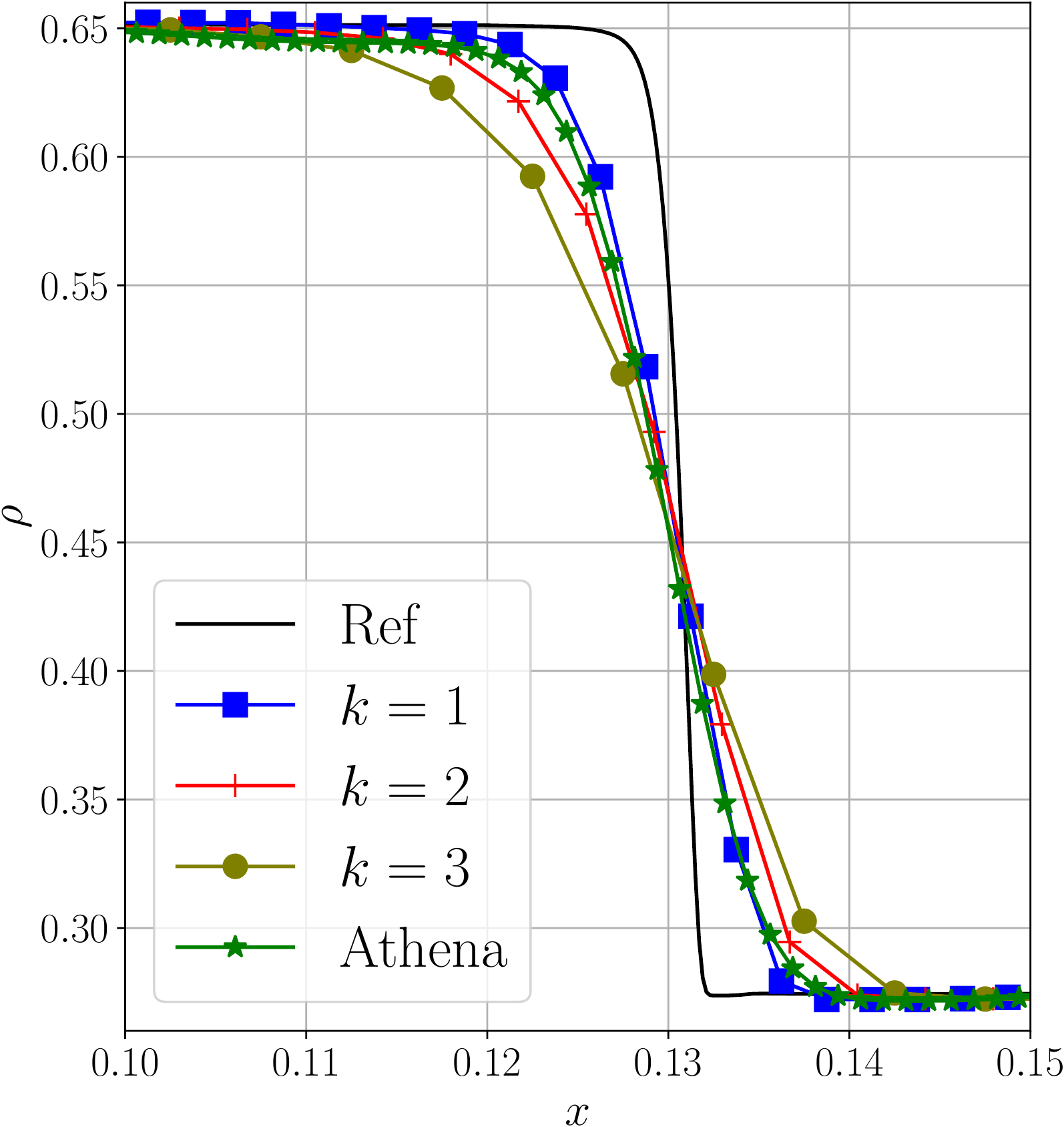} &
\includegraphics[width=0.33\textwidth]{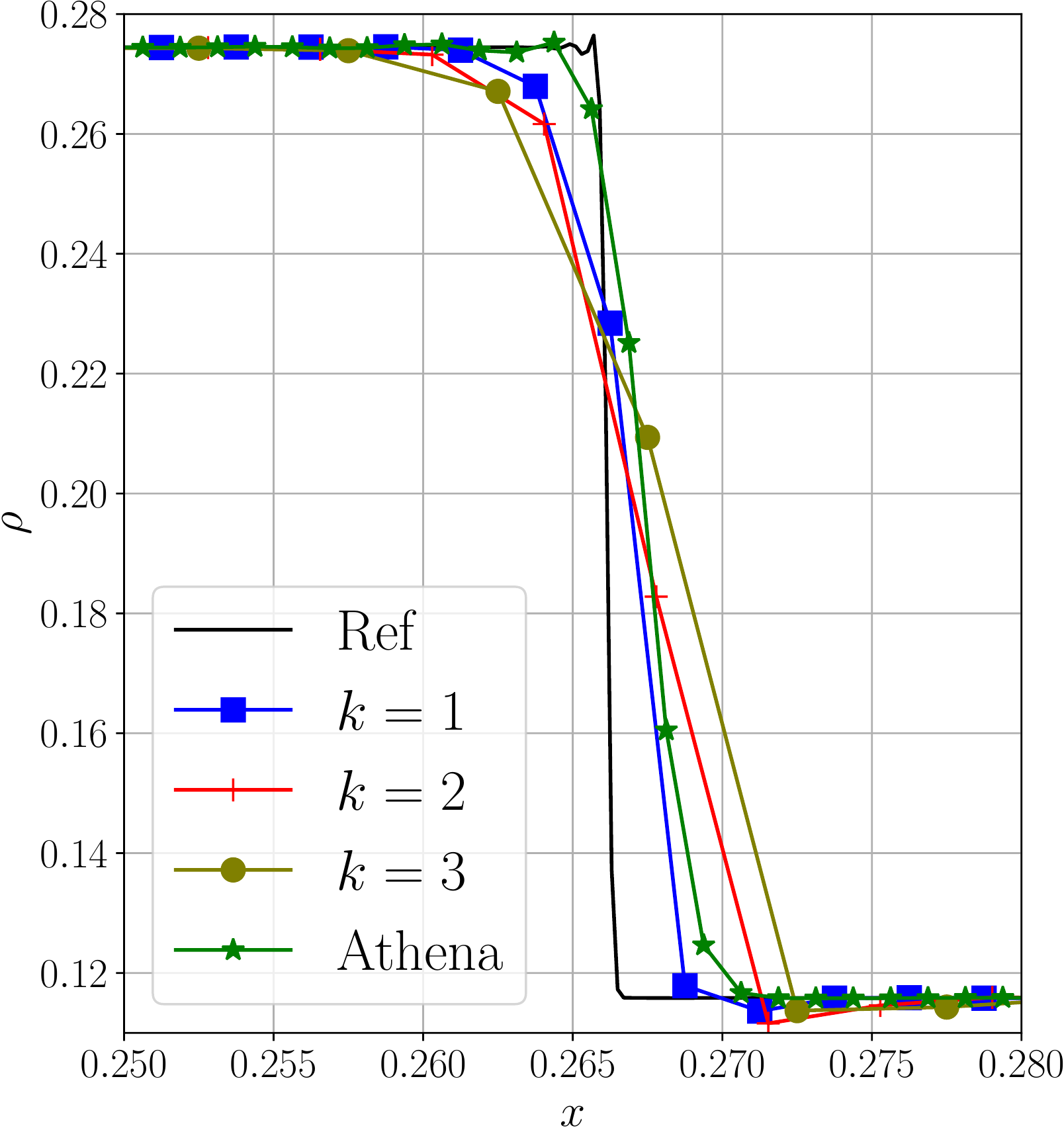} \\
(a) & (b) & (c)
\end{tabular}
\end{center}
\caption{Brio-Wu test case with same number of degrees of freedom and HLLC flux. Athena (1600 cells), $k=1, 800$ cells, $k=2, 534$ cells, $k=3, 400$ cells. Comparison of $\rho$   (a) around compound wave, (b) around contact wave, (c) around shock region.}
\label{fig:bw3}
\end{figure}

This is a classical shock tube problem for MHD \cite{Brio1988} and solution of this problem contains the fast rarefaction wave, the intermediate shock followed by a
slow rarefaction wave, the contact discontinuity, a slow shock, and a fast rarefaction wave. The initial condition has a discontinuity; for $x < 0$, the state is given by
\[
\rho = 1, \quad p = 1, \quad \vel = (0,0,0), \quad \bthree = (0.75, 1, 0)
\]
and for $x > 0$, it is given by
\[
\rho = 0.125, \quad p = 0.1, \quad \vel = (0,0,0), \quad \bthree = (0.75, -1, 0)
\]
which corresponds to Sod test case for hydrodynamics. The value of $\gamma$ is taken to be 5/3. We have computed the numerical solution for degree $k=1,2,3$  and LxF, HLL, HLLC fluxes at  time $T=0.2$ using $800$ cells.  In Figure~\ref{fig:bw1}, we have compared the numerical solutions obtained using different fluxes  with the reference solution computed using Athena code\footnote{Code taken from \url{https://github.com/PrincetonUniversity/athena-public-version} at git version \tt{273e451e16d3a5af594dd0a}} with 10000 cells.  The numerical results show that all the waves have been captured crisply and with very little or no oscillations. In Figure~\ref{fig:bw2}, we show zoomed view of the density plot around the compound wave, contact and shock waves, where we also compare with the results from Athena code using 800 cells. We can observe that the present method yields very good results. The resolution of the shock wave is quite good from all the solvers and the HLLC solver yields the best resolution, especially of the contact wave since only this solver explicitly includes the contact wave in its construction. \resb{In Figure~\ref{fig:bw3}, we perform computations at different degree and meshes so that the total number of degrees of freedom is similar in each case. We compare also with the finite volume results from Athena at the same resolution. While the degree $k=1$ results compare well with Athena, the higher order results are slightly diffused in the contact region due to use of coarser mesh and a TVD-type limiter.  A more sophisticated limiter with good sub-cell resolution is required to  achieve accurate results with high order schemes in the presence of discontinuities.}

\subsection{Consistency of 2-D Riemann solver}
\begin{figure}
\begin{center}
\begin{tabular}{cc}
\includegraphics[width=0.49\textwidth]{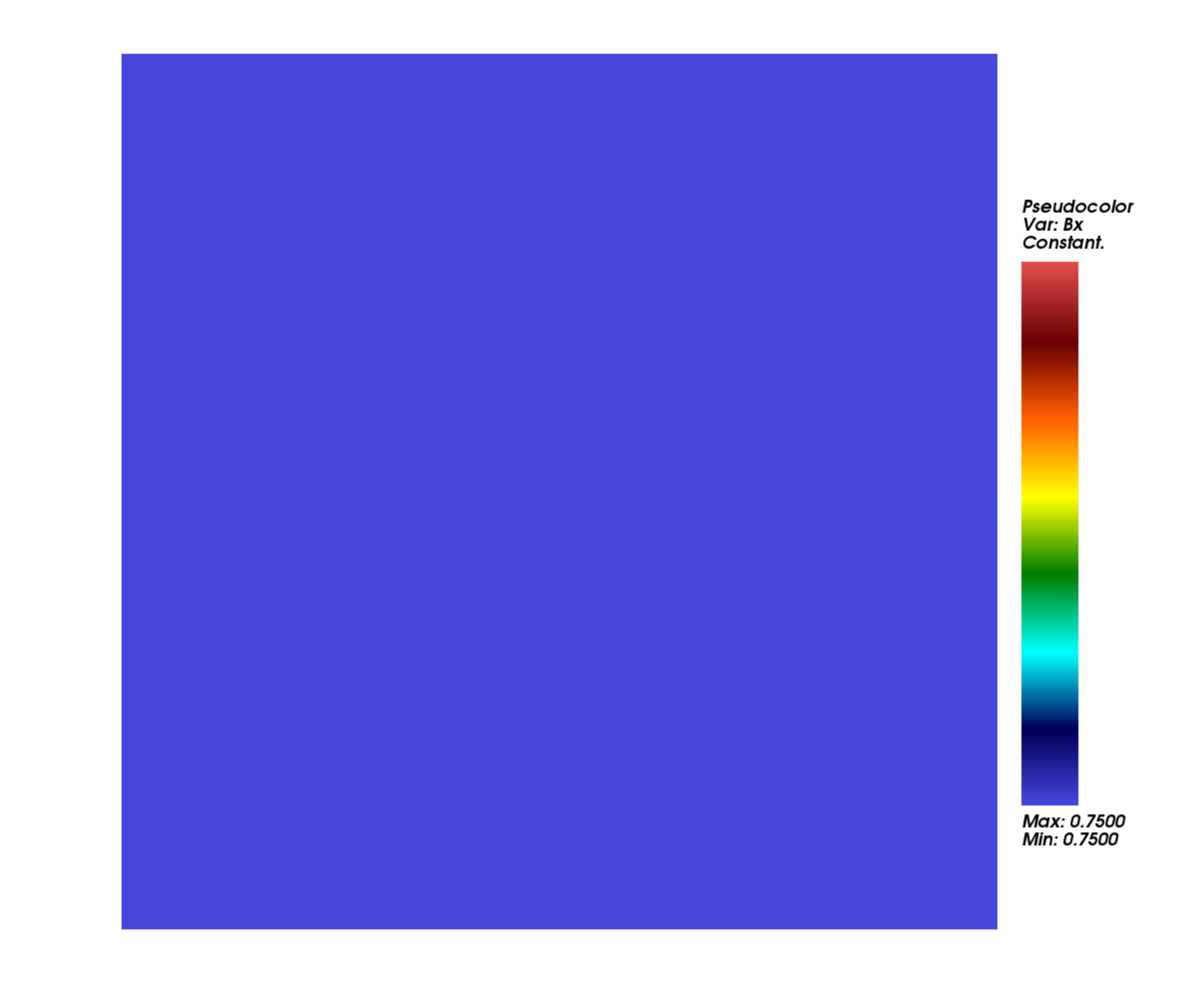} &
\includegraphics[width=0.49\textwidth]{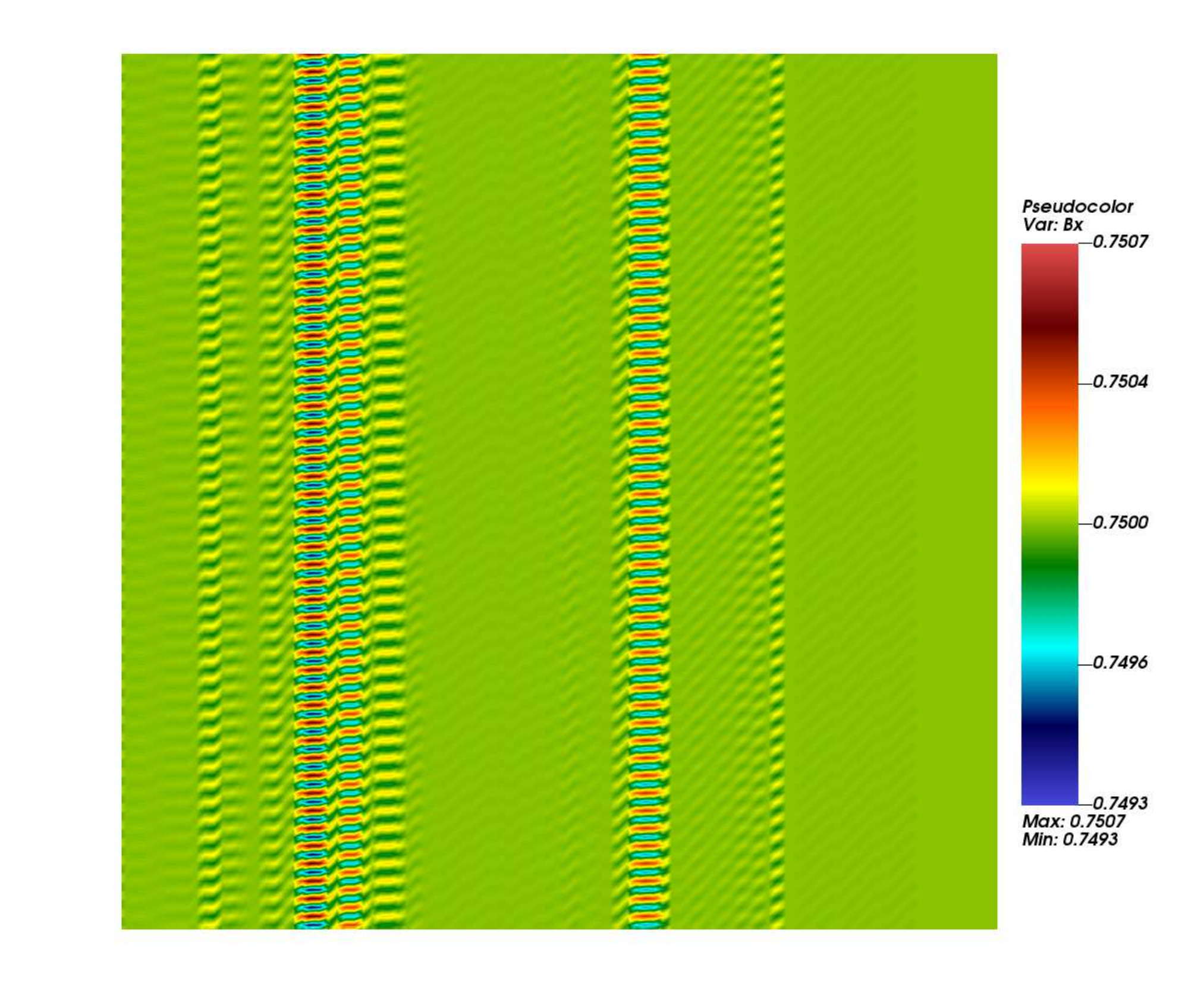} \\
(a) & (b)
\end{tabular}
\end{center}
\caption{Brio-Wu test case with HLL flux and $100 \times 100$ mesh. Color plot of $B_x$ using (a)~consistent 2-D Riemann solver, (b) inconsistent 2-D Riemann solver.}
\label{fig:consist}
\end{figure}

We take the Brio-Wu Riemann data to create a 2-D Riemann problem with $\con^{sw} = \con^{nw} = \con^L$ and $\con^{se} = \con^{ne} = \con^R$. The estimate of $E_z$ from both the 1-D and 2-D HLL Riemann solvers is same and equal to $-2.5400697250351683$ whereas if we use only the first term in \eqref{eq:Ezss}, we get a value of $-1.2700348625175841$, which has a very different magnitude. We run the Brio-Wu computation using both the consistent and inconsistent versions of the 2-D HLL Riemann solver on a 2-D domain $[-1,1]\times[-1,1]$ with a mesh of $100 \times 100$ cells. The resulting magnetic field component $B_x$ is shown in Figure~\ref{fig:consist}. We see that the consistent solver keeps the constancy of $B_x$ whereas the inconsistent version is not able to do so.

\subsection{Rotated shock tube}
\begin{figure}
\begin{center}
\newcolumntype{C}{ >{\centering\arraybackslash} m{1em} }
\newcolumntype{D}{ >{\centering\arraybackslash} m{5cm} }
\begin{tabular}{C D D D}
\rotatebox{90}{$k=1$} & \includegraphics[width=0.33\textwidth]{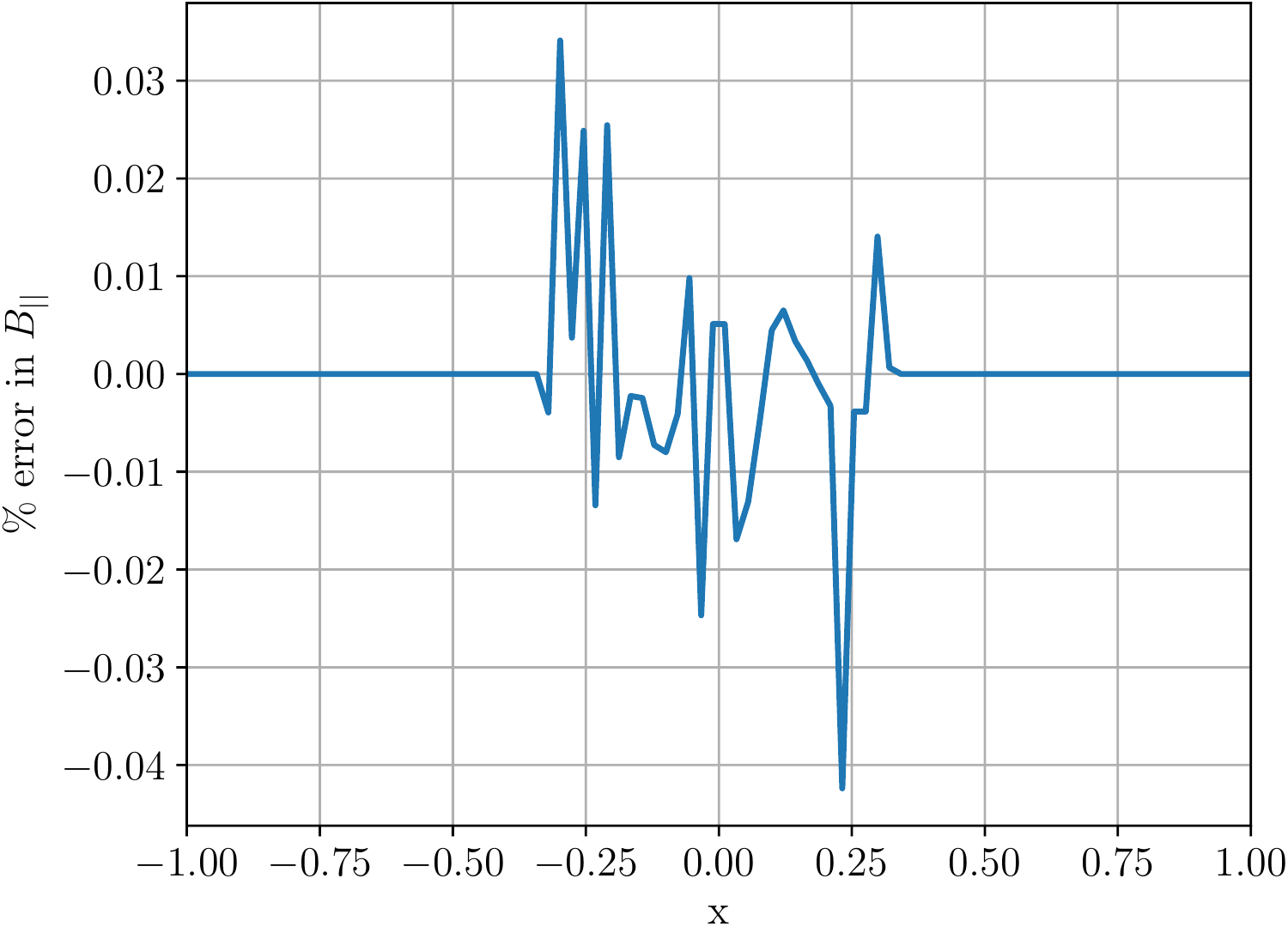} &
\includegraphics[width=0.33\textwidth]{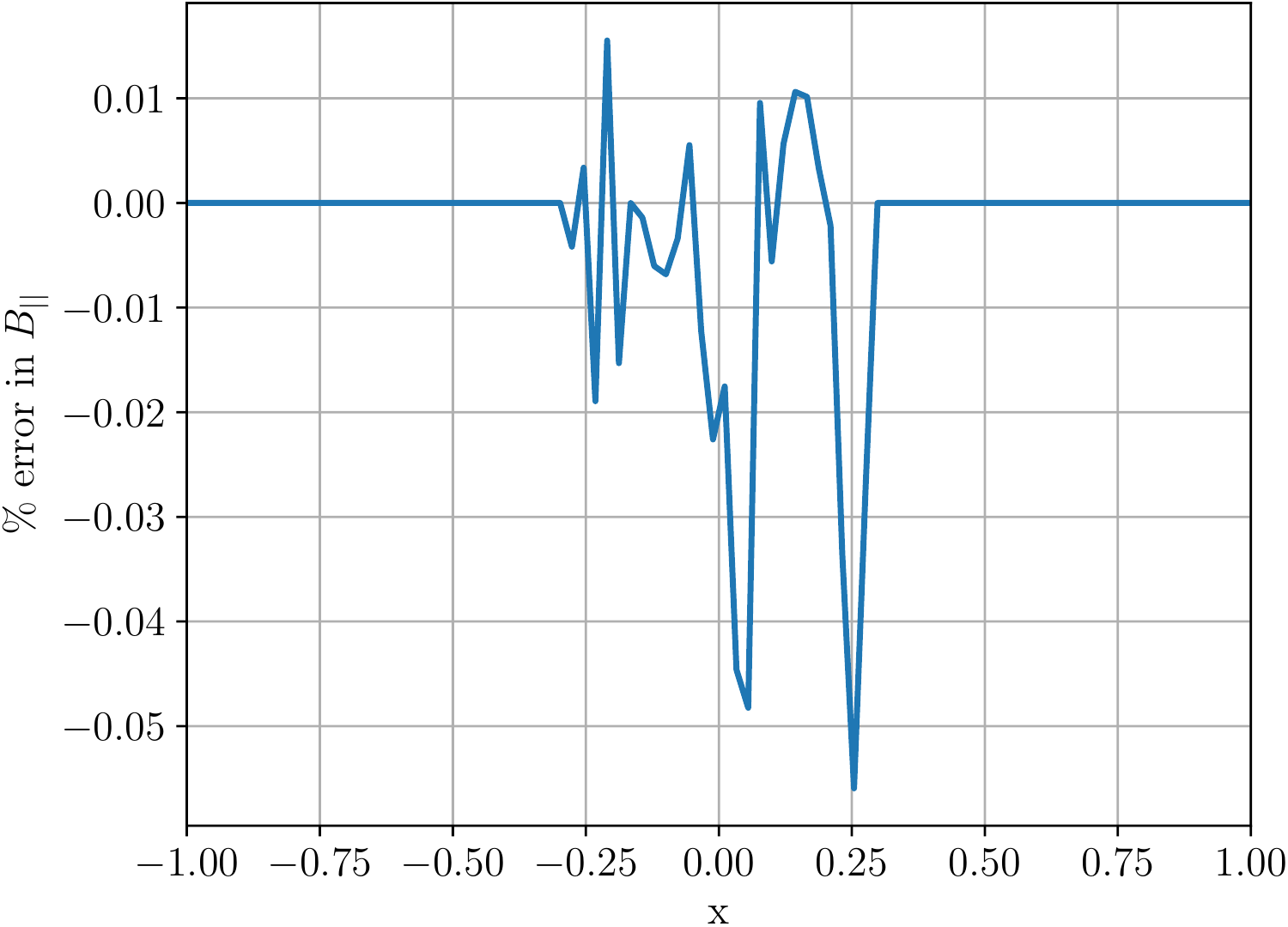} &
\includegraphics[width=0.33\textwidth]{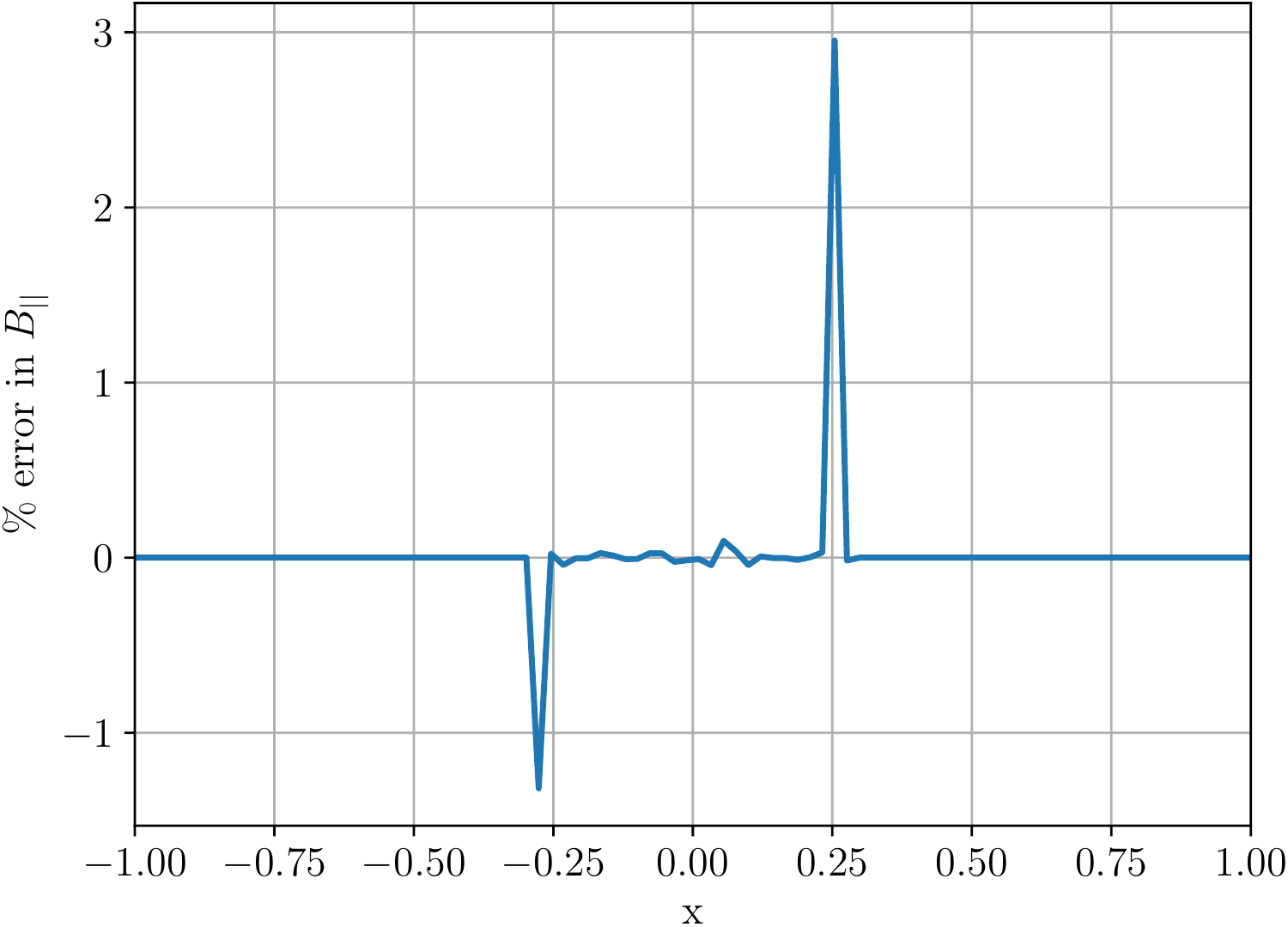}\\
\rotatebox{90}{$k=2$} & \includegraphics[width=0.33\textwidth]{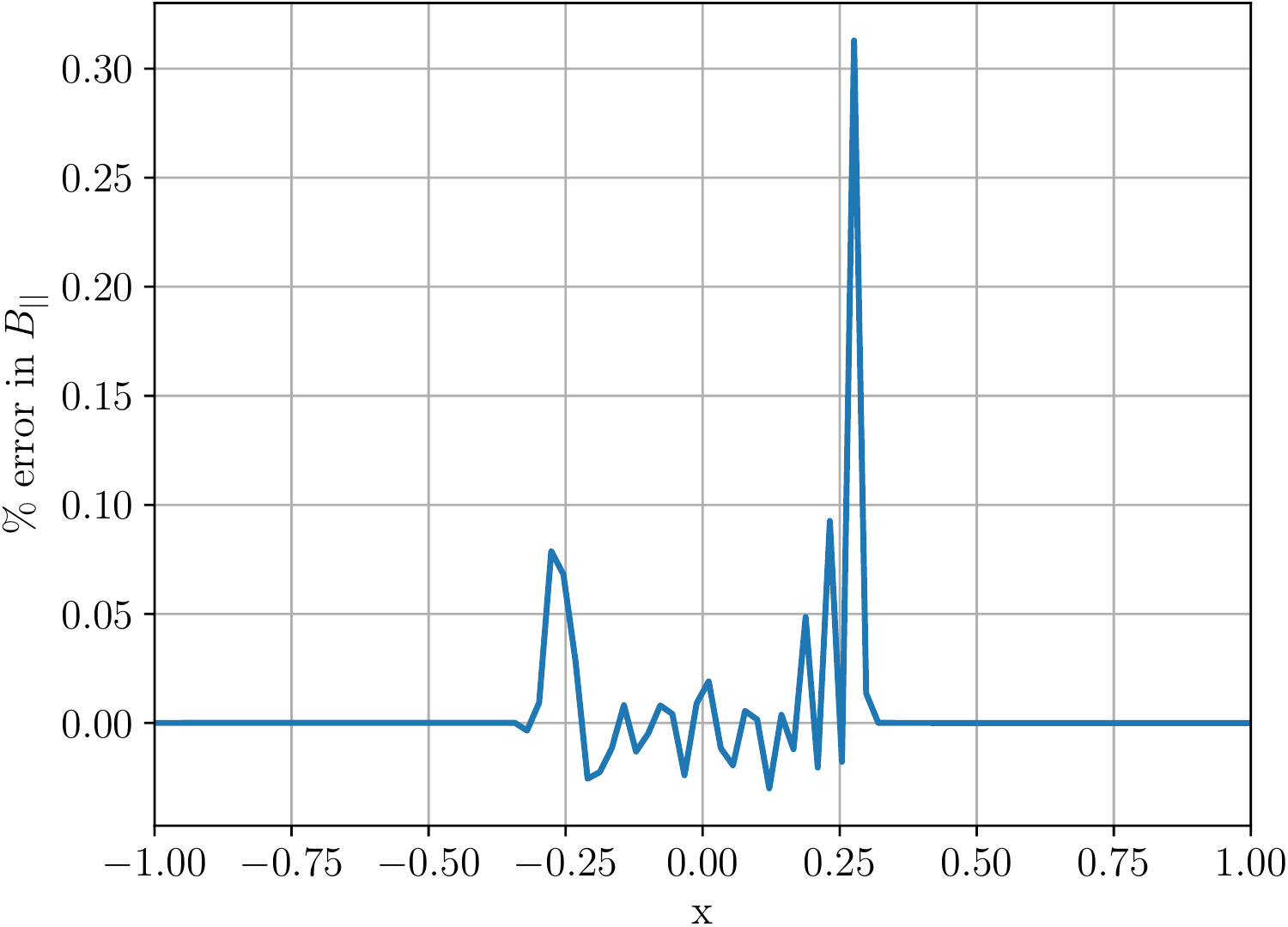} &
\includegraphics[width=0.33\textwidth]{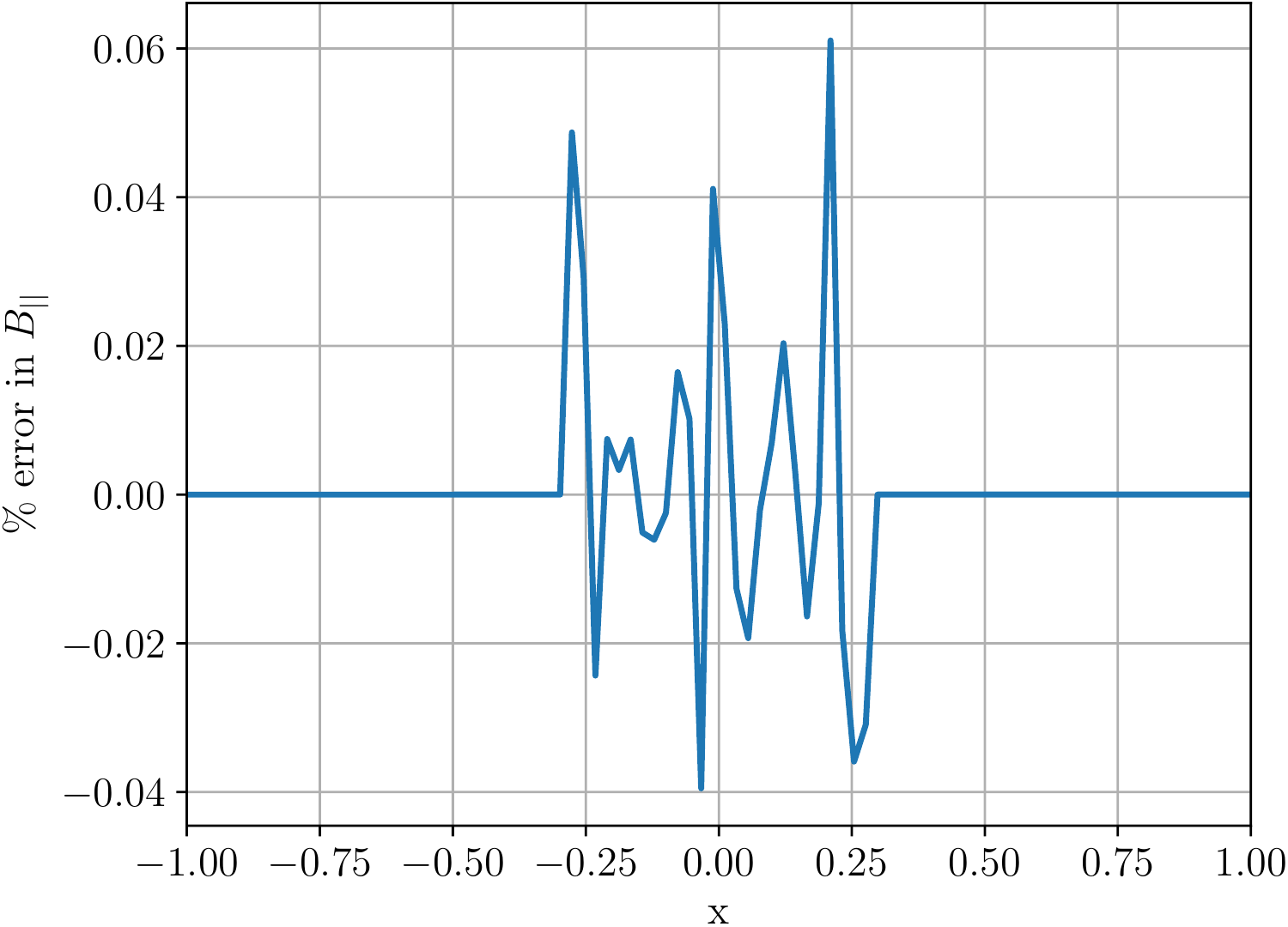} &
\includegraphics[width=0.33\textwidth]{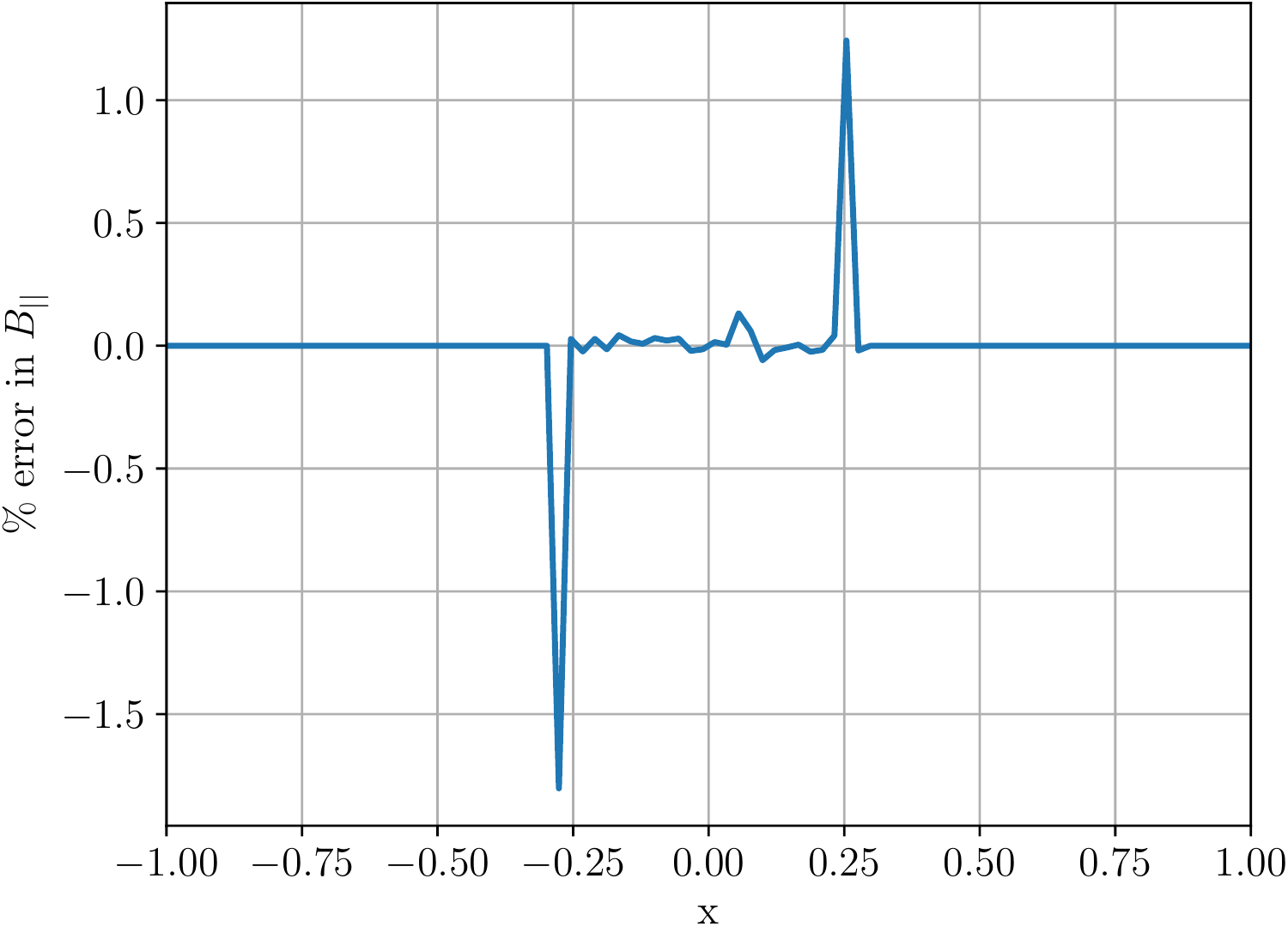}\\
\rotatebox{90}{$k=3$} & \includegraphics[width=0.33\textwidth]{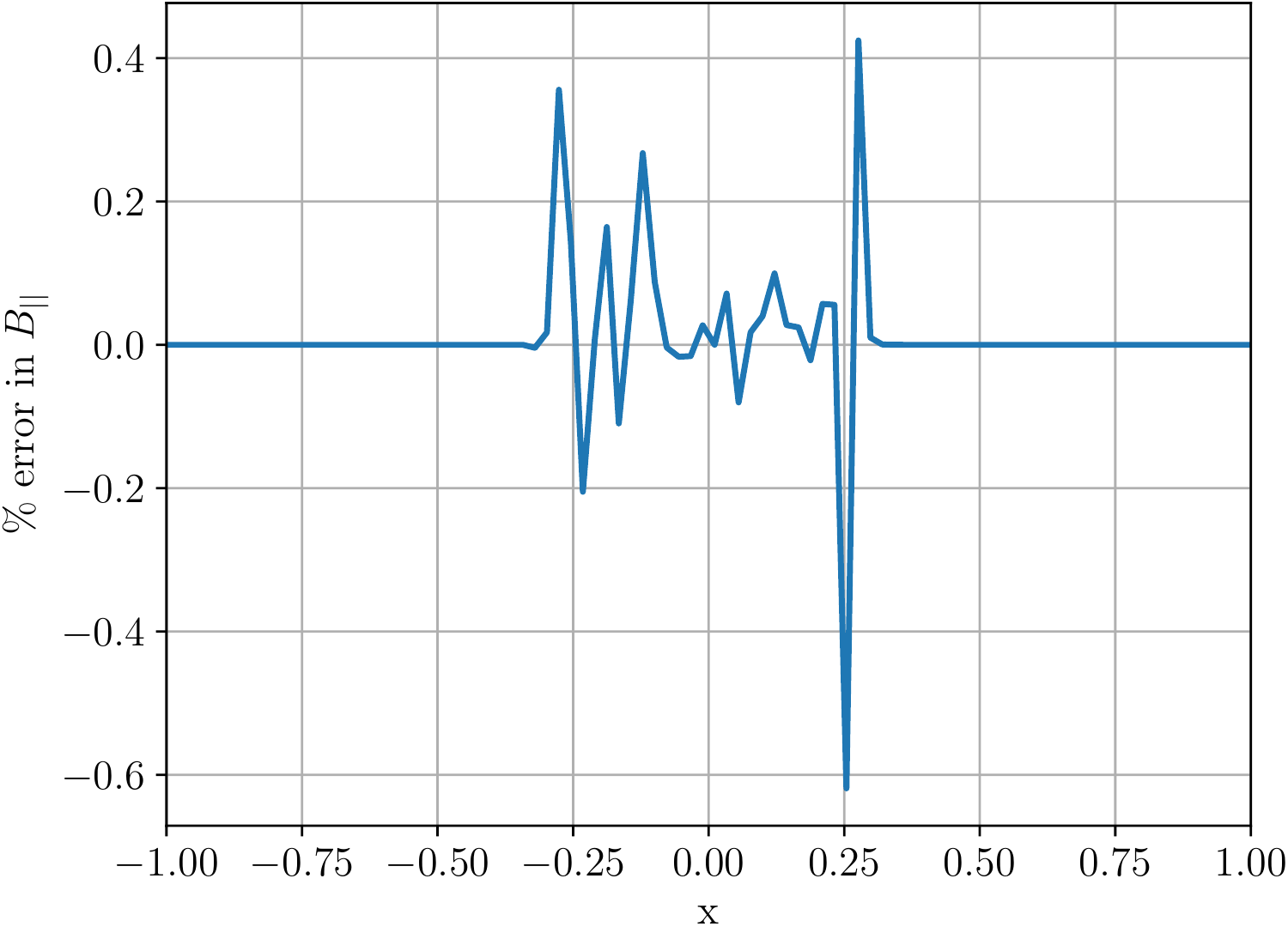} &
\includegraphics[width=0.33\textwidth]{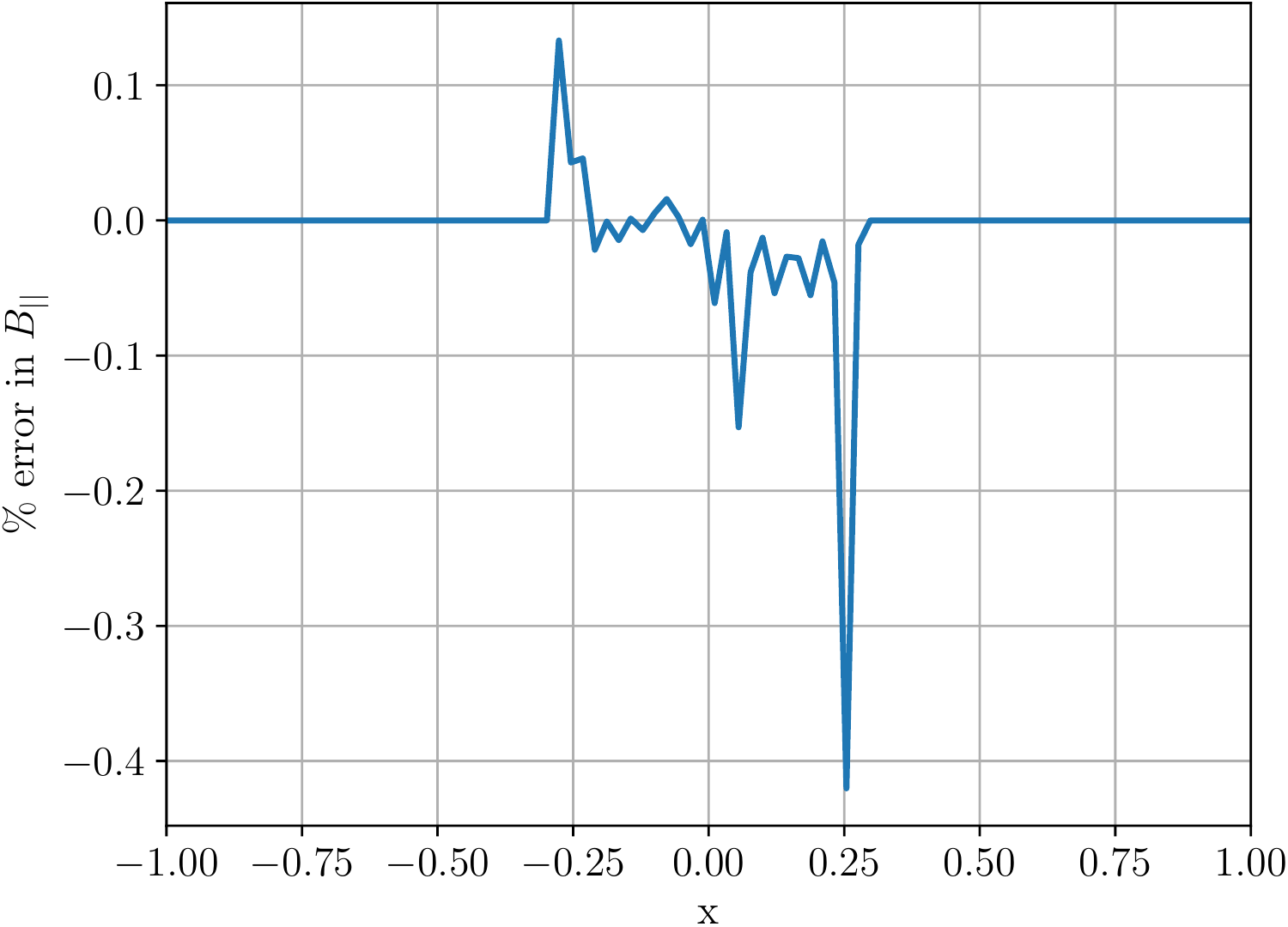} &
\includegraphics[width=0.33\textwidth]{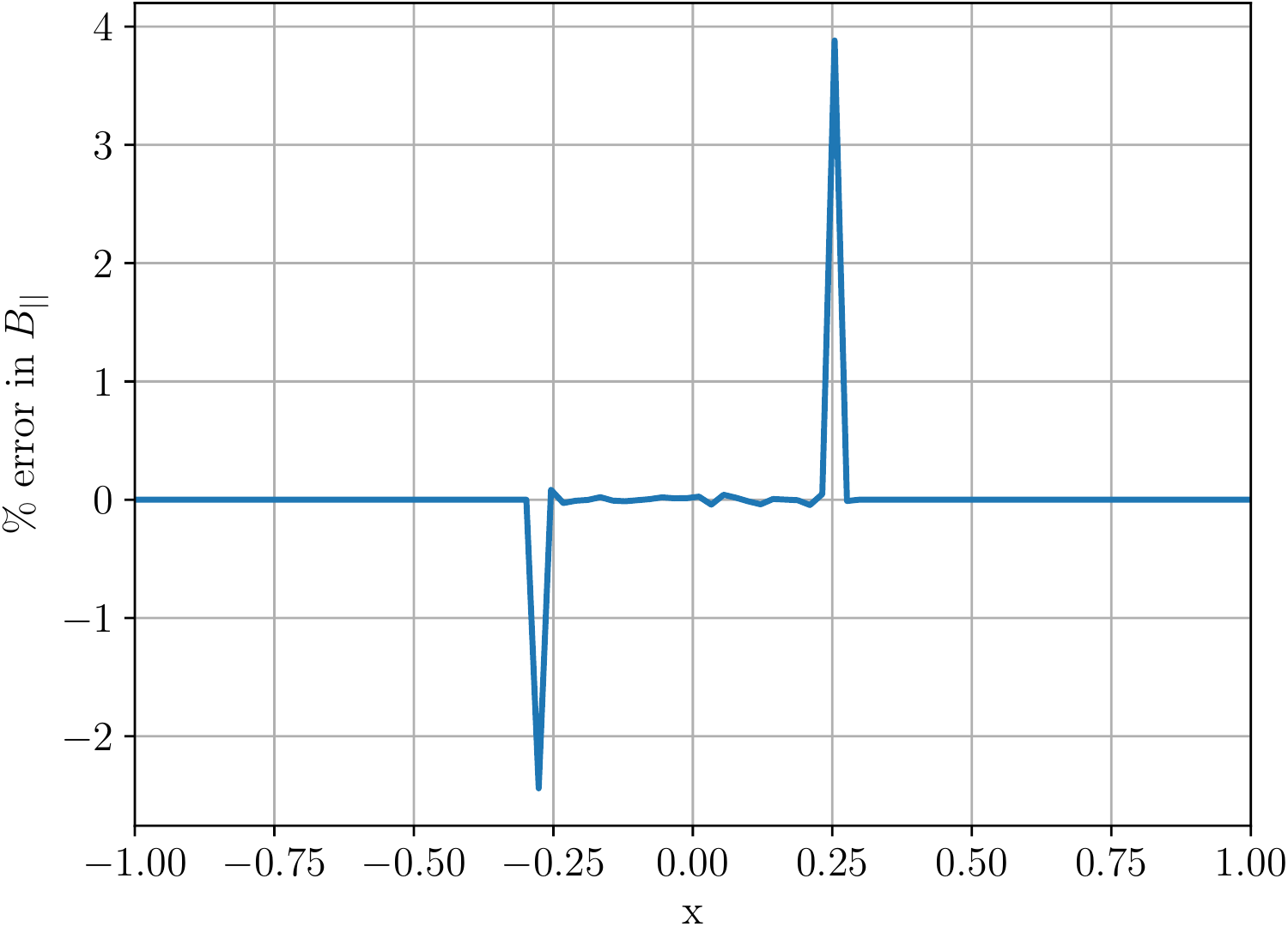}\\
 & LxF  & HLL &  HLLC 
\end{tabular}
\caption{Comparison of the percentage error on the parallel magnetic field $B_\|$ for rotated shock tube test over a grid of size $128\times 128$ using LxF, HLL, and HLLC fluxes. First row degree $k=1$, second row degree $k=2$, and third row degree $k=3$.}
\label{fig:rstube0}
\end{center}
\end{figure}

The initial condition is a Riemann problem which is aligned at an angle to the mesh~\cite{Toth2000}. We take the domain to be $[-1,+1] \times [-1,+1]$ with Neumann boundary conditions, and the initial discontinuity is across the line $x+y = 0$. Throughout the domain, the density $\rho = 1$ and the magnetic field is
\[
B_x = B_\perp \cos\alpha - B_\| \sin\alpha, \qquad B_y = B_\perp \sin\alpha + B_\| \cos\alpha, \qquad B_z = 0
\]
In the region $x + y < 0$, the remaining quantities are given by
\[
p = 20, \qquad \vel =10(\cos\alpha,~\sin\alpha,~0)
\]
and for $x + y > 0$
\[
p = 1, \qquad \vel = -10 (\cos\alpha,~\sin\alpha,~0)
\]
where $\alpha = \pi/4$ is the orientation of the initial discontinuity. Note that this represents a one dimensional Riemann problem when viewed along the line $x=y$. The solution is computed up to the time $T = 0.08/\cos\alpha$ and we plot the solution along the line $x = y$. At this final time, the solution along this line is not affected by the flow features that develop due to boundary conditions. The exact solution should have $B_\| = -B_x \sin\alpha + B_y \cos\alpha = 5/\sqrt{4\pi}$ and $B_\perp = B_x \cos\alpha + B_y \sin\alpha = 5/\sqrt{4\pi}$ along the line $x=y$. The constancy of $B_\perp$ is difficult to obtain in numerical schemes which do not satisfy the divergence-free constraint. Numerical solutions are computed over a grid of size $128\times128$ using LxF, HLL, and HLLC fluxes. In Figure~\ref{fig:rstube0}, we have compared the relative percentage error on the magnetic field $B_\|$ for the considered fluxes and $k=1,2,3$. Though we cannot clearly classify which flux is best, the LxF and HLL fluxes yield the smallest and very similar levels of error, while the errors for HLLC are highest. However, even the largest observed error is similar to what is observed with other standard constraint preserving schemes~\cite{Toth2000}. The large errors are observed only at the location of the discontinuity and the error is small in other regions, which is a benefit obtained due to the divergence-free methods.

\subsection{Orszag-Tang vortex}
\begin{figure}
\begin{center}
\begin{tabular}{ccc}
\includegraphics[width=0.32\textwidth]{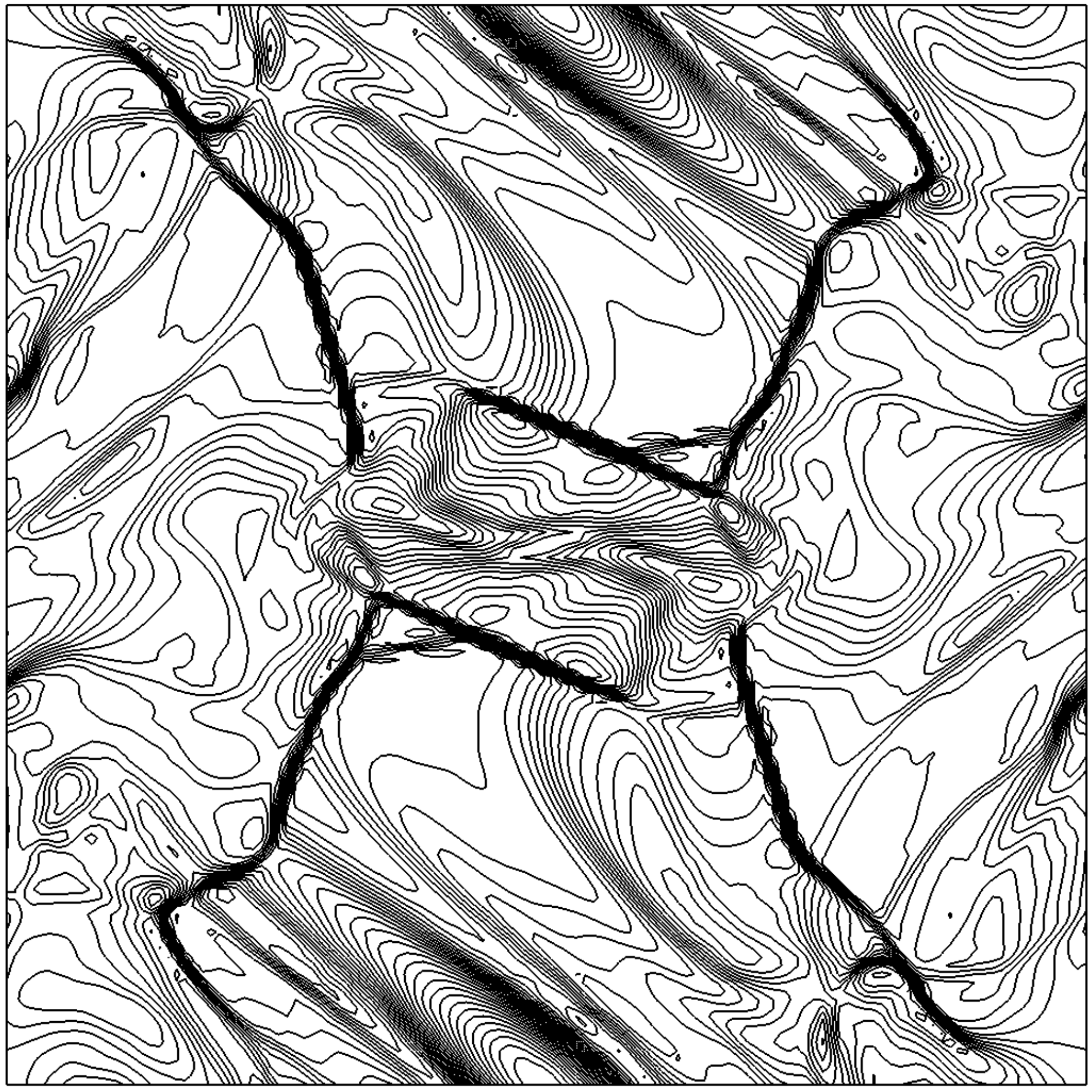} &
\includegraphics[width=0.32\textwidth]{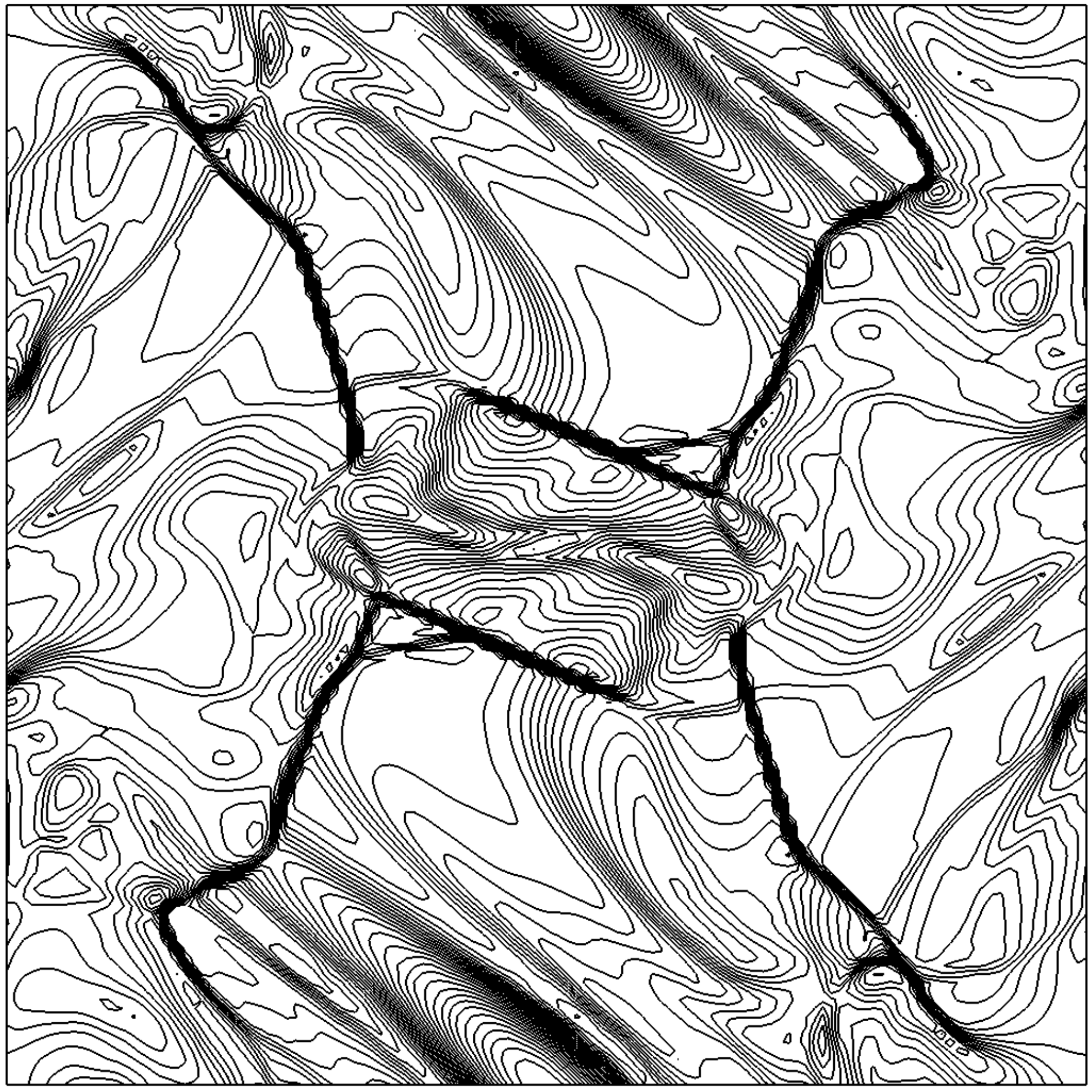} &
\includegraphics[width=0.32\textwidth]{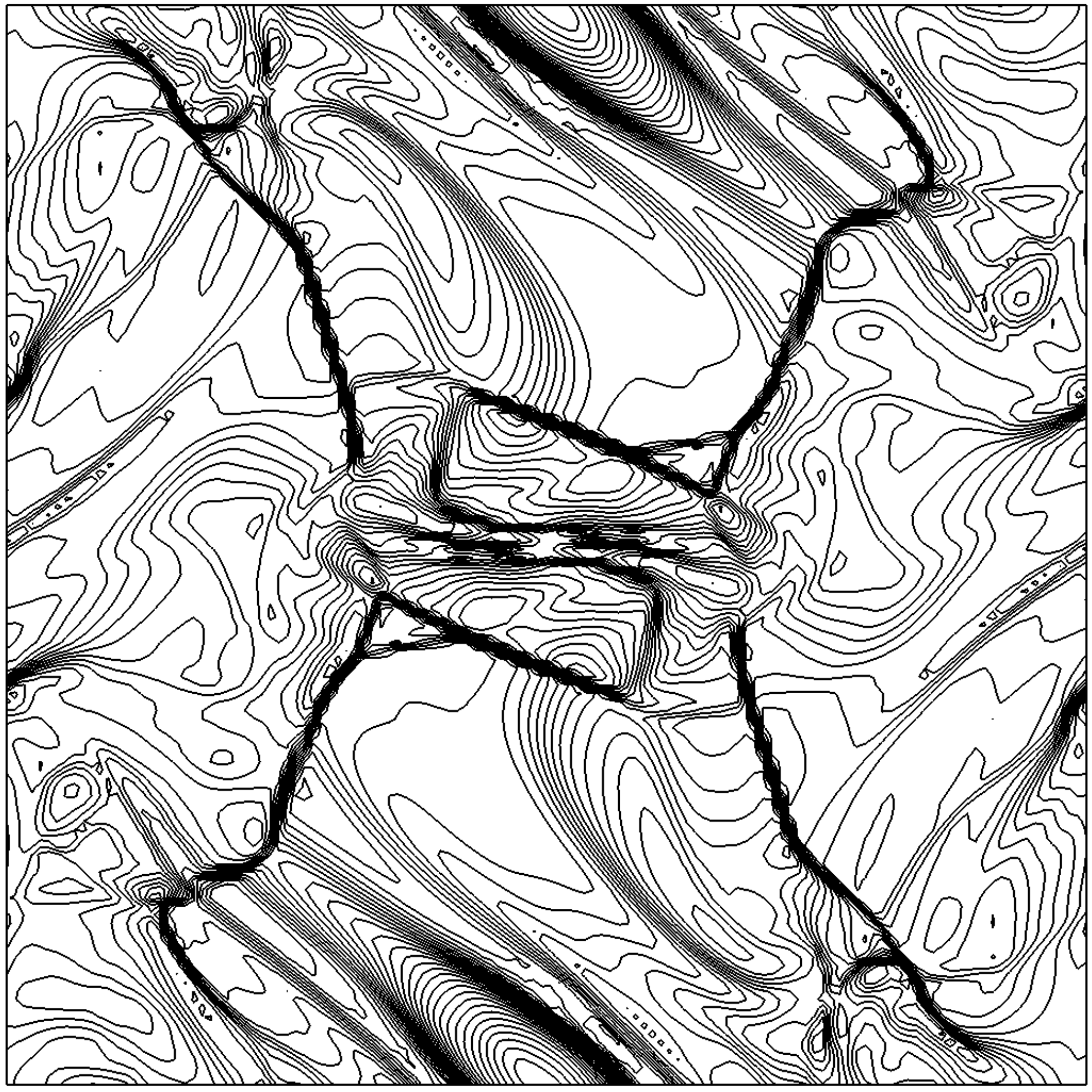} \\
\includegraphics[width=0.32\textwidth]{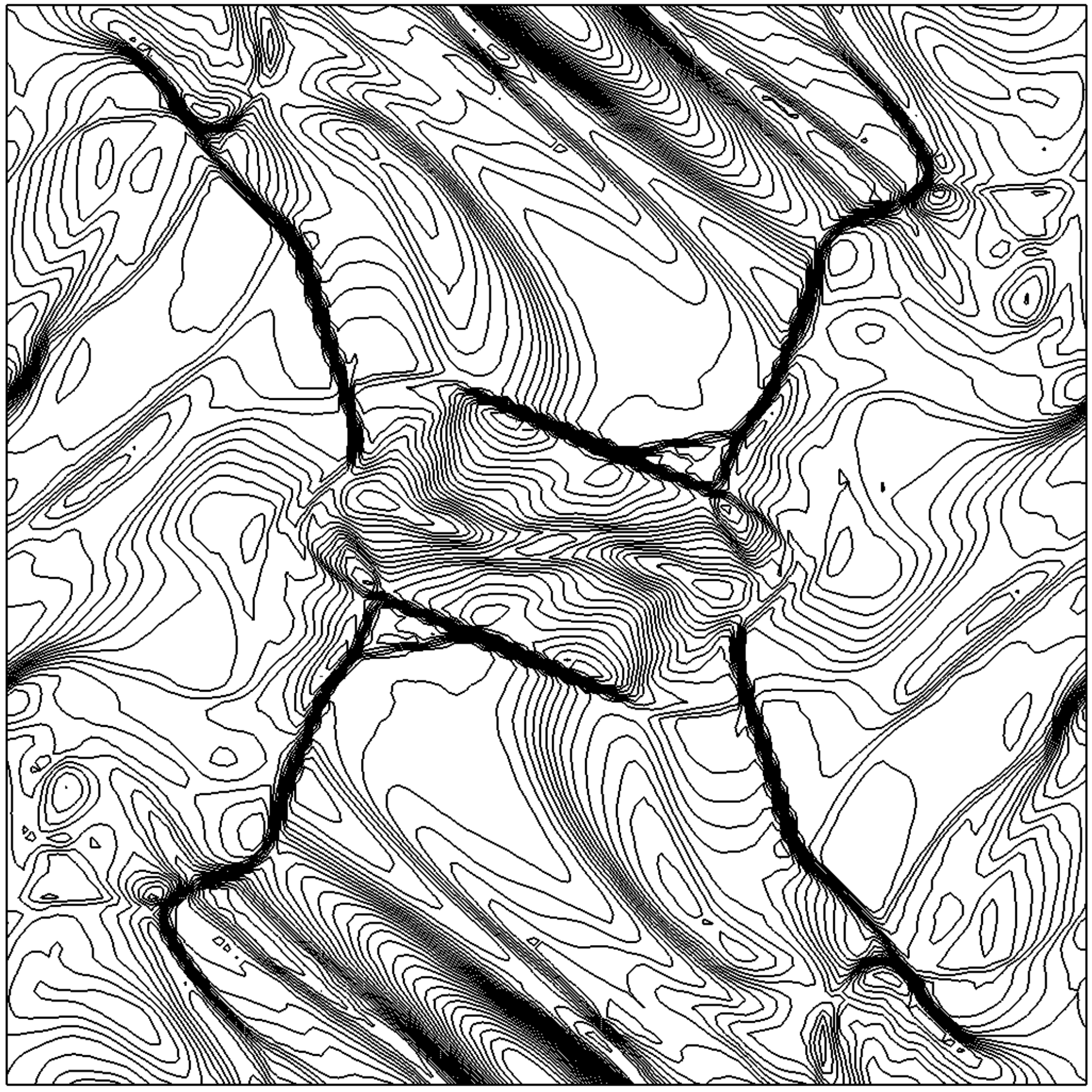} &
\includegraphics[width=0.32\textwidth]{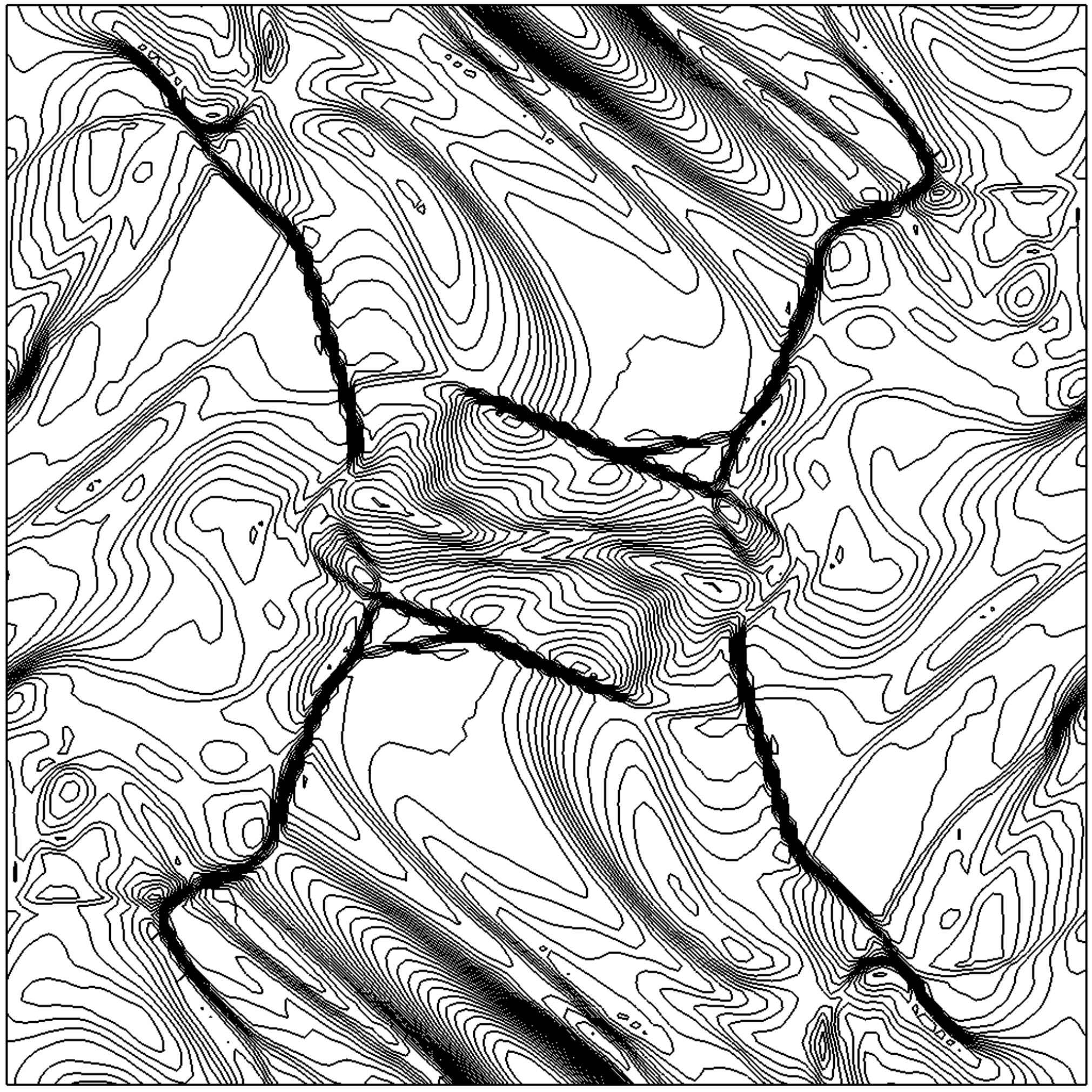} &
\includegraphics[width=0.32\textwidth]{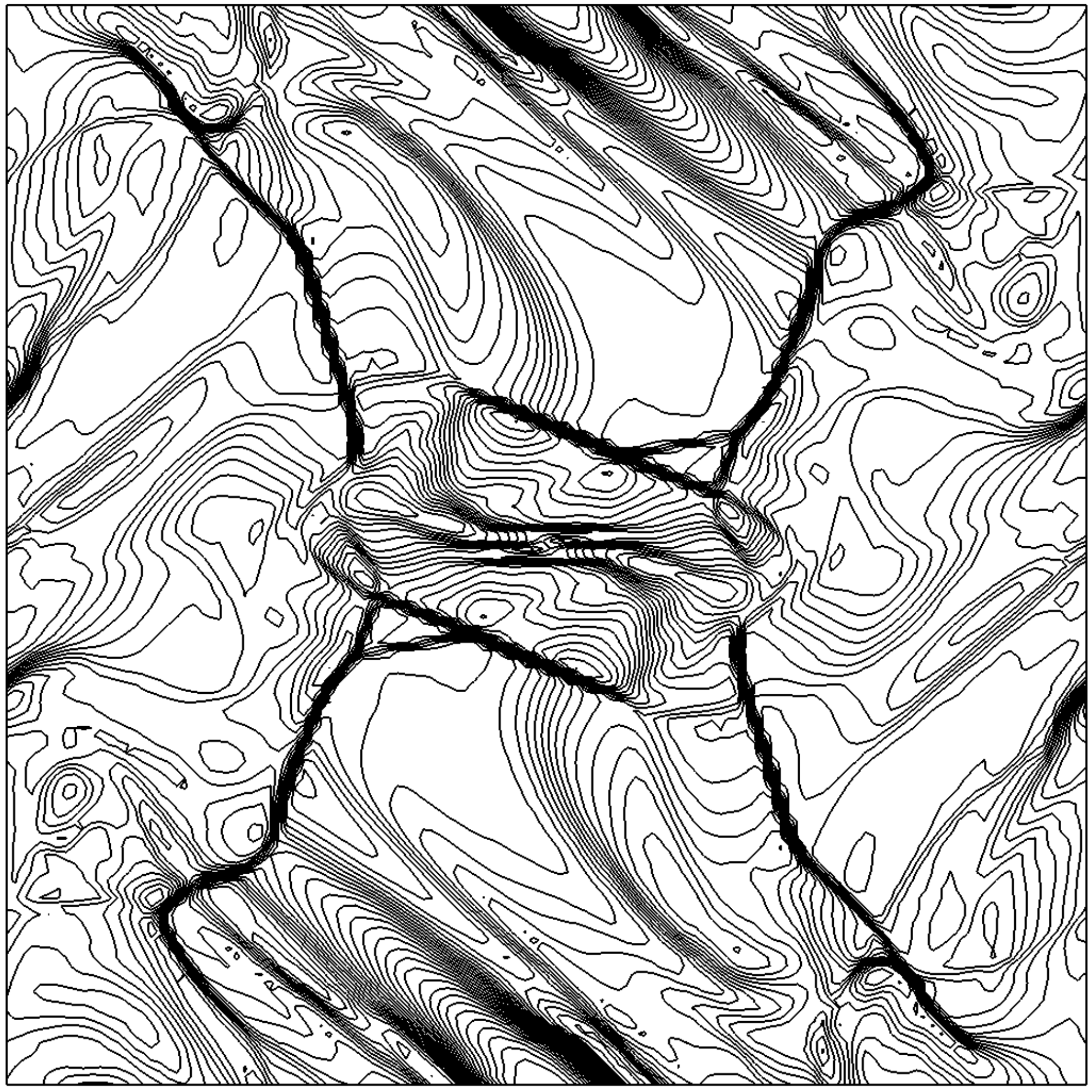} \\
\includegraphics[width=0.32\textwidth]{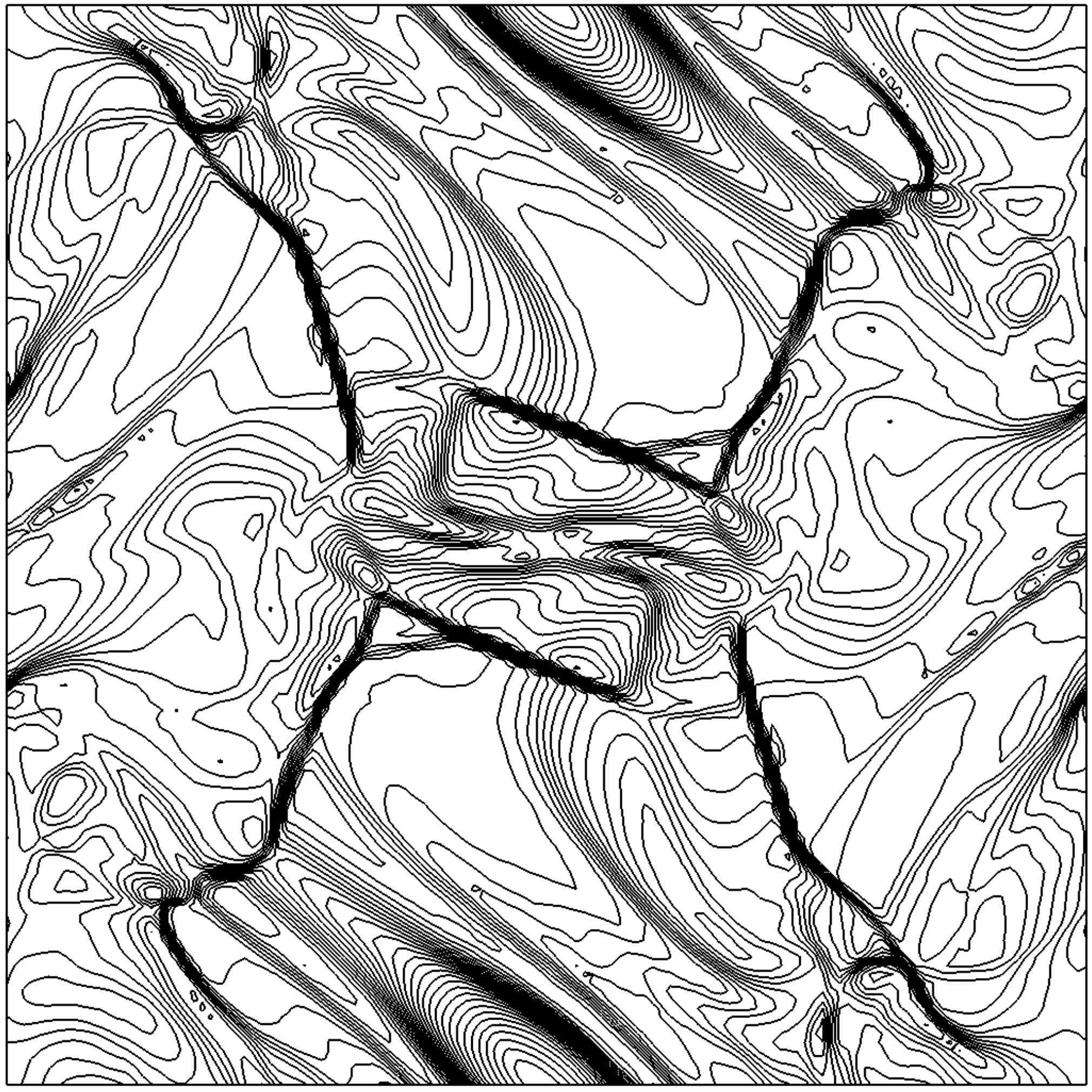} &
\includegraphics[width=0.32\textwidth]{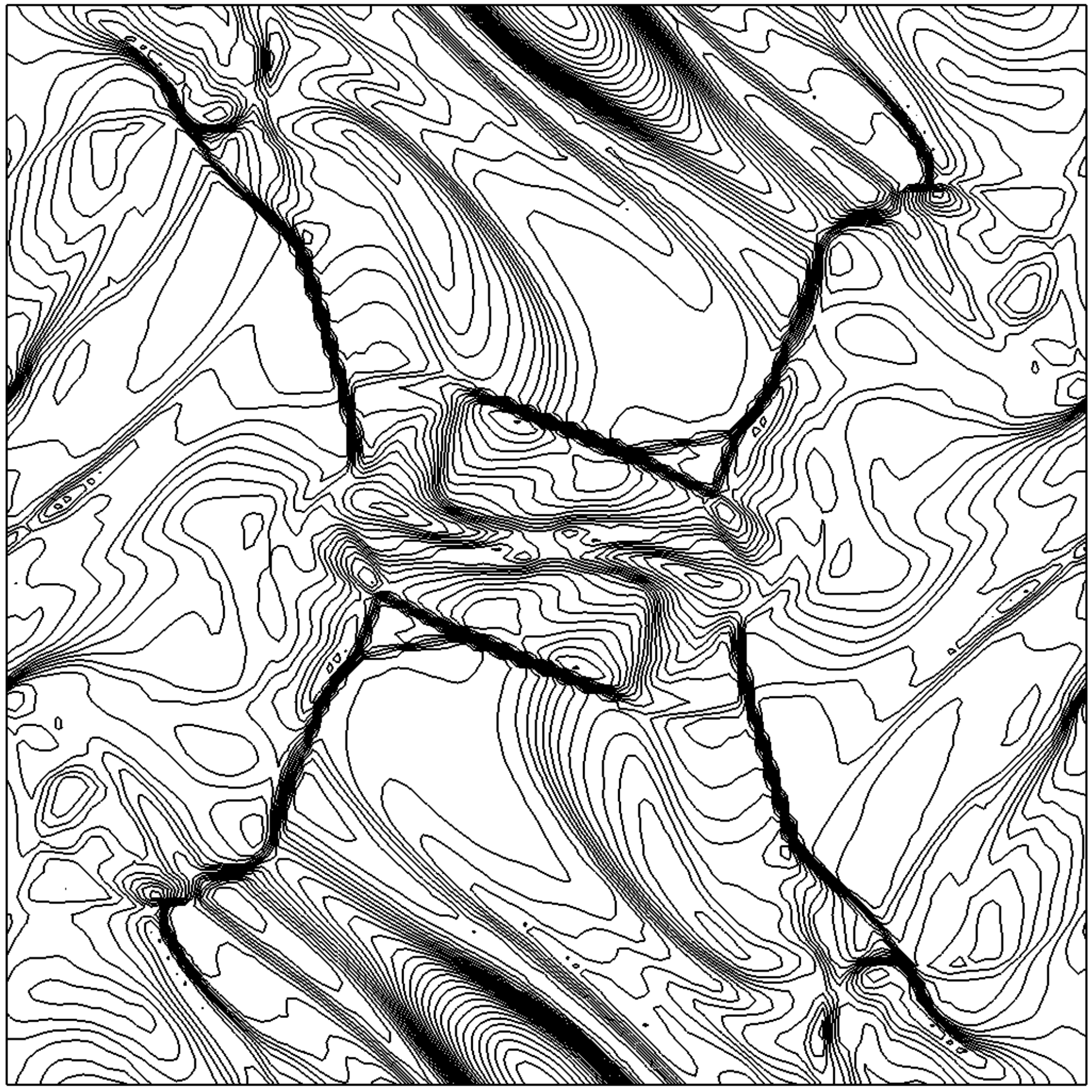} &
\includegraphics[width=0.32\textwidth]{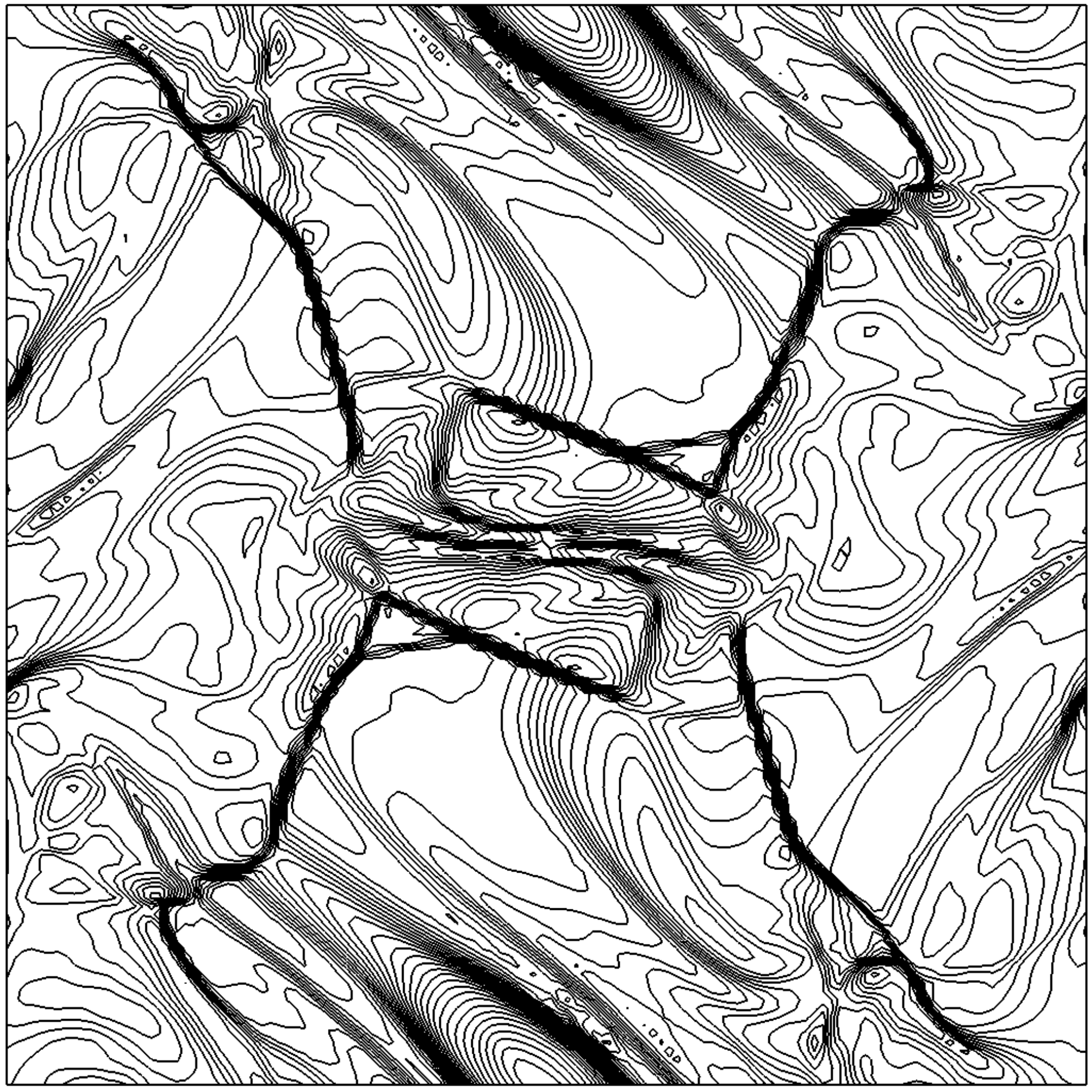} \\
LxF & HLL & HLLC \\
\end{tabular}
\caption{ Orszag-Tang test using LxF, HLL, HLLC fluxes on  $128\times 128$ mesh. 30 density contours in $(0.08,0.5)$. Top row: $k=1$, middle row: $k=2$, bottom row: $k=3$}
\label{fig:ot1}
\end{center}
\end{figure}

\begin{figure}
\begin{center}
\begin{tabular}{ccc}
\includegraphics[width=0.32\textwidth]{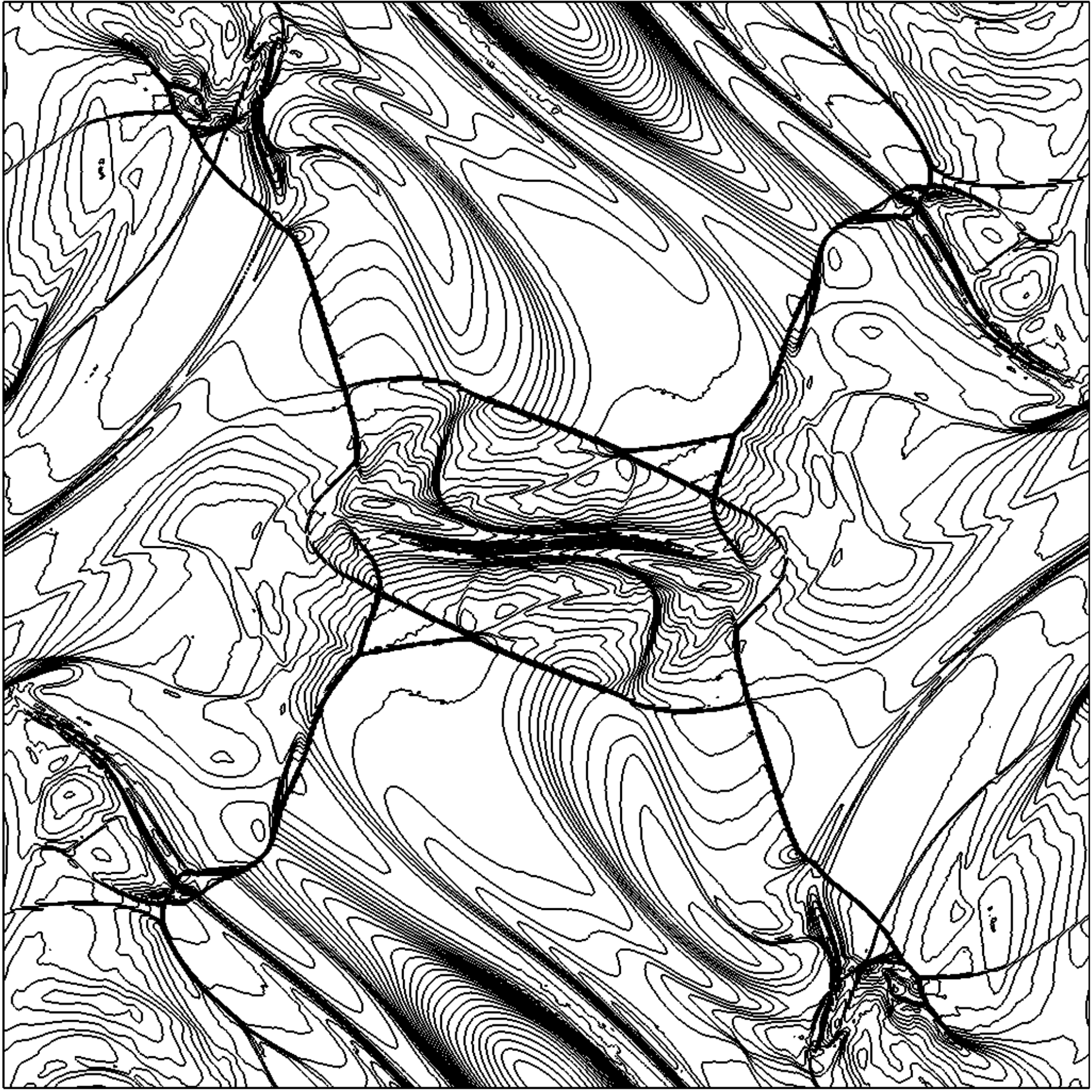} &
\includegraphics[width=0.32\textwidth]{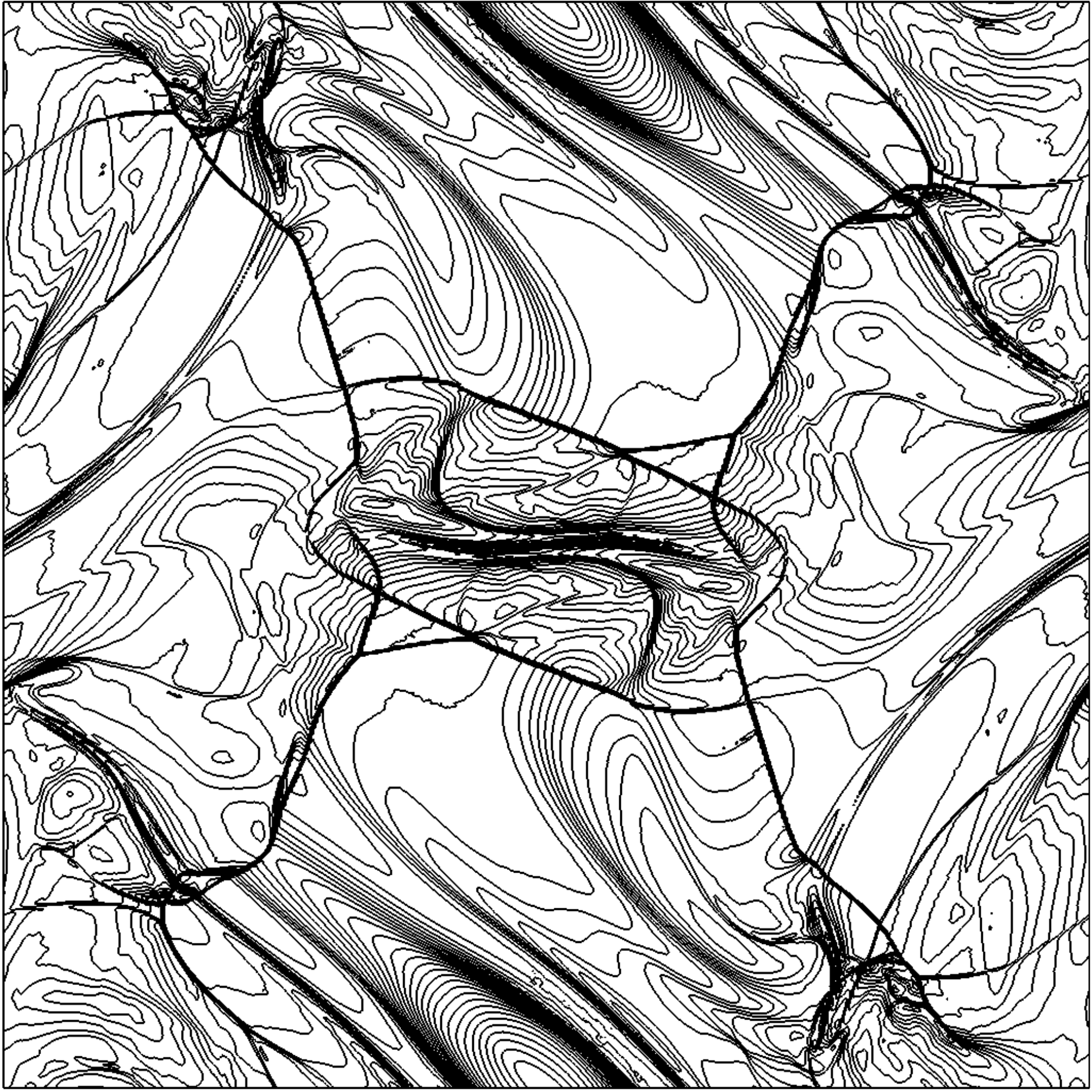} &
\includegraphics[width=0.32\textwidth]{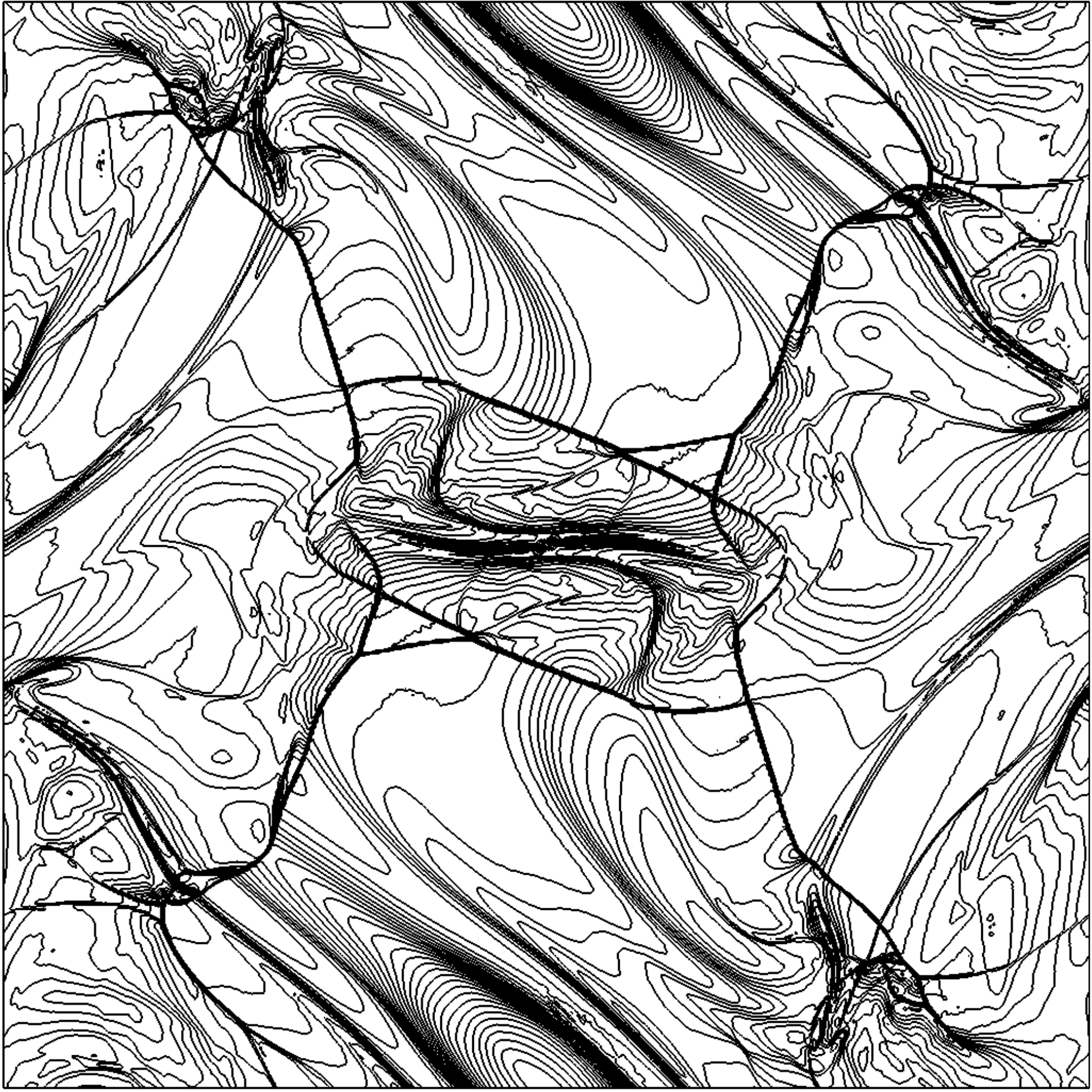} \\
LxF & HLL & HLLC \\
\end{tabular}
\caption{ Orszag-Tang test using LxF, HLL, HLLC fluxes on  $512\times 512$ mesh and degree $k=3$. 30 density contours in the interval $(0.08,0.5)$.}
\label{fig:ot2}
\end{center}
\end{figure}

\begin{figure}
\begin{center}
\begin{tabular}{ccc}
\includegraphics[width=0.32\textwidth]{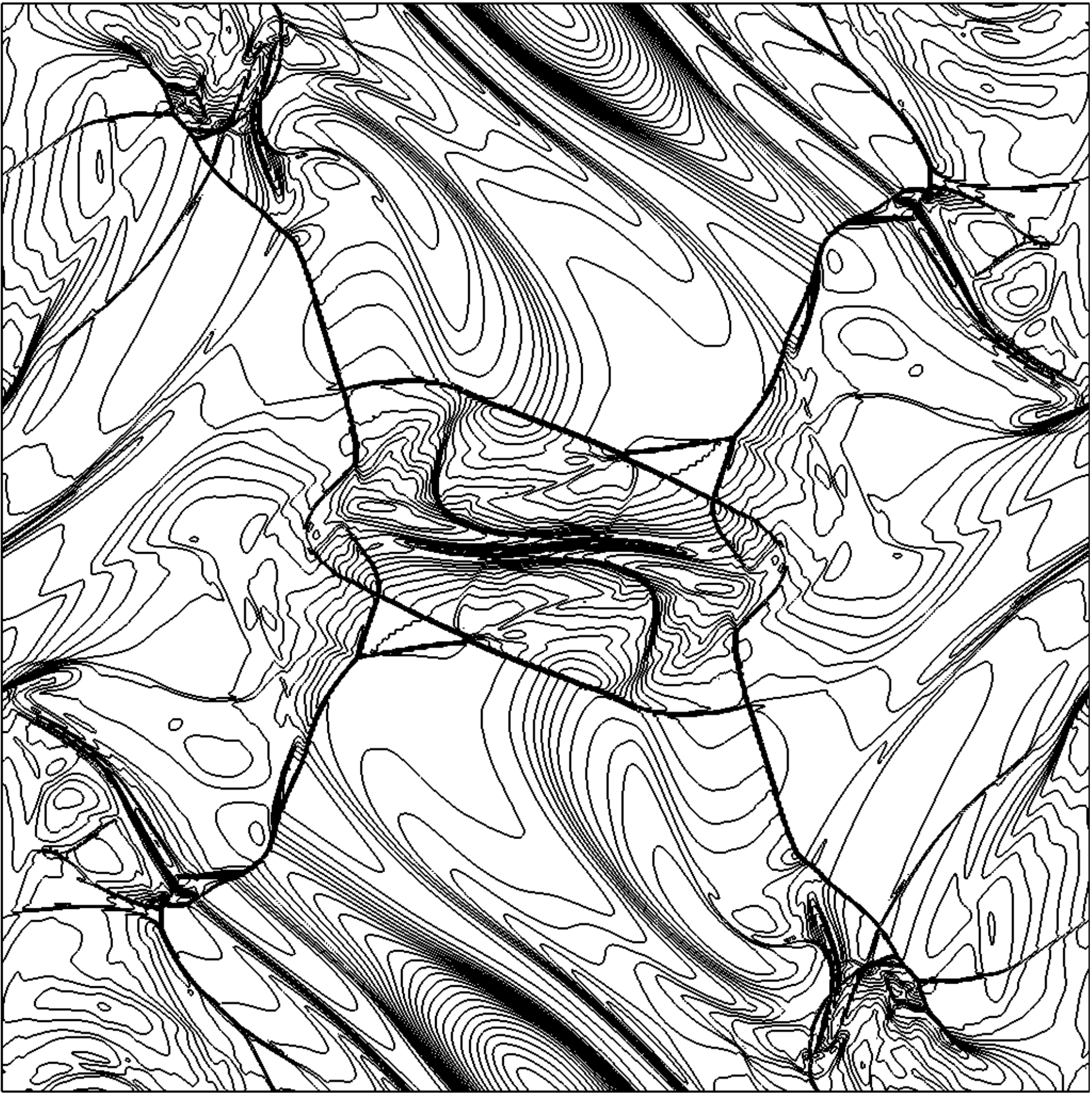} &
\includegraphics[width=0.32\textwidth]{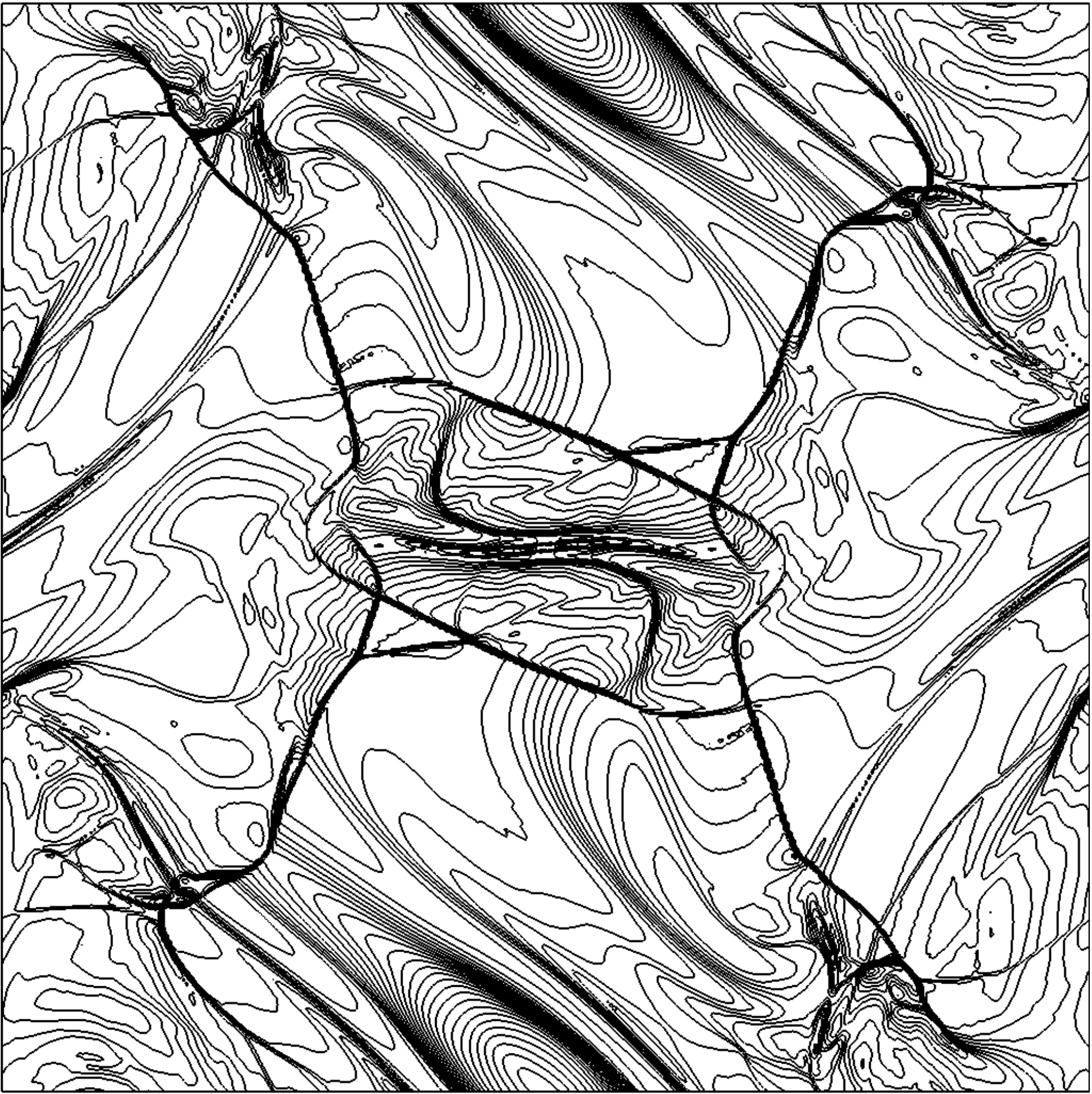} &
\includegraphics[width=0.32\textwidth]{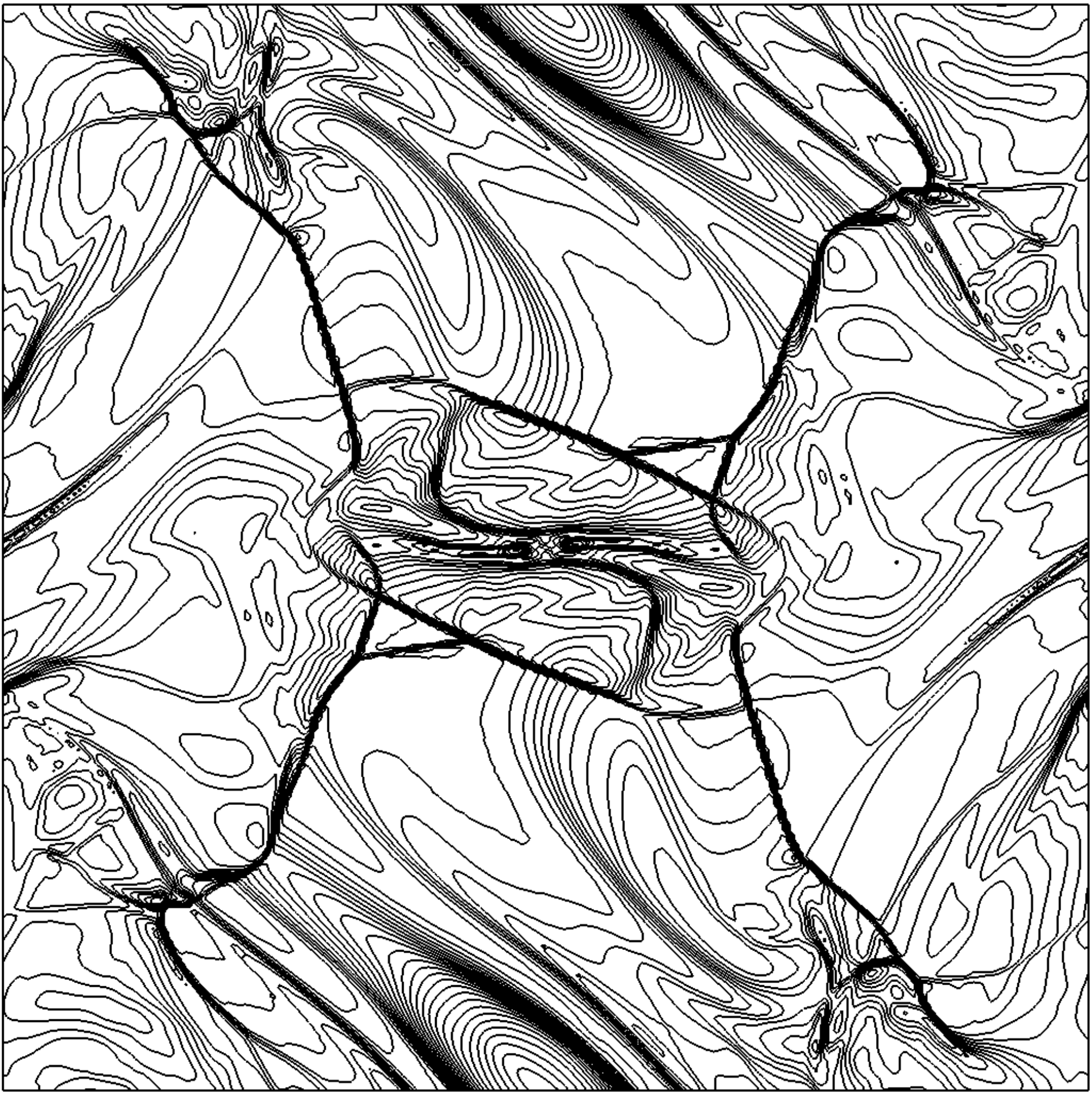} \\
(a)  & (b) & (c) \\
\end{tabular}
\caption{Orszag-Tang test using HLL flux at time $0.5$ units.  (a) $k=1$, $512\times512$ cells, (b) $k=2$, $342\times342$ cells, (c) $k=3$, $256\times 256$ cells. 30 density contours in the interval $(0.08,0.5)$. The degree and mesh size are chosen so that all three cases have nearly the same number of degrees of freedom.}
\label{fig:ot3}
\end{center}
\end{figure}

\begin{figure}
\begin{center}
\begin{tabular}{ccc}
\includegraphics[width=0.32\textwidth]{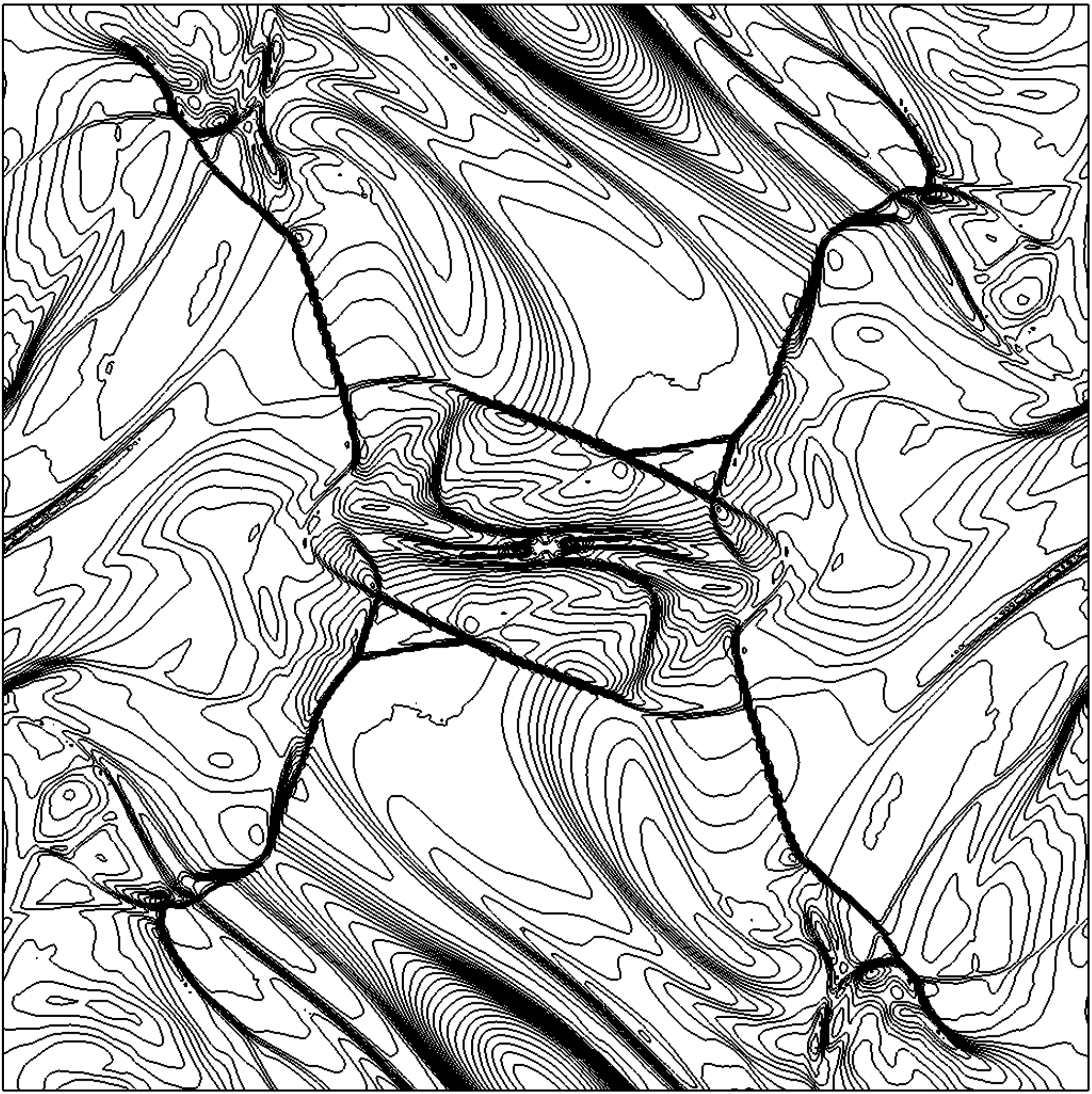} &
\includegraphics[width=0.32\textwidth]{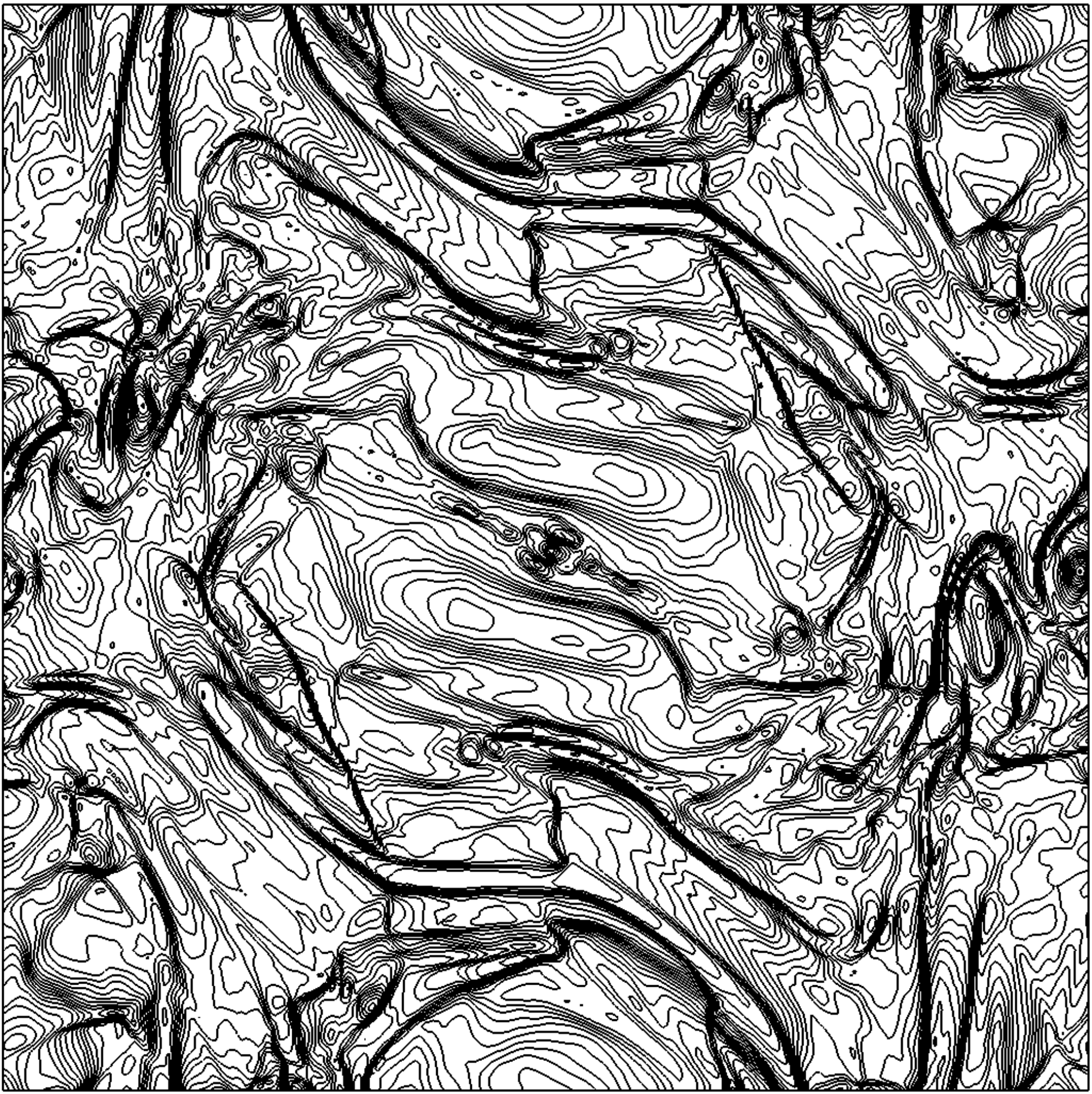} &
\includegraphics[width=0.32\textwidth]{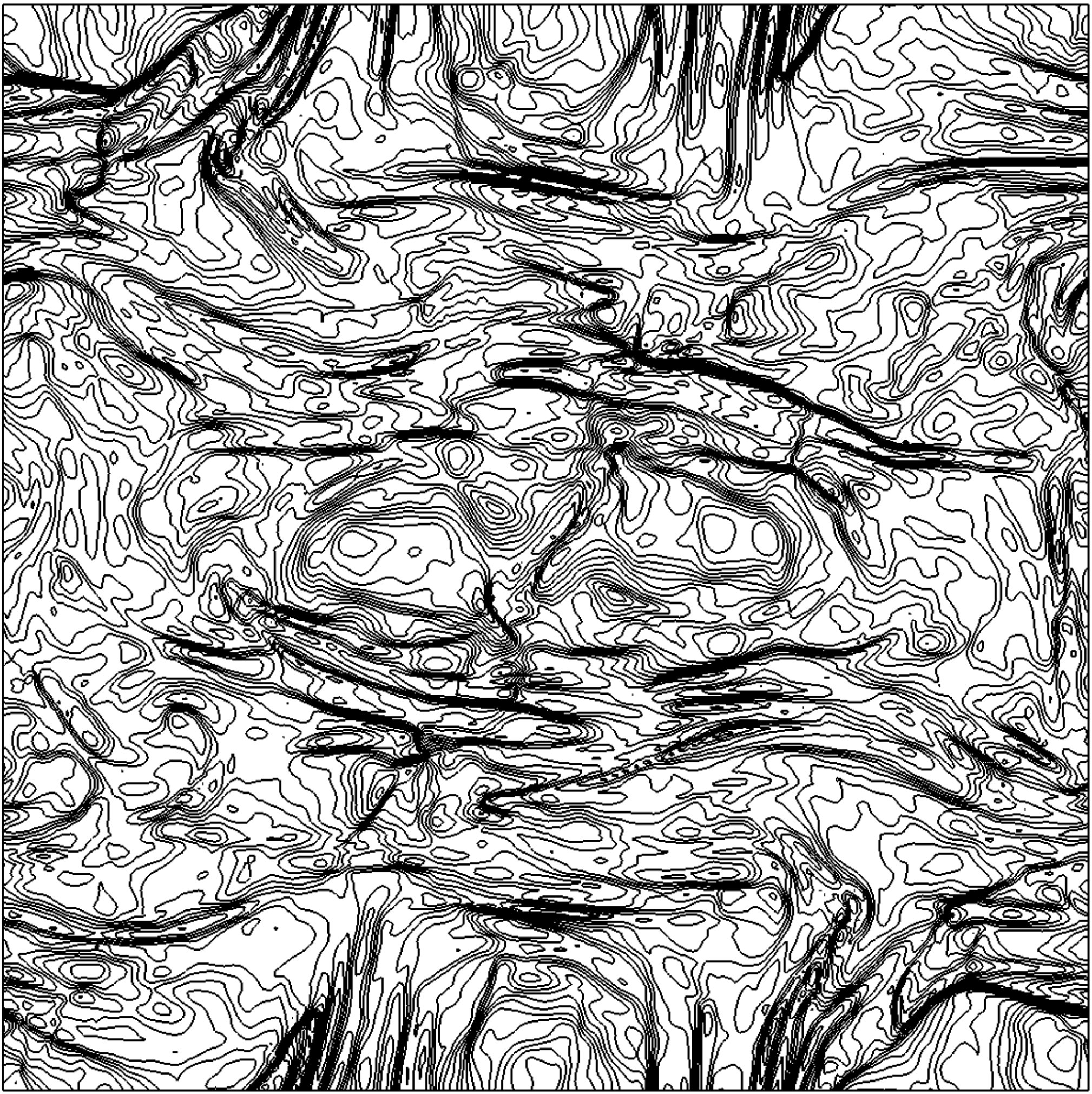} \\
$t=0.5$ & $t=1$ & $t=2$ \\
\includegraphics[width=0.32\textwidth]{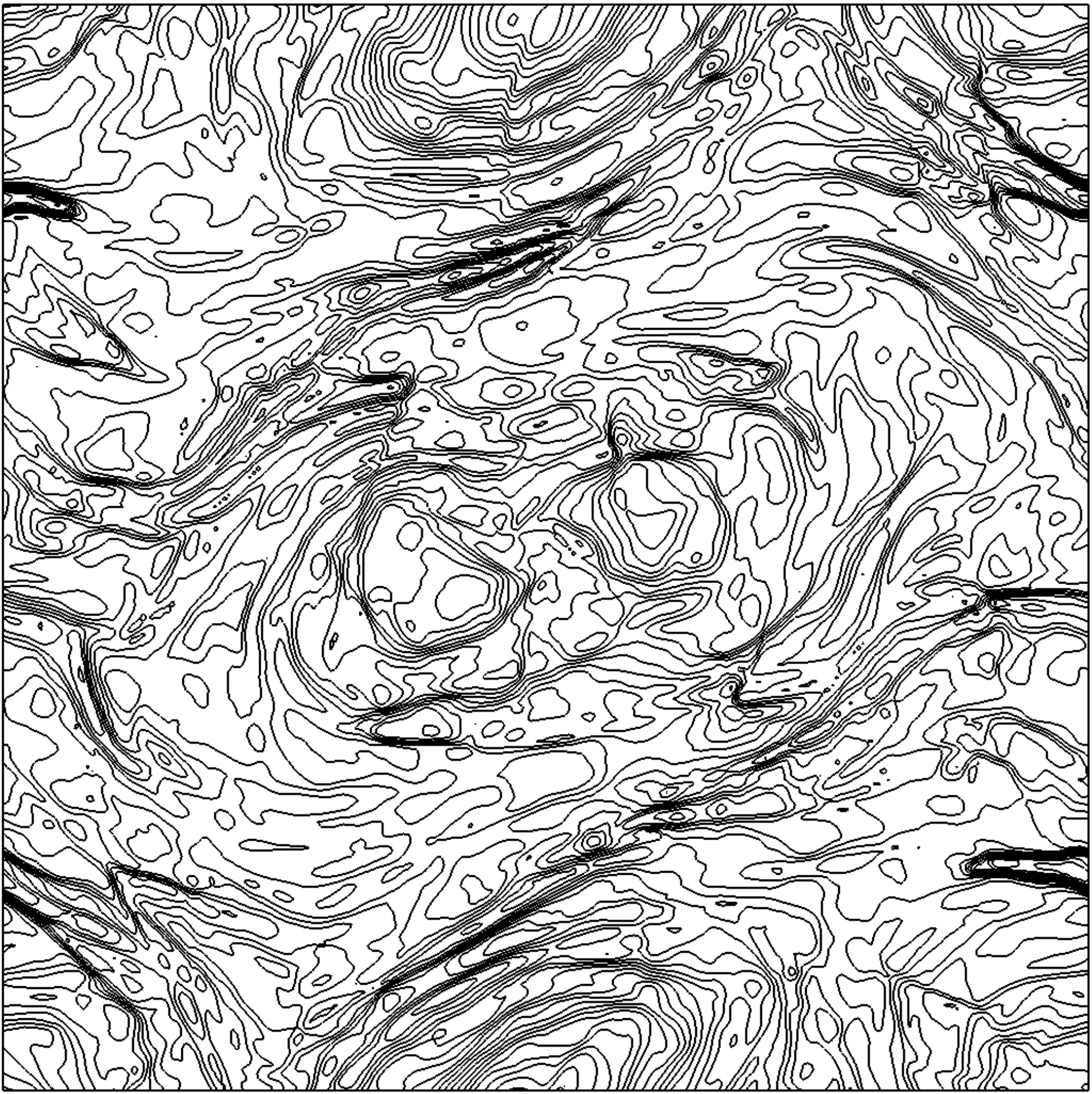} &
\includegraphics[width=0.32\textwidth]{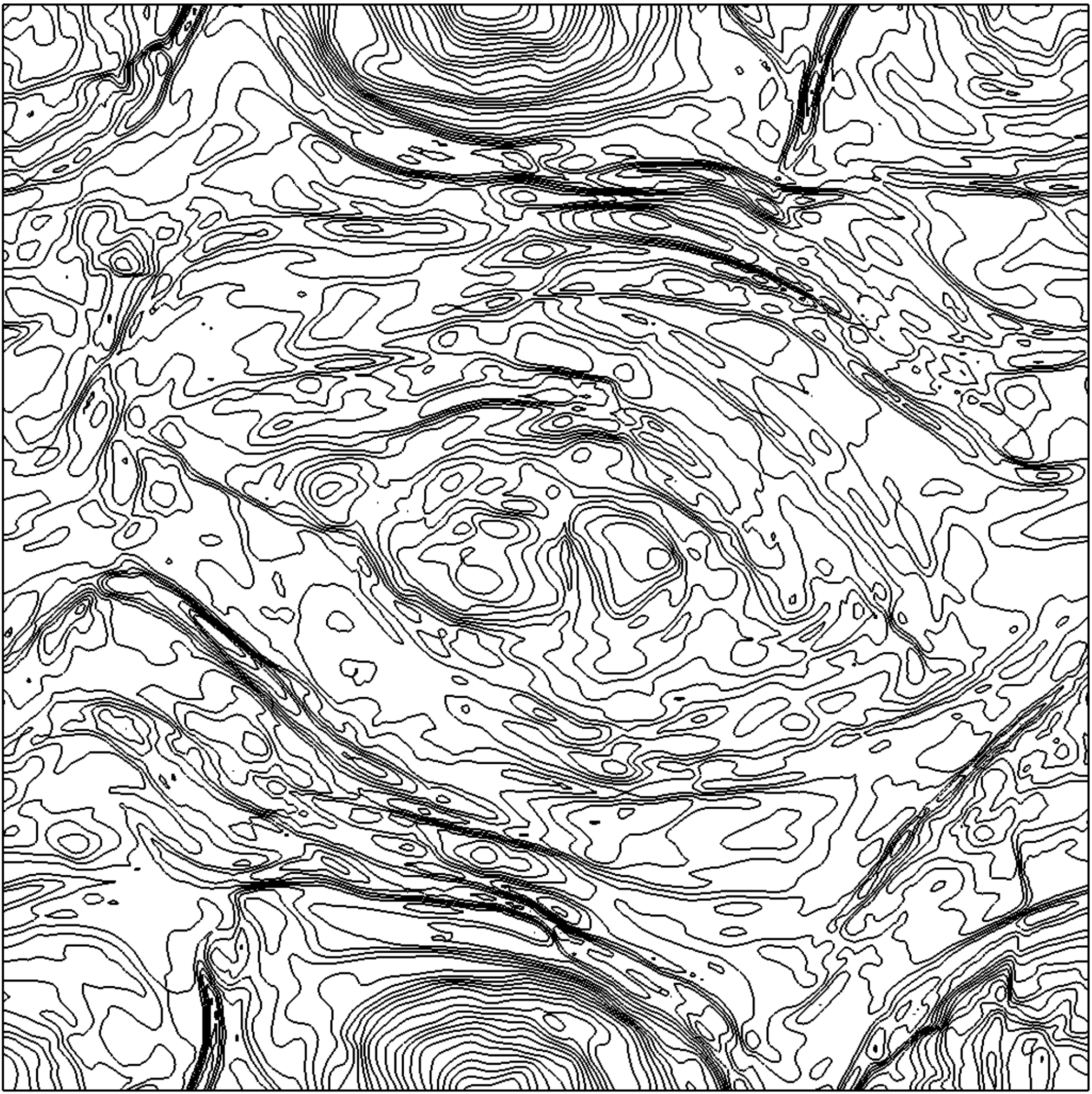} &
\includegraphics[width=0.32\textwidth]{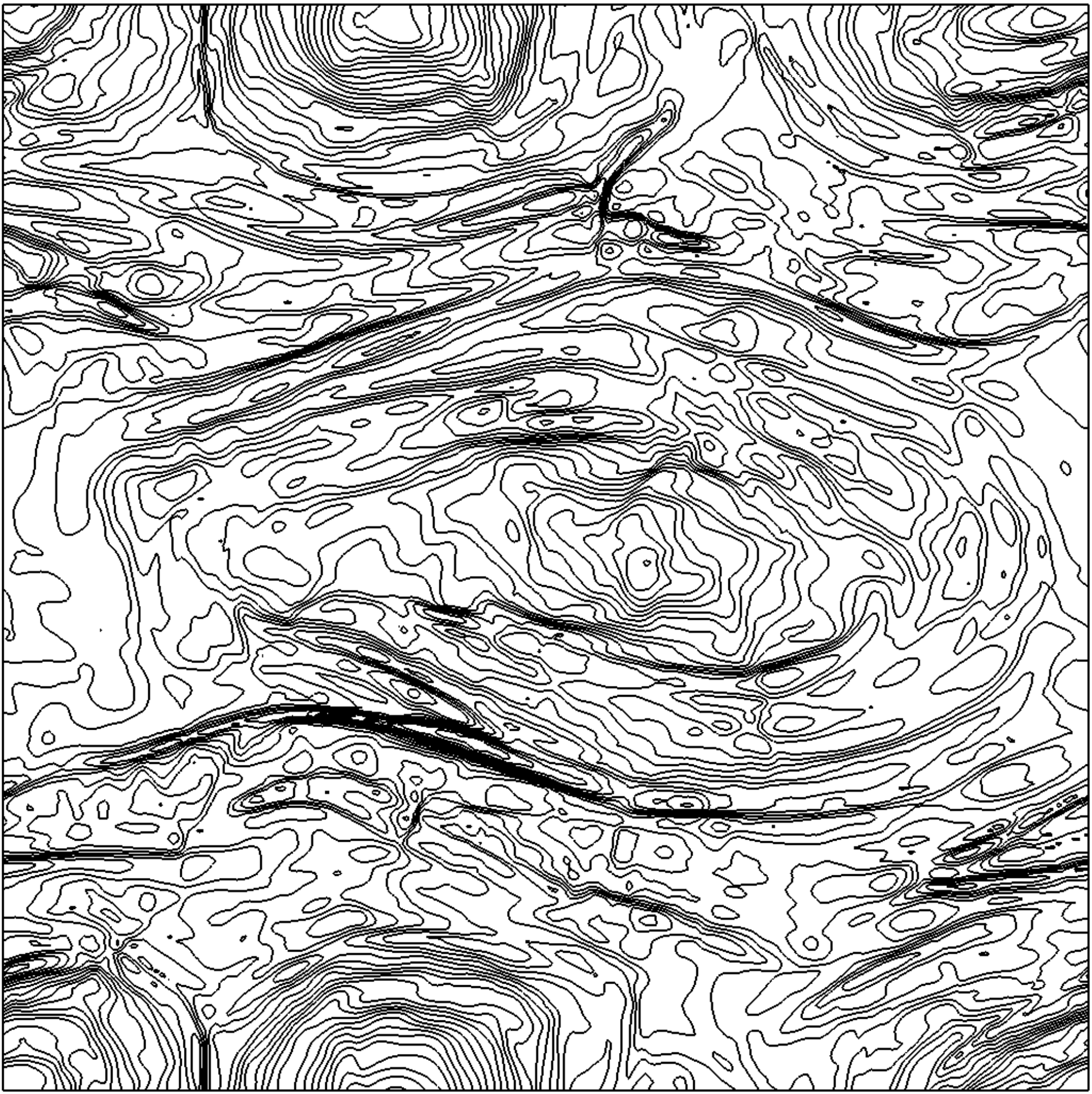}\\
$t=3$ & $t=4$ & $t=5$ \\
\end{tabular}
\caption{Long time simulation of Orszag-Tang test using HLL flux, $k=3$ on  $256\times 256$ mesh. Contours of density are shown.}
\label{fig:ot4}
\end{center}
\end{figure}

\begin{figure}
\begin{center}
\begin{tabular}{cc}
\includegraphics[width=0.48\textwidth]{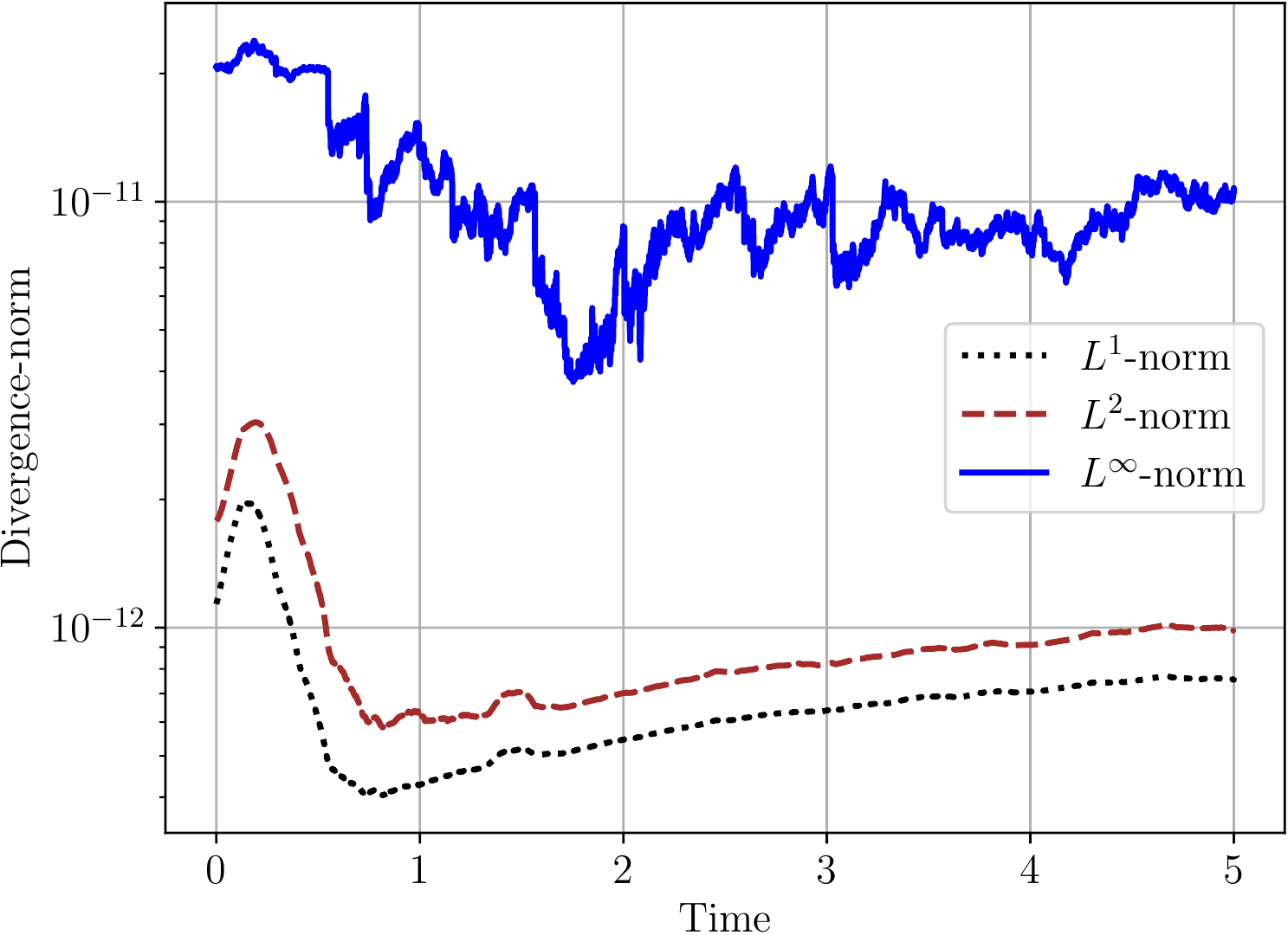} &
\includegraphics[width=0.48\textwidth]{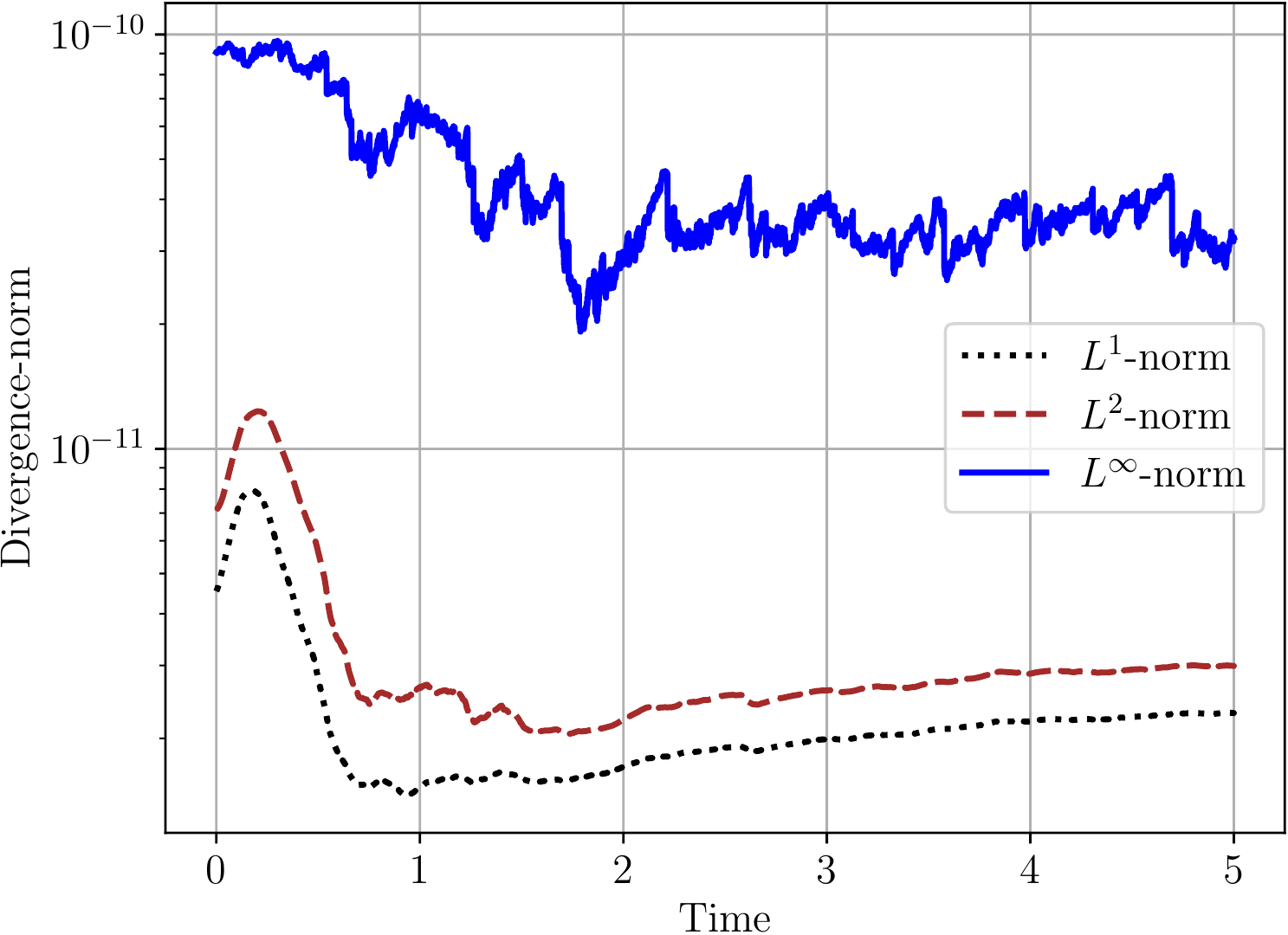} \\
(a) & (b)
\end{tabular}
\caption{Divergence norm of Orszag-Tang test using HLL flux, $k=3$ on  (a) $128 \times 128$ mesh, (b) $256\times 256$ mesh.}
\label{fig:ot5}
\end{center}
\end{figure}

This test case is first proposed in~\cite{Orszag1979} and used as a benchmark test case for many numerical algorithms for MHD. We initialize the problem with smooth initial data which later leads to the formation of more complex flow having many discontinuities as the non-linear system evolves forward in time. If the divergence error is not controlled sufficiently during simulations then numerical schemes may show instability~\cite{Li2005},~\cite{Li2011}. Even higher order local divergence free DG schemes also shows instability with time \cite{Li2005}.  The initial condition is given by
\[
\rho = \frac{25}{36\pi}, \quad p=\frac{5}{12\pi}, \quad \vel= (-\sin(2\pi y), \ \sin(2\pi x), \ 0)
\]
\[
\bthree = \frac{1}{\sqrt{4\pi}}(-\sin(2\pi y), \ \sin(4\pi x), \ 0)
\]
We have performed the numerical computations over the domain $[0,1]\times [0,1]$ with periodic boundary conditions on all sides.
The numerical solutions are computed up to the time $T=0.5$. In Figure~\ref{fig:ot1}, we have compared the LxF, HLL, and HLLC flux for density variable and $k=1,2,3$ over a mesh of size $128\times 128$. We can observe from the Figure~\ref{fig:ot1} that HLLC flux resolves features more sharply in comparison to HLL and LxF flux, e.g., in the central part of the domain. The solution on a finer mesh of $512\times 512$ cells is shown in Figure~\ref{fig:ot2} and we observe that all the features are now resolved more sharply. \resb{To study the effectiveness of high order methods, we also perform computations on different resolutions by matching the number of degrees of freedom as shown in Figure~\ref{fig:ot3}. As the degree is increased from 1 to 3, the mesh size is reduced so that all three cases have approximately the same number of degrees of freedom. We can see from the figures that the case $k=3$ which has a smaller mesh size is still able to capture the solution features.} \resa{In Figure~\ref{fig:ot4}, we show results of long time simulation upto time $t=5$ units. The solution becomes turbulent at long times but the computations remain stable.} \resb{We also monitor the divergence norm as a function of time, as shown in Figure~\ref{fig:ot5}. Theoretically, the numerical scheme must preserve the divergence at all times. In practice, due to round off errors, we see that it is not exactly preserved and is not exactly zero, but the values remain small and do not increase with time.}

\subsection{Rotor test}
\begin{figure}
\begin{center}
\begin{tabular}{ccc}
\includegraphics[width=0.32\textwidth]{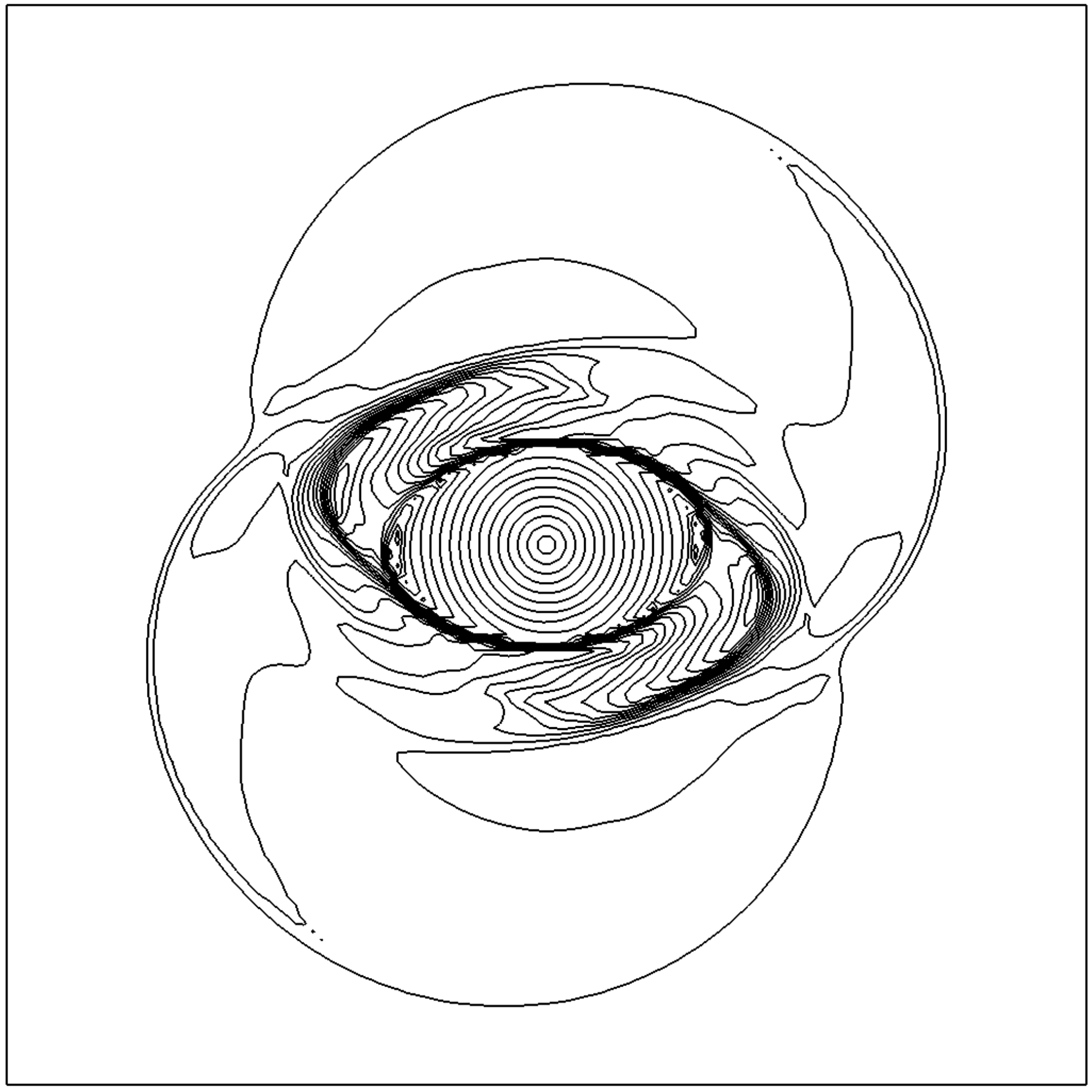} &
\includegraphics[width=0.32\textwidth]{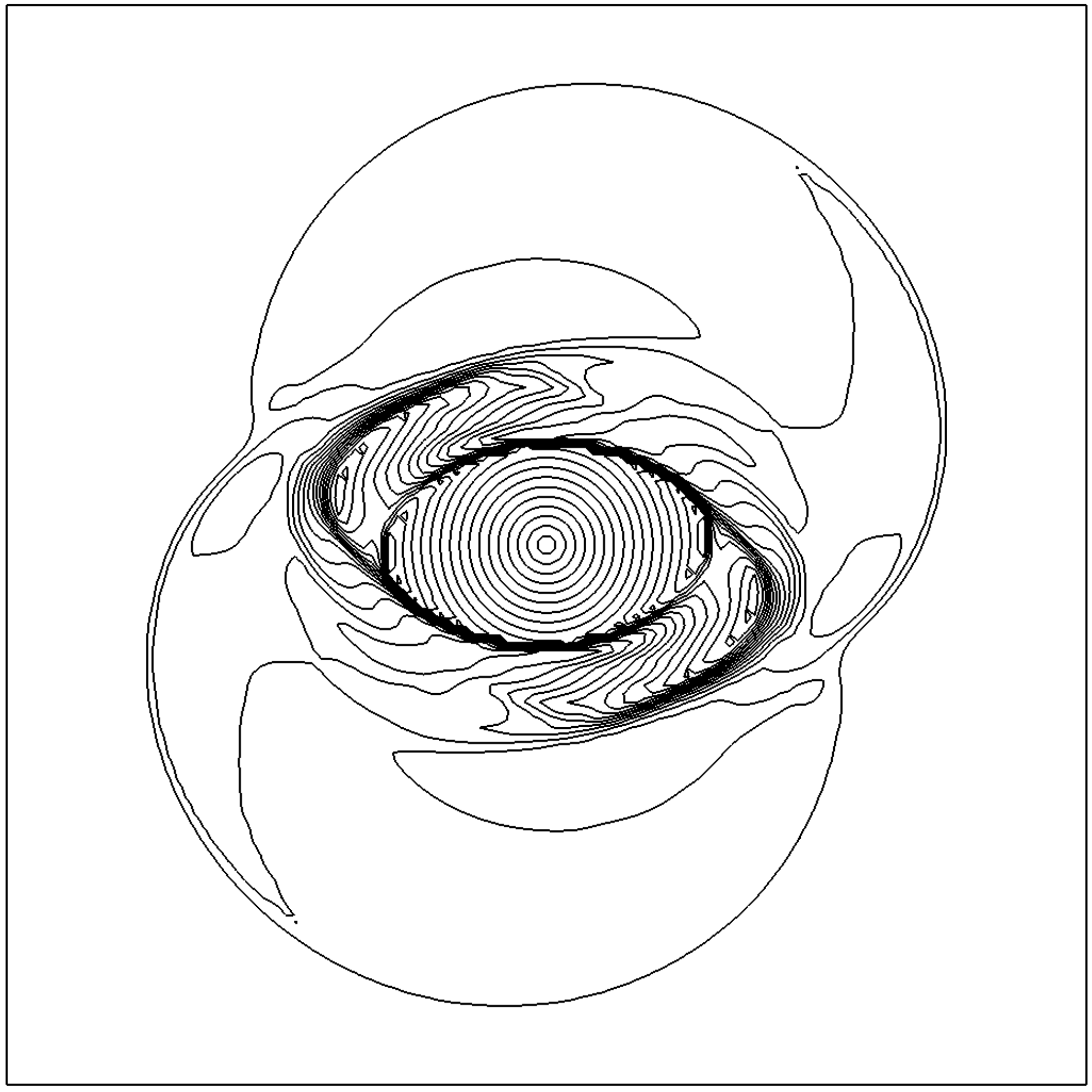} &
\includegraphics[width=0.32\textwidth]{rotor/k1_rotor_hll_128} \\
\includegraphics[width=0.32\textwidth]{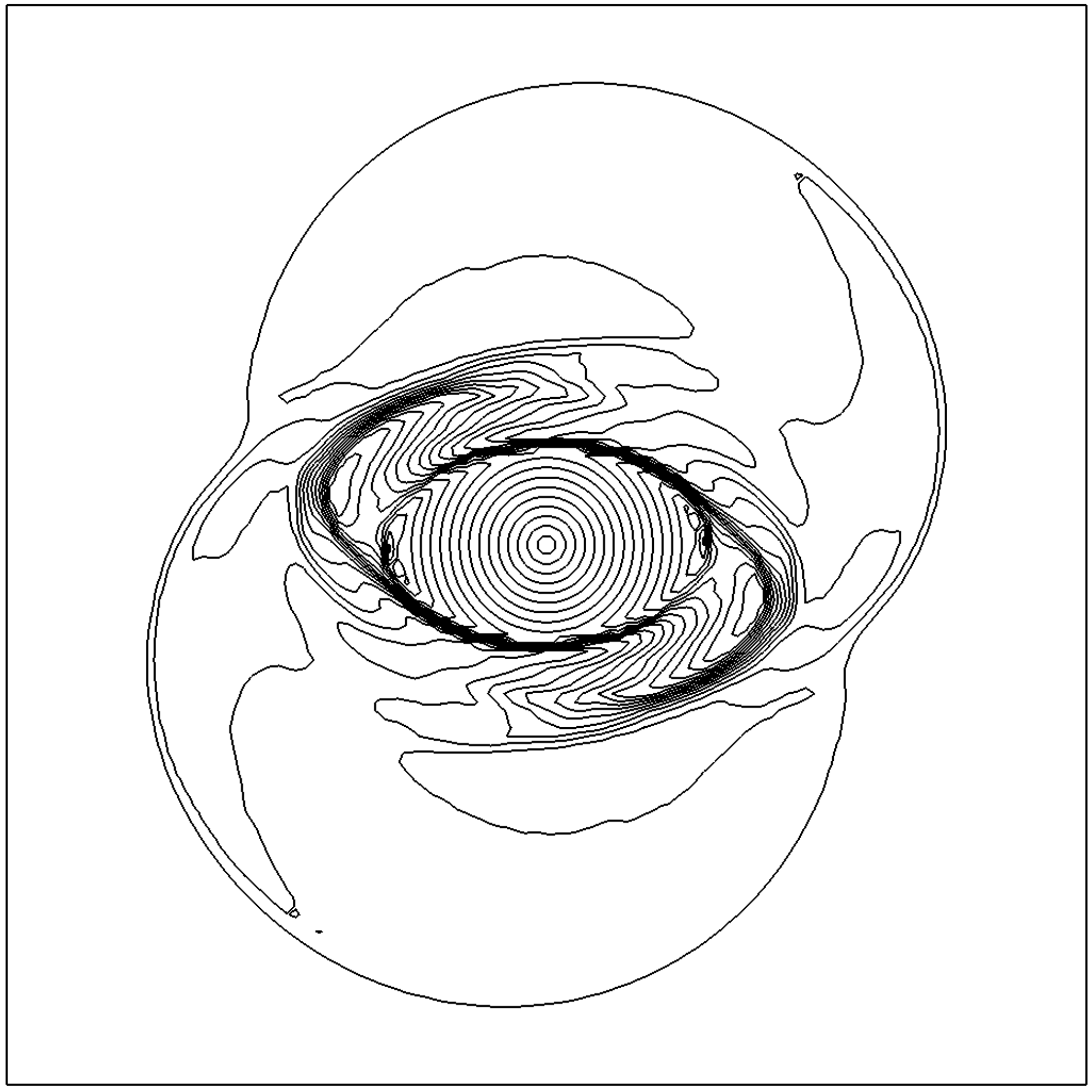} &
\includegraphics[width=0.32\textwidth]{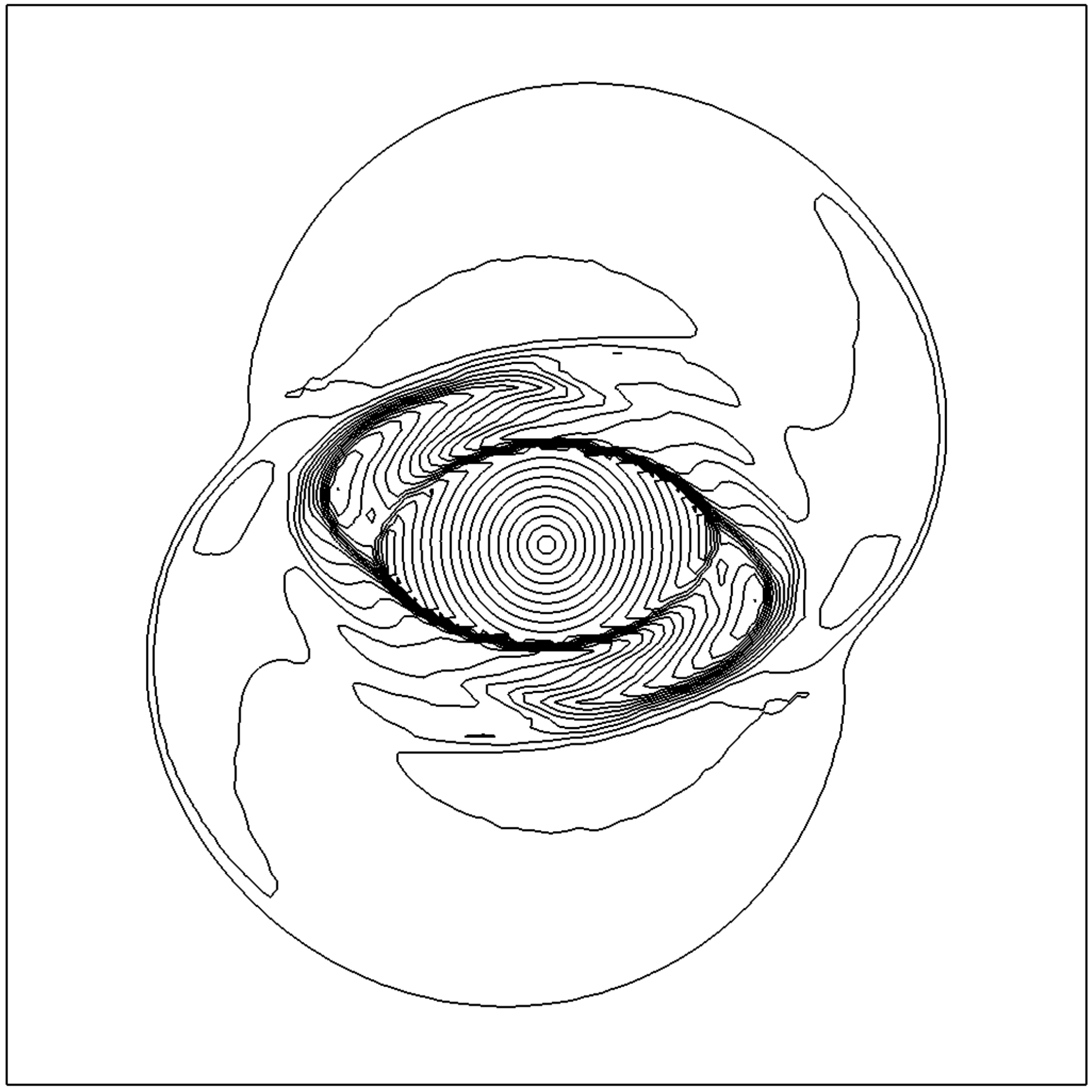} &
\includegraphics[width=0.32\textwidth]{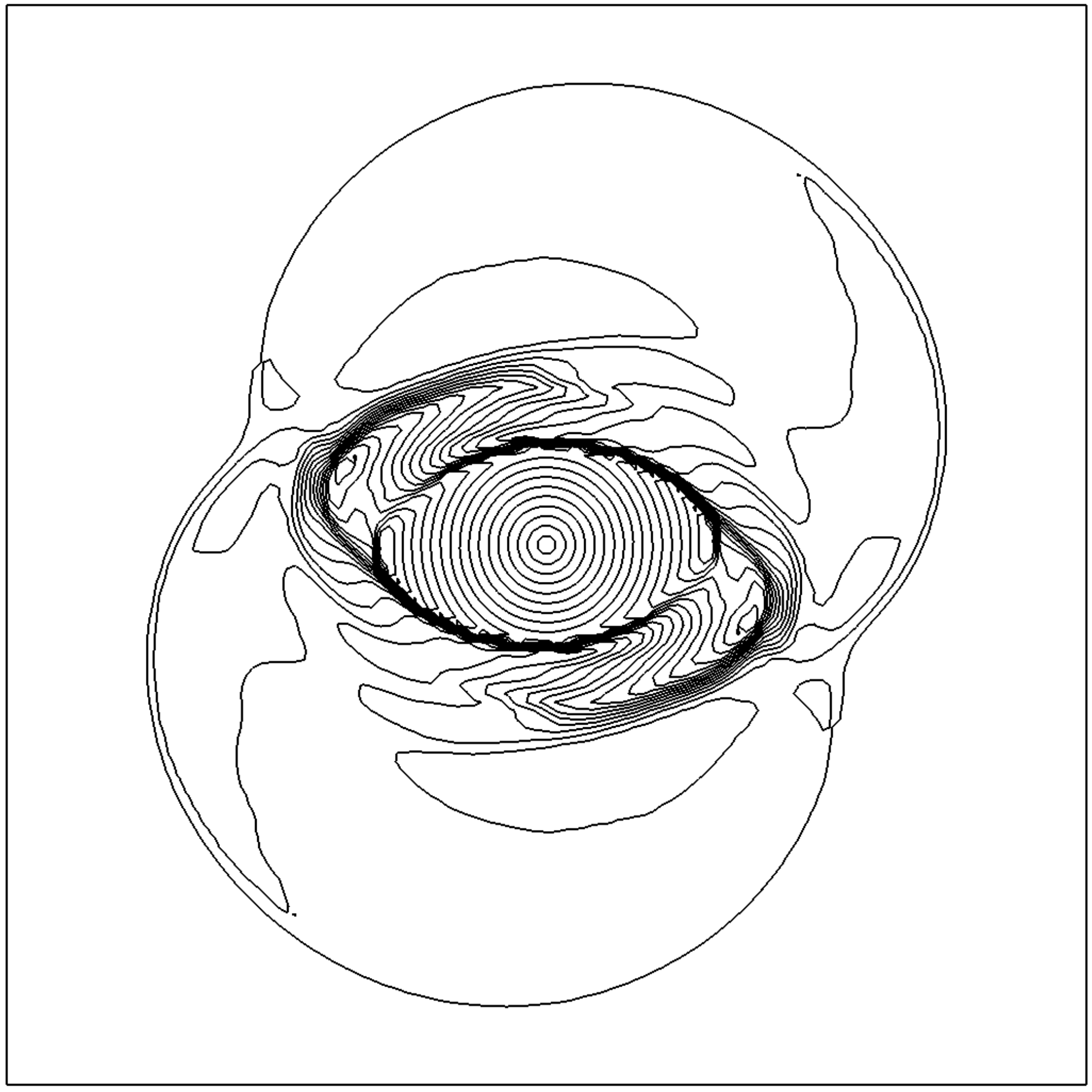} \\
\includegraphics[width=0.32\textwidth]{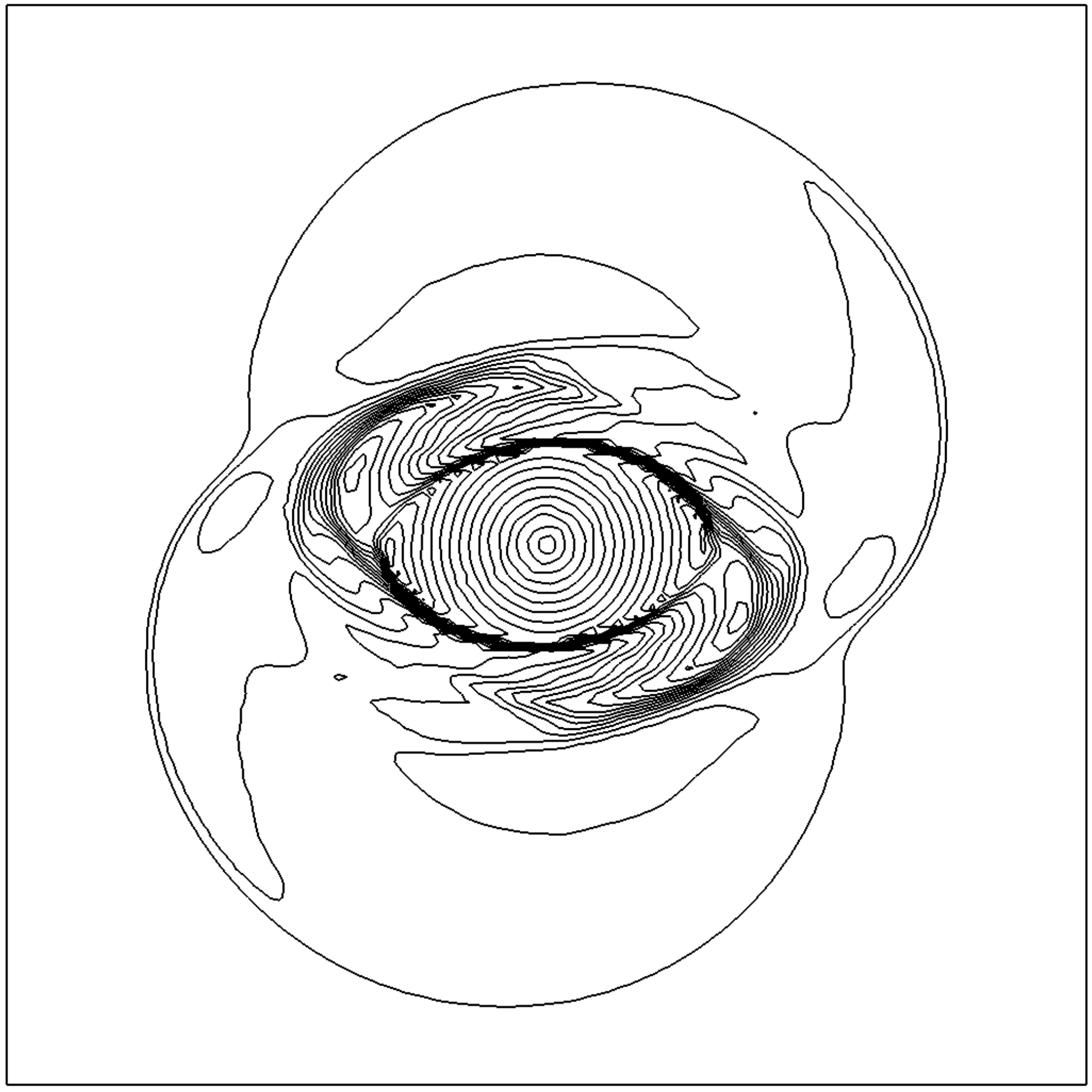} &
\includegraphics[width=0.32\textwidth]{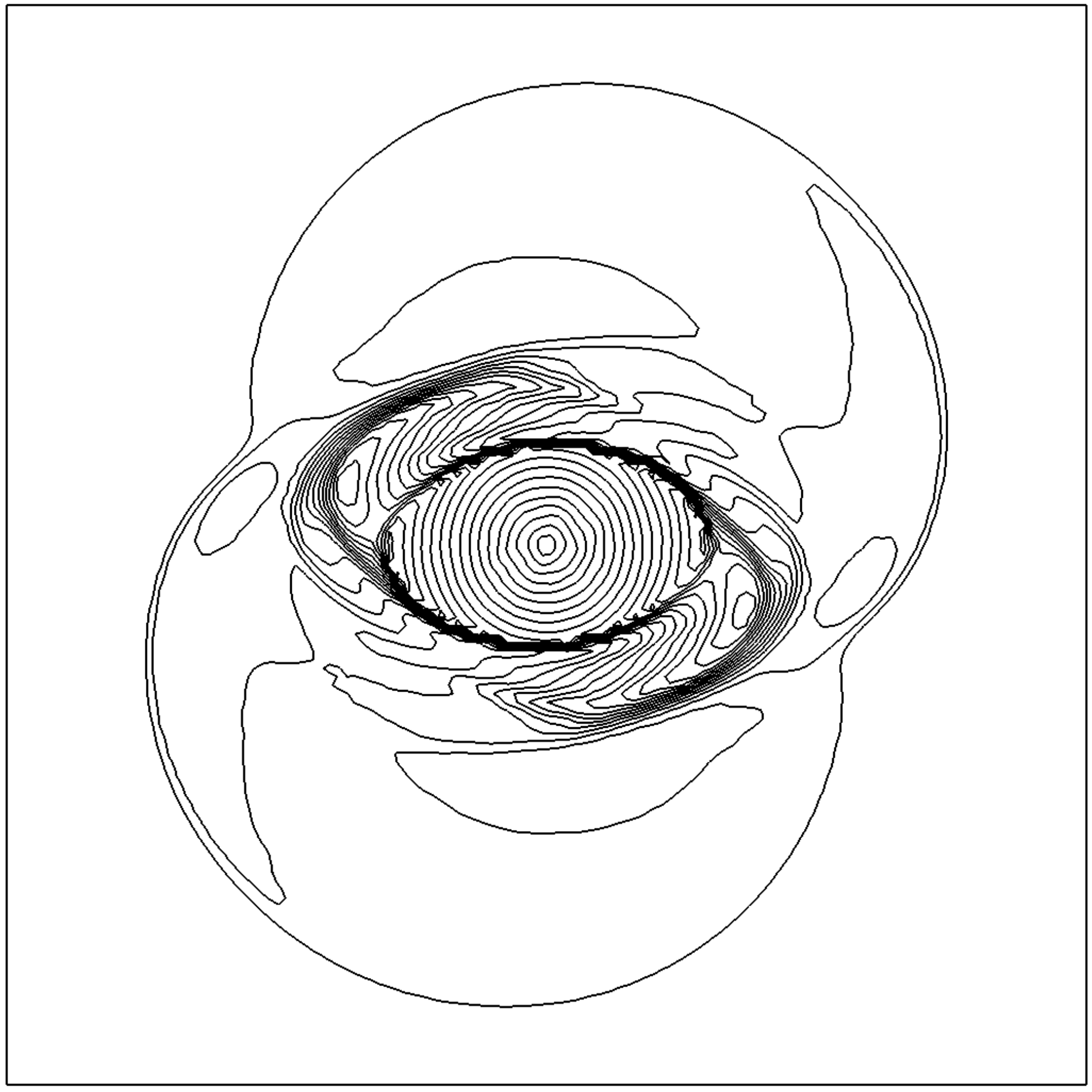} &
\includegraphics[width=0.32\textwidth]{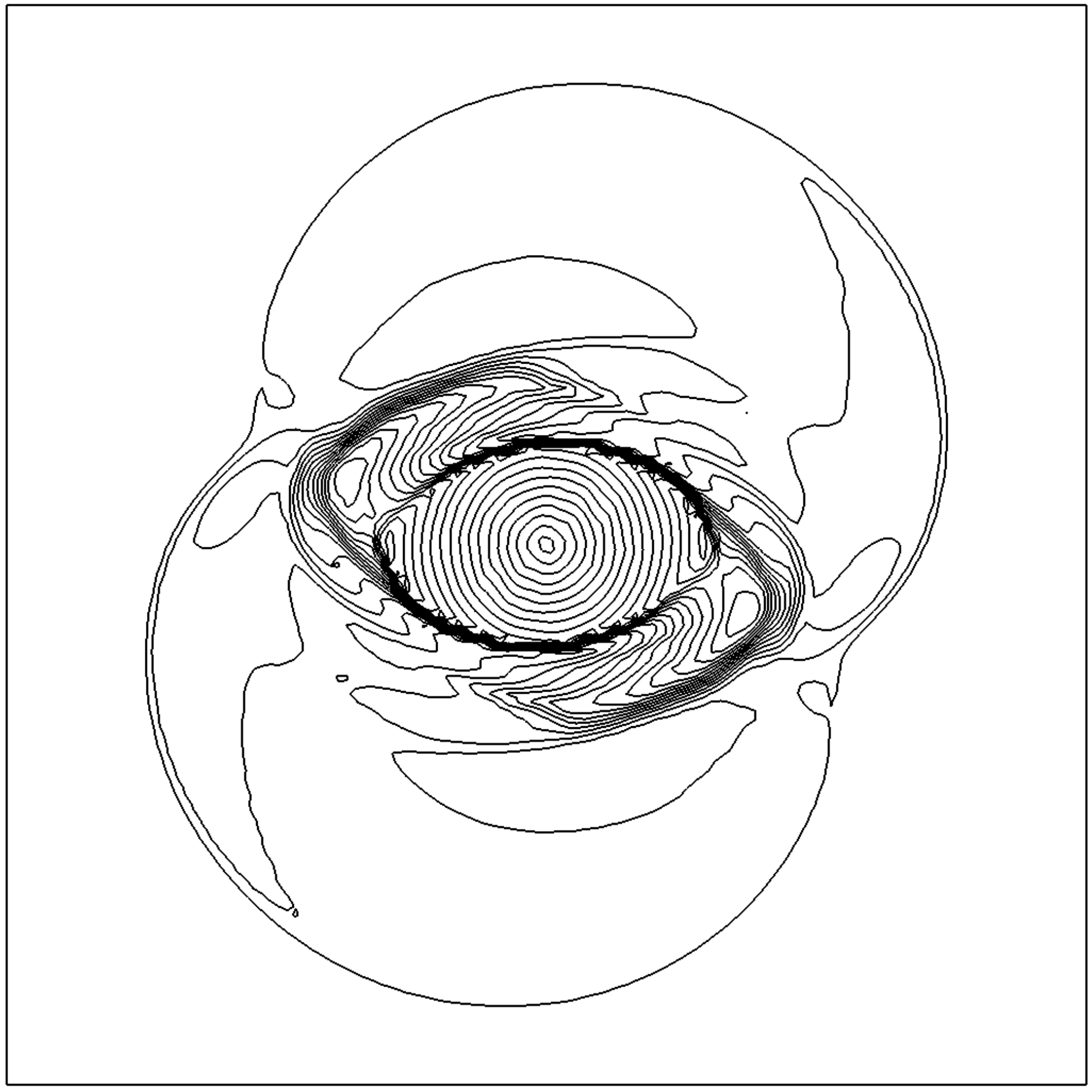} \\
LxF & HLL & HLLC \\
\end{tabular}
\caption{Rotor test using degree $k=1$ (top row), degree $k=2$   (middle row) and degree $k=3$ (bottom row) over $128 \times 128$ mesh. 20   Mach contours in $(0,4.5)$.}
\label{fig:rotor}
\end{center}
\end{figure}

\begin{figure}
\begin{center}
\begin{tabular}{ccc}
\includegraphics[width=0.32\textwidth]{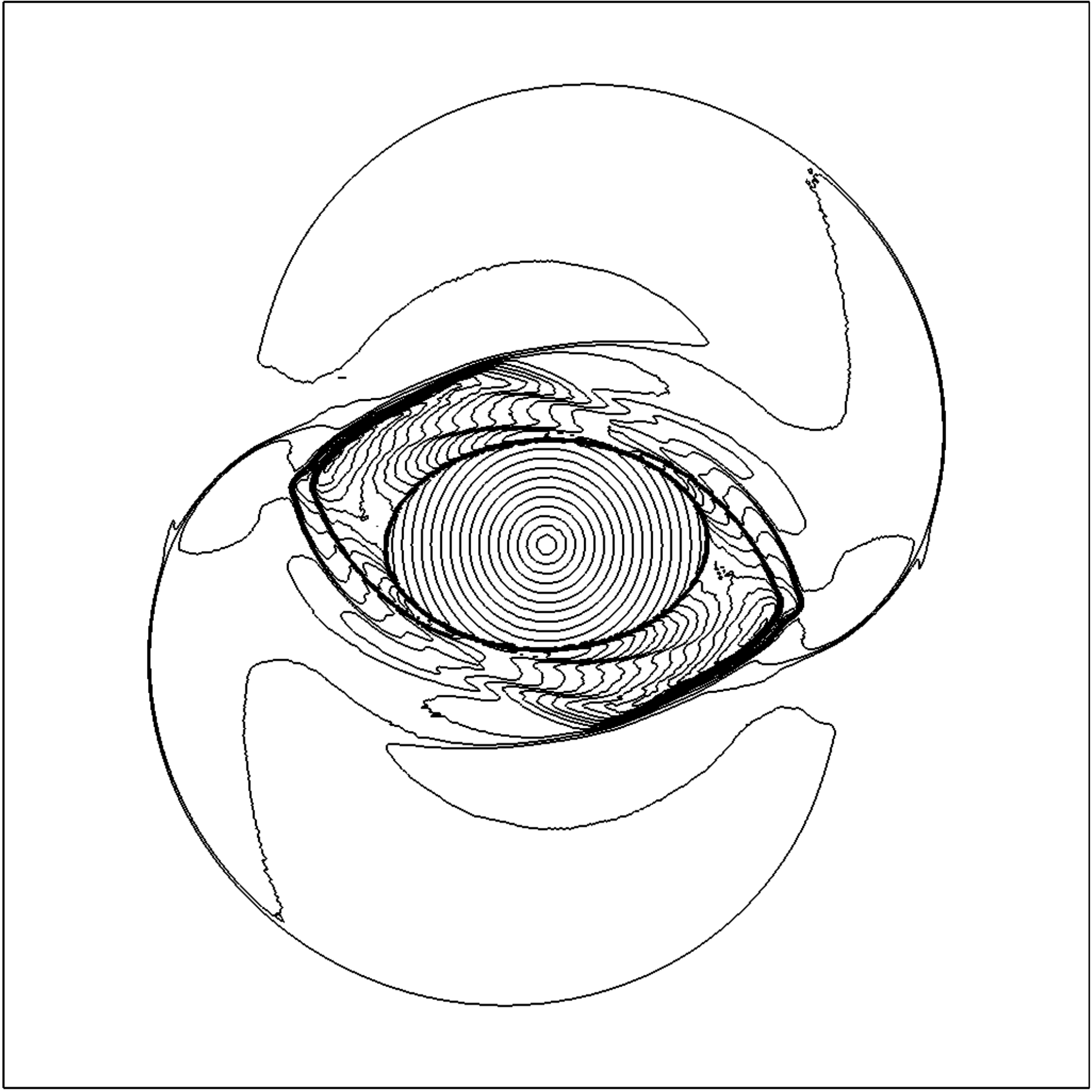} &
\includegraphics[width=0.32\textwidth]{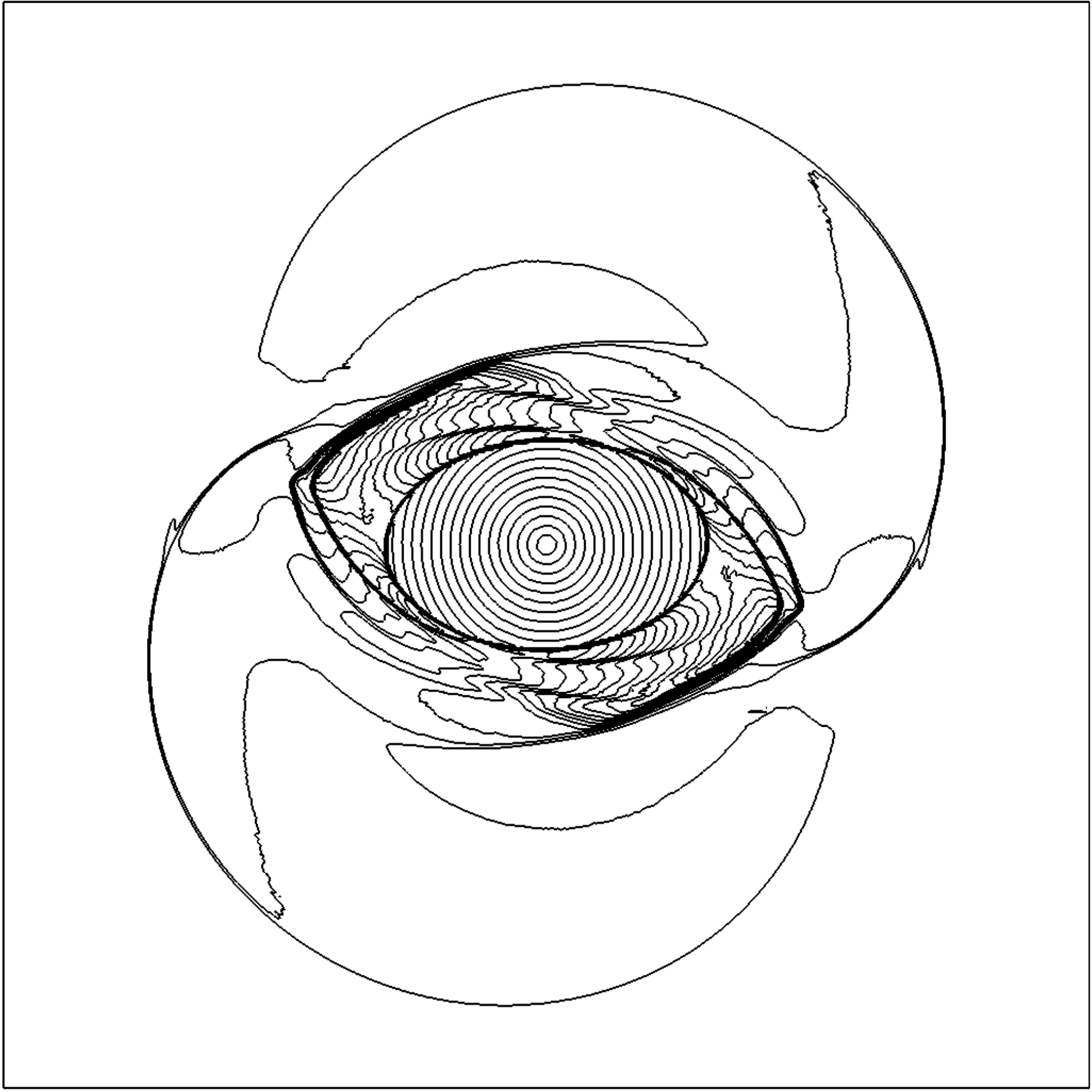} &
\includegraphics[width=0.32\textwidth]{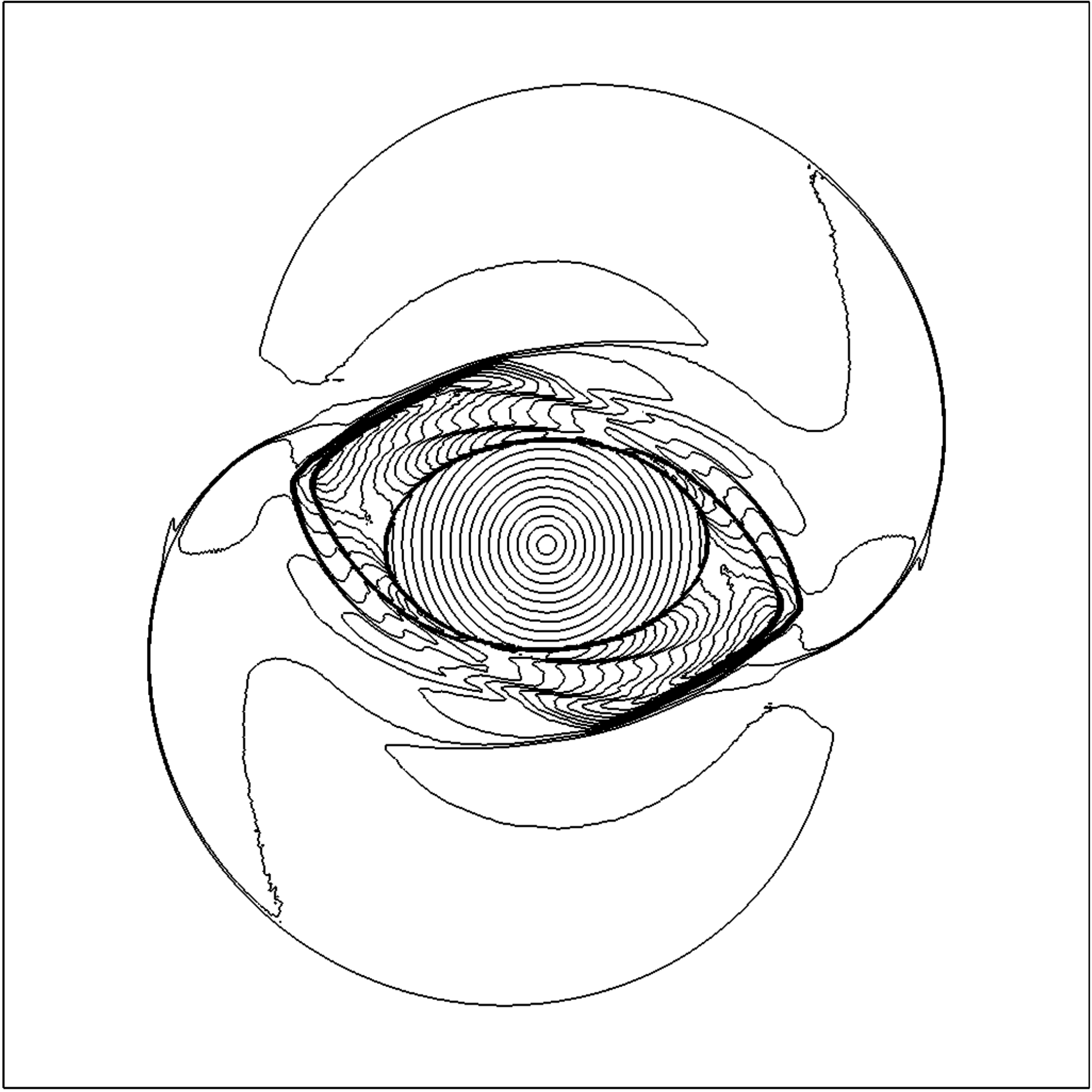} \\
LxF & HLL & HLLC \\
\end{tabular}
\caption{Rotor test using degree $k=3$  and $512 \times 512$ mesh. 20 Mach contours in the interval $(0,4.5)$.}
\label{fig:rotor1}
\end{center}
\end{figure}

This test case was first proposed in~\cite{Balsara1999}, but we use the version given in~\cite{Toth2000}. This problem describes the spinning of a dense rotating disc of fluid in the center while ambient fluid are at rest. The magnetic field wraps around the rotating dense fluid turned it into an oblate shape.  If the numerical scheme is not sufficiently control the divergence-error in the magnetic field, distortion can be observed in Mach number~\cite{Li2005}. The computational domain is $[0,1] \times [0,1]$ with periodic  boundary conditions on all sides, and the initial condition is given as follows. For $r < r_0$,
\[
\rho = 10, \qquad
\vel = \frac{u_0}{r_0}( -(y-\shalf), \ (x-\shalf), \ 0)
\]
and for $r_0 < r < r_1$
\[
\rho = 1 + 9f, \qquad
\vel = \frac{f u_0}{r}( -(y-\shalf), \ (x-\shalf), \ 0), \qquad f = \frac{r_1 - r}{r_1 - r_0}
\]
and for $r > r_1$
\[
\rho = 1, \qquad \vel = (0,~0,~0)
\]
with $r_0 = 0.1$, $r_1 = 0.115$ and $u_0=2$. The rest of the quantities are     constant and given by
\[
p = 1, \qquad \bthree = \frac{1}{\sqrt{4\pi}}(5, \ 0, \ 0)
\]
We set $\gamma = 1.4$ and domain is discretized with $512\times 512$ mesh points. The numerical solutions are computed using LxF, HLL, HLLC flux up to the time $T = 0.15$ units. In Figure~\ref{fig:rotor}, we have shown the Mach number for considered fluxes and degree 1 to 3 over a mesh of size $128\times 128$. We can observe from the figures that in all cases, circularly rotating velocity field in the central part is captured well and the solutions remain stable. The  results on a finer mesh of  $512\times 512$ cells is shown in  Figure~\ref{fig:rotor1}, and we can observe that all the solution features are now resolved very sharply.



\subsection{Magnetic field loop test}
This test case \cite{Gardiner2005} involves the advection of magnetic field loop over a periodic domain. The numerical simulations are performed over the computational domain $[-1,+1] \times [-0.5,+0.5]$ with periodic boundary conditions in both directions. The initial density, pressure and velocity are uniform in the domain and given by
\[
\rho = 1, \qquad p = 1, \qquad \vel = (2,~1,~0)
\]
while the magnetic field is given by
\[
\bthree = \begin{cases}
A_0(-y/r,~x/r,~0) & r < r_0 \\
(0,~0,~0) & \textrm{otherwise}
\end{cases}
\]
The parameters in the initial condition are $A_0 = 10^{-3}$ and $r_0=0.3$, and the solution is computed up to a time of $T=1$ units.
\begin{figure}
\begin{center}
\begin{tabular}{ccc}
\includegraphics[width=0.3\textwidth]{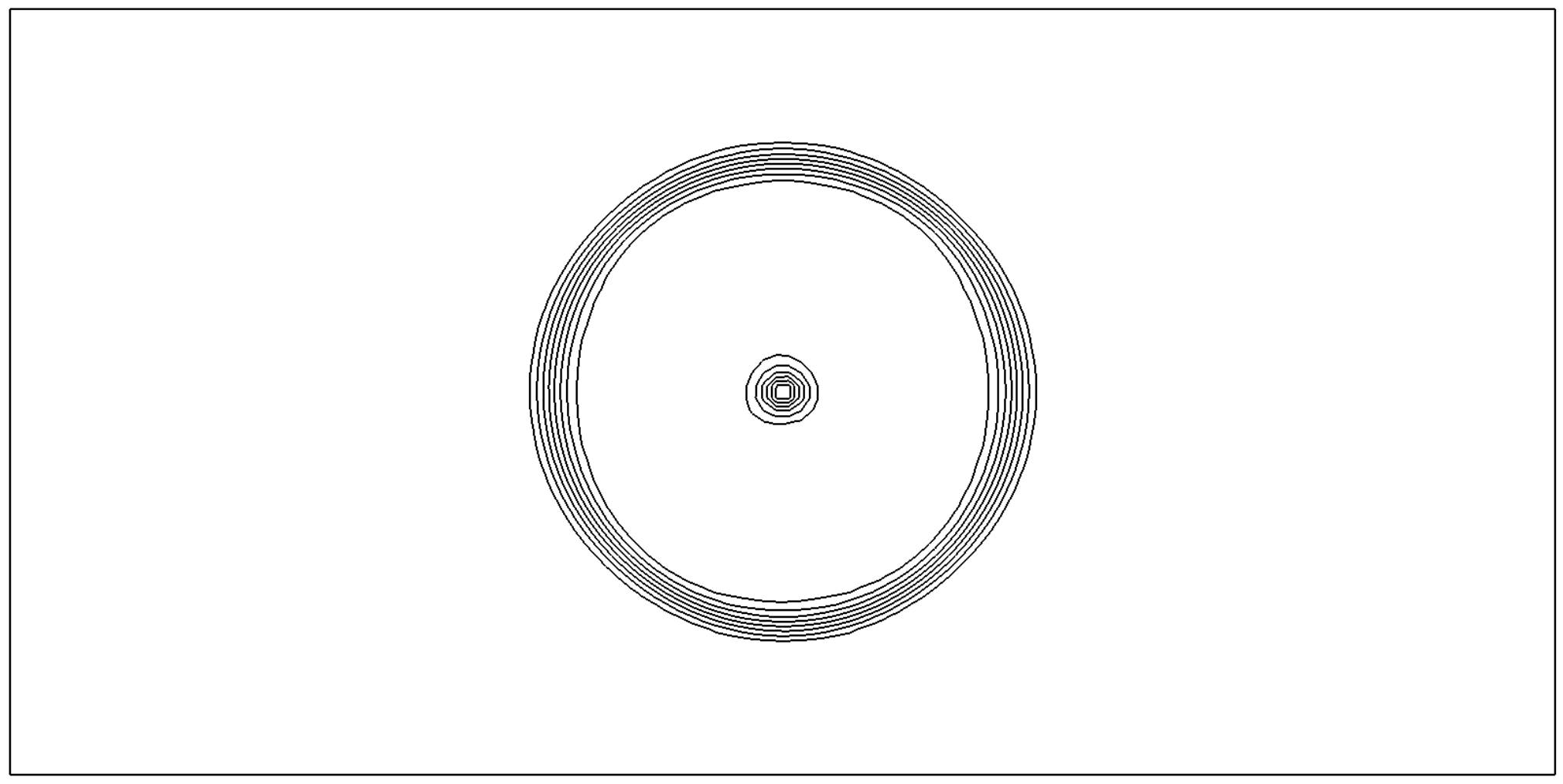}&
\includegraphics[width=0.3\textwidth]{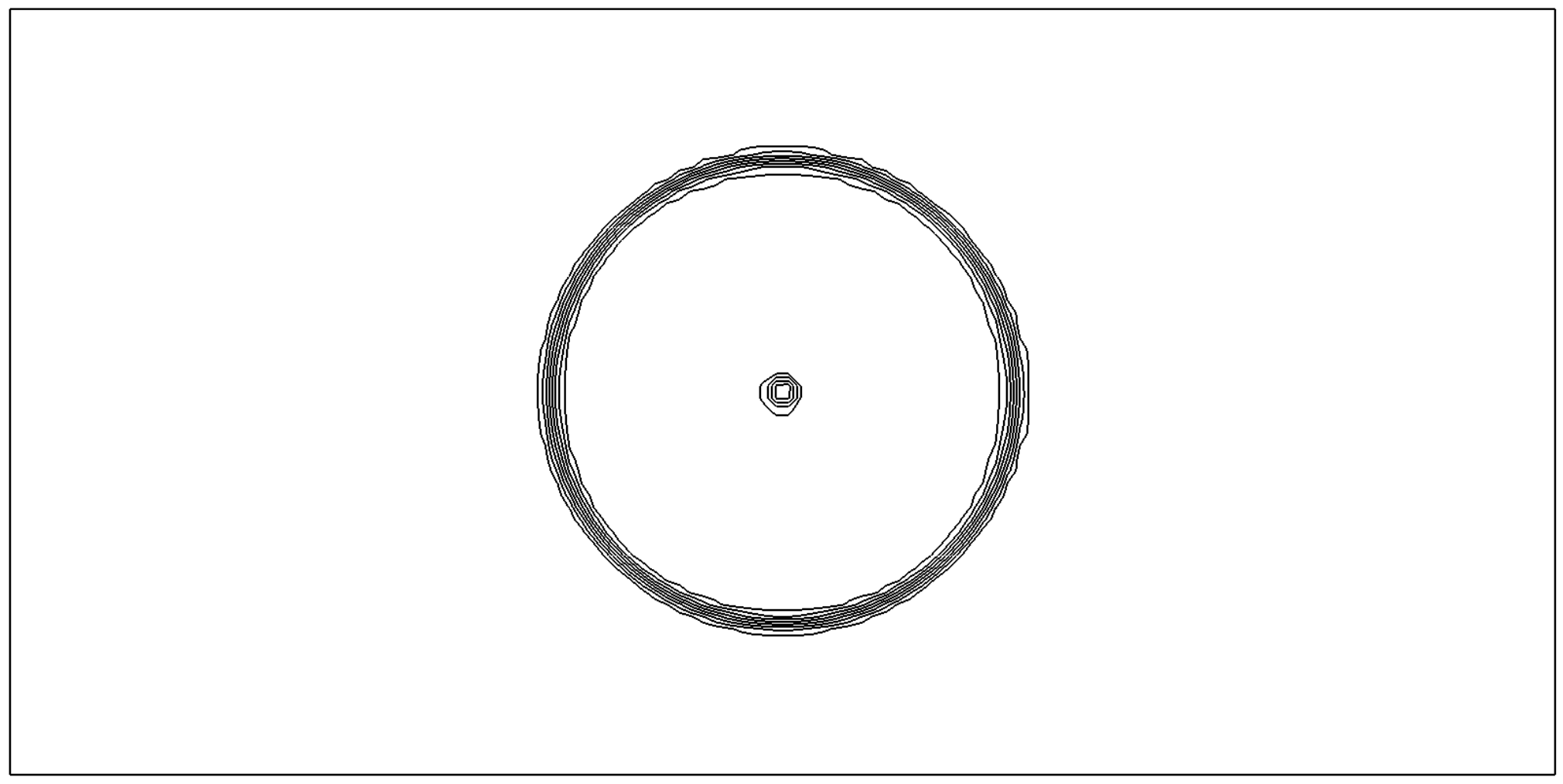}&
\includegraphics[width=0.3\textwidth]{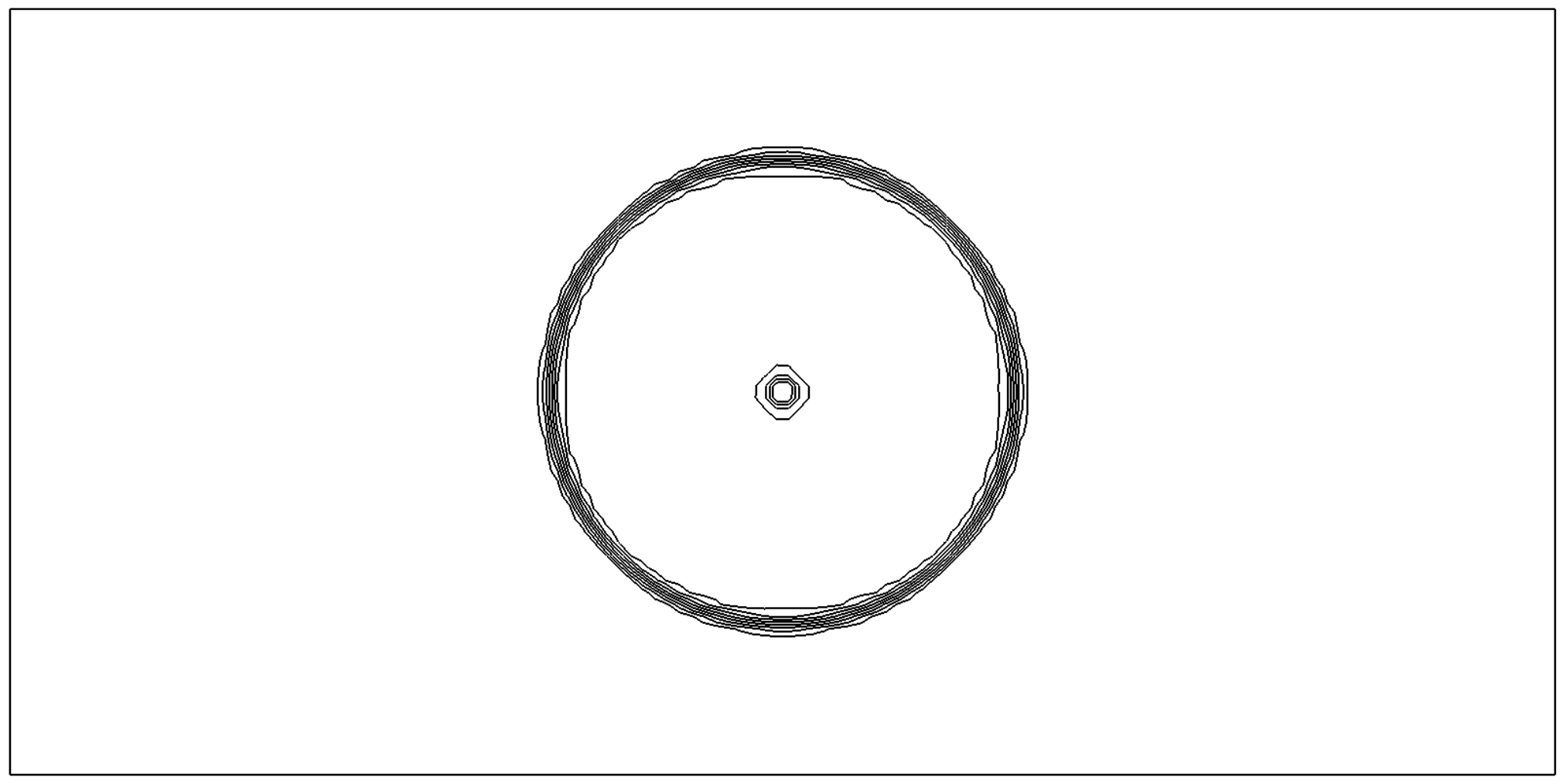} \\
\includegraphics[width=0.3\textwidth]{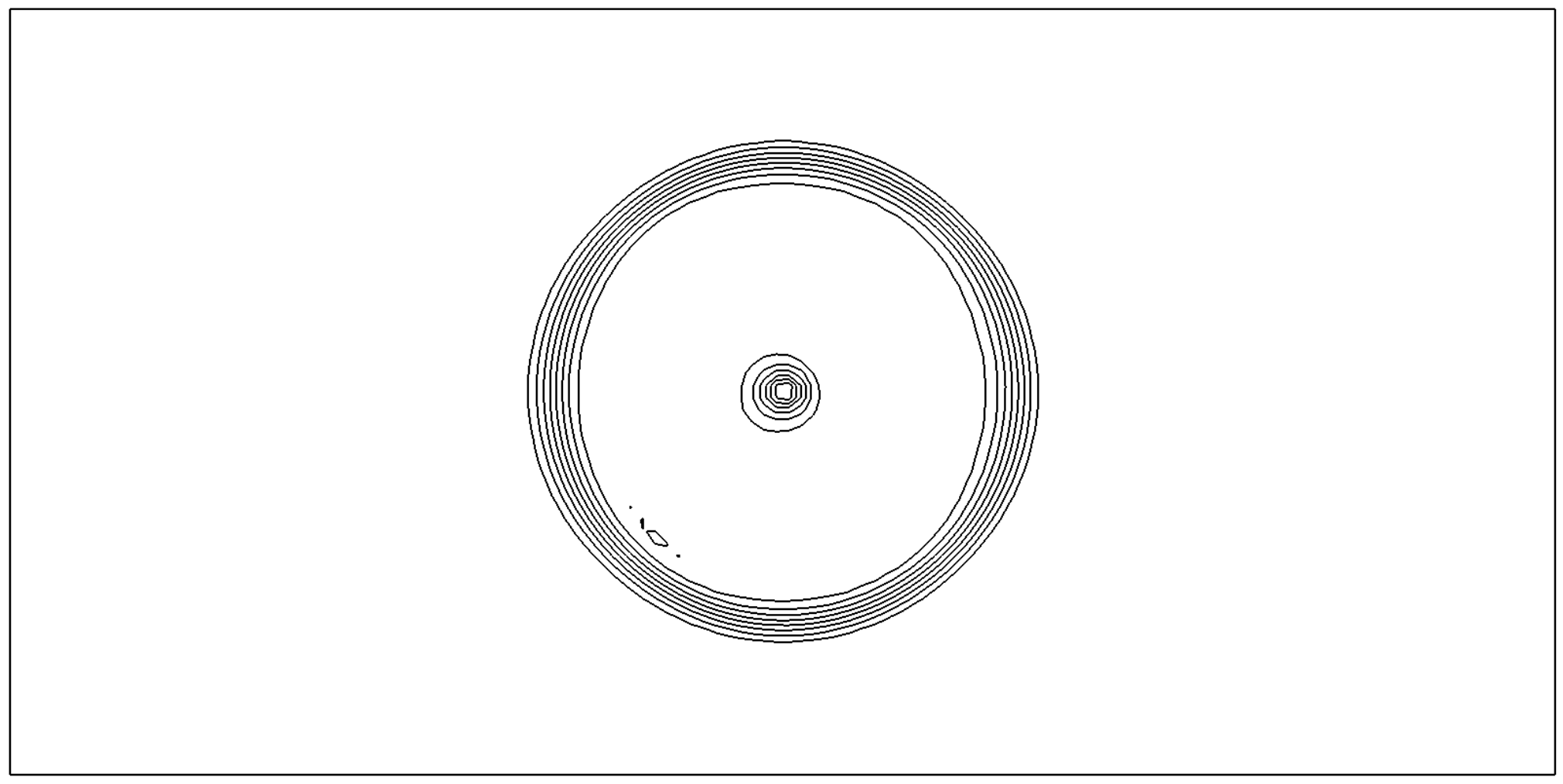}&
\includegraphics[width=0.3\textwidth]{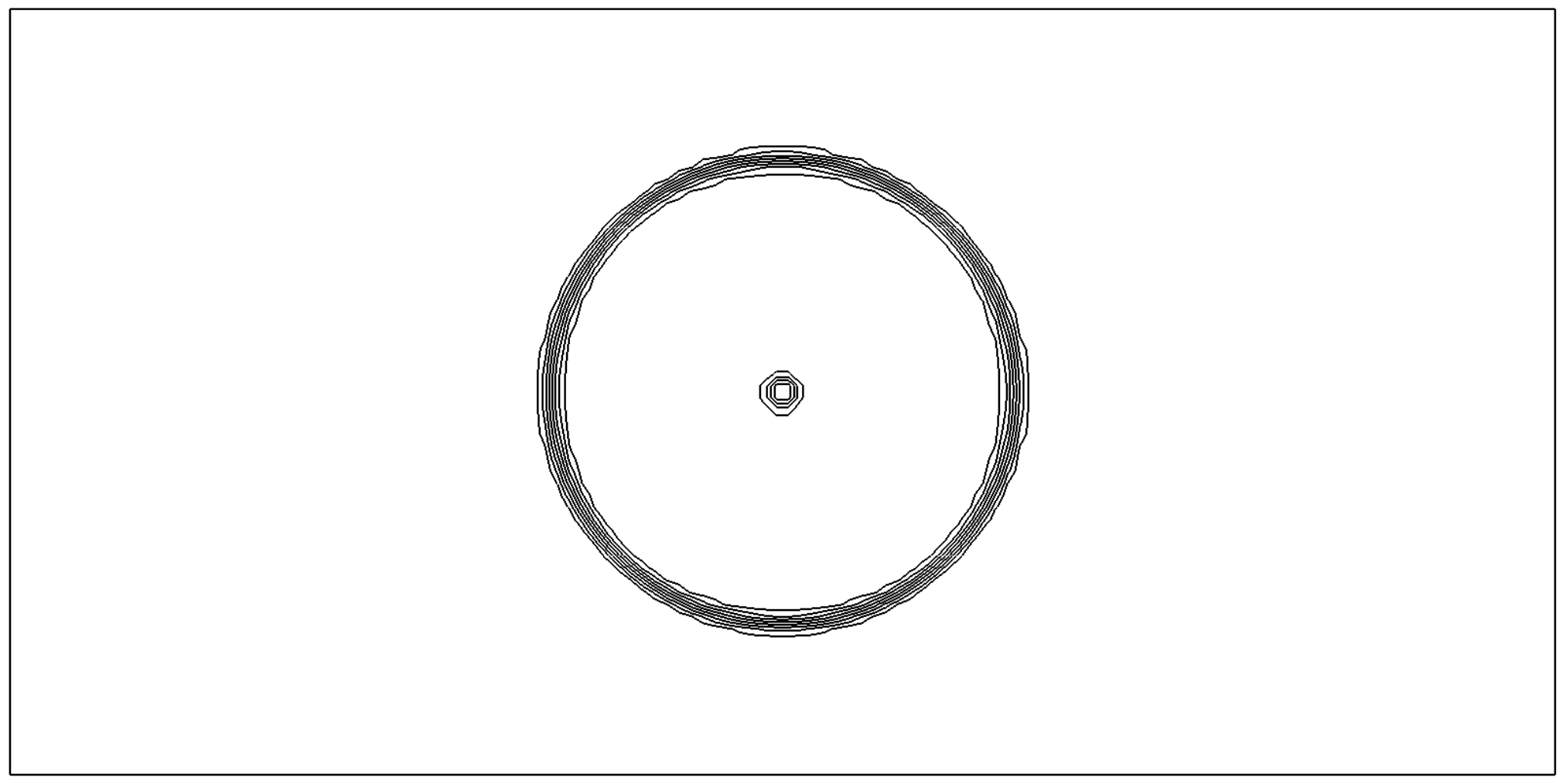}&
\includegraphics[width=0.3\textwidth]{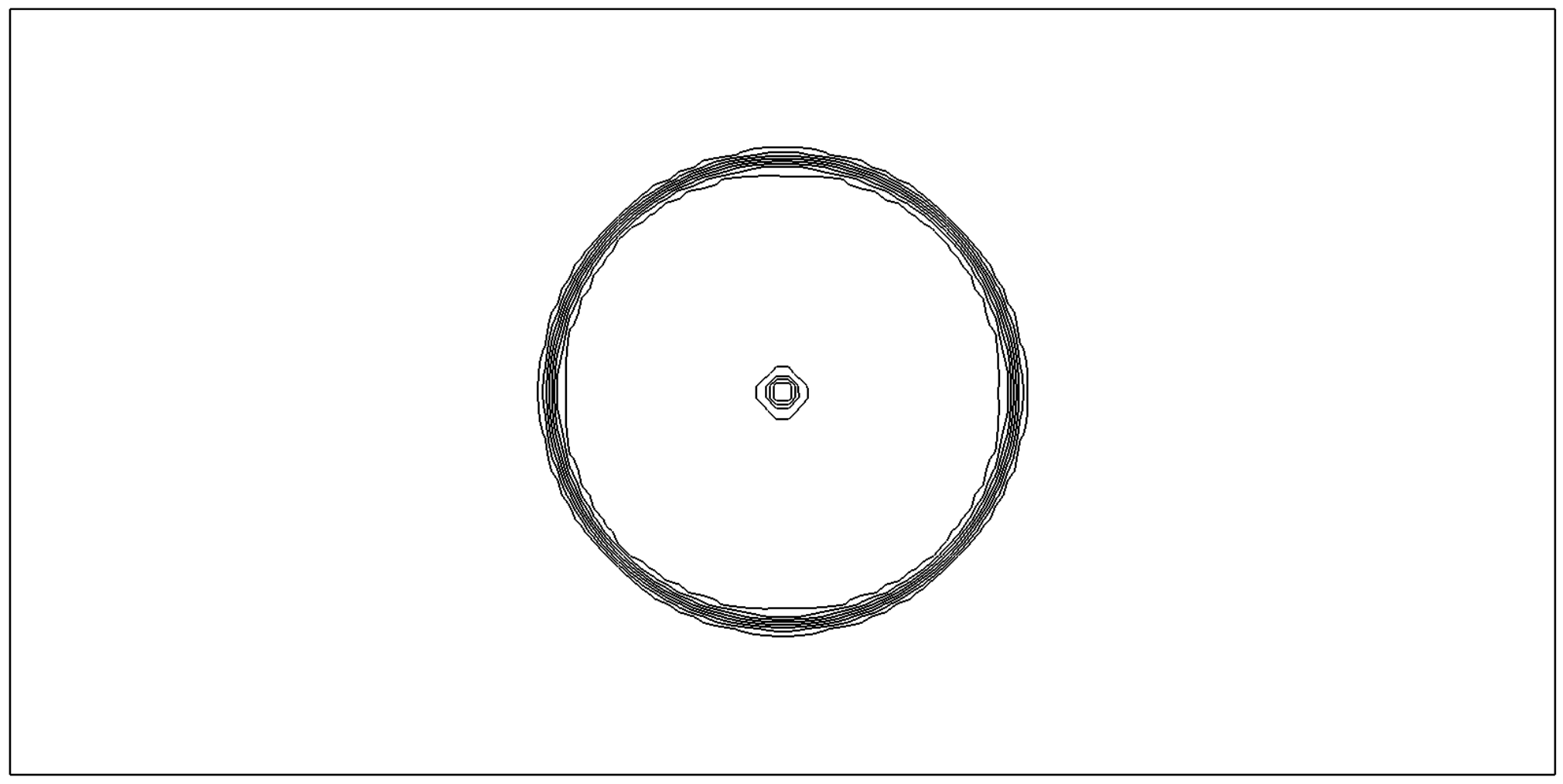} \\
\includegraphics[width=0.3\textwidth]{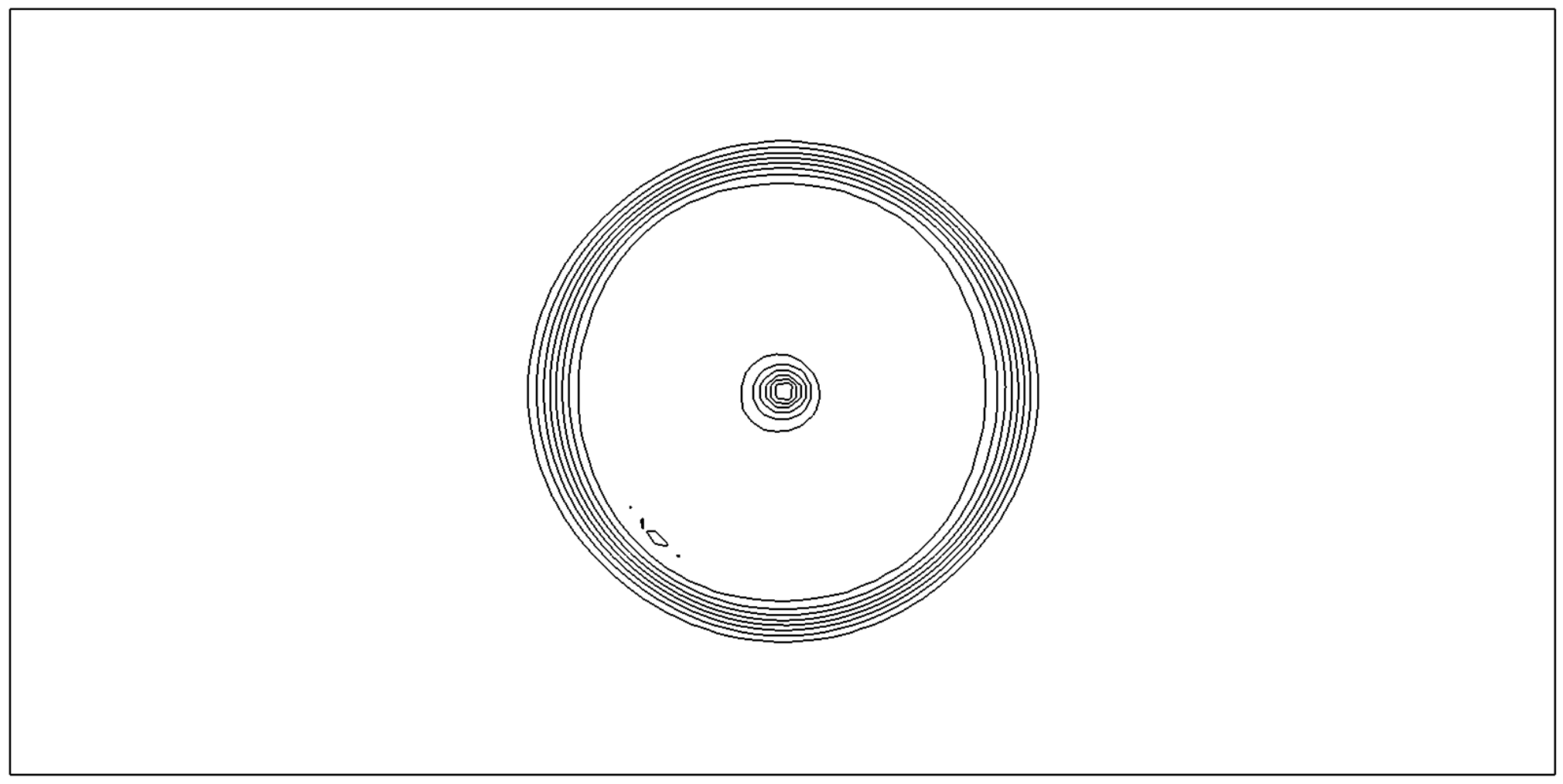}&
\includegraphics[width=0.3\textwidth]{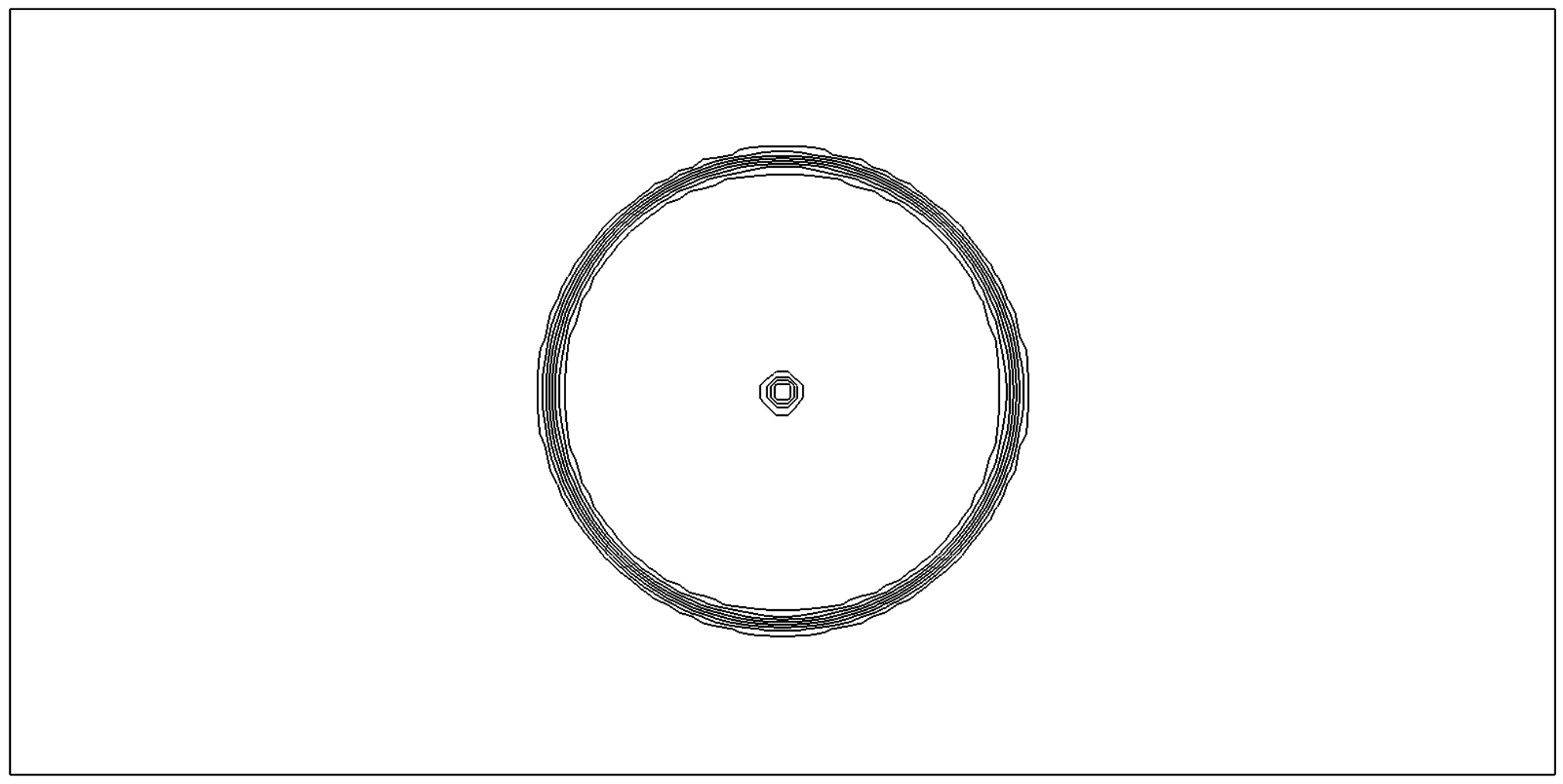}&
\includegraphics[width=0.3\textwidth]{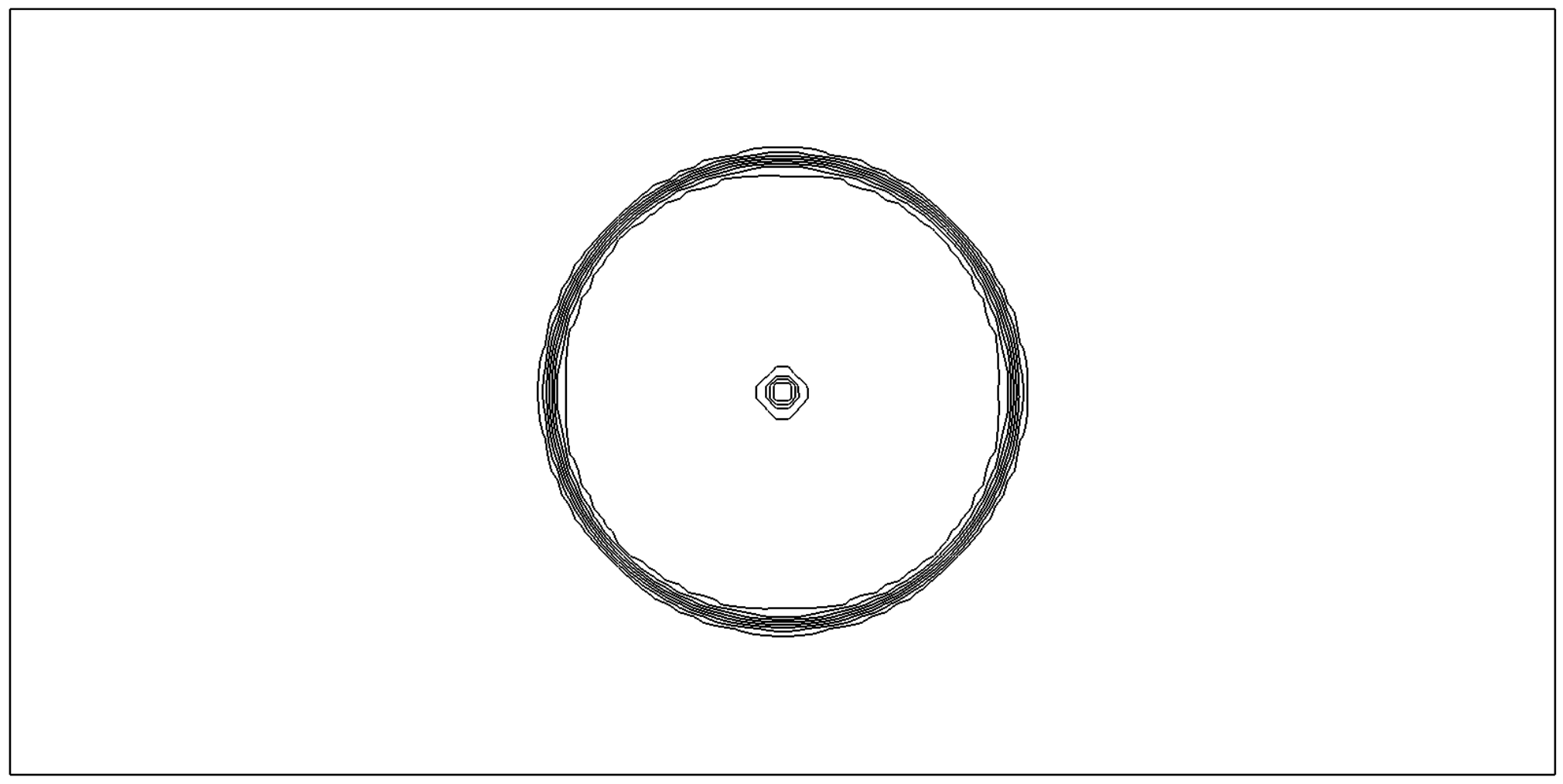} \\
  $k=1$ &  $k=2$ &   $k=3$\\
\end{tabular}
\end{center}
\caption{Contour plots of $\sqrt{B_x^2 + B_y^2}$ for loop advection test using $128 \times 64$ mesh at time $t=1$; 10 contours are shown in the range $(0, 0.00109)$. Top row: LxF flux, middle row: HLL flux, bottom row: HLLC flux.}
\label{fig:loop1}
\end{figure}
\begin{figure}
\begin{center}
\begin{tabular}{ccc}
\includegraphics[width=0.3\textwidth]{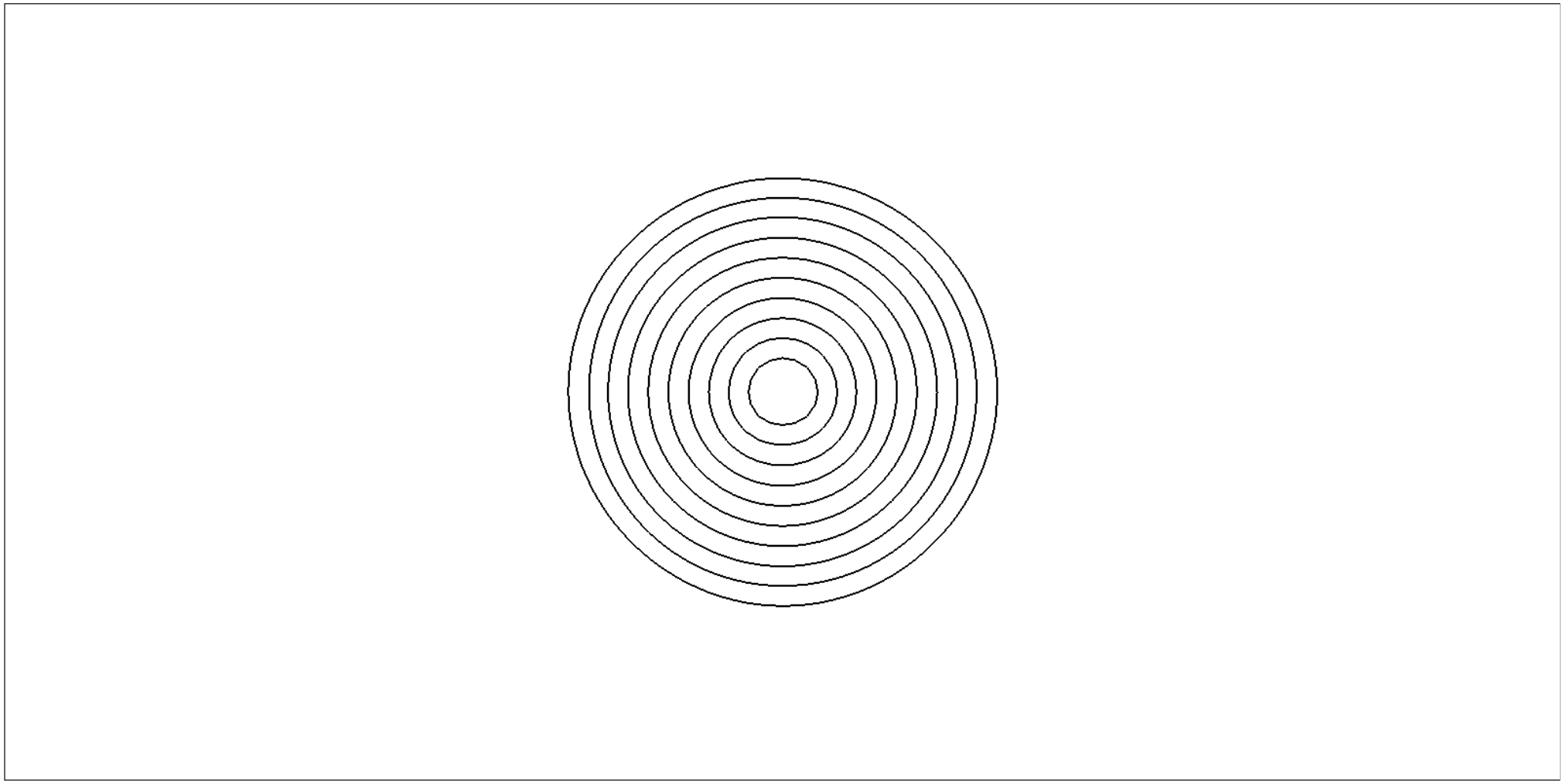}&
\includegraphics[width=0.3\textwidth]{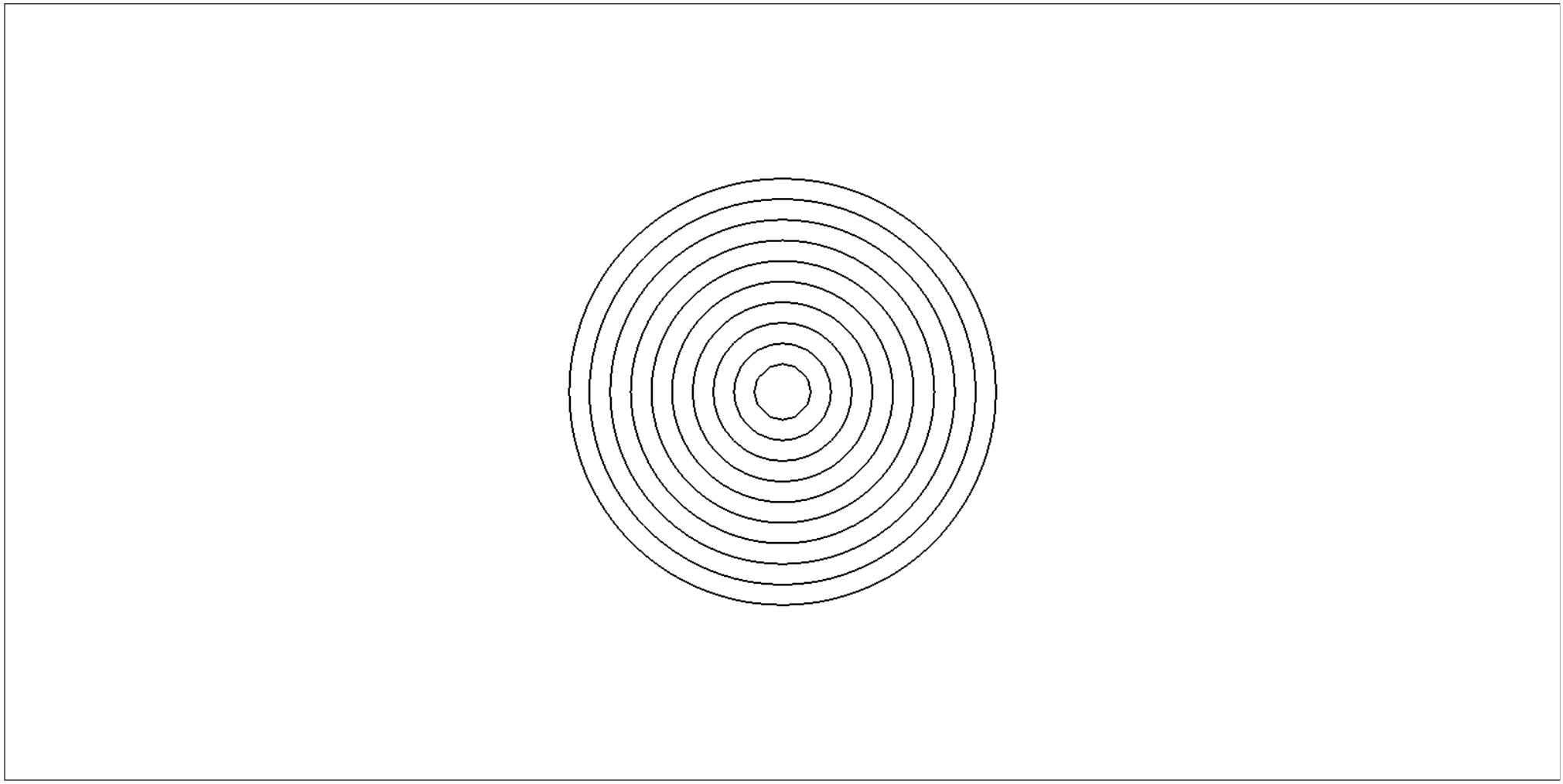}&
\includegraphics[width=0.3\textwidth]{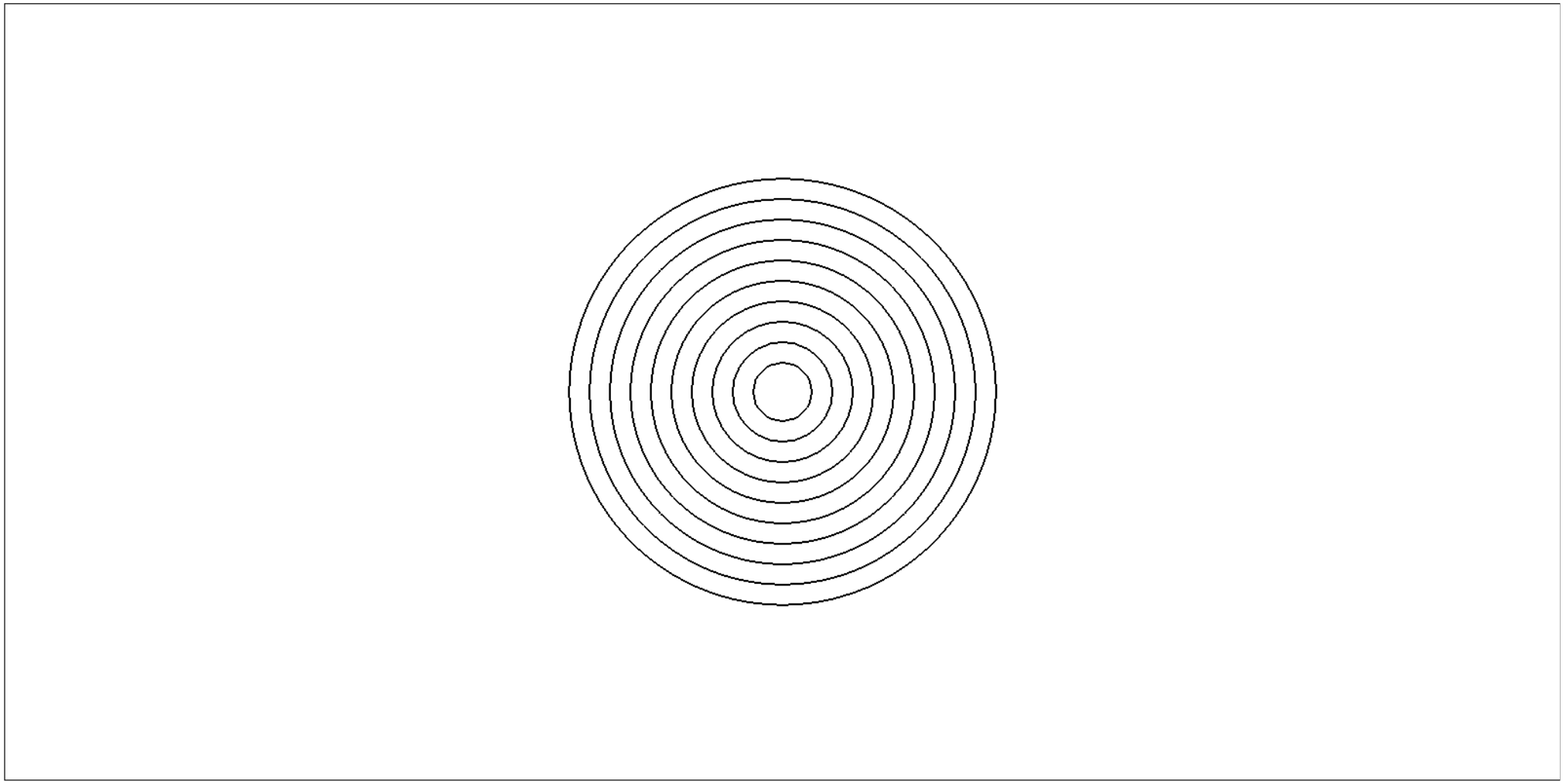} \\
\includegraphics[width=0.3\textwidth]{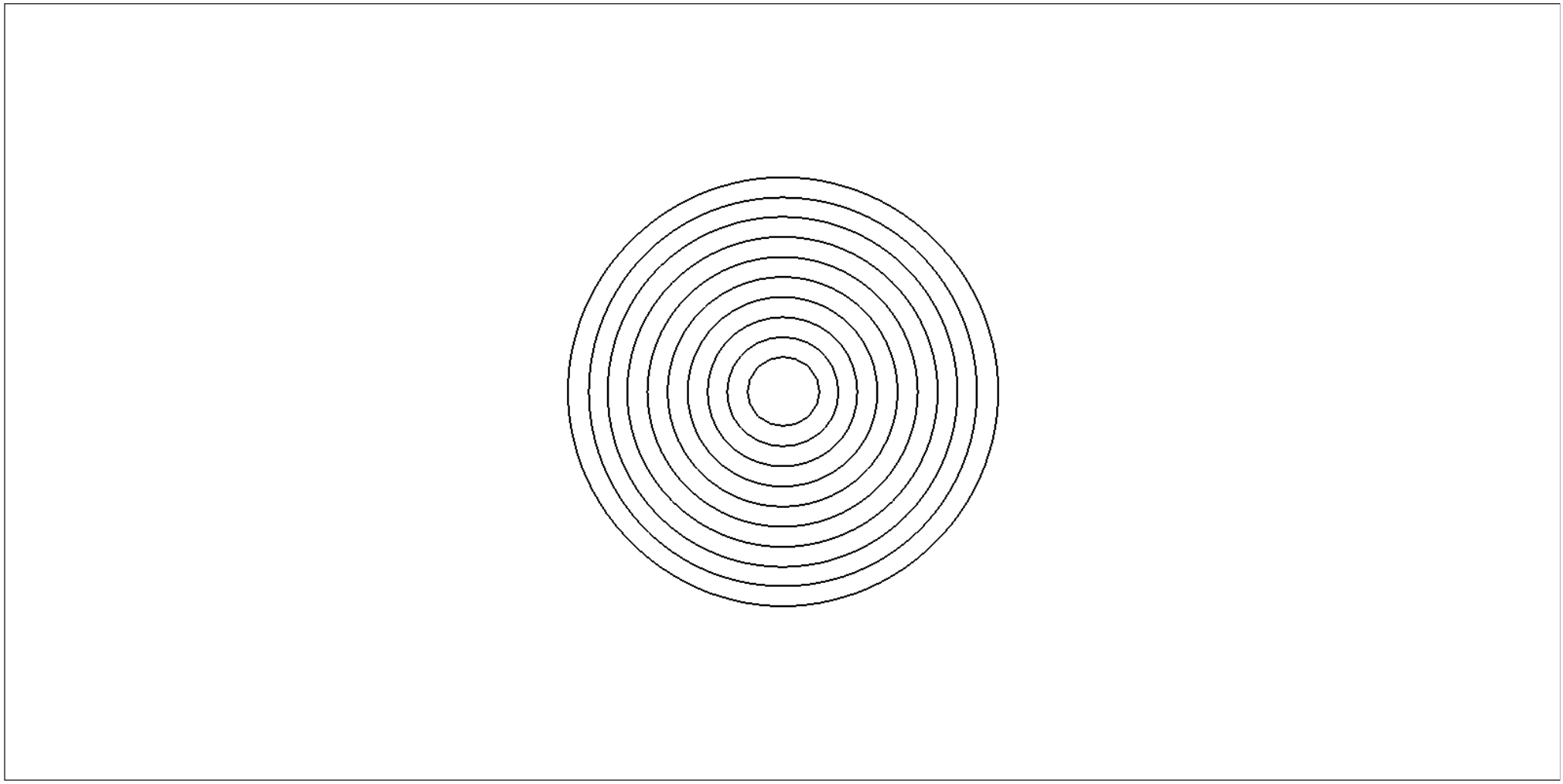}&
\includegraphics[width=0.3\textwidth]{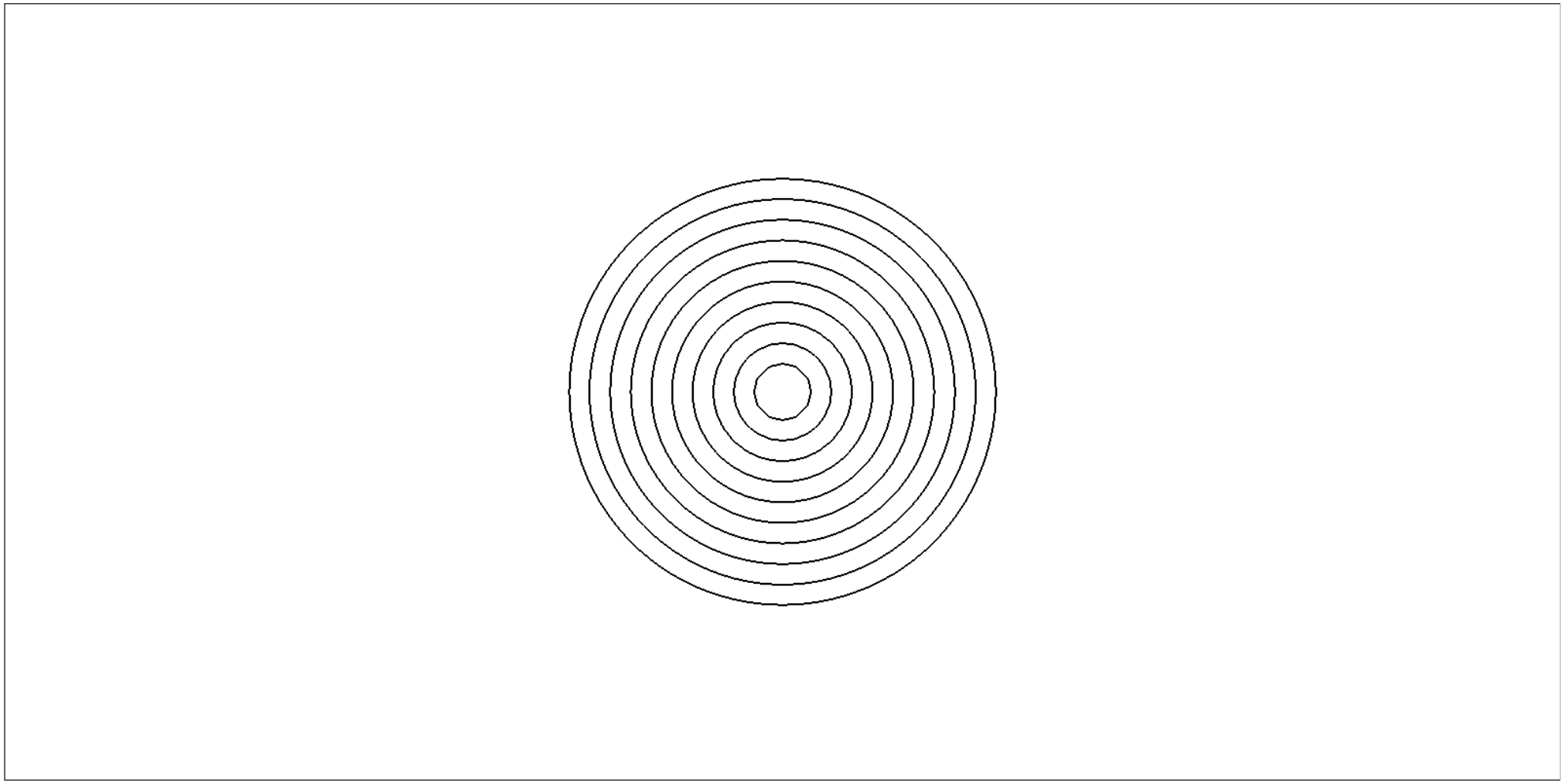}&
\includegraphics[width=0.3\textwidth]{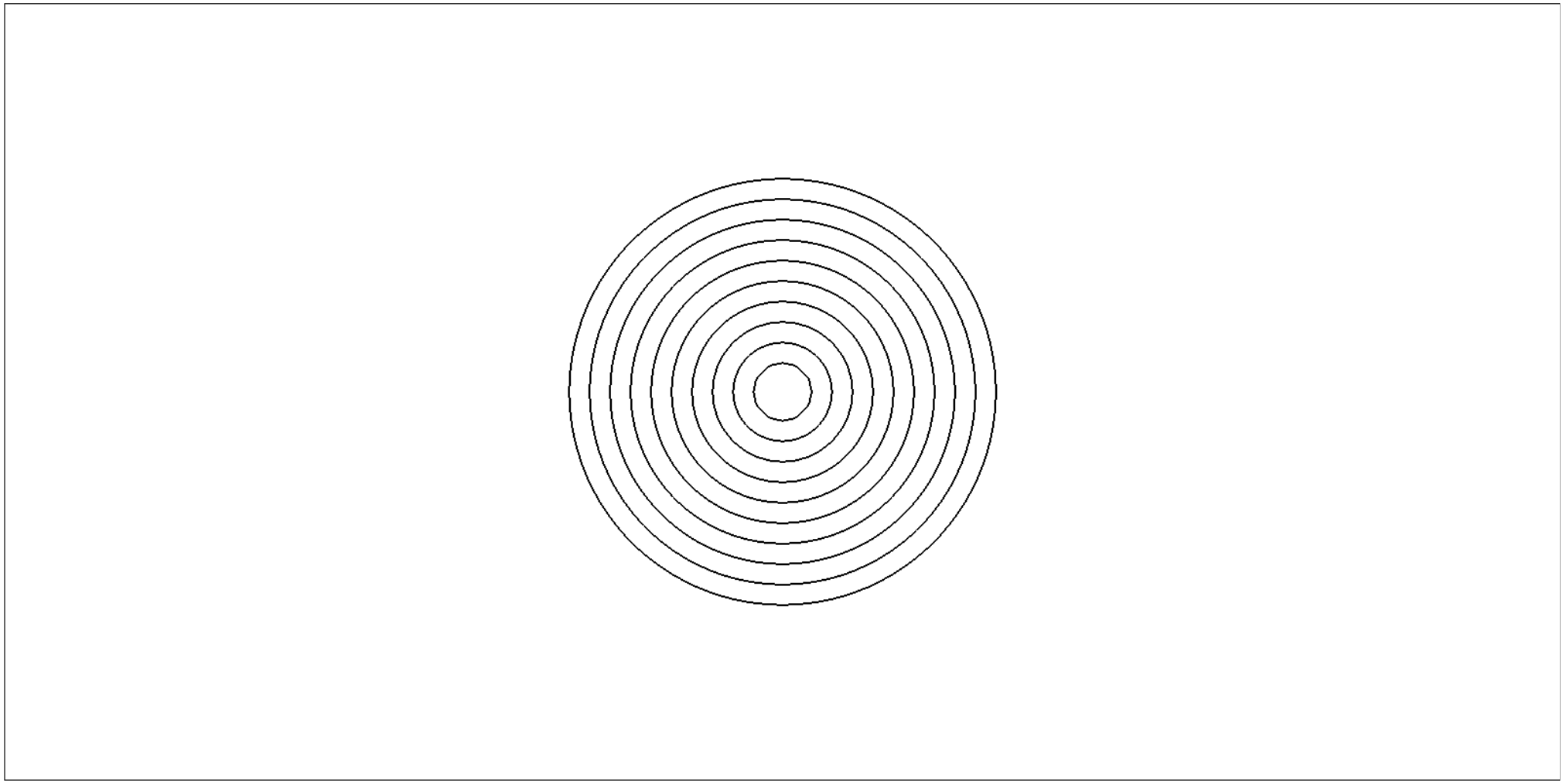} \\
\includegraphics[width=0.3\textwidth]{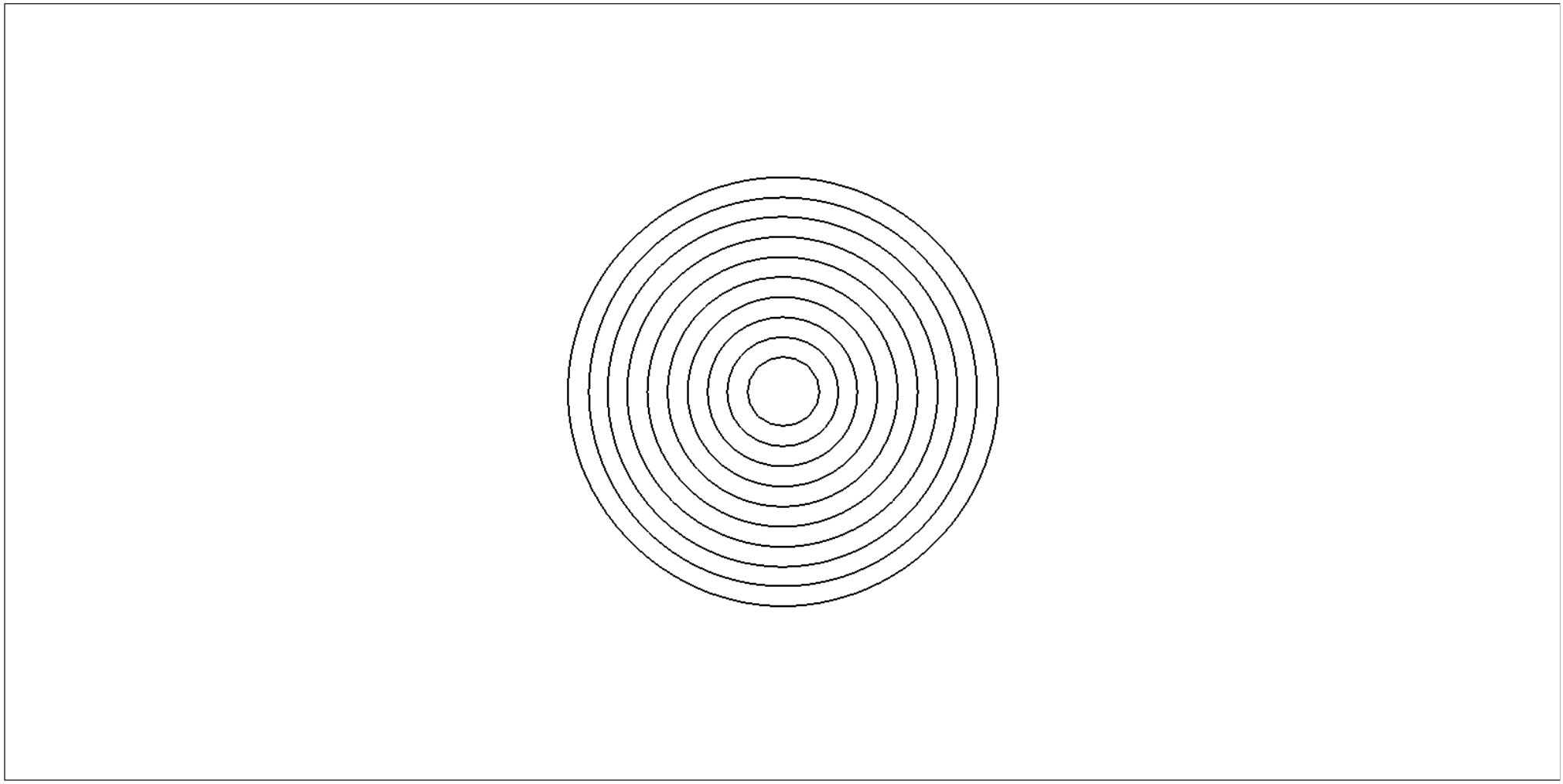}&
\includegraphics[width=0.3\textwidth]{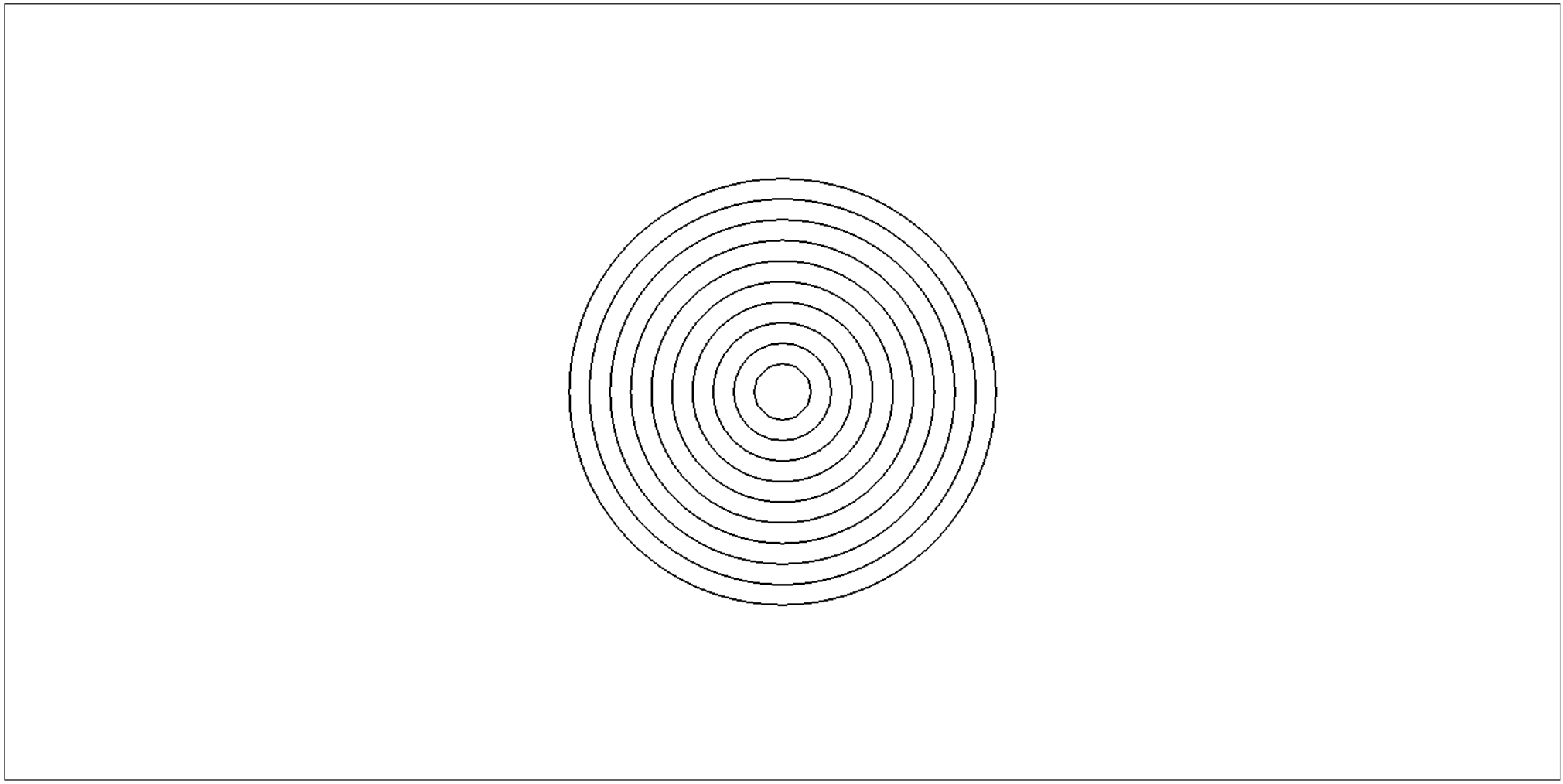}&
\includegraphics[width=0.3\textwidth]{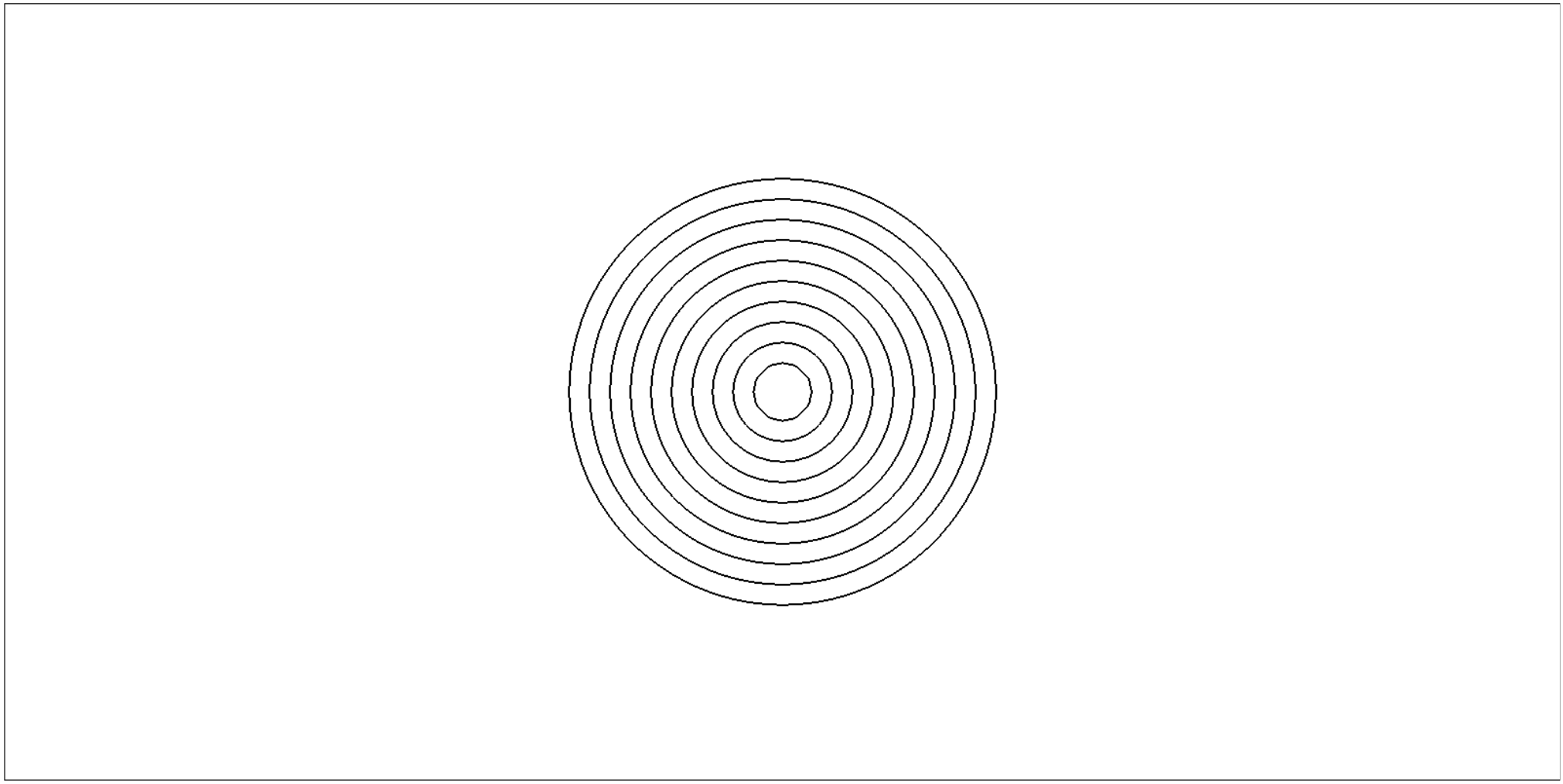} \\
  $k=1$ &  $k=2$ &  $k=3$\\
\end{tabular}
\end{center}
\caption{Contour plots of magnetic potential for loop advection test using $128 \times 64$ mesh at time $t=1$; 10 contours. First row LxF flux, second row HLL flux, third row HLLC flux.}
\label{fig:loop2}
\end{figure} 
\begin{figure}
\begin{center}
\begin{tabular}{cc}
 \includegraphics[width=0.48\textwidth]{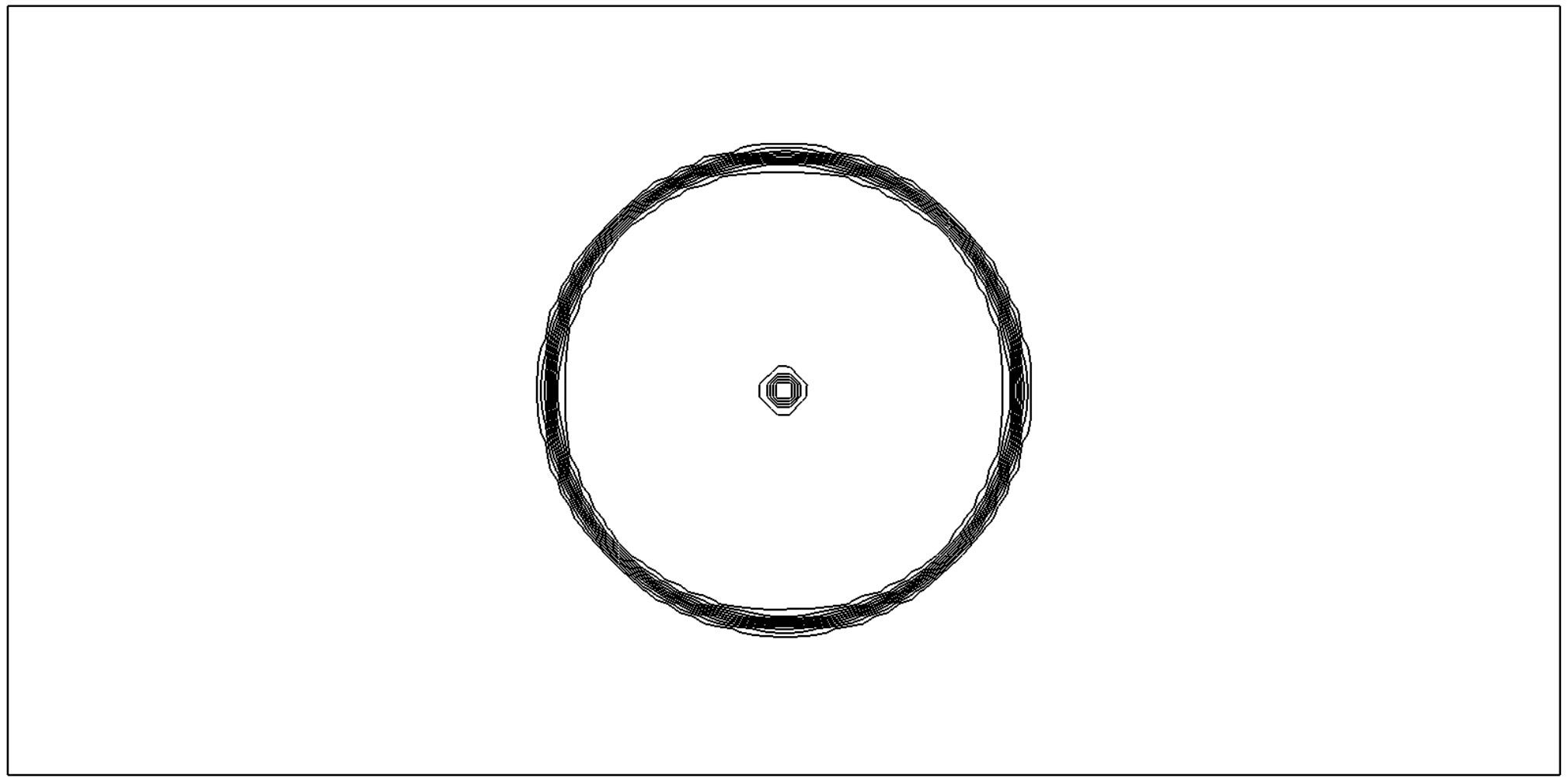} &
\includegraphics[width=0.48\textwidth]{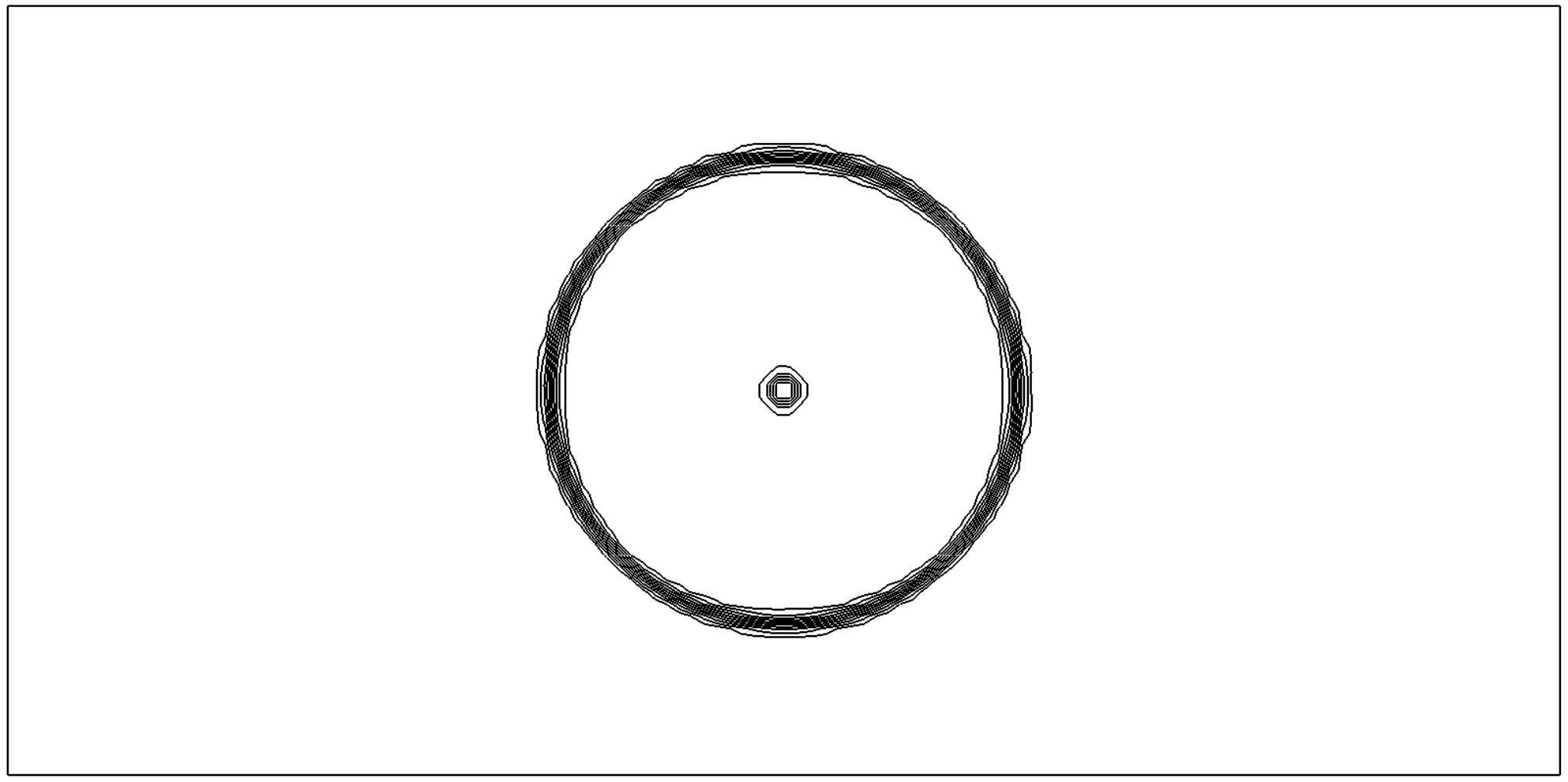} \\
(a) $t=1$ & (b) $t=10$ 
\end{tabular}
\end{center}
\caption{Contour plots of $\sqrt{B_x^2 + B_y^2}$ for loop advection test at time $t=1$ and $t=10$ using $128 \times 64$ mesh, degree $k=3$ and HLL flux.}
\label{fig:loop3}
\end{figure}
In Figure~\ref{fig:loop1}, we have depicted the magnitude of magnetic pressure  $\sqrt{B_x^2 + B_y^2}$ obtained using Lxf, HLL, and HLLC flux for second to fourth order schemes over a grid of size $128\times 64$. The magnetic field loop advects over the domain and returns to its initial position. Since this solution is essentially linear advection of $\B$, the use of shock indicator as described in~\cite{Fu2018} is very critical to reduce the dissipation from limiters. We can observe the numerical dissipation around the center and boundary of the advected loop where the solution is less smooth. We can also observe that numerical dissipation is reduced as we move from second to fourth order scheme. In Figure~\ref{fig:loop2}, we have depicted the contour plots of magnetic potential for different fluxes and second to fourth order schemes using ten contour lines. We can observe from Figure~\ref{fig:loop2} that the proposed schemes are able to preserve the shape and symmetry of magnetic field lines during simulations. \resa{Finally, we also show a long time simulation result in Figure~\ref{fig:loop3} using the fourth order scheme. At time $t=10$, the loop has advected through the domain for 10 times, and the scheme is still able to capture the features quite accurately.}

\subsection{Blast wave test}
\label{sec:blast}
\begin{figure}
\begin{center}
\begin{tabular}{ccc}
\includegraphics[width=0.33\textwidth]{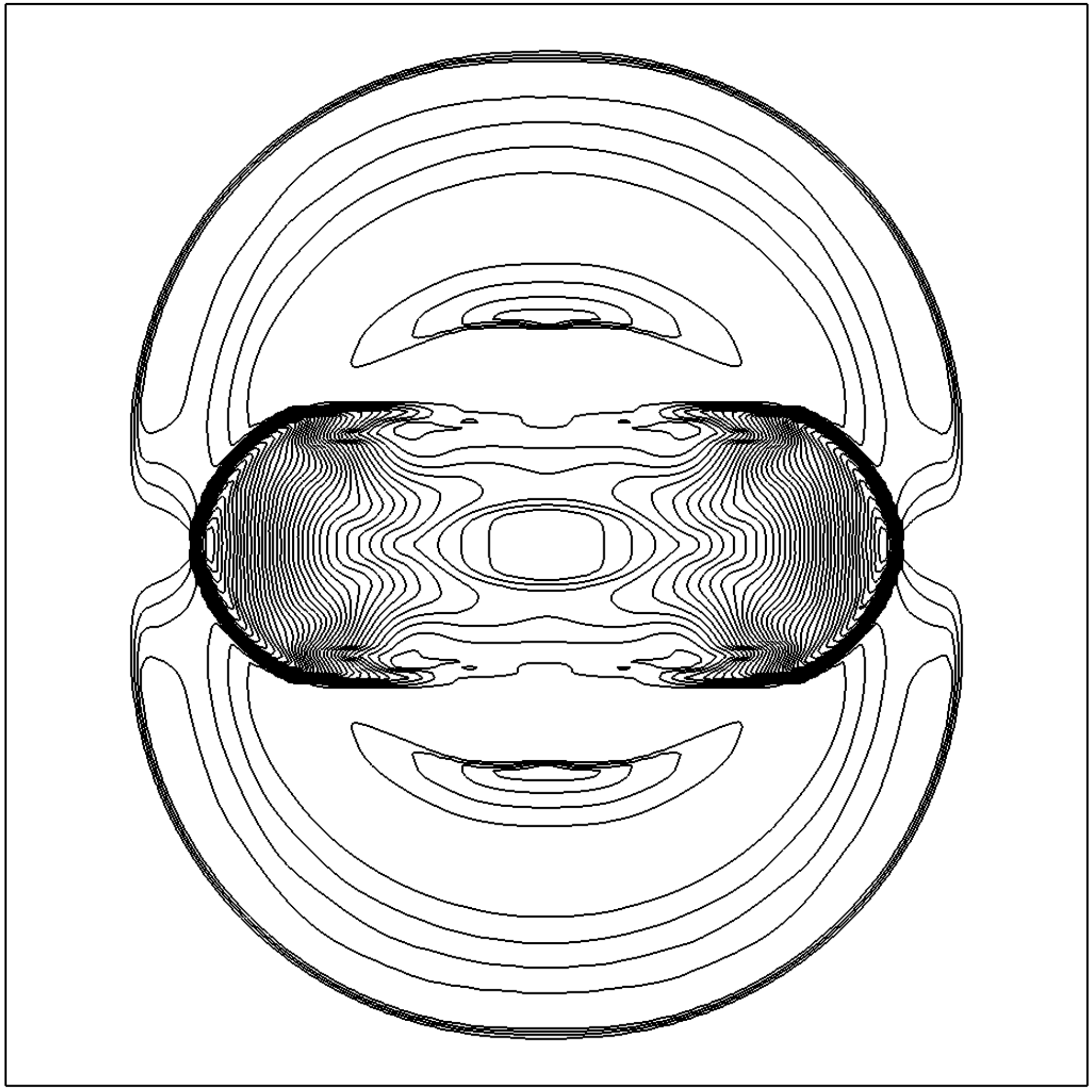} &
\includegraphics[width=0.33\textwidth]{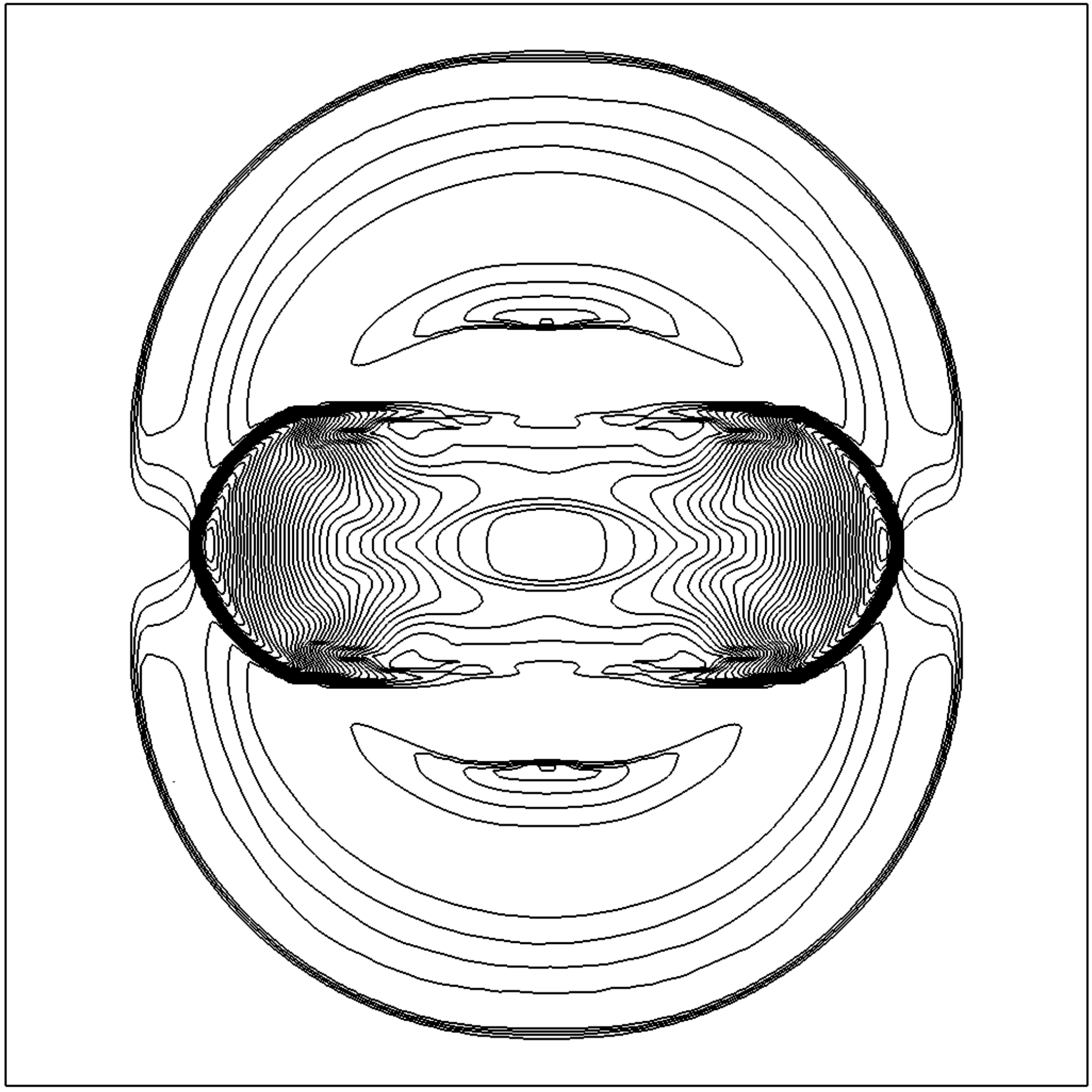} &
\includegraphics[width=0.33\textwidth]{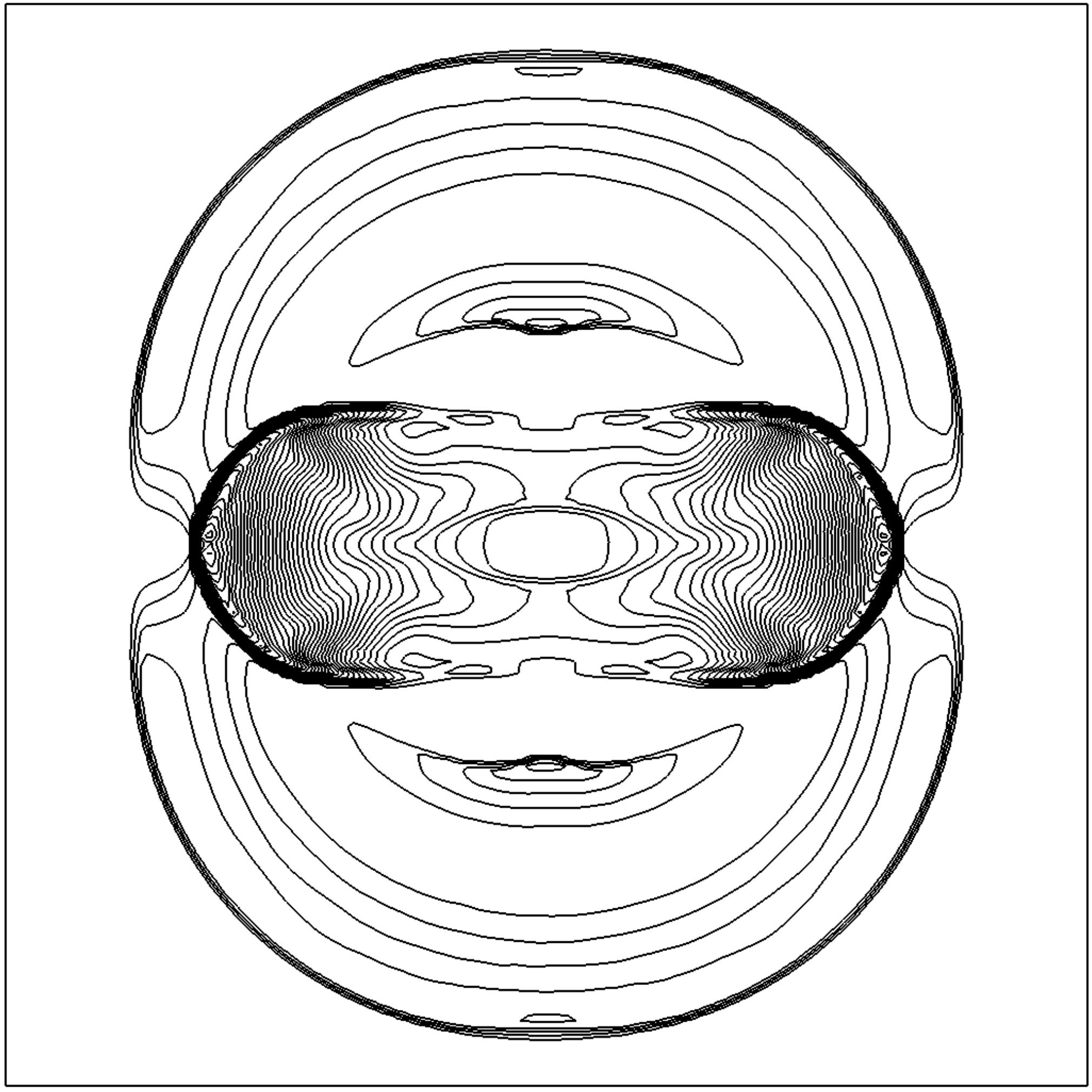} \\
(a) $v_x^2 + v_y^2$ (LxF) & (b) $v_x^2 + v_y^2$ (HLL) & (c) $v_x^2 + v_y^2$ (HLLC) \\
\includegraphics[width=0.33\textwidth]{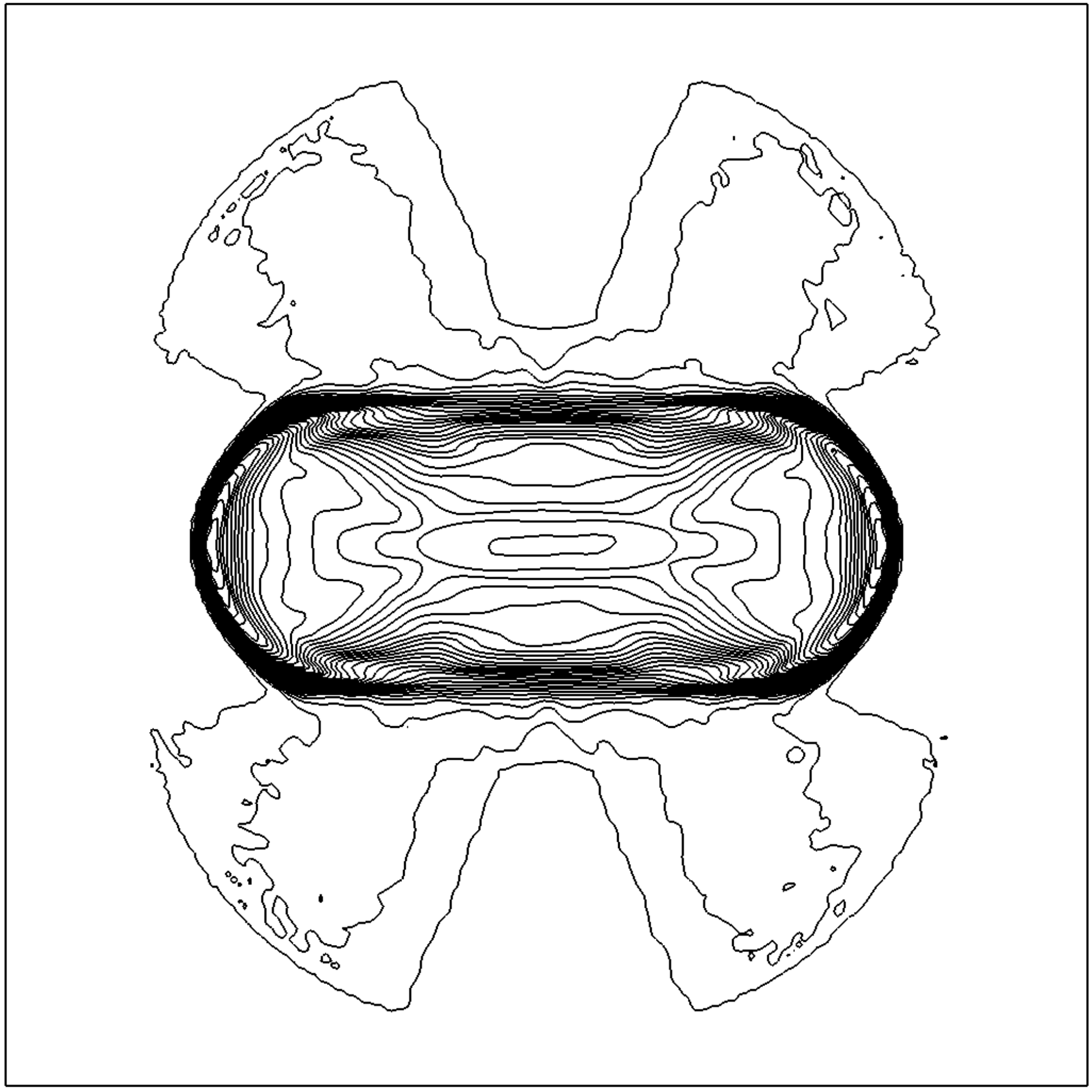} &
\includegraphics[width=0.33\textwidth]{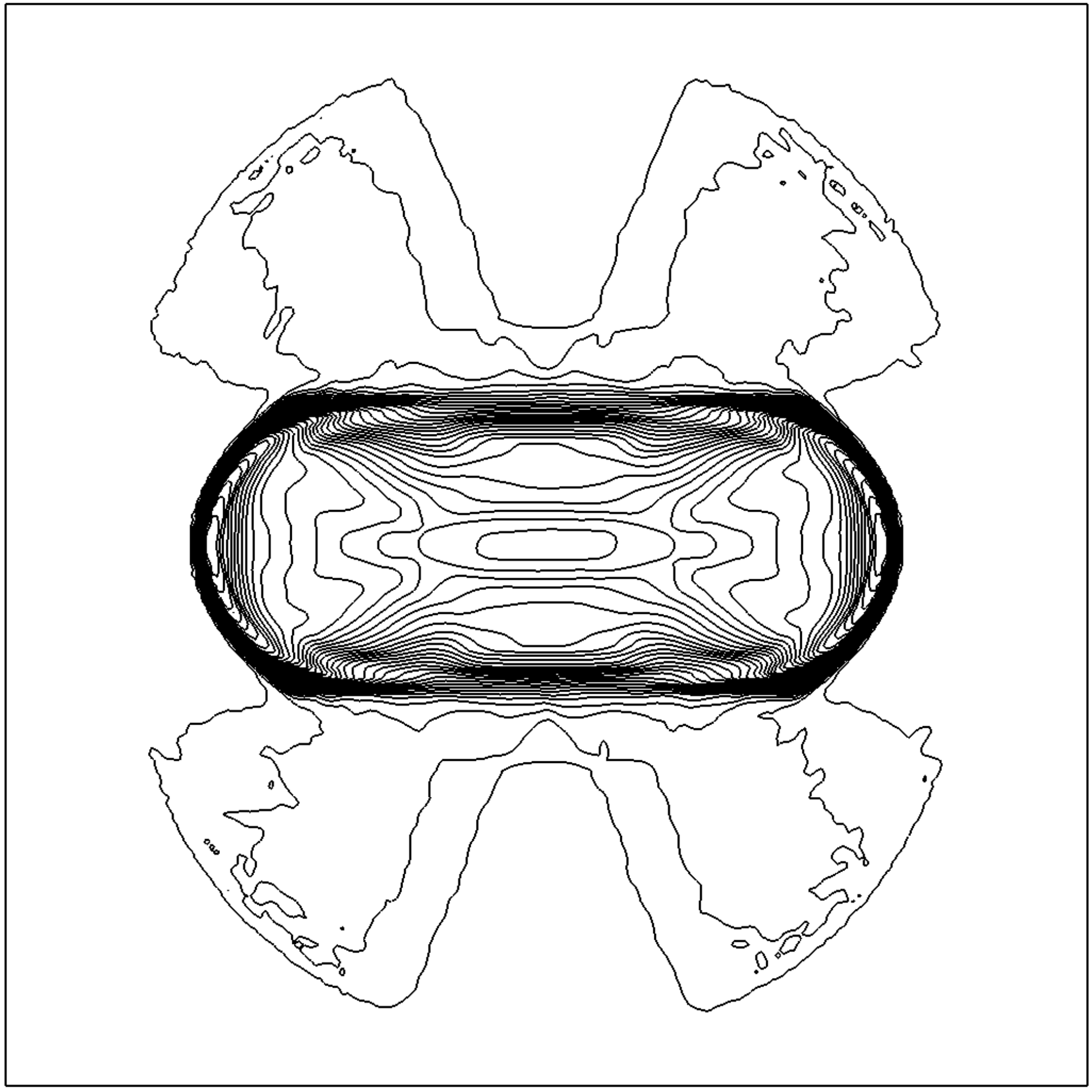} &
\includegraphics[width=0.33\textwidth]{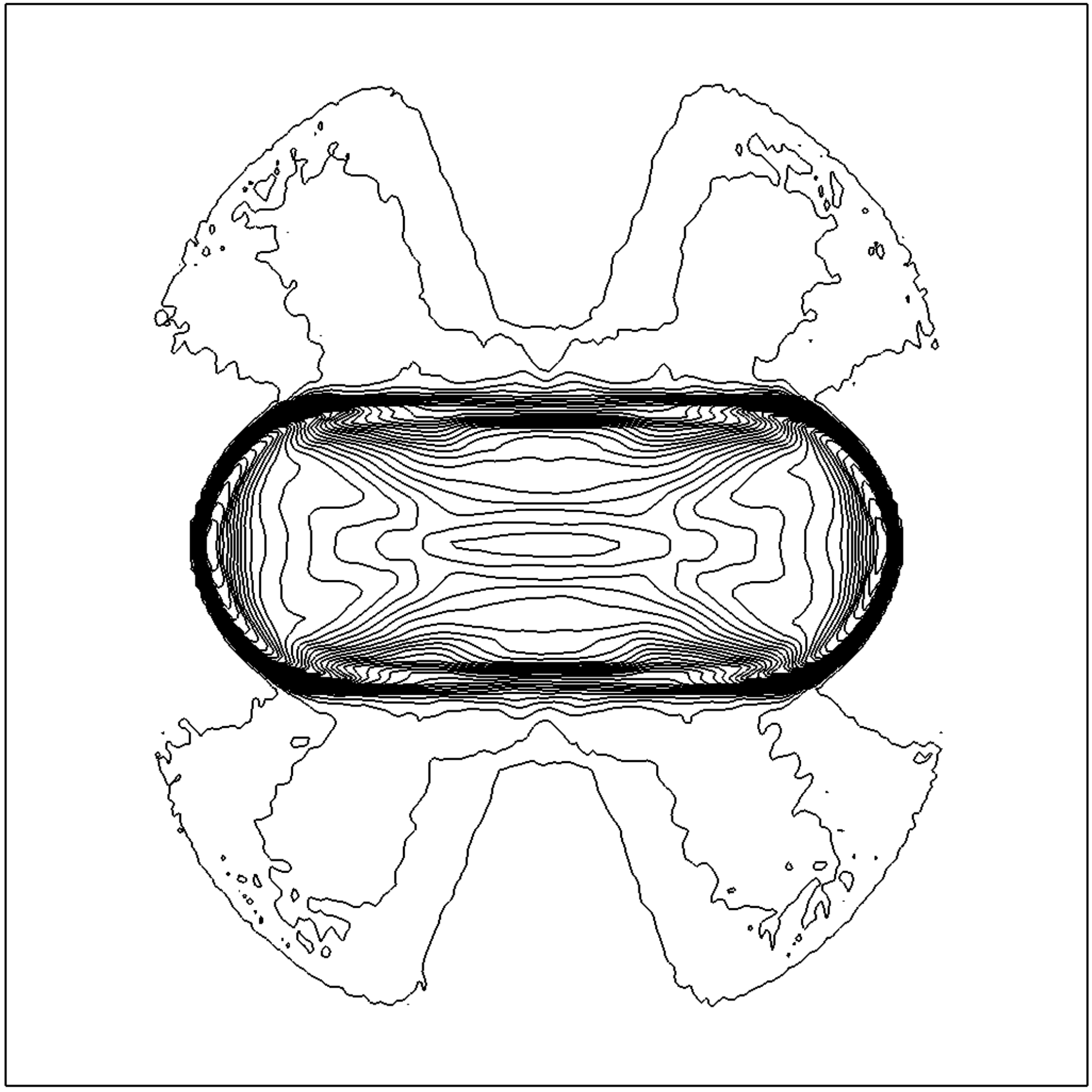} \\
(d) $p$ (LxF) & (e) $p$ (HLL) & (f) $p$ (HLLC) \\
\includegraphics[width=0.33\textwidth]{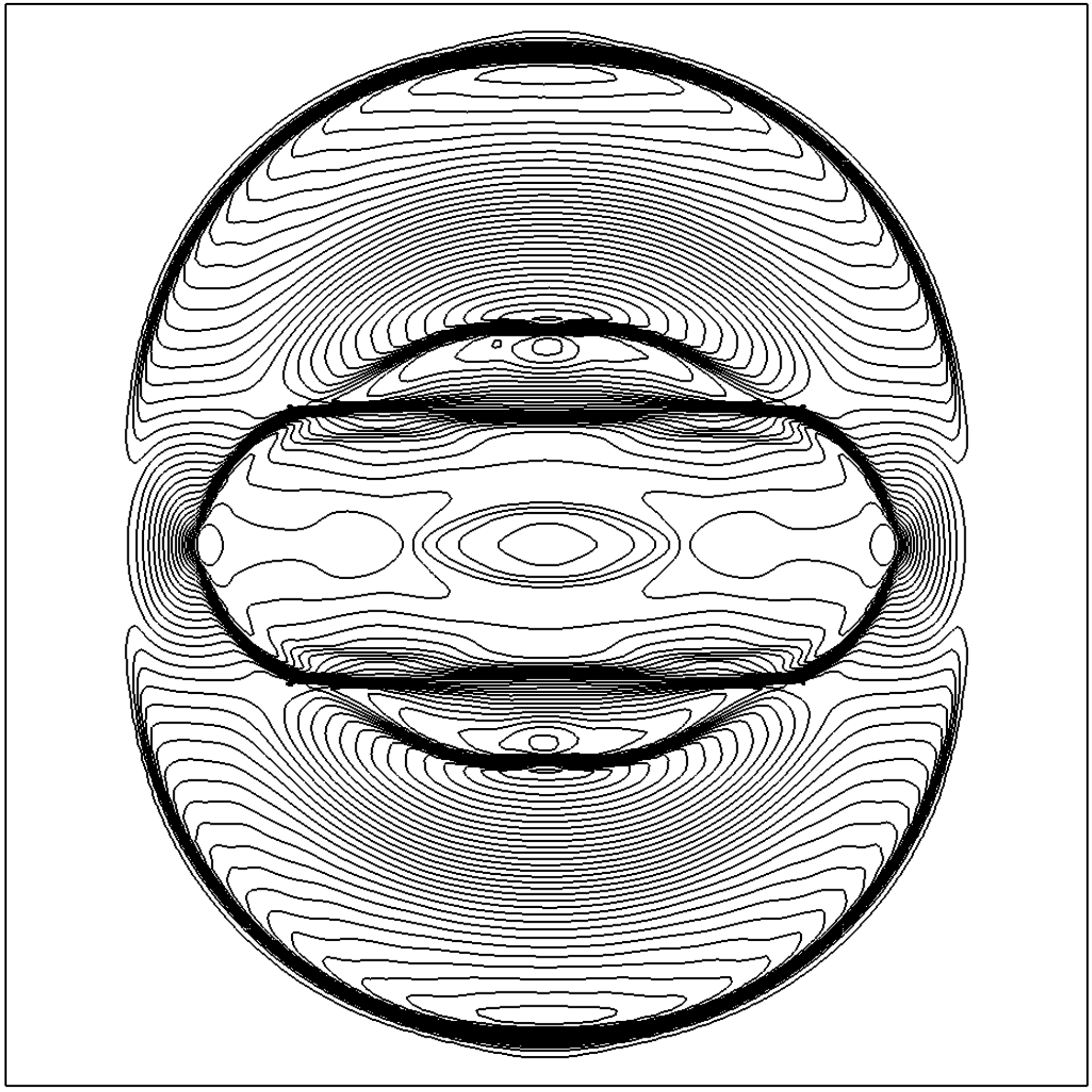}  &
\includegraphics[width=0.33\textwidth]{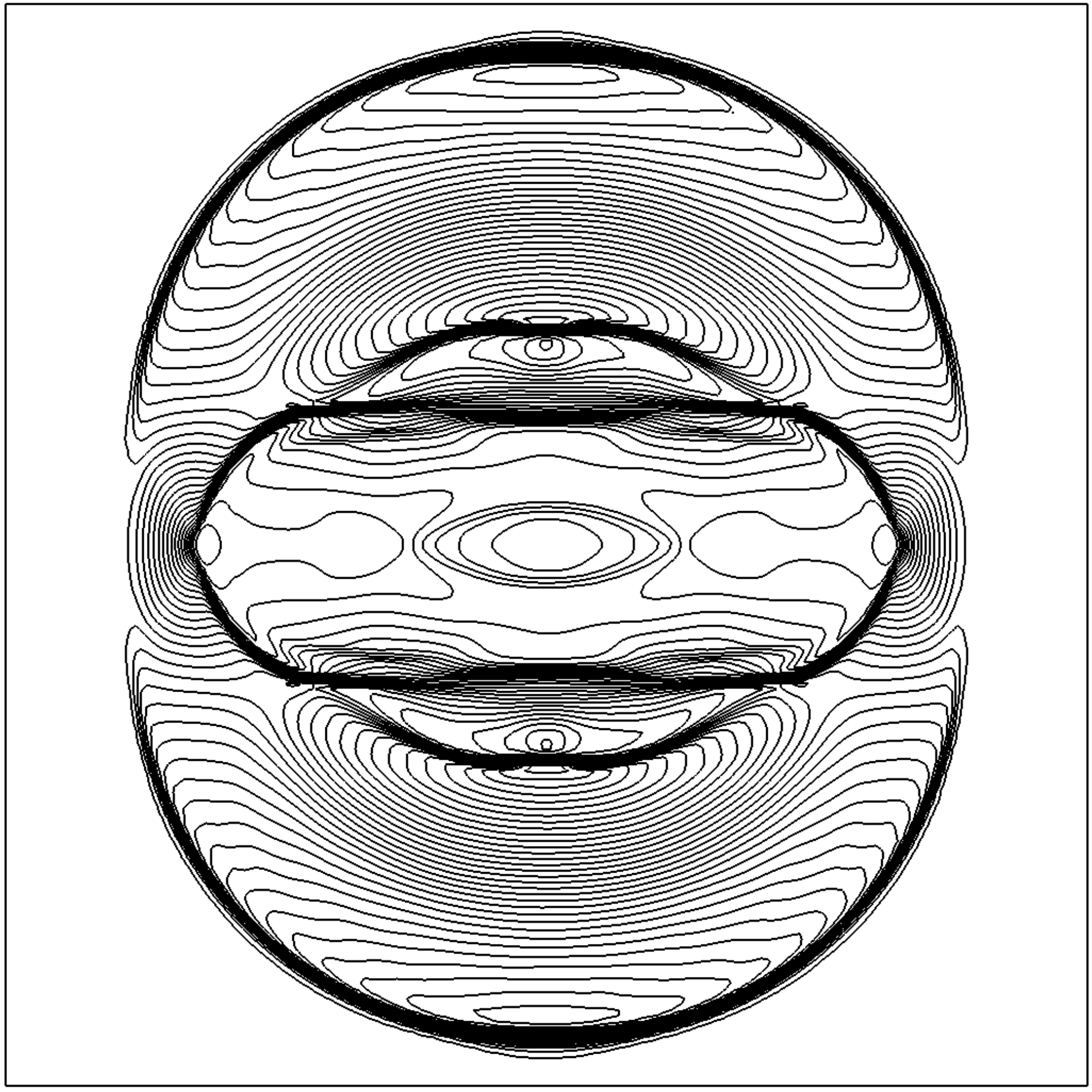}  &
\includegraphics[width=0.33\textwidth]{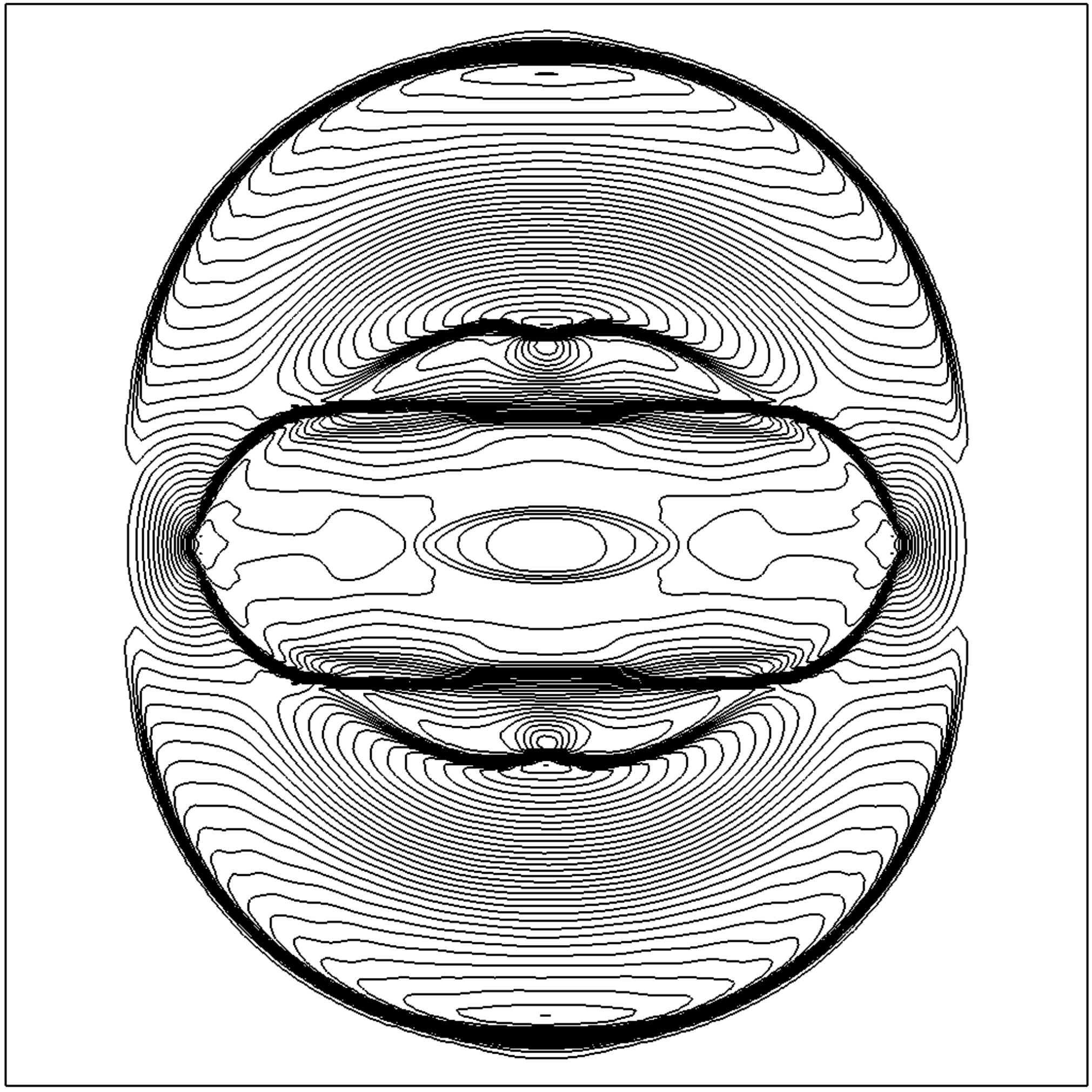}  \\
(g) $\frac{1}{2}(B_x^2 + B_y^2)$ (LxF) & (h) $\frac{1}{2}(B_x^2 + B_y^2)$ (HLL) & (i) $\frac{1}{2}(B_x^2 + B_y^2)$ (HLLC)
\end{tabular}
\caption{Blast test case using LxF, HLL and HLLC flux for degree $k=1$ on $200\times 200$ mesh with 40 contours. Top row: velocity square, middle row: pressure, bottom row: magnetic pressure}
\label{fig:blast1}
\end{center}
\end{figure}

The MHD blast wave test  introduced by Balsara \& Spicer  \cite{Balsara1999}, is a challenging test problem and often used  as a benchmark  for testing the robustness of the numerical algorithms in terms of maintaining positivity of solutions. The problem is initialized with constant density, velocity, and magnetic field except the pressure. The initial condition is given by
 \[
\rho = 1, \quad \vel= (\ 0, \ 0, \ 0),  \quad
\bthree = \frac{1}{\sqrt{4\pi}}(100, \ 0, \ 0), \qquad 
p = \begin{cases}
1000 & r < 0.1 \\
0.1 & r > 0.1 \end{cases}
\]
where $r^2=(x-0.5)^2+(y-0.5)^2$. The computational domain is $[0,1] \times [0,1]$ with periodic  boundary conditions on all sides. The numerical experiments are performed over a grid of size $200\times200$  up to the time $T=0.01$.  As discussed earlier in the paper, the positivity of solutions cannot be guaranteed by constraint preserving schemes and this becomes an issue especially when we have low values of plasma beta where $\beta = 2 p /|\bthree|^2$. In the present test case, we have $\beta = O(10^{-4})$ and in a few cells, the pressure can become negative in which case it is reset to a small value. This happens in at most one or two cells, also infrequently during the time iterations. In Figures~\ref{fig:blast1}, we have shown numerical solution for degree $k=1$ for squared velocity, pressure and magnetic pressure. The results at higher degree are not shown as they look similar to the case of $k=1$. 

\section{Summary and conclusions}
The paper develops an arbitrary order discontinuous Galerkin method for the compressible ideal MHD equations which naturally preserves the divergence-free condition on the magnetic field. This is known to be an important structural property of the solutions whose satisfaction is directly related to the accuracy and robustness of the method. The magnetic field is approximated in terms of Raviart-Thomas polynomials which  automatically ensures that the normal component of the magnetic field is continuous across the cell faces. The DG scheme evolves the degrees of freedom using a combination of face-based and cell-based DG schemes which automatically preserves the divergence of the magnetic field. Being a DG method, it requires numerical fluxes which are supplied via an approximate Riemann solver. We have proposed simple HLL-type multi-dimensional Riemann solvers which are consistent with their 1-D counterparts. Since we deal with non-linear flows, the solutions can develop discontinuities which requires some form of non-linear limiter but this can destroy the condition on the divergence. We can recover the divergence-free property by performing a local divergence-free reconstruction which makes use of information on the divergence and curl of the magnetic field. Many numerical tests presented here show the accuracy and robustness of the method. The positivity property of the scheme is however not possible to prove at present within the framework of divergence-free DG schemes, but we show that a heuristic application of scaling limiter can yield stable computations.
\section*{Acknowledgments}
Praveen Chandrashekar would like to acknowledge support from SERB-DST, India, under the MATRICS grant (MTR/2018/000006) and Department of Atomic Energy,  Government of India, under  project no.  12-R\&D-TFR-5.01-0520. Rakesh Kumar would like to acknowledge funding support from the  National Post-doctoral Fellowship (PDF/2018/002621)  administered by SERB-DST,  India.
\appendix
\section{Eigenvectors of the MHD system}
This section lists the right and left eigenvectors of the flux Jacobian in the $x$ direction which are taken from~\cite{Jiang1999},~\cite{Brio1988}. The eigenvector formulae correspond to the following ordering of the conserved variables: $[\rho,\rho v_x, \rho v_y, \rho v_x, B_x, B_y, B_z, \tote]$. Define $b = (b_x,b_y,b_z) = (B_x,B_y,B_z)/\sqrt{\rho}$ and $b^2 = b_x^2 + b_y^2 + b_z^2$. The sound speed, fast and slow speeds are given by
\[
a^2 = \frac{\gamma p}{\rho}, \qquad c_{f,s}^2 = \half\left[ a^2 + b^2 \pm \sqrt{(a^2 + b^2)^2 - 4 a^2 b_x^2} \right]
\]
Define
\[
 (\beta_y, \beta_z) = \begin{cases}
                       \frac{(B_y,B_z)}{\sqrt{B_y^2+B_z^2}} & \mbox{If}~  B_y^2+B_z^2\neq 0,\\
                       \left(\frac{1}{\sqrt{2}},\frac{1}{\sqrt{2}}\right) & \mbox{otherwise}
                      \end{cases}
\]

\[
 (\alpha_f, \alpha_s)=\begin{cases}
                       \frac{(\sqrt{a^2-c_s^2},\sqrt{c_f^2-a^2})}{\sqrt{c_f^2-c_s^2}} & \mbox{If}~  B_y^2+B_z^2\neq 0~\mbox{or}~ \gamma p\neq B_x^2,\\
                       \left(\frac{1}{\sqrt{2}},\frac{1}{\sqrt{2}}\right)             & \mbox{otherwise},
                      \end{cases}
\]

\begin{align*}
 \Gamma_f & = \alpha_f c_f v_x -\alpha_s c_s \mbox{sgn}(B_x)(\beta_y v_y +\beta_z v_z),\\
 \Gamma_a & = \mbox{sgn}(B_x) (\beta_z v_y -\beta_y v_z),\\
 \Gamma_s & = \alpha_s c_s v_x +\alpha_f c_f \mbox{sgn}(B_x) (\beta_y v_y +\beta_z v_z).
\end{align*}

The right eigenvectors are given by
\[
r_1= \begin{bmatrix}
\alpha_f\\
\alpha_f(v_x-c_f) \\
\alpha_f v_y+\alpha_s c_s \beta_y \mbox{sgn}{(B_x)} \\
\alpha_f v_z +\alpha_s c_s \beta_z \mbox{sgn}{(B_x)} \\
0\\
\frac{a\alpha_s \beta_y}{\sqrt{\rho}}\\
\frac{a\alpha_s \beta_z}{\sqrt{\rho}}\\
\alpha_f(\frac{v^2}{2}+c_f^2 -\gamma_2 a^2)-\Gamma_f
\end{bmatrix}, \qquad
r_2= \begin{bmatrix}
0\\
0 \\
-\beta_z \mbox{sgn}{(B_x)} \\
\beta_y \mbox{sgn}{(B_x)} \\
0\\
-\frac{\beta_z}{\sqrt{\rho}}\\
\frac{\beta_y}{\sqrt{\rho}}\\
-\Gamma_a
\end{bmatrix}, \qquad
r_3= \begin{bmatrix}
\alpha_s\\
\alpha_s(v_x-c_s) \\
\alpha_s v_y-\alpha_f c_f \beta_y \mbox{sgn}{(B_x)} \\
\alpha_s v_z -\alpha_f c_f \beta_z \mbox{sgn}{(B_x)} \\
0\\
-\frac{a\alpha_f \beta_y}{\sqrt{\rho}}\\
-\frac{a\alpha_f \beta_z}{\sqrt{\rho}}\\
\alpha_s(\frac{v^2}{2}+c_s^2 -\gamma_2 a^2)-\Gamma_s
\end{bmatrix}
\]

\[
r_4= \begin{bmatrix}
1\\
v_x\\
v_y \\
v_z\\
0\\
0\\
0\\
\frac{v^2}{2}
\end{bmatrix}, \qquad
r_5= \begin{bmatrix}
0\\
0 \\
0\\
0\\
1\\
0\\
0\\
B_x
\end{bmatrix}, \qquad
r_6= \begin{bmatrix}
\alpha_s\\
\alpha_s(v_x+c_s) \\
\alpha_s v_y+\alpha_f c_f \beta_y \mbox{sgn}{(B_x)} \\
\alpha_s v_z +\alpha_f c_f \beta_z \mbox{sgn}{(B_x)}\\
0\\
-\frac{a\alpha_f \beta_y}{\sqrt{\rho}}\\
-\frac{a\alpha_f \beta_z}{\sqrt{\rho}}\\
\alpha_s(\frac{v^2}{2}+c_s^2 -\gamma_2 a^2)+\Gamma_s
\end{bmatrix}
\]

\[
r_7= \begin{bmatrix}
0\\
0\\
-\beta_z \mbox{sgn}{(B_x)} \\
\beta_y \mbox{sgn}{(B_x)}  \\
0\\
\frac{\beta_z}{\sqrt{\rho}}\\
-\frac{\beta_y}{\sqrt{\rho}} \\
-\Gamma_a
\end{bmatrix}, \qquad
r_8= \begin{bmatrix}
\alpha_f\\
\alpha_f(v_x+c_f) \\
\alpha_f v_y-\alpha_s c_s \beta_y \mbox{sgn}{(B_x)} \\
\alpha_f v_z -\alpha_s c_s \beta_z \mbox{sgn}{(B_x)}\\
0\\
\frac{a\alpha_s \beta_y}{\sqrt{\rho}}\\
\frac{a\alpha_s \beta_z}{\sqrt{\rho}}\\
\alpha_f(\frac{v^2}{2}+c_f^2 -\gamma_2 a^2)+\Gamma_f
\end{bmatrix}
\]
The left eigenvectors are given by
\[
l_1= \frac{1}{2a^2}\begin{bmatrix}
\gamma_1 \alpha_f v^2 +\Gamma_f, ~
(1-\gamma)\alpha_f v_x -\alpha_f c_f,~
(1-\gamma)\alpha_f v_y +c_s \alpha_s \beta_y \mbox{sgn}{(B_x)},~\hdots\\
(1-\gamma)\alpha_f v_z +c_s \alpha_s \beta_z \mbox{sgn}{(B_x)},
-{(\gamma-1) \alpha_f B_x},~
(1-\gamma)\alpha_f B_y+a \alpha_s \sqrt{\rho}\beta_y,~\hdots\\
(1-\gamma)\alpha_f B_z+a \alpha_s \sqrt{\rho}\beta_z,~ 
(\gamma-1) \alpha_f 
\end{bmatrix}, 
\]

\[
l_2= \frac{1}{2}\begin{bmatrix}
\Gamma_a, ~
0,~
-\beta_z \mbox{sgn}{(B_x)},~
\beta_y \mbox{sgn}{(B_x)},~ 
0,~
-\sqrt{\rho} \beta_z,~
 \sqrt{\rho} \beta_y,~
0
\end{bmatrix}, 
\]

\[
l_3= \frac{1}{2a^2}\begin{bmatrix}
\gamma_1 \alpha_s v^2 +\Gamma_s, ~
(1-\gamma)\alpha_s v_x -\alpha_s c_s,~ 
(1-\gamma)\alpha_s v_y -c_f\alpha_f \beta_y \mbox{sgn}{(B_x)},~\hdots\\
(1-\gamma)\alpha_s v_z -c_f\alpha_f \beta_z \mbox{sgn}{(B_x)},~
-B_x (\gamma-1) \alpha_s,~
(1-\gamma)\alpha_s B_y - a \alpha_f \sqrt{\rho}\beta_y,~\hdots\\
(1-\gamma)\alpha_s B_z - a \alpha_f \sqrt{\rho}\beta_z,~ 
(\gamma-1) \alpha_s
\end{bmatrix}, 
\]

\[
l_4= \begin{bmatrix}
1-\frac{1}{2} \tau v^2, ~
\tau v_x,~ 
\tau v_y,~
\tau v_z,~
\tau B_x,~
\tau B_y,~
\tau B_z,~ 
-\tau
\end{bmatrix}, 
\]

\[
l_5= \begin{bmatrix}
0, ~
0,~ 
0,~
0,~
1,~
0,~
0,~ 
0
\end{bmatrix}, 
\]

\[
l_6= \frac{1}{2a^2}\begin{bmatrix}
\gamma_1 \alpha_s v^2 -\Gamma_s, ~
(1-\gamma)\alpha_s v_x +\alpha_s c_s,~ 
(1-\gamma)\alpha_s v_y +c_f\alpha_f \beta_y \mbox{sgn}{(B_x)},~\hdots\\
(1-\gamma)\alpha_s v_z +c_f\alpha_f \beta_z \mbox{sgn}{(B_x)},~
-B_x (\gamma-1) \alpha_s,~
(1-\gamma)\alpha_s B_y - a \alpha_f \sqrt{\rho}\beta_y,~\hdots\\
(1-\gamma)\alpha_s B_z - a \alpha_f \sqrt{\rho}\beta_z,~ 
(\gamma-1) \alpha_s
\end{bmatrix}, 
\]

\[
l_7= \frac{1}{2}\begin{bmatrix}
g_a, ~
0,~ 
-\beta_z \mbox{sgn}{(B_x)},~
 \beta_y \mbox{sgn}{(B_x)},~
0,~
\sqrt{\rho} \beta_z,~
-\sqrt{\rho} \beta_y,~ 
0
\end{bmatrix}, 
\]

\[
l_8= \frac{1}{2a^2}\begin{bmatrix}
\gamma_1 \alpha_f v^2 -\Gamma_f, ~
(1-\gamma)\alpha_f v_x +\alpha_f c_f,~
(1-\gamma)\alpha_f v_y -c_s \alpha_s \beta_y \mbox{sgn}{(B_x)},~~\hdots\\
(1-\gamma)\alpha_f v_z -c_s \alpha_s \beta_z \mbox{sgn}{(B_x)},
-{(\gamma-1) \alpha_f B_x},~
(1-\gamma)\alpha_f B_y+a \alpha_s \sqrt{\rho}\beta_y,~~\hdots\\
(1-\gamma)\alpha_f B_z+a \alpha_s \sqrt{\rho}\beta_z,~ 
(\gamma-1) \alpha_f 
\end{bmatrix}, 
\]

\section{Limiting and divergence-free reconstruction}
\label{sec:dfr}
When the solution on the faces $b_x$, $b_y$ is limited as explained in Section~(\ref{sec:lim}), we lose the divergence-free property of the magnetic field. To recover this property, we have to perform a divergence-free reconstruction step. We explain this reconstruction process for second, third and fourth order accuracy. The fifth order version is given in~\cite{Hazra2019} together with more details on the reconstruction idea. The resulting polynomial has the structure of the BDM polynomial on rectangles, see~\cite{Brezzi1991}, Equation (3.29). Here we explain how the RT polynomial can be modified to recover divergence-free property. \resa{In two dimensions, the BDM polynomial has $(k+1)(k+2)+2$ degrees of freedom while its divergence has $\half k(k+1)$ coefficients.} Up to third order accuracy, the reconstruction can be performed using only the solution on the faces $(b_x, b_y)$, but at fourth order and higher, we need additional information which is supplied in the form of the curl of the magnetic field. This additional information is available to us via the cell moments $\alpha, \beta$.

\resa{For example, at fourth order, the BDM polynomial has $(3+1)(3+2)+2 =22$ degrees of freedom, while the face solution $(b_x,b_y)$, which are polynomials of degree 3, provide $(4+4+4+4)-1=15$ degrees of freedom, where one piece of information is redundant since the face solution satisfies $\int_{\partial K} \B \cdot \un \ud s = 0$ on each cell $K$. The divergence-free condition on the BDM polynomial yields $\half (3)(3+1)=6$ conditions. So we have a total of $15+6 = 21$ equations but 22 coefficients to be determined. Hence we need to supply one additional piece of information to completely determine the BDM polynomial.}

\resa{We have the following inclusions $\poly_k^2 \subset \bdm(k) \subset \rt(k)$ and the RT polynomial has many more basis functions than the BDM polynomial}.  In the reconstruction process, we set some of the coefficients $\{a,b\}$ in the RT polynomial to zero but this does not affect the accuracy since only coefficients $a_{ij}, b_{ij}$ with $i+j > k$ are set to zero, and we  retain the $\poly_k$ part of the solution.
\subsection{Degree $k=1$}
The divergence of the vector field $\B \in \rt(1)$ is given by
\[
\begin{aligned}
\nabla\cdot\B = \frac{a_{10}}{\Delta x} + \frac{b_{01}}{\Delta y} & + \left( \frac{2 a_{20}}{\Delta x} + \frac{b_{11}}{\Delta y} \right) \phi_1(\xi) + \left( \frac{a_{11}}{\Delta x} + \frac{2 b_{02}}{\Delta y} \right) \phi_1(\eta) \\
& + 2 \left( \frac{a_{21}}{\Delta x} + \frac{b_{12}}{\Delta y} \right) \phi_1(\xi) \phi_1(\eta)
\end{aligned}
\]
The coefficients $a_{ij}, b_{ij}$ are related to the face solution and cell moments according to Table~\ref{tab:rec1}. The constant term is already zero. The linear terms can be made zero by setting
\[
\alpha_{00} = \half (a_0^- + a_0^+)  + \frac{1}{12} (b_1^+ - b_1^-) \frac{\dx}{\dy}
\]
\[
\beta_{00} = \half (b_0^- + b_0^+) + \frac{1}{12} (a_1^+ - a_1^-) \frac{\dy}{\dx}
\]
This will however destroy the conservation property since $\alpha_{00}$, $\beta_{00}$ are cell averages of $B_x$, $B_y$ respectively. The bilinear term can be made zero by individually setting $a_{21} = b_{12} = 0$ which yields
\[
\alpha_{01} = \half(a_1^- + a_1^+), \qquad \beta_{10} = \half (b_1^- + b_1^+)
\]
By this process we would have modified all the cell moments and the resulting reconstruction coincides with that of Balsara.
\subsection{Degree $k=2$}
The divergence of the vector field $\B \in \rt(2)$ is given by
\[
\begin{aligned}
\nabla\cdot\B = \ & \left[ \frac{1}{\Delta x}\left(a_{10}+\frac{a_{30}}{10}\right) + \frac{1}{\Delta y}\left(b_{01} + \frac{b_{03}}{10}\right) \right] + \left[\frac{2}{\Delta x}a_{20} +\frac{1}{\Delta y}\left(b_{11}+\frac{b_{13}}{10}\right)\right] \phi_1(\xi) \\
&+ \left[\frac{1}{\Delta x}\left(a_{11}+\frac{a_{31}}{10}\right)+\frac{2}{\Delta y}b_{02}\right] \phi_1(\eta) + \left[ \frac{3}{\Delta x}a_{30}+\frac{1}{\Delta y}\left(b_{21}+\frac{b_{23}}{10}\right)\right] \phi_2(\xi) \\
&+ 2\left[ \frac{a_{21}}{\Delta x}+\frac{b_{12}}{\Delta y}\right] \phi_1(\xi) \phi_1(\eta) + \left[ \frac{1}{\Delta x}\left(a_{12}+\frac{a_{32}}{10}\right)+\frac{3}{\Delta y}b_{03}\right] \phi_2(\eta) \\
&+ \left[ \frac{2}{\Delta x} a_{22} +\frac{3}{\Delta y}b_{13}\right] \phi_1(\xi)\phi_2(\eta)+\left[ \frac{3}{\Delta x}a_{31}+\frac{2}{\Delta y}b_{22} \right]\phi_2(\xi)\phi_1(\eta)\\
&+ 3\left[ \frac{a_{32}}{\Delta x} + \frac{b_{23}}{\Delta y}\right] \phi_2(\xi)\phi_2(\eta)
\end{aligned}
\]
The constant term is already zero. The linear terms can be made zero by setting
\[
\alpha_{00} = \half (a_0^- + a_0^+)  + \frac{1}{12} (b_1^+ - b_1^-) \frac{\dx}{\dy}, \qquad
\beta_{00} = \half (b_0^- + b_0^+) + \frac{1}{12} (a_1^+ - a_1^-) \frac{\dy}{\dx}
\]
The quadratic terms are zero by choosing
\[
\alpha_{10} = a_0^+ - a_0^- + \frac{1}{30}(b_2^+ - b_2^-) \frac{\dx}{\dy}, \qquad \beta_{01} = b_0^+ - b_0^- + \frac{1}{30}(a_2^+ - a_2^-) \frac{\dy}{\dx}
\]
and also setting $a_{21} = b_{12} = 0$ which yields
\[
\alpha_{01} = \shalf (a_1^- + a_1^+), \qquad \beta_{10} = \shalf (b_1^- + b_1^+)
\]
In the cubic terms, we set each coefficient to zero, $a_{22} = a_{31} = a_{32} = b_{22} = b_{13} = b_{23} = 0$, which yields \\
\begin{minipage}{0.5\textwidth}
\begin{eqnarray*}
\alpha_{02} &=& \shalf(a_2^- + a_2^+)  \\
\alpha_{11} &=& a_1^+ - a_1^-  \\
\end{eqnarray*}
\end{minipage}
\begin{minipage}{0.5\textwidth}
\begin{eqnarray*}
\beta_{20} &=& \shalf (b_2^- + b_2^+) \\
\beta_{11} &=& b_1^+ - b_1^- \\
\end{eqnarray*}

\end{minipage}

\noindent
Finally, in the biquadratic term, we set $a_{32} = b_{23} = 0$ to obtain
\[
\alpha_{12} = a_2^+ - a_2^-, \qquad \beta_{21} = b_2^+ - b_2^-
\]
\subsection{Degree $k=3$}
The divergence of the vector field $\B \in \rt(3)$ is given by
\[
\begin{aligned}
\nabla\cdot\B = \ & \left[ \frac{1}{\Delta x}\left(a_{10}+\frac{a_{30}}{10}\right) + \frac{1}{\Delta y}\left(b_{01} + \frac{b_{03}}{10}\right) \right]\\
+& \left[\frac{1}{\Delta x}\left(2a_{20}+\frac{6}{35}a_{40}\right) +\frac{1}{\Delta y}\left(b_{11}+\frac{b_{13}}{10}\right)\right] \phi_1(\xi) \\
+& \left[\frac{1}{\Delta x}\left(a_{11}+\frac{a_{31}}{10}\right)+\frac{1}{\Delta y}\left(2b_{02}+\frac{6}{35}b_{04}\right)\right] \phi_1(\eta)\\
+& \left[ \frac{3}{\Delta x}a_{30}+\frac{1}{\Delta y}\left(b_{21}+\frac{b_{23}}{10}\right)\right] \phi_2(\xi)+ \left[ \frac{1}{\Delta x}\left(a_{12}+\frac{a_{32}}{10}\right)+\frac{3}{\Delta y}b_{03}\right] \phi_2(\eta) \\
+& \left[ \frac{1}{\Delta x}\left(2a_{21}+\frac{6}{35}a_{41}\right)+\frac{1}{\Delta y}\left(2b_{12}+\frac{6}{35}b_{14}\right)\right] \phi_1(\xi) \phi_1(\eta) \\
+ & \left[\frac{4}{\Delta x}a_{40} + \frac{1}{\Delta y}\left(b_{31}+\frac{b_{33}}{10}\right)\right]\phi_3(\xi)+\left[ \frac{1}{\Delta x}\left(2 a_{22}+\frac{6}{35}a_{42}\right) +\frac{3}{\Delta y}b_{13}\right] \phi_1(\xi)\phi_2(\eta)\\
+&\left[ \frac{3}{\Delta x}a_{31}+\frac{1}{\Delta y}\left(2b_{22}+\frac{6}{35}b_{24}\right) \right]\phi_2(\xi)\phi_1(\eta)
 + \left[\frac{1}{\Delta x}\left(a_{13}+\frac{a_{33}}{10}\right) +\frac{4}{\Delta y} b_{04}\right]\phi_3(\eta)\\
+& 3\left[ \frac{a_{32}}{\Delta x} + \frac{b_{23}}{\Delta y}\right] \phi_2(\xi)\phi_2(\eta) + \left[\frac{1}{\Delta x}\left(2 a_{23}+\frac{6}{35}a_{43}\right)+\frac{4}{\Delta y} b_{14}\right]\phi_1(\xi) \phi_3(\eta) \\
+&\left[\frac{4}{\Delta x}a_{41} +\frac{1}{\Delta y}\left(2b_{32} +\frac{6}{35}b_{34} \right) \right] \phi_3(\xi)\phi_1(\eta) \\
+& \left[\frac{3}{\Delta x} a_{33} +\frac{4}{\Delta y}b_{24}\right]\phi_2(\xi) \phi_3(\eta)
+\left[ \frac{4}{\Delta x}a_{42} + \frac{3}{\Delta y}b_{33}\right] \phi_3(\xi) \phi_2(\eta) \\
+& \left[ \frac{4}{\Delta x} a_{43}+\frac{4}{\Delta y}b_{34} \right]\phi_3(\xi) \phi_3(\eta)
\end{aligned}
\]
The constant term is already zero. The linear terms can be made zero by setting
\[
\alpha_{00} = \half (a_0^- + a_0^+)  + \frac{1}{12} (b_1^+ - b_1^-) \frac{\dx}{\dy}, \qquad
\beta_{00} = \half (b_0^- + b_0^+) + \frac{1}{12} (a_1^+ - a_1^-) \frac{\dy}{\dx}
\]
The quadratic terms which are coefficients of $\phi_2(\xi)$, $\phi_2(\eta)$ become zero by choosing
\[
\alpha_{10} = a_0^+ - a_0^- + \frac{1}{30}(b_2^+ - b_2^-) \frac{\dx}{\dy}, \qquad \beta_{01} = b_0^+ - b_0^- + \frac{1}{30}(a_2^+ - a_2^-) \frac{\dy}{\dx}
\]
The coefficient of $\phi_1(\xi)\phi_1(\eta)$ gives only one equation but there are two unknowns; adding an extra equation $\omega= \beta_{10}-\alpha_{01}$, we can solve for the two coefficients
\[
\alpha_{01} = \frac{1}{\left(1+\frac{\dy}{\dx}\right)}\left(r_2-\omega +r_1\frac{\dy}{\dx}\right), \qquad \beta_{10} = \omega +\alpha_{01}
\]
where $r_1=\half(a_1^{-}+a_1^{+})$ and $r_2=\half(b_1^{-}+b_1^{+})$. The cubic terms are zeros by choosing
\[
 \alpha_{20} = -\half(b_1^{+}-b_1^{-})\frac{\dx}{\dy}+\frac{3}{140}(b_3^{+}-b_3^{-})\frac{\dx}{\dy}, \qquad
 \beta_{02}=-\frac{1}{2}(a_1^{+}-a_1^{-})\frac{\dy}{\dx}+\frac{3}{140}(a_3^{+}-a_3^{-})\frac{\dy}{\dx}
\]
and setting $2 a_{22}+\frac{6}{35}a_{42} = 2b_{22}+\frac{6}{35}b_{24} = a_{31} = b_{13} = 0$, which yields\\
\begin{minipage}{0.5\textwidth}
\begin{eqnarray*}
\alpha_{02} &=& \shalf(a_2^- + a_2^+)  \\
\alpha_{11} &=& a_1^+ - a_1^-  \\
\end{eqnarray*}
\end{minipage}
\begin{minipage}{0.5\textwidth}
\begin{eqnarray*}
\beta_{20} &=& \shalf (b_2^- + b_2^+) \\
\beta_{11} &=& b_1^+ - b_1^- \\
\end{eqnarray*}
\end{minipage}

\noindent
In the higher order term greater then three, we set each coefficient to zero,
$a_{32}=b_{23} = 2a_{23}+\frac{6}{35}a_{43} = b_{14}=a_{41} = 2b_{32}+\frac{6}{35}b_{34} = a_{33}=b_{24}=a_{42}=b_{33}=a_{43}=b_{34}=0$, which yields\\
\begin{minipage}{0.5\textwidth}
\begin{eqnarray*}
 \alpha_{12} &=& a_2^+ - a_2^-\\
 \alpha_{03} &=& \shalf(a_3^+ +a_3^-)\\
 \alpha_{21} &=& 6(r_1-\alpha_{01})\\
 \alpha_{13} &=& \alpha_3^+ -\alpha_3^-\\
 \alpha_{22} &=& 0\\
 \alpha_{23} &=& 0
\end{eqnarray*}
\end{minipage}
\begin{minipage}{0.5\textwidth}
\begin{eqnarray*}
 \beta_{21} &=& b_2^+ - b_2^-\\
 \beta_{12} &=& 6(r_2-\beta_{10})\\
 \beta_{30} &=& \shalf(b_3^+ +b_3^-)\\
 \beta_{22} &=& 0\\
 \beta_{31} &=& b_3^+-b_3^-\\
 \beta_{32} &=& 0
\end{eqnarray*}
\end{minipage}\\

\noindent
We see that at fourth order, we need an extra information which we took in the form of the quantity $\omega$ in order to complete the divergence-free reconstruction. Note $\omega$ gives us information about the curl of the magnetic field.
\section{Setting the initial condition}
\label{sec:ic}
Let $\psi_h \in \tpoly_{k+1,k+1}$ be a continuous interpolation of the magnetic potential $\psi$ which can be achieved using $(k+2)\times(k+2)$ GLL nodes. Then we can set the magnetic field as $B_x = \df{\psi_h}{y}$, $B_y = -\df{\psi_h}{x}$ which will be exactly divergence-free. But in our work, we want to set the initial condition in terms of the polynomials $b_x,b_y$ and the moments $\alpha,\beta$. We can perform an $L^2$ projection of $\nabla\times(\psi_h e_z)$ to initialize $b_x,b_y$ which will be exact, and the moments can be computed using the same GLL nodes for quadrature as are used to define $\psi_h$. Let $\xi_i, i=1,2,\ldots,k+2$ denote the GLL nodes and let $\ell_i(\xi), i=1,2,\ldots,k+2$ be the Lagrange polynomials. Define the barycentric weights
\[
w_j = \frac{1}{\prod_{i=1, i\ne j}^{k+2} (\xi_j - \xi_i)}
\]
Then the derivatives of Lagrange polynomials at the GLL nodes are given by
\[
D_{ij}  = \ell_j'(\xi_i) = \frac{w_j}{w_i} \frac{1}{\xi_i - \xi_j}, \quad i \ne j, \qquad D_{ii} = \ell_i'(\xi_i) = -\sum_{j=1, j \ne i}^{k+2} D_{ij}
\]
The derivatives of the potential at the GLL nodes are given by
\[
\df{\psi_h}{x}(\xi_i,\xi_j) = \frac{1}{\dx} [D \cdot \psi(:,j)]_i, \qquad \df{\psi_h}{y}(\xi_i,\xi_j) = \frac{1}{\dy} [D \cdot \psi(i,:)]_j
\]
The cell moments are initialized as $\alpha = \alpha^h(\partial_y \psi_h)$ and $\beta = \beta^h(-\partial_x\psi_h)$ where the superscript $h$ denotes that we compute the integrals using $(k+2)^2$-point GLL quadrature which is exact for the integrands involved in the cell moments.
\section{Running the code}
The code is written in Fortan90 and works only in serial. Some OpenMP has been implemented but this is not properly tested and may have some bugs. Each test case must be implemented in a header files like \texttt{alfven.h} and the test case is selected while compiling the code along with some other options. The way to compile the code is
\begin{lstlisting}
make <problem> \
     DEGREE=0|1|2|3 \
     NX=<integer> \
     NY=<integer> \
     FLUX=lxf|hll|hllc \
     LIMIT=none|tvd|weno|mdl \
     INDICATOR=no|yes \
     POSLIM=no|yes \
     CHECKPOS=no|yes
\end{lstlisting}
\begin{itemize}
\item \texttt{<problem>} can be \texttt{alfven, vortex, ot, rotor, rstube, loop, briowu, blast}
\item \texttt{NX} and \texttt{NY} are grid sizes in $x$ and $y$ directions.
\item \texttt{LIMIT=none} is default; if you don't want limiter, this parameter need not be specified.
\end{itemize}
For example to run the vortex test which does not require any limiter, compile and run like this
\begin{lstlisting}
make vortex DEGREE=3 NX=100 NY=100 FLUX=hll
./mhd > log.txt &
\end{lstlisting}
The solution is saved in Tecplot format in files named \texttt{avg\#\#\#\#.plt} which can be viewed using VisIt. These files contain the cell average solution. A more detailed solution with sub-sampling is also written at initial and final times in files \texttt{sol000.plt} and \texttt{sol0001.plt} respectively. Some test cases write specialized files also which is shown at the end of the code run. There are many Python scripts available for making plots, e.g.,
\begin{lstlisting}
visit -nowin -cli -s <path to>/contour_ot.py Rho
\end{lstlisting}
will generate contour plots of density, while
\begin{lstlisting}
visit -nowin -cli -s <path to>/pseudo_ot.py Rho
\end{lstlisting}
will generate color plots of density.


\bibliographystyle{siam}
\bibliography{bbibtex}

\end{document}